\documentclass[11pt,twoside,openany]{book}
\usepackage{amsmath}
\usepackage{amssymb}
\usepackage{latexsym}
\usepackage{graphicx}
\usepackage{xcolor}
\usepackage{array}
\usepackage{lscape}
\usepackage{rotating}
\usepackage{hyperref}
\usepackage{fullpage}
\usepackage{amsrefs}
\usepackage{tikz}
\usetikzlibrary{matrix,arrows,decorations.pathmorphing, cd}
\usepackage[UKenglish]{babel}
\usepackage[UKenglish]{isodate}
\usepackage{enumitem}
\usepackage{tcolorbox}

\usepackage{multicol}
\usepackage{fancyhdr}
\usepackage{titlesec}
\titleformat{\section}[block]{\bfseries\filcenter}%
{{\upshape\S\,\thesection\enspace}}{.5em}{}
\titleformat{\subsection}[block]{\scshape\filcenter}%
{{\upshape\thesubsection\enspace}}{.5em}{}

\usepackage{fancyvrb}
\usepackage{algorithmic}
\usepackage[pointlessenum]{paralist}

\usepackage{wasysym}
\usepackage{amsthm}

\newcommand{\Hom}       {\operatorname{Hom}}
\newcommand{\Mod}	    {\operatorname{Mod}}
\newcommand{\Ext}		{\operatorname{Ext}}
\newcommand{\Tor}		{\operatorname{Tor}}
\newcommand{\Tr}			{\operatorname{Tr}}
\newcommand{\Arr}		{\operatorname{Arr}}
\newcommand{\Fr}			{\operatorname{Fr}}
\newcommand{\St}			{\operatorname{St}}

\newcommand{\ann}       {\operatorname{ann}}
\newcommand{\cok}       {\operatorname{cok}}

\newcommand{\im}       {\operatorname{im}}
\newcommand{\supp}	   {\operatorname{supp}}
\newcommand{\id}			{\operatorname{id}}
\newcommand{\loc}		{\operatorname{loc}}

\newcommand{\F}         {{\mathbb{F}}}          
\newcommand{\Z}         {{\mathbb{Z}}}
\newcommand{\Q}         {{\mathbb{Q}}}
\newcommand{\R}         {{\mathbb{R}}}

\newcommand{\mcF}		{{\mathcal{F}}}

\newcommand{\mcM}		{{\mathcal{M}}}
\newcommand{\mcZ}		{{\mathcal{Z}}}

\newcommand{\al}        {\alpha}
\newcommand{\bt}        {\beta} 
\newcommand{\gm}        {\gamma}
\newcommand{\dl}        {\delta}

\newcommand{\io}			{\iota}

\newcommand{\lmb}		{\lambda}
\newcommand{\sg}        {\sigma}
\newcommand{\ta}			{\theta}
\newcommand{\vphi}		{\varphi}

\newcommand{\Ta}       {\Theta}
\newcommand{\Sg}			{\Sigma}
\newcommand{\Lmb}		{\Lambda}
\newcommand{\Dl}			{\Delta}

\newcommand{\tns}        {\otimes}
\newcommand{\tk}        {\setminus}
\newcommand{\subs}      {\subseteq}

\newcommand{\x}         {\times}

\newcommand{\w}			{\wedge}

\newcommand{\bbm}       {\left[\begin{matrix}}
\newcommand{\ebm}       {\end{matrix}\right]}

\newcommand{\mcA}		{\mathcal{A}}
\newcommand{\mcB}		{\mathcal{B}}
\newcommand{\mcC}	    {\mathcal{C}}
\newcommand{\mcD}		{\mathcal{D}}

\newcommand{\mcP}		{\mathcal{P}}
\newcommand{\mcS}		{\mathcal{S}}
\newcommand{\mcW}		{\mathcal{W}}

\newcommand{\mfm}		{\mathfrak{m}}

\newcommand{\incl}		{\hookrightarrow}

\newcommand{\td}			{\widetilde}

\newtheorem{theorem}{Theorem}[section]
\newtheorem{conj}[theorem]{Conjecture}
\newtheorem{lem}[theorem]{Lemma}
\newtheorem{prop}[theorem]{Proposition}
\newtheorem{cor}[theorem]{Corollary}
\theoremstyle{definition}
\newtheorem{rmk}[theorem]{Remark}
\newtheorem{definition}[theorem]{Definition}
\newtheorem{exmp}[theorem]{Example}

\cleanlookdateon
\date{\today}
\begin{document}
\frontmatter

\begin{titlepage}
\vspace*{\stretch{1}}
\begin{center}
\includegraphics{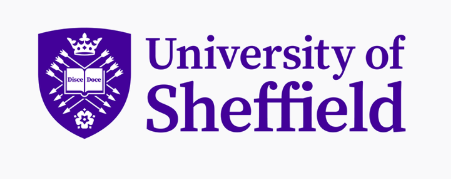}
\end{center}
\begin{center}\bfseries
{\LARGE Reflexive Modules, the Infinite Root Algebra and the Generating Hypothesis}\\
\vspace{1.2cm}
{\LARGE Oliver House}\\
\vspace{1cm}
{\Large Submitted for the degree of Doctor of Philosophy}\\
\bigskip
{\Large School of Mathematical and Physical Sciences}\\
\bigskip
{\large April 2025}\\
\bigskip\bigskip\bigskip
{\large Supervisor: Neil Strickland}
\end{center}
\vspace*{\stretch{3}}
\begin{center}
\textbf{\large University of Sheffield}
\end{center}
\end{titlepage}

\thispagestyle{empty}
\vspace*{\stretch{3}}
\begin{center}\scshape
Abstract
\end{center}
\vspace*{\stretch{1}}
\begin{quote}

This thesis concerns the algebraic consequences of Freyd's Generating Hypothesis, and explores the question of whether there exists a self-injective ring $R$ that can be constructed purely algebraically that exhibits some of the known properties of the stable homotopy ring, including some conjectured properties that follow from Freyd's Generating Hypothesis. As an example, we investigate the infinite root algebra of Hahn series $P$, firstly by establishing results for the related Hahn ring $A$. In particular, we prove that the $\Ta$-reflexive $A$-modules and the multibasic $A$-modules are the same.

\end{quote}
\vspace*{\stretch{7}}
\newpage

\thispagestyle{empty}
\vspace*{\stretch{1}}
\begin{center}\itshape
To Carien
\end{center}
\vspace*{\stretch{3}}

\newpage

\thispagestyle{empty}
\vspace*{\stretch{3}}
\begin{center}\scshape
Acknowledgements 
\end{center}
\vspace*{\stretch{1}}
\begin{quote}

I would like to thank the EPSRC for my PhD scholarship; my supervisor, Neil Strickland, for his invaluable mathematical insights; Fiona Maisey, for support and guidance through what has been an extremely difficult few years; Moty Katzman, for a number of fruitful conversations about mathematics, programming and finance; the users of the Stack Exchange network for helping me resolve all manner of queries. I would also like to thank my partner Carien, my parents, my friends, and my cat Basil, without whom I doubt I would have made it through to the end. 

\end{quote}
\vspace*{\stretch{7}}

\newpage

\thispagestyle{empty}
\vspace*{\stretch{3}}
\begin{center}\scshape
Declaration 
\end{center}
\vspace*{\stretch{1}}
\begin{quote}

I, the author, confirm that this thesis is my own work. It has not been presented for any other any other degree. Where information has been derived from other sources, I confirm that this has been indicated in the work. 

\end{quote}
\vspace*{\stretch{7}}

\tableofcontents

\cleardoublepage

\chapter{Introduction}

The Spanier-Whitehead category $\mcF$ of finite spectra is a triangulated category of central importance to stable homotopy theory \cite{Mar83}. Its objects, the finite spectra, can be thought of as pairs $(X,n)$, where $X$ is a finite pointed CW complex and $n\in\Z$. It is equivalent to the full subcategory of ‘finite' objects in a larger category $\mcB$, namely Boardman's stable homotopy category of spectra \cite{Str20}. From $\mcF$, we can construct a graded-commutative ring $R$ and a functor $\pi_*:\mcF\to\Mod_R$. Freyd's Generating Hypothesis claims that this functor $\pi_*$ is faithful; this claim has some interesting algebraic consequences, many of which are outlined in \cite{Hov07}, \cite{ShS14} and \cite{Fre65}. For example, if Freyd's Generating Hypothesis holds, then $\pi_*$ must be fully faithful, and therefore $\mcF$ must be equivalent to a full subcategory of $\Mod_R$. Another consequence is that $\pi_*(X)$ is injective for all $X\in\mcF$, and in particular that $R=\pi_*(S)$ is self-injective.

It is known that $R_i=0$ for $i<0$ and $R_0=\Z$, and $R_i$ has been calculated up to around $i=60$ (see \cite{IWX20}). However, $R_i$ in general is not known. Following \cite{ShS14}, we shall focus on a variant $\mcF_p$ of $\mcF$ for a fixed prime $p$. This is obtained by tensoring the hom-sets of $\mcF$ with $\Z_p$, and we briefly discuss in \S\ref{p completion} how $\mcF_p$ sits inside a different ambient category to $\mcB$. It is established in \cite{Hov07} that $\pi_*$ is faithful in the $\mcF$ case if and only if it is for the $\mcF_p$ case. 

Given that $R_i$ is difficult to compute beyond a certain point, it is natural to ask if there might be purely algebraic examples of graded-commutative rings exhibiting similar properties to the stable homotopy ring, including some of those properties which follow as consequences of Freyd's Generating Hypothesis. Since such purely algebraic examples may be easier to work with in many respects than the stable homotopy ring itself, they might shed some light on the Generating Hypothesis by mimicking the conjectured behaviour of the stable homotopy ring. Such an algebraic example $R$ should exhibit a triangulated category of injective graded $R$-modules containing $R$ itself, representing the embedding of $\mcF_p$ into $\Mod_R$. A number of interesting candidates are identified in \cite{ShS14}, one of which is the infinite root algebra, which will be the main focus of this thesis. 

In Chapters \ref{hahn ring chapter}, \ref{reflexivity} and \ref{infinite root algebra chapter}, we investigate the properties of the Hahn ring $A$ and the infinite root algebra $P$. Despite the self-injective ring $P$ showing promising signs, we establish that $\Mod_P$ cannot contain an appropriate triangulated subcategory for it be effectively compared with the stable homotopy ring.

In Chapter \ref{freyd cats chapter}, we present an axiomatic generalisation of the situation for the $p$-complete Spanier-Whitehead category $\mcF_p$, following \cite{HPS97}, leading to the definition of Freyd and semi-Freyd categories, where a Freyd category is a semi-Freyd category for which the Generating Hypothesis holds.

In Chapter 2, we develop the theory of multibasic modules over the Hahn ring $A$, obtained from the field of Hahn series with value group $\R$ and residue field $\F$. Then in Chapter 3, we develop the theory of $\Ta$-reflexive $A$-modules, and show that these are the same as the multibasic modules. In Chapter 4, we use the results in Chapters 2 and 3 to establish corresponding results about the infinite root algebra $P$, a quotient of $A$. Finally, we include an appendix covering the case of $\Ta$-reflexive and multibasic modules over a complete discrete valuation ring, which is essentially a simplified account of the contents of Chapters 2 and 3. 

\newpage

\mainmatter

\chapter{Freyd Categories and the Generating Hypothesis} \label{freyd cats chapter}

We start this chapter by fixing an additive category $\mcC$ and an additive automorphism $[1]$ on $\mcC$ sending $f:X\to Y$ to $f[1]:X[1]\to Y[1]$. We shall introduce additional structure and properties on $\mcC$ bit by bit to demonstrate how the results depend on various assumptions. We denote the set of morphisms $X\to Y$ in $\mcC$ by $[X,Y]$ and the coproduct by $\bigvee$, mirroring the notation for a wedge sum in topology. We denote the category of abelian groups by Ab.

In Section \ref{triang cat section}, we provide a brief exposition of the fundamentals of triangulated categories, following \cite{HoJ10} and then \cite{HPS97} for some later definitions. This involves introducing the notion of a diagram \begin{center}
    \begin{tikzcd} [sep=2cm]
     X \arrow[r,"u"] & 
     Y \arrow[r,"v"] & 
     Z \arrow[r,"w"] & 
     X[1]\;\;\;.
    \end{tikzcd}
 \end{center} being ‘exact’ relative to the automorphism $[1]$. Note that $\mcC$ is not necessarily abelian; we instead choose a class of such diagrams (called ‘triangles') and insist on various axioms being satisfied to give meaning to a triangle being exact. Of particular interest will be the full subcategory $\mcF$ of ‘small' objects (those $T$ for which the functor $[T,-]$ is ‘homological', meaning that in addition to being exact it also preserves arbitrary coproducts). We shall see that $\mcF$ is itself a triangulated category in a natural way.
 
In Section \ref{freyd comp}, we follow \cite{Fre65} to construct the Freyd completion of $\mcF$, whose objects are the morphisms in $\mcF$ and whose morphisms are equivalence classes of commuting squares in $\mcF$. We can think of the Freyd completion as an extension of $\mcF$ to an abelian category, the crucial change being that ‘weak' kernels and cokernels in $\mcF$ (where the requirement of uniqueness in the universal property is relaxed) become genuine kernels and cokernels. We prove that $\Fr(\mcF)$ is an abelian category with enough projectives and enough injectives. Adding in some further assumptions about $\mcC$, we use \cite{Nee01} to establish that the projectives and injectives in $\Fr(\mcF)$ coincide (this is known as an abelian Frobenius category), and that the projective-injective objects are precisely those isomorphic to an identity morphism in $\mcF$, so that we can regard $\mcF$ as the full subcategory of projective-injective objects in $\Fr(\mcF)$ \cite{Fre65}. Finally, we establish a periodicity result relating to Ext groups in $\Fr(\mcF)$ \cite{Fre65}.

Having just constructed an (abelian) Frobenius category $\Fr(\mcF)$ from $\mcC$, in Section \ref{stable cat section} we turn our attention to an example of a triangulated category that can be constructed from a Frobenius category, following \cite{HoZ12}. Given a Frobenius category $\mcA$, we quotient out any morphisms that factor through a projective-injective object to form its ‘stable category' $\St(\mcA)$. The projective-injective objects in $\mcA$ become zero objects in $\St(\mcA)$. Connecting this with Section 1.2, we can thus form another triangulated category $\St(\Fr(\mcF))$ from $\mcC$. However, as we will see, the functor $\mcF\to\St(\Fr(\mcF))$ is simply the zero functor. 

In Section \ref{tensor triang}, we follow \cite{HPS97} to introduce a notion of tensor product to $\mcC$ with a generalisation of the tensor-hom adjunction, leading to the definition of a tensor triangulated category. We use this tensor product to construct an associated graded-commutative ring $R$ and a functor $\pi_*:\mcC\to\Mod_R$, where $\Mod_R$ is the category of $\Z$-graded $R$-modules.

In Section \ref{reflexive graded modules}, we fix a graded-commutative ring $R$ and an injective graded $R$-module $\Ta$. We then construct an abelian full subcategory $\mcM$ of $\Mod_R$ with a duality functor $D:\mcM\to\mcM$, and establish various properties of $\mcM$. These results will apply to the stable homotopy ring $R$ associated with our fixed category $\mcC$ in Section \ref{freyd semi freyd}, with similar arguments being used in Chapters \ref{reflexivity} and \ref{infinite root algebra chapter} relating to the Hahn ring $A$ and infinite root algebra $P$.

In Section \ref{monogenic}, we insist on some further properties relating to the closed symmetric monoidal structure, leading to the definition of a monogenic stable homotopy category, as described in \cite{HPS97}. We also discuss four equivalent interpretations of the subcategory $\mcF$ established there. 

In Section \ref{freyd semi freyd}, we impose various further conditions on $\mcC$ relating to $R$, leading to the definition of a semi-Freyd category. A Freyd category is then a semi-Freyd category for which the restricted functor $\pi_*:\mcF\to\Mod_R$ is faithful, a much stronger assumption. We then follow \cite{Hov07} to explore various algebraic consequences stemming from this assumption. For example, if $\mcC$ is a Freyd category, then $\pi_*$ is also full, and therefore an embedding of categories. We shall also see that $\pi_*(X)$ is injective in $\Mod_R$ for all $X\in\mcF$, and in particular that $R$ is self-injective as a graded-commutative ring. Finally, we shall see that $R$ has no non-trivial finitely presented graded ideals \cite{ShS14}.

In Section \ref{sw cat section}, we construct the Spanier-Whitehead category $\mcF$. In Section \ref{p completion}, we fix a prime $p$ and construct the $p$-completed Spanier-Whitehead category, starting with Boardman's stable homotopy category $\mcB$ \cite{Boa65}. Following \cite{Str24}, we see that this is equivalent to the full subcategory of small objects in some semi-Freyd category $\mcC$. The conjecture that we can assume here that $\mcC$ is in fact a Freyd category (rather than just semi-Freyd) is a variant of the Generating Hypothesis, attributed to Peter Freyd in the 1960s \cite{Fre65}.  

\section{Triangulated Categories} \label{triang cat section}

This section is largely based on \cite{HoJ10}, with some later definitions coming from \cite{HPS97}. We start by fixing an additive category $\mcC$ and an additive automorphism $[1]$ sending $f:X\to Y$ to $f[1]:X[1]\to Y[1]$, which we call the \textit{suspension} of $\mcC$. For each $i\in\Z$, we write $[1]^i$ as $[i]$ and call $[-1]$ the \textit{desuspension} of $\mcC$.  

\begin{definition} \cite[\S 3]{HoJ10} We define the additive category $\Tr(\mcC)$ as follows. 

A \textit{triangle} in $\mcC$ consists of objects $X$, $Y$ and $Z$, and morphisms $u:X\to Y$, $v:Y\to Z$ and $w:Z\to X[1]$. We represent this via the following diagram. 

\begin{center}
    \begin{tikzcd} [sep=2cm]
     X \arrow[r,"u"] & 
     Y \arrow[r,"v"] & 
     Z \arrow[r,"w"] & 
     X[1]\;\;\;.
    \end{tikzcd}
 \end{center} The objects in $\Tr(\mcC)$ are precisely the triangles in $\mcC$. The morphisms from \begin{center}
    \begin{tikzcd} [sep=2cm]
     X \arrow[r,"u"] & 
     Y \arrow[r,"v"] & 
     Z \arrow[r,"w"] & 
     X[1]
    \end{tikzcd}
 \end{center} to \begin{center}
    \begin{tikzcd} [sep=2cm]
     X' \arrow[r,"u'"] & 
     Y' \arrow[r,"v'"] & 
     Z' \arrow[r,"w'"] & 
     X'[1]
    \end{tikzcd}
 \end{center} are the triples $(f,g,h)$, where $f:X\to X'$, $g:Y\to Y'$, $h:Z\to Z'$ and the following diagram commutes. \begin{center}
    \begin{tikzcd} [sep=2cm]
     X \arrow[r,"u"] \arrow[d,"f"] & 
     Y \arrow[r,"v"] \arrow[d,"g"] & 
     Z \arrow[r,"w"] \arrow[d,"h"] & 
     X[1] \arrow[d,"f{[1]}"] \\
     X' \arrow[r,"u'"] &
     Y' \arrow[r,"v'"] &
     Z' \arrow[r,"w'"] &
     X'[1] 
    \end{tikzcd}
 \end{center} It is routine to check that these morphisms form an abelian group under componentwise addition. Composition of morphisms is defined componentwise and is clearly biadditive. The identity on \begin{center}
    \begin{tikzcd} [sep=2cm]
     X \arrow[r,"u"] & 
     Y \arrow[r,"v"] & 
     Z \arrow[r,"w"] & 
     X[1]
    \end{tikzcd}
 \end{center} is $(\id_X,\id_Y,\id_Z)$. It is clear that $\Tr(\mcC)$ has zero object \begin{center}
    \begin{tikzcd} [sep=2cm]
     0 \arrow[r] & 
     0 \arrow[r] & 
     0 \arrow[r] & 
     0[1]\;\;\;.
    \end{tikzcd}
 \end{center} The biproduct of the triangles \begin{center}
    \begin{tikzcd} [sep=2cm]
     X \arrow[r,"u"] & 
     Y \arrow[r,"v"] & 
     Z \arrow[r,"w"] & 
     X[1]
    \end{tikzcd}
 \end{center} and \begin{center}
    \begin{tikzcd} [sep=2cm]
     X' \arrow[r,"u'"] & 
     Y' \arrow[r,"v'"] & 
     Z' \arrow[r,"w'"] & 
     X'[1]
    \end{tikzcd}
 \end{center} is the triangle \begin{center}
    \begin{tikzcd} [sep=2cm]
     X\vee X' \arrow[r,"u\vee u'"] & 
     Y\vee Y' \arrow[r,"v\vee v'"] & 
     Z\vee Z' \arrow[r,"w\vee w'"] & 
    X[1]\vee X'[1]
    \end{tikzcd}
 \end{center} with embeddings \begin{center}
    \begin{tikzcd} [sep=2cm]
     X \arrow[r,"u"] \arrow[d,"i_X"] & 
     Y \arrow[r,"v"] \arrow[d,"i_Y"] & 
     Z \arrow[r,"w"] \arrow[d,"i_Z"] & 
     X[1] \arrow[d,"i_X{[1]}"] \\
     X\vee X' \arrow[r,"u\vee u'"] &
     Y\vee Y' \arrow[r,"v\vee v'"] &
     Z\vee Z' \arrow[r,"w\vee w'"] &
     X[1]\vee X'[1] 
    \end{tikzcd}
 \end{center} 
 and
 \begin{center}
    \begin{tikzcd} [sep=2cm]
     X' \arrow[r,"u'"] \arrow[d,"i_{X'}"] & 
     Y' \arrow[r,"v'"] \arrow[d,"i_{Y'}"] & 
     Z' \arrow[r,"w'"] \arrow[d,"i_{Z'}"] & 
     X'[1] \arrow[d,"i_{X'}{[1]}"] \\
     X\vee X' \arrow[r,"u\vee u'"] &
     Y\vee Y' \arrow[r,"v\vee v'"] &
     Z\vee Z' \arrow[r,"w\vee w'"] &
     X[1]\vee X'[1]\;\;\;,
    \end{tikzcd}
 \end{center} and projections \begin{center}
    \begin{tikzcd} [sep=2cm]
     X\vee X' \arrow[r,"u\vee u'"] \arrow[d,"\pi_X"] & 
     Y\vee Y' \arrow[r,"v\vee v'"] \arrow[d,"\pi_Y"] & 
     Z\vee Z' \arrow[r,"w\vee w'"] \arrow[d,"\pi_Z"] & 
     X[1]\vee X'[1] \arrow[d,"\pi_X{[1]}"] \\
     X \arrow[r,"u"] &
     Y \arrow[r,"v"] &
     Z \arrow[r,"w"] &
     X[1] 
    \end{tikzcd}
 \end{center} and \begin{center}
    \begin{tikzcd} [sep=2cm]
     X\vee X' \arrow[r,"u\vee u'"] \arrow[d,"\pi_{X'}"] & 
     Y\vee Y' \arrow[r,"v\vee v'"] \arrow[d,"\pi_{Y'}"] & 
     Z\vee Z' \arrow[r,"w\vee w'"] \arrow[d,"\pi_{Z'}"] & 
     X[1]\vee X'[1] \arrow[d,"\pi_{X'}{[1]}"] \\
     X' \arrow[r,"u'"] &
     Y' \arrow[r,"v'"] &
     Z' \arrow[r,"w'"] &
     X'[1]\;\;\;.
    \end{tikzcd}
 \end{center} We conclude that $\Tr(\mcC)$ is an additive category. To each object $X$ in $\mcC$ we associate a triangle \begin{center}
    \begin{tikzcd} [sep=2cm]
     X \arrow[r,"\id_X"] & 
     X \arrow[r] & 
     0 \arrow[r] & 
     X[1]\;\;\;,
    \end{tikzcd} \end{center} and to each morphism $f:X\to Y$ in $\mcC$ we associate a morphism \begin{center}
    \begin{tikzcd} [sep=2cm]
     X \arrow[r,"\id_X"] \arrow[d,"f"] & 
     X \arrow[r] \arrow[d,"f"] & 
     0 \arrow[r] \arrow[d] & 
     X[1] \arrow[d,"f{[1]}"] \\
     Y \arrow[r,"\id_Y"] &
     Y \arrow[r] &
     0 \arrow[r] &
     Y[1]\;\;\;.
    \end{tikzcd}
 \end{center} This defines a fully faithful additive functor $\mcC\to\Tr(\mcC)$ that is injective on objects. We may thus regard $\mcC$ as a full additive subcategory of $\Tr(\mcC)$. \end{definition}
 
\begin{definition} \cite[Definition 3.1]{HoJ10} \label{definition of triangulated structure} A \textit{triangulation} of $\mcC$ is a full subcategory of $\Tr(\mcC)$ containing $\mcC$ that is closed under isomorphisms, and subject to the following axioms A1-A4 (we shall call the triangles in this subcategory \textit{exact}): \\

(A1) For each morphism $u:X\to Y$ in $\mcC$, there exists an object $Z$ and morphisms $v:Y\to Z$ and $w:Z\to X[1]$ for which \begin{center}
    \begin{tikzcd} [sep=2cm]
     X \arrow[r,"u"] & 
     Y \arrow[r,"v"] & 
     Z \arrow[r,"w"] & 
     X[1]
    \end{tikzcd}
 \end{center} is exact. We call such an exact triangle a \textit{completion} of the morphism $u$. We also call $v$ a \textit{cofibre} of $u$ and $-w[-1]$ a \textit{fibre} of $u$. \\
 
(A2) If the triangles \begin{center}
    \begin{tikzcd} [sep=2cm]
     X \arrow[r,"u"] & 
     Y \arrow[r,"v"] & 
     Z \arrow[r,"w"] & 
     X[1]
    \end{tikzcd}
 \end{center} and \begin{center}
    \begin{tikzcd} [sep=2cm]
     X' \arrow[r,"u'"] & 
     Y' \arrow[r,"v'"] & 
     Z' \arrow[r,"w'"] & 
     X'[1]
    \end{tikzcd}
 \end{center} are both exact, and $f:X\to X'$, $g:Y\to Y'$ satisfy $u'f=gu$, then there exists a morphism $h:Z\to Z'$ for which $(f,g,h)$ is a morphism from the first triangle to the second triangle (i.e. the below diagram commutes). \begin{center}
    \begin{tikzcd} [sep=2cm]
     X \arrow[r,"u"] \arrow[d,"f"] & 
     Y \arrow[r,"v"] \arrow[d,"g"] & 
     Z \arrow[r,"w"] \arrow[d,"h",dashed] & 
     A[1] \arrow[d,"f{[1]}"] \\
     X' \arrow[r,"u'"] &
     Y' \arrow[r,"v'"] &
     Z' \arrow[r,"w'"] &
     X'[1]
    \end{tikzcd}
 \end{center} We call such a morphism $(f,g,h)$ a \textit{completion} of the pair $(f,g)$.  \\
 
(A3) The triangle \begin{center}
    \begin{tikzcd} [sep=2cm]
     X \arrow[r,"u"] & 
     Y \arrow[r,"v"] & 
     Z \arrow[r,"w"] & 
     X[1]
    \end{tikzcd}
 \end{center} is exact if and only if the triangle \begin{center}
    \begin{tikzcd} [sep=2cm]
     Y \arrow[r,"v"] & 
     Z \arrow[r,"w"] & 
     X[1] \arrow[r,"-u{[1]}"] & 
     Y[1]
    \end{tikzcd}
 \end{center} is exact. \\ 
 
(A4) (Octahedral Axiom) Let $u:X\to Y$ and $v:Y\to Z$ be morphisms in $\mcC$. Given completions \begin{center}
    \begin{tikzcd} [sep=2cm]
     X \arrow[r,"u"] & 
     Y \arrow[r] & 
     Z' \arrow[r] & 
     X[1]\;\;\;,
    \end{tikzcd}
 \end{center} \begin{center}
    \begin{tikzcd} [sep=2cm]
     Y \arrow[r,"v"] & 
     Z \arrow[r] & 
     X' \arrow[r] & 
     Y[1]\;\;\;,
    \end{tikzcd}
 \end{center} and \begin{center}
    \begin{tikzcd} [sep=2cm]
     X \arrow[r,"vu"] & 
     Z \arrow[r] & 
     Y' \arrow[r] & 
     X[1]\;\;\;,
    \end{tikzcd}
 \end{center} of $u$, $v$ and $vu$, there exist morphisms making an exact triangle \begin{center}
    \begin{tikzcd} [sep=2cm]
     Z' \arrow[r] & 
     Y' \arrow[r] & 
     X' \arrow[r] & 
     Z'[1]
    \end{tikzcd}
 \end{center} and for which the following diagram commutes. \begin{center}
    \begin{tikzcd} [sep=2cm]
     X \arrow[r,"u"] \arrow[d,equal] & 
     Y \arrow[r] \arrow[d,"v"] & 
     Z' \arrow[r] \arrow[d] & 
     X[1] \arrow[d,equal] \\
     X \arrow[r,"vu"] \arrow[d,"u"] &
     Z \arrow[r] \arrow[d,equal] &
     Y' \arrow[r] \arrow[d] &
     X[1] \arrow[d,"u{[1]}"] \\
     Y \arrow[r,"v"] \arrow[d] &
     Z \arrow[r] \arrow[d] &
     X' \arrow[r] \arrow[d,equal] &
     Y[1] \arrow[d] \\
     Z' \arrow[r] &
     Y' \arrow[r] &
     X' \arrow[r] &
     Z'[1]\;\;\;.
    \end{tikzcd}
 \end{center} \end{definition}
 
\begin{definition} A \textit{triangulated category} is an additive category $\mcC$ equipped with an additive automorphism $[1]$ and a triangulation. We shall henceforth assume that $\mcC$ comes equipped with a fixed triangulation, so that it is a triangulated category. \end{definition}

\begin{prop} \cite[Prop 4.1]{HoJ10} \label{composite of morphisms in exact triangle is zero} Let \begin{center}
    \begin{tikzcd} [sep=2cm]
     X \arrow[r,"u"] & 
     Y \arrow[r,"v"] & 
     Z \arrow[r,"w"] & 
     X[1]
    \end{tikzcd}
 \end{center} be exact. Then $vu=0$, $wv=0$ and $u[1]w=0$. \end{prop}
 
\begin{proof} We know that \begin{center}
    \begin{tikzcd} [sep=2cm]
     X \arrow[r,"\id_X"] & 
     X \arrow[r] & 
     0 \arrow[r] & 
     X[1]
    \end{tikzcd}
 \end{center} is exact. It then follows from A2 that there exists a morphism $0\to Z$ for which \begin{center}
    \begin{tikzcd} [sep=2cm]
     X \arrow[r,"\id_X"] \arrow[d,equal] & 
     X \arrow[r] \arrow[d,"u"] & 
     0 \arrow[r] \arrow[d,dashed] & 
     X[1] \arrow[d,equal] \\
     X \arrow[r,"u"] &
     Y \arrow[r,"v"] &
     Z \arrow[r,"w"] &
     X[1]
    \end{tikzcd}
 \end{center} commutes. Clearly this morphism is just the zero morphism, and we see that $vu=0$. Then \begin{center}
    \begin{tikzcd} [sep=2cm]
     Y \arrow[r,"v"] & 
     Z \arrow[r,"w"] & 
     X[1] \arrow[r,"-u{[1]}"] & 
     Y[1]
    \end{tikzcd}
 \end{center} is exact by A3, so that $wv=0$ and finally \begin{center}
    \begin{tikzcd} [sep=2cm]
     Z \arrow[r,"w"] & 
     X[1] \arrow[r,"-u{[1]}"] & 
     Y[1] \arrow[r,"-v{[1]}"] & 
     Z[1]
    \end{tikzcd}
 \end{center} is exact, so that $u[1]w=0$. \end{proof}
 
\begin{definition} \cite[Remark 1.1.5]{Nee01} $\mcC^{\text{op}}$ is an additive category, $[-1]$ is an additive automorphism of $\mcC^{\text{op}}$ and $\Tr(\mcC^{\text{op}})=\Tr(\mcC)^{\text{op}}$ (where $\mcC^{\text{op}}$ has the suspension and desuspension the opposite way round). Then the exact triangles in $\mcC$ also satisfy the axioms for a triangulated category with respect to $\mcC^{\text{op}}$, so that $\mcC^{\text{op}}$ is also a triangulated category. \end{definition}
 
\begin{definition} \cite[Definition 1.1.3]{HPS97} \cite[Remark 1.1.9]{Nee01} We say that an additive functor $H:\mcC\to\text{Ab}$ is \textit{exact} if whenever \begin{center}
    \begin{tikzcd} [sep=2cm]
     X \arrow[r,"u"] & 
     Y \arrow[r,"v"] & 
     Z \arrow[r,"w"] & 
     X[1]
    \end{tikzcd}
 \end{center} is an exact triangle in $\mcC$, we have that \begin{center}
    \begin{tikzcd} [sep=2cm]
     H(X) \arrow[r,"H(u)"] & 
     H(Y) \arrow[r,"H(v)"] & 
     H(Z)
    \end{tikzcd}
 \end{center} is an exact sequence of abelian groups. A \textit{homology functor} on $\mcC$ is an exact functor $H:\mcC\to\text{Ab}$ where if $\{X_i\}_{i\in I}$ is a family of objects in $\mcC$ for which the coproduct exists, then the map $$\bigoplus_{i\in I}H(X_i)\to H\left(\bigvee_{i\in I}X_i\right)$$ is an isomorphism. If $H$ is contravariant, then we say that $H$ is \textit{exact} if \begin{center}
    \begin{tikzcd} [sep=2cm]
     H(Z) \arrow[r,"H(v)"] & 
     H(Y) \arrow[r,"H(u)"] & 
     H(X)
    \end{tikzcd}
 \end{center} is always exact, or equivalently that $H$ is exact on $\mcC^{\text{op}}$. We also say that $H$ is a \textit{cohomology functor} on $\mcC$ if $$H\left(\bigvee_{i\in I}X_i\right)\to\prod_{i\in I}H(X_i)$$ is always an isomorphism whenever the coproduct for $\{X_i\}_{i\in I}$ exists, or equivalently that $H$ is a homology functor on $\mcC^{\text{op}}$. Using A3, we can for any exact functor $H:\mcC\to\text{Ab}$ extend the exact sequence \begin{center}
    \begin{tikzcd} [sep=2cm]
     H(Z) \arrow[r,"H(v)"] & 
     H(Y) \arrow[r,"H(u)"] & 
     H(X)
    \end{tikzcd}
 \end{center} to see that \begin{center}
    \begin{tikzcd} [sep=1.5cm]
     \cdots \arrow[r,"H(w{[i-1]})"] & H(X{[i]}) \arrow[r,"H(u{[i]})"] & 
     H(Y{[i]}) \arrow[r,"H(v{[i]})"] & 
     H(Z{[i]}) \arrow[r,"H(w{[i]})"] & \cdots
    \end{tikzcd}
 \end{center} is exact, and similarly in the contravariant case. Later we will introduce to $\mcC$ the assumption that the coproduct $\bigvee_{i\in I}X_i$ of any set-indexed family $\{X_i\}_{i\in I}$ always exists, simplifying the above definitions. \end{definition}
 
\begin{prop} \cite[Proposition 4.2]{HoJ10} \cite[Remark 1.1.11]{Nee01} \label{long exact sequence induced by exact triangle covariant} $[T,-]$ is exact for each object $T$, and $[-,T]$ is a cohomology functor on $\mcC$. \end{prop}
 
\begin{proof} Let \begin{center}
    \begin{tikzcd} [sep=2cm]
     X \arrow[r,"u"] & 
     Y \arrow[r,"v"] & 
     Z \arrow[r,"w"] & 
     X[1]
    \end{tikzcd}
 \end{center} be exact. We know from Proposition \ref{composite of morphisms in exact triangle is zero} that $vu=0$. Since the functor $[T,-]$ is additive, it follows that $[T,v][T,u]=0$. It is clear that $\im([T,u])$ is contained in $\ker([T,v])$. Suppose that $\al\in\ker([T,v])$. Then $v\al=0$. Using A3, we see that \begin{center}
    \begin{tikzcd} [sep=2cm]
     T \arrow[r] & 
     0 \arrow[r] & 
     T[1] \arrow[r,"-\id_{T{[1]}}"] & 
     T[1]
    \end{tikzcd}
 \end{center} and \begin{center}
    \begin{tikzcd} [sep=2cm]
     Y \arrow[r,"v"] & 
     Z \arrow[r,"w"] & 
     X[1] \arrow[r,"-u{[1]}"] & 
     Y[1]
    \end{tikzcd}
 \end{center} are both exact. Then since $v\al=0$, the left hand square in \begin{center}
    \begin{tikzcd} [sep=2cm]
     T \arrow[r] \arrow[d,"\al"] & 
     0 \arrow[r] \arrow[d] & 
     T[1] \arrow[r,"-\id_{T{[1]}}"] \arrow[d,"\ta",dashed] & 
     T[1] \arrow[d,"\al{[1]}"] \\
     Y \arrow[r,"v"] &
     Z \arrow[r,"w"] &
     X[1] \arrow[r,"-u{[1]}"] &
     Y[1]
    \end{tikzcd}
 \end{center} commutes, so by A2 there exists a morphism $\ta:T[1]\to X[1]$ with $u[1]\ta=\al[1]$. Applying $[-1]$ then shows that $\al=u\ta[-1]\in\im([T,u])$, so that $\im([T,u])=\ker([T,v])$ and therefore \begin{center}
    \begin{tikzcd} [sep=2cm]
     {[T,X]} \arrow[r,"{[T,u]}"] & 
     {[T,Y]} \arrow[r,"{[T,v]}"] & 
     {[T,Z]}
    \end{tikzcd}
 \end{center} is exact and we conclude that $[T,-]:\mcC\to\text{Ab}$ is exact. Since $\mcC^{\text{op}}$ is also a triangulated category with $T\in\mcC^{\text{op}}$, we see that $[-,T]=\Hom_{\mcC^{\text{op}}}(T,-):\mcC^{\text{op}}\to\text{Ab}$ is exact, so that $[-,T]:\mcC\to\text{Ab}$ is contravariant exact. If $\{X_i\}_{i\in I}$ is a family of objects in $\mcC$ with coproduct $X$, then the map $$[X,T]\to\prod_{i\in I}[X_i,T]$$ sending each map $X\to T$ to the family of composites $X_i\to X\to T$ is an isomorphism, so $[-,T]$ is a cohomology functor on $\mcC$, and we are done. \end{proof}
 
\begin{definition} \cite[Definition 1.1.2]{HPS97} We say that $T\in\mcC$ is \textit{small} if $[T,-]$ is a homology functor on $\mcC$. \end{definition}

\begin{exmp} Suppose that all infinite coproducts in $\mcC$ exist, and let $S\in\mcC$ be small with $[S,S]=\Z$. Let $X=\bigvee_{i=0}^\infty S$ and suppose that $X$ is small. Then since $[-,S]$ and $[-,X]$ are cohomological, and $[X,-]$ and $[S,-]$ are homological, we see that $$\bigoplus_i\left(\prod_j\Z\right)\simeq\bigoplus_i[X,S]\simeq[X,X]\simeq\prod_i[S,X]\simeq\prod_i\left(\bigoplus_j\Z\right),$$ with the composite being the inclusion. The right hand side can be identified with the abelian group of infinite matrices with integer entries, where each row has only finitely many non-zero entries; the left hand side can be identified with the subgroup of those where only finitely many columns are non-zero. Since we have shown that the inclusion is surjective, the identity matrix must have only finitely many non-zero columns, a contradiction. So $X$ cannot be small. \end{exmp}

\begin{definition} \cite[Definition 1.1.3(c)]{HPS97} We say that a cohomology functor on $\mcC$ is \textit{representable} if it is naturally isomorphic to $[-,T]$ for some $T\in\mcC$. We will later add in the assumption that all cohomology functors on $\mcC$ are representable. \end{definition}

\begin{prop} \cite[Proposition 4.3]{HoJ10} \label{triangulated five lemma} (Triangulated Five Lemma) Let \begin{center}
    \begin{tikzcd} [sep=2cm]
     X \arrow[r,"u"] \arrow[d,"f"] & 
     Y \arrow[r,"v"] \arrow[d,"g"] & 
     Z \arrow[r,"w"] \arrow[d,"h"] & 
     X[1] \arrow[d,"f{[1]}"] \\
     X' \arrow[r,"u'"] &
     Y' \arrow[r,"v'"] &
     Z' \arrow[r,"w'"] &
     X'[1]
    \end{tikzcd}
 \end{center} be a morphism of exact triangles. If any two of $f,g,h$ are isomorphisms, then all of them are, so that $(f,g,h)$ is an isomorphism in $\Tr(\mcC)$. \end{prop}
 
\begin{proof} First assume that $f$ and $g$ are isomorphisms. Applying Proposition \ref{long exact sequence induced by exact triangle covariant} to both exact triangles, we see that \begin{center}
    \begin{tikzcd} [sep=1.5cm]
     {[Z',X]} \arrow[r,"{[Z',u]}"] & 
     {[Z',Y]} \arrow[r,"{[Z',v]}"] & 
     {[Z',Z]} \arrow[r,"{[Z',w]}"] & 
     {[Z',X[1]]} \arrow[r,"{[Z',u{[1]}]}"] &
     {[Z',Y[1]]}
    \end{tikzcd}
 \end{center} and \begin{center}
    \begin{tikzcd} [sep=1.5cm]
     {[Z',X']} \arrow[r,"{[Z',u']}"] & 
     {[Z',Y']} \arrow[r,"{[Z',v']}"] & 
     {[Z',Z']} \arrow[r,"{[Z',w']}"] & 
     {[Z',X'[1]]} \arrow[r,"{[Z',u'{[1]}]}"] &
     {[Z',Y'[1]]}
    \end{tikzcd}
 \end{center} are both exact. Since $f$ and $g$ are isomorphisms and functors preserve isomorphisms, we see that $[Z',f]$, $[Z',g]$, $[Z',f[1]]$ and $[Z',g[1]]$ are isomorphisms. Since $(f,g,h)$ is a morphism of triangles, we see that the diagram \begin{center}
    \begin{tikzcd} [sep=1.4cm]
     {[Z',X]} \arrow[r,"{[Z',u]}"] \arrow[d,"{[Z',f]}"] & 
     {[Z',Y]} \arrow[r,"{[Z',v]}"] \arrow[d,"{[Z',g]}"] & 
     {[Z',Z]} \arrow[r,"{[Z',w]}"] \arrow[d,"{[Z',h]}"] & 
     {[Z',X[1]]} \arrow[r,"{[Z',u{[1]}]}"] \arrow[d,"{[Z',f{[1]}]}"] &
     {[Z',Y[1]]} \arrow[d,"{[Z',g{[1]}]}"] \\
     {[Z',X']} \arrow[r,"{[Z',u']}"] & 
     {[Z',Y']} \arrow[r,"{[Z',v']}"] & 
     {[Z',Z']} \arrow[r,"{[Z',w']}"] & 
     {[Z',X'[1]]} \arrow[r,"{[Z',u'{[1]}]}"] &
     {[Z',Y'[1]]}
    \end{tikzcd}
 \end{center} commutes and it follows from the Five Lemma for abelian groups that $[Z',h]$ is an isomorphism. Then set $q=[Z',h]^{-1}(\id_{Z'}):Z'\to Z$. Then $hq=\id_{Z'}$ so that $h$ is a split epimorphism. Using A3, we see from this that if any two of $f,g,h$ are isomorphisms, then the remaining one is a split epimorphism. By applying this result to $\mcC^{\text{op}}$, we see that the remaining one is also a split monomorphism, and therefore an isomorphism. \end{proof}
 
\begin{prop} \cite[Lemma 13.4.10]{StP} \label{direct sum of exact trgs is exact} Let \begin{center}
    \begin{tikzcd} [sep=2cm]
     X_i \arrow[r,"u_i"] & 
     Y_i \arrow[r,"v_i"] & 
     Z_i \arrow[r,"w_i"] & 
     X_i[1]
    \end{tikzcd}
 \end{center} be triangles for $i=0,1$. These are both exact if and only if \begin{center}
    \begin{tikzcd} [sep=2cm]
     X_0\vee X_1 \arrow[r,"u_0\vee u_1"] & 
     Y_0\vee Y_1 \arrow[r,"v_0\vee v_1"] & 
     Z_0\vee Z_1 \arrow[r,"w_0\vee w_1"] & 
     (X_0\vee X_1)[1]
    \end{tikzcd}
 \end{center} is exact. \end{prop}
 
\begin{proof} See \cite[Lemma 13.4.10]{StP} for a proof. \end{proof}
 
\begin{prop} \cite[Proposition 4.4]{HoJ10} \cite[Lemma 12.4]{ShS14} \label{zeros in exact tr split monos epis} Let \begin{center}
    \begin{tikzcd} [sep=2cm]
     X \arrow[r,"u"] & 
     Y \arrow[r,"v"] & 
     Z \arrow[r,"w"] & 
     X[1]
    \end{tikzcd}
 \end{center} be exact. Then $u$ is a split monomorphism if and only if $v$ is a split epimorphism if and only if $w=0$. Under these equivalent conditions there exists an isomorphism $\ta:Y\to X\vee Z$ for which the following diagram commutes, and therefore is an isomorphism of triangles. \begin{center}
    \begin{tikzcd} [sep=2cm] \label{split exact triangle diagram}
     X \arrow[r,"u"] \arrow[d,equal] & 
     Y \arrow[r,"v"] \arrow[d,"\ta",dashed] & 
     Z \arrow[r,"0"] \arrow[d,equal] & 
     X[1] \arrow[d,equal] \\
     X \arrow[r] &
     X\vee Z \arrow[r] &
     Z \arrow[r,"0"] &
     X[1] 
    \end{tikzcd}
 \end{center} Furthermore, every monomorphism and every epimorphism in $\mcC$ splits. \end{prop}
 
\begin{proof} If $u$ is a split monomorphism, then since \begin{center}
    \begin{tikzcd} [sep=2cm]
     Z[-1] \arrow[r,"-w{[-1]}"] & 
     X \arrow[r,"u"] & 
     Y \arrow[r,"v"] & 
     Z
    \end{tikzcd}
 \end{center} is exact, we see that $-w[-1]=0$ and therefore $w=0$. Conversely, if $w=0$, then since \begin{center}
    \begin{tikzcd} [sep=2cm]
     0 \arrow[r] & 
     X \arrow[r,"\id_X"] & 
     X \arrow[r] & 
     0
    \end{tikzcd}
 \end{center} is exact, we use A2 to see that the following diagram can be completed to a morphism of triangles. \begin{center}
    \begin{tikzcd} [sep=2cm]
     Z[-1] \arrow[r,"0"] \arrow[d] & 
     X \arrow[r,"u"] \arrow[d,equal] & 
     Y \arrow[r,"v"] \arrow[d,"\ta",dashed] & 
     Z \arrow[d] \\
     0 \arrow[r] &
     X \arrow[r,"\id_X"] &
     X \arrow[r] &
     0
    \end{tikzcd}
 \end{center} So $\ta u=\id_X$ and therefore $u$ is a split monomorphism. A similar argument shows that $v$ is a split epimorphism if and only if $w=0$. 
 
If the stated equivalent conditions hold, then we combine $v$ with a left inverse for $u$ to obtain a map $\ta:Y\to X\vee Z$ for which \eqref{split exact triangle diagram} commutes, and use Proposition \ref{direct sum of exact trgs is exact} to see that the bottom row of this diagram is exact. Then we use the Triangulated Five Lemma to see that $\ta$ is an isomorphism, so that the diagram is an isomorphism of exact triangles. 

Finally, if $u$ is a monomorphism, then we can complete it to an exact triangle \begin{center}
    \begin{tikzcd} [sep=2cm]
     X \arrow[r,"u"] & 
     Y \arrow[r,"v"] & 
     Z \arrow[r,"w"] & 
     X[1] \;\;\;.
    \end{tikzcd}
 \end{center} Since $u$ is a monomorphism, we use the same argument as for the split monomorphism case above to see that $w=0$ and therefore $u$ is a split monomorphism. A similar argument shows that every epimorphism splits. \end{proof}

\begin{prop} \cite{San11} Let $u:X\to Y$ be a morphism in $\mcC$, and let \begin{center}
    \begin{tikzcd} [sep=2cm]
     X \arrow[r,"u"] & 
     Y \arrow[r,"v"] & 
     Z \arrow[r,"w"] & 
     X[1]
    \end{tikzcd}
 \end{center} and \begin{center}
    \begin{tikzcd} [sep=2cm]
     X \arrow[r,"u"] & 
     Y \arrow[r,"v'"] & 
     Z' \arrow[r,"w'"] & 
     X[1]
    \end{tikzcd}
 \end{center} be completions of $u$ to exact triangles. Then there exists an isomorphism $\ta:Z\to Z'$ with $\ta v=v'$ and $w'\ta=w$, so that the fibre, cofibre and completion are unique up to (non-unique) isomorphism. \begin{center}
    \begin{tikzcd} [sep=2cm]
     X \arrow[r,"u"] \arrow[d,equal] & 
     Y \arrow[r,"v"] \arrow[d,equal] & 
     Z \arrow[r,"w"] \arrow[d,"\ta",dashed] & 
     X[1] \arrow[d,equal] \\
     X \arrow[r,"u"] &
     Y \arrow[r,"v'"] &
     Z' \arrow[r,"w'"] &
     X[1]
    \end{tikzcd}
 \end{center} \end{prop}
 
\begin{proof} It follows easily from A2 that there exists a morphism $\ta:Z\to Z'$ as indicated in the above diagram. It is then an immediate consequence of the Triangulated Five Lemma (Proposition \ref{triangulated five lemma}) that $\ta$ is in fact an isomorphism. \end{proof}
 
\begin{definition} A \textit{weak kernel} of a morphism $f:X\to Y$ in an additive category consists of an object $K$ and a morphism $\phi:K\to X$ for which $f\phi=0$ and possessing the following property:

For each object $K'$ and morphism $\psi:K'\to X$ with $f\psi=0$, there exists a morphism $\al:K'\to K$ with $\phi\al=\psi$.

\begin{center}
    \begin{tikzcd} [sep=2cm]
     K' \arrow[d,dashed,"\al"] \arrow[rd,"\psi"] & {} & {} \\
     K \arrow[r,"\phi"] & X \arrow[r,"f"] & Y \\ 
     \end{tikzcd}\end{center} 

Similarly, a \textit{weak cokernel} of $f$ consists of an object $C$ and a morphism $\phi:Y\to C$ for which $\phi f=0$ and possessing the following property:

For each object $C'$ and morphism $\psi:Y\to C'$ with $\psi f=0$, there exists a morphism $\al:C\to C'$ with $\al\phi=\psi$.

\begin{center}
    \begin{tikzcd} [sep=2cm]
     {} & {} & C' \\
     X \arrow[r,"f"] & Y \arrow[r,"\phi"] \arrow[ru,"\psi"] & C \arrow[u,dashed,"\al",swap] \\ 
     \end{tikzcd}\end{center} \end{definition}

\begin{definition} We say that an additive category is \textit{weakly abelian} if every morphism $f:X\to Y$ has a weak kernel and a weak cokernel, and there exists a morphism $\vphi:K\to X$ where $f$ is a weak cokernel of $\vphi$, and also there exists a morphism $\psi:Y\to C$ where $f$ is a weak kernel of $\psi$. \end{definition}

We now establish that $\mcC$ (and, by the same argument, all triangulated categories) is weakly abelian. We will use this later on to show that the full subcategory $\mcF$ of small objects in $\mcC$ can be extended to an abelian category (its Freyd completion, $\Fr(\mcF)$) with some interesting properties. 

\begin{lem} \cite{San11} \label{weak kernels and cokernels in exact triangle} Fibres are weak kernels and cofibres are weak cokernels. \end{lem}
 
\begin{proof} Let $v:Y\to Z$ have a fibre $u:X\to Y$ and a cofibre $w:Z\to X[1]$. Then \begin{center}
    \begin{tikzcd} [sep=2cm]
     X \arrow[r,"u"] & 
     Y \arrow[r,"v"] & 
     Z \arrow[r,"w"] & 
     X[1]
    \end{tikzcd}
 \end{center} is exact, so that $vu=0$. If $X'$ is an object and $u':X'\to Y$ with $vu'=0$, then since $[X',-]$ is exact, we have that \begin{center}
    \begin{tikzcd} [sep=2cm]
     {[X',X]} \arrow[r,"{[X',u]}"] & 
     {[X',Y]} \arrow[r,"{[X',v]}"] & 
     {[X',Z]} 
    \end{tikzcd}
 \end{center} is exact. So $u'\in\ker([X',v])=\im([X',u])$ and then $u'=u\ta$ for some $\ta:X'\to X$ so that $u$ is a weak kernel for $v$. The dual argument shows that $w$ is a weak cokernel of $v$. \begin{center}
    \begin{tikzcd} [sep=2cm]
     X' \arrow[d,dashed,"\ta"] \arrow[rd,"u'"] & {} & {} \\
     X \arrow[r,"u"] & Y \arrow[r,"v"] & Z \\ 
     \end{tikzcd}\end{center} \end{proof}

\begin{theorem} \cite{San11} $\mcC$ is weakly abelian. \end{theorem}

\begin{proof} Let $v:Y\to Z$ be a morphism in $\mcC$. Since $v$ has a fibre $u$ and a cofibre $w$, and these are a weak kernel and a weak cokernel for $v$ respectively by Lemma \ref{weak kernels and cokernels in exact triangle}, we know that $v$ has a weak kernel and a weak cokernel. Then $v$ is a cofibre and therefore a weak cokernel of $u$, and $v$ is a fibre and therefore a weak kernel of $w$. \end{proof}

\section{The Freyd Completion} \label{freyd comp}

So far as I am aware, the Freyd completion of an additive category was first described in \cite[\S 3]{Fre65}, upon which this section is largely based. We first define thick and localising subcategories of $\mcC$, and establish that the full subcategory of small objects in $\mcC$ is thick. We then construct the Freyd completion of the full subcategory of small objects in $\mcC$, and establish various properties of this new category. 

\begin{definition} \cite[Definition 1.4.3(a)]{HPS97} \label{thick localising subcat def} A full subcategory $\mcD$ of $\mcC$ is said to be \textit{thick} if

1. Whenever \begin{center}
    \begin{tikzcd} [sep=2cm]
     X \arrow[r,"u"] & 
     Y \arrow[r,"v"] & 
     Z \arrow[r,"w"] & 
     X[1]
    \end{tikzcd}
 \end{center} is an exact triangle in $\mcC$ and any two of $X$, $Y$, $Z$ are contained in $\mcD$, then they are all contained in $\mcD$. 

2. If $Y\in\mcD$ and $X$ is a retract of $Y$, then $X\in\mcD$ (i.e. $\mcD$ is closed under retracts). 

We say that $\mcD$ is \textit{localising} if it is thick, and also closed under arbitrary coproducts i.e. if $\{X_i\}_{i\in I}$ is a family of objects in $\mcD$, then $\bigvee_{i\in I}X_i$ is also in $\mcD$. 

Let $\mcS$ be a set of objects in $\mcC$. The intersection of all thick (localising) subcategories of $\mcC$ containing $\mcS$ is itself a thick (localising) subcategory containing $\mcS$, and is the smallest one. We call this the thick (localising) subcategory of $\mcC$ \textit{generated} by $\mcS$, and write it as $\text{thick}(\mcS)$ ($\loc(\mcS)$). Thick subcategories of $\mcC$ inherit a triangulated structure from $\mcC$. \end{definition}

\begin{definition} Denote by $\mcF$ the full subcategory of $\mcC$ consisting of the small objects. \end{definition}

\begin{prop} \cite[Lemma 13.37.2]{StP} $\mcF$ is a thick subcategory of $\mcC$. \end{prop}

\begin{proof} Let $X\in\mcF$, $n\in\Z$, and $\{A_i\}_{i\in I}$ be a family of objects in $\mcC$ with wedge sum $A$. Then since $X\in\mcF$, we see that $$\bigoplus_i[X[n],A_i]\simeq\bigoplus_i[X,A_i[-n]]\simeq[X,\bigvee_i(A_i[-n])]\simeq[X,A[-n]]\simeq[X[n],A],$$ so that $X[n]\in\mcF$. Now let \begin{center}
    \begin{tikzcd} [sep=2cm]
     X \arrow[r,"u"] & 
     Y \arrow[r,"v"] & 
     Z \arrow[r,"w"] & 
     X[1]
    \end{tikzcd}
 \end{center} be exact, and $\{A_i\}_{i\in I}$ be a family of objects with wedge sum $A$. We then use Proposition \ref{long exact sequence induced by exact triangle covariant} to see that the following diagram commutes, with exact rows. \begin{center}
    \begin{tikzcd} [sep=1.2cm]
     \cdots \arrow[r] &
     \bigoplus_i[X{[1]},A_i] \arrow[r] \arrow[d] & 
     \bigoplus_i[Z,A_i] \arrow[r] \arrow[d] & 
     \bigoplus_i[Y,A_i] \arrow[r] \arrow[d] & 
     \bigoplus_i[X,A_i] \arrow[d] \arrow[r] &
     \cdots \\
     \cdots \arrow[r] &
     {[X{[1]},A]} \arrow[r] &
     {[Z,A]} \arrow[r] &
     {[Y,A]} \arrow[r] &
     {[X,A]} \arrow[r] &
     \cdots
    \end{tikzcd}
 \end{center} If $X,Y\in\mcF$, then also $X[1],Y[1]\in\mcF$, so we can use the Five Lemma to see that the map $\bigoplus_i[Z,A_i]\to[Z,A]$ is an isomorphism, and therefore $Z\in\mcF$. A similar argument shows that $\mcF$ is closed under retracts, so that $\mcF$ is thick. \end{proof}

\begin{definition} \cite[Section II.4]{Mac78} We define the \textit{arrow category} of $\mcF$ as follows. The objects in $\Arr(\mcF)$ are precisely the morphisms in $\mcF$. For objects $u:X\to Y$ and $v:Z\to W$, we define the abelian group $$\Hom_{\Arr(\mcF)}(u,v)=\{(f,g)\in[X,Z]\x[Y,W]:vf=gu\}.$$ \begin{center}
    \begin{tikzcd}[sep=2cm]
     X \arrow[r,"u"] \arrow[d,"f",swap] &
     Y \arrow[d,"g"] \\
     Z \arrow[r,"v"] & 
     W
    \end{tikzcd}
 \end{center} Composition is given by $$(f',g')(f,g)=(f'f,g'g).$$ For each object $u:X\to Y$, we set $\id_u=(\id_X,\id_Y)\in\Hom_{\Arr(\mcF)}(u,u)$. \begin{center}
    \begin{tikzcd}[sep=2cm]
     X \arrow[r,"u"] \arrow[d,equal,swap] &
     Y \arrow[d,equal] \\
     X \arrow[r,"u"] & 
     Y
    \end{tikzcd}
 \end{center}If $0$ is a zero object in $\mcF$, then the unique morphism $0\to 0$ is a zero object in $\Arr(\mcF)$. 

Let $u:X\to Y$ and $v:Z\to W$ be objects in $\Arr(\mcF)$. We form an object $u\vee v:X\vee Z\to Y\vee W$. We then define morphisms $i_u=(i_X,i_Y):u\to u\vee v$, $i_v=(i_Z,i_W):v\to u\vee v$, $\pi_u=(\pi_X,\pi_Y):u\vee v\to u$ and finally $\pi_v=(\pi_Z,\pi_W):u\vee v\to v$. It is routine to check that $(u\vee v,i_u,i_v,\pi_u,\pi_v)$ is a biproduct of $u$ and $v$; we conclude that $\Arr(\mcF)$ is an additive category. Note that the functor $$(f:X\to Y)\mapsto((f,f):\id_X\to\id_Y):\mcF\to\Arr(\mcF)$$ is fully faithful. \end{definition}

\begin{definition} \cite[\S 3]{Fre65} We now construct $\Fr(\mcF)$ as a quotient category of $\Arr(\mcF)$. For objects $u:X\to Y$ and $v:Z\to W$ in $\Arr(\mcF)$, set $$\Hom_0(u,v)=\{(f,g)\in[X,Z]\x[Y,W]:vf=gu=0\}.$$ We write $f\sim g$ to indicate that $f-g\in\Hom_0(u,v)$. This defines an equivalence relation on $\Hom(u,v)$ since $\Hom_0(u,v)$ is a subgroup. Since $\sim$ respects composition and additivity, we can define the \textit{Freyd completion} of $\mcF$ to be the quotient category $\Fr(\mcF)=\Arr(\mcF)/\sim$. By first applying the functor $$f\mapsto(f,f):\mcF\to\Arr(\mcF),$$ (sending each object to its identity map) and then the quotient functor $\Arr(\mcF)\to\Fr(\mcF)$, we obtain a fully faithful functor $\mcF\to\Fr(\mcF)$ and thus regard $\mcF$ as a full subcategory of $\Fr(\mcF)$ with this functor being the inclusion functor. \end{definition}

\begin{definition} We extend the suspension functor to $\Fr(\mcF)$ as follows. It sends each object $u:X\to Y$ in $\Fr(\mcF)$ to the object $u[1]:X[1]\to Y[1]$ in $\Fr(\mcF)$. Represent a morphism $\chi:u\to v$ in $\Fr(\mcF)$ by the following commuting square. \begin{center}
    \begin{tikzcd}[sep=2cm]
     X \arrow[r,"u"] \arrow[d,"f",swap] &
     Y \arrow[d,"g"] \\
     Z \arrow[r,"v"] & 
     W
    \end{tikzcd}
 \end{center} Then \begin{center}
    \begin{tikzcd}[sep=2cm]
     X[1] \arrow[r,"u{[1]}"] \arrow[d,"f{[1]}",swap] &
     Y[1] \arrow[d,"g{[1]}"] \\
     Z[1] \arrow[r,"v{[1]}"] & 
     W[1]
    \end{tikzcd}
 \end{center} represents a morphism, and one can check that this process respects equivalence. We label this morphism $\chi[1]:u[1]\to v[1]$. One can easily check that this defines an additive automorphism on $\Fr(\mcF)$, and is exact. \end{definition}

\begin{lem} \cite[Lemma 3.1.2]{Fre65} \label{monos and epis preservation in freyd completion} Consider the following commuting square in $\mcF$. \begin{center}
    \begin{tikzcd}[sep=2cm]
     X \arrow[r,"u"] \arrow[d,"f",swap] &
     Y \arrow[d,"g"] \\
     Z \arrow[r,"v"] & 
     W
    \end{tikzcd}
 \end{center} If $f$ is an epimorphism in $\mcF$, then $[(f,g)]:u\to v$ is an epimorphism in $\Fr(\mcF)$. If $g$ is a monomorphism in $\mcF$, then $[(f,g)]:u\to v$ is a monomorphism in $\Fr(\mcF)$. \end{lem}

\begin{proof} Let $\vphi:A\to B$ be an object in $\Fr(\mcF)$ and $\al:v\to\vphi$ be a morphism with $\al[(f,g)]=0$. Then $\al=[(a,b)]$ for some $a:Z\to A$ and $b:W\to B$ with $bv=\vphi a$. \begin{center}
    \begin{tikzcd}[sep=2cm]
     X \arrow[r,"u"] \arrow[d,"f",swap] &
     Y \arrow[d,"g"] \\
     Z \arrow[r,"v"] \arrow[d,"a",swap] & W \arrow[d,"b"] \\ 
     A \arrow[r,"\vphi"] & B
    \end{tikzcd}
 \end{center} This shows that $\vphi af=0$. If $f$ is an epimorphism, then $\vphi a=0$ so that $\al=0$. It follows that $[(f,g)]$ is an epimorphism in $\Fr(\mcF)$. A similar proof shows that if $g$ is a monomorphism in $\mcF$, then $[(f,g)]$ is a monomorphism in $\Fr(\mcF)$. \end{proof}

\begin{lem} \cite[Lemma 3.1.4]{Fre65} \label{forming kernels and cokernels in freyd completion} Consider the following commuting square in $\mcF$. \begin{center}
    \begin{tikzcd}[sep=2cm]
     X \arrow[r,"u"] \arrow[d,"f",swap] &
     Y \arrow[d,"g"] \\
     Z \arrow[r,"v"] & 
     W
    \end{tikzcd}
 \end{center} If $Y=W$ and $g=\id_Y$, then $[(f,g)]:u\to v$ is a kernel in $\Fr(\mcF)$. Similarly, if $X=Z$ and $f=\id_X$, then $[(f,g)]$ is a cokernel in $\Fr(\mcF)$. \end{lem}

\begin{proof} Let $\chi=[(f,g)]$. First suppose that $Y=W$ and $g=\id_Y$. Let $\pi:Y\to P$ be a cofibre for $u$, so that $\pi u=0$. Then $[(\id_Z,\pi)]:v\to\pi v$ is a morphism in $\Fr(\mcF)$ with $[(\id_Z,\pi)]\chi=0$. \begin{center}
    \begin{tikzcd}[sep=2cm]
     X \arrow[r,"u"] \arrow[d,"f",swap] &
     Y \arrow[d,equal] \\
     Z \arrow[r,"v"] \arrow[d,equal,swap] & Y \arrow[d,"\pi"] \\
     Z \arrow[r,"\pi v"] & P 
    \end{tikzcd}
 \end{center} Let $u':X'\to Y'$ be an object in $\Fr(\mcF)$ and $\al:u'\to v$ be a morphism with $[(\id_Z,\pi)]\al=0$. Then $\al=[(x,y)]$ for some $x:X'\to Z$, $y:Y'\to Y$ with $yu'=vx$. \begin{center}
    \begin{tikzcd}[sep=2cm]
     X' \arrow[r,"u'"] \arrow[d,"x",swap] &
     Y' \arrow[d,"y"] \\
     Z \arrow[r,"v"] \arrow[d,equal,swap] & Y \arrow[d,"\pi"] \\
     Z \arrow[r,"\pi v"] & P 
    \end{tikzcd}
 \end{center} Since $u$ is a fibre and therefore a weak kernel of $\pi$, there exists a morphism $\ta:X'\to X$ with $u\ta=yu'=vx$. \begin{center}
    \begin{tikzcd} [sep=2cm]
     X' \arrow[d,dashed,"\ta",swap] \arrow[rd,"vx=yu'"] & {} & {} \\
     X \arrow[r,"u"] & Y \arrow[r,"\pi"] & P \\ 
     \end{tikzcd}\end{center} This gives a morphism $[(\ta,y)]:u'\to u$ with $\chi[(\ta,y)]=\al$. \begin{center}
    \begin{tikzcd}[sep=2cm]
     X' \arrow[r,"u'"] \arrow[d,"\ta",swap] &
     Y' \arrow[d,"y"] \\
     X \arrow[r,"u"] \arrow[d,"f",swap] & Y \arrow[d,equal] \\
     Z \arrow[r,"v"] & Y 
    \end{tikzcd}
 \end{center} If $\rho:u'\to u$ is a morphism with $\chi\rho=\al$, then $\rho=[(\ta,y)]$ by Lemma \ref{monos and epis preservation in freyd completion} since $\chi$ is a monomorphism. We conclude that $\chi:u\to v$ is a kernel of $[(\id_Z,\pi)]:v\to\pi v$. 
 
\begin{center}
    \begin{tikzcd} [sep=2cm]
     u' \arrow[d,dashed,"{[(\ta,y)]}",swap] \arrow[rd,"\al"] & {} & {} \\
     u \arrow[r,"\chi"] & v \arrow[r,"{[(\id_Z,\pi)]}"] & \pi v \\ 
     \end{tikzcd}\end{center} 

Similarly, if $X=Z$ and $f=\id_X$, then let $\sg:P\to X$ be a fibre of $v$. Then $[(\sg,u)]:\sg\to u$ is a morphism in $\Fr(\mcF)$ for which $\chi:u\to v$ is a cokernel. \end{proof}

\begin{lem} \cite[Lemma 3.1.3]{Fre65} \label{weak kernels implies kernels in freyd comp} Let a morphism $\chi$ in $\Fr(\mcF)$ be represented by the following commuting square. \begin{center}
    \begin{tikzcd}[sep=2cm]
     X \arrow[r,"u"] \arrow[d,"f",swap] &
     Y \arrow[d,"g"] \\
     Z \arrow[r,"v"] & 
     W
    \end{tikzcd}
 \end{center} Let $k:K\to X$ and $\pi:W\to P$ be a fibre and a cofibre of $vf=gu:X\to W$, respectively. Then $[(k,\id_Y)]$ and $[(\id_Z,\pi)]$ are a kernel and a cokernel of $\chi$, respectively. Therefore every morphism in $\Fr(\mcF)$ has both a kernel and a cokernel. \end{lem}

\begin{proof} We have a morphism $[(k,\id_Y)]:uk\to u$ with $\chi[(k,\id_Y)]=0$. \begin{center}
    \begin{tikzcd}[sep=2cm]
     K \arrow[r,"uk"] \arrow[d,"k",swap] & Y \arrow[d,equal] \\
     X \arrow[r,"u"] \arrow[d,"f",swap] & Y \arrow[d,"g"] \\
     Z \arrow[r,"v"] & W
    \end{tikzcd}
 \end{center} Let $\vphi:A\to B$ be an object in $\Fr(\mcF)$ and $\al:\vphi\to u$ be a morphism with $\chi\al=0$. Then $\al=[(x,y)]$ for some $x:A\to X$ and $y:B\to Y$ with $y\vphi=ux$ and we see that $vfx=gy\vphi=0$. \begin{center}
    \begin{tikzcd}[sep=2cm]
     A \arrow[r,"\vphi"] \arrow[d,"x",swap] & B \arrow[d,"y"] \\
     X \arrow[r,"u"] \arrow[d,"f",swap] & Y \arrow[d,"g"] \\
     Z \arrow[r,"v"] & W
    \end{tikzcd}
 \end{center} Since $k$ is a weak kernel of $vf$, there must exist a morphism $\ta:A\to K$ with $k\ta=x$. \begin{center}
    \begin{tikzcd} [sep=2cm]
     A \arrow[d,dashed,"\ta",swap] \arrow[rd,"x"] & {} & {} \\
     K \arrow[r,"k"] & X \arrow[r,"vf"] & W \\ 
     \end{tikzcd}\end{center} This gives a morphism $[(\ta,y)]:\vphi\to uk$ with $[(k,\id_Y)][(\ta,y)]=\al$. \begin{center}
    \begin{tikzcd}[sep=2cm]
     A \arrow[r,"\vphi"] \arrow[d,"\ta",swap] & B \arrow[d,"y"] \\
     K \arrow[r,"uk"] \arrow[d,"k",swap] & Y \arrow[d,equal] \\
     X \arrow[r,"u"] & Y
    \end{tikzcd}
 \end{center} 

Since $\id_Y$ is a monomorphism in $\mcF$, we know from Lemma \ref{monos and epis preservation in freyd completion} that $[(k,\id_Y)]$ is a monomorphism in $\Fr(\mcF)$. It follows that if $\rho:\vphi\to uk$ with $[(k,\id_Y)]\rho=\al$, then $\rho=[(\ta,y)]$, so that $[(k,\id_Y)]$ is a kernel of $\chi$. Using a similar argument, we see that $[(\id_Z,\pi)]:v\to\pi v$ is a cokernel of $\chi$. \begin{center}
    \begin{tikzcd} [sep=2cm]
     \vphi \arrow[d,dashed,"{[(\ta,y)]}",swap] \arrow[rd,"\al"] & {} & {} \\
     uk \arrow[r,"{[(k,\id_Y)]}"] & u \arrow[r,"\chi"] & v \\ 
     \end{tikzcd}\end{center}\end{proof}

\begin{lem} \cite[Theorem 3.1]{Fre65} \label{morphism in freyd comp splits kernel cokernel} Every morphism in $\Fr(\mcF)$ is the composite of a cokernel followed by a kernel. \end{lem}

\begin{proof} Let $\chi$ be a morphism in $\Fr(\mcF)$ represented by the commuting square \begin{center}
    \begin{tikzcd}[sep=2cm]
     X \arrow[r,"u"] \arrow[d,"f",swap] &
     Y \arrow[d,"g"] \\
     Z \arrow[r,"v"] & 
     W
    \end{tikzcd}
 \end{center} Then $[(\id_X,g)]:u\to vf$ and $[(f,\id_W)]:vf\to v$ are morphisms with $[(f,\id_W)][(\id_X,g)]=\chi$. \begin{center}
    \begin{tikzcd}[sep=2cm]
     X \arrow[r,"u"] \arrow[d,equal,swap] & Y \arrow[d,"g"] \\
     X \arrow[r,"vf"] \arrow[d,"f",swap] & W \arrow[d,equal] \\
     Z \arrow[r,"v"] & W
    \end{tikzcd}
 \end{center} It follows from Lemma \ref{forming kernels and cokernels in freyd completion} that $[(\id_X,g)]$ is a cokernel and $[(f,\id_W)]$ is a kernel. \end{proof}

\begin{lem} \label{mono and cokernel implies iso} Let $f:X\to Y$ be a monomorphism (epimorphism) in an additive category $\mcD$ that is also a cokernel (kernel). Then $f$ is an isomorphism. \end{lem}

\begin{proof} Suppose that $f:X\to Y$ is a monomorphism that is also a cokernel. Then $f$ is a cokernel for some map $k:K\to X$. So $fk=0=f0$ with $f$ a monomorphism, giving that $k=0$. Then the universal property of the cokernel shows that $f$ must have a left inverse and therefore be a split monomorphism. Since $f$ is a cokernel, it must also be an epimorphism, and therefore an isomorphism. A similar argument shows that an epimorphism that is also a kernel must be an isomorphism. \end{proof}

\begin{lem} \cite[Theorem 3.1]{Fre65} \label{monos in freyd comp are kernels} Each monomorphism in $\Fr(\mcF)$ is a kernel and each epimorphism in $\Fr(\mcF)$ is a cokernel. \end{lem}

\begin{proof} Let $\chi:u\to v$ be a morphism in $\Fr(\mcF)$. Then there exists an object $\vphi$ and morphisms $\al:u\to\vphi$, $\bt:\vphi\to v$ for which $\chi=\bt\al$, where $\al$ is a cokernel and $\bt$ is a kernel by Lemma \ref{morphism in freyd comp splits kernel cokernel}. \begin{center}
    \begin{tikzcd} [sep=2cm]
     {} & \vphi \arrow[rd,"\bt"] & {} \\
     u \arrow[rr,"\chi"] \arrow[ru,"\al"] & {} & v \\ 
     \end{tikzcd}\end{center} If $\chi$ is a monomorphism, then $\al$ is a monomorphism and a cokernel, and therefore an isomorphism by Lemma \ref{mono and cokernel implies iso}. Since $\al$ is an isomorphism and $\bt$ is a kernel, we know that $\chi$ is a kernel. Similarly, if $\chi$ is an epimorphism, then $\bt$ is an epimorphism and a kernel, and therefore an isomorphism. Since $\al$ is a cokernel and $\bt$ is an isomorphism, we know that $\chi$ is a cokernel. \end{proof}

\begin{theorem} \cite[Theorem 3.1]{Fre65} $\Fr(\mcF)$ is an abelian category. \end{theorem}

\begin{proof} We already know that $\Fr(\mcF)$ is an additive category. We know from Lemma \ref{weak kernels implies kernels in freyd comp} that every morphism in $\Fr(\mcF)$ has both a kernel and a cokernel. We also know from Lemma \ref{monos in freyd comp are kernels} that monomorphisms in $\Fr(\mcF)$ are kernels and epimorphisms in $\Fr(\mcF)$ are cokernels. We conclude that $\Fr(\mcF)$ is an abelian category. \end{proof}

\begin{prop} \cite[Lemma 3.1.5]{Fre65} \label{inj proj obj in freyd comp} Let $T$ be an object in $\mcF$. Then $\id_T$ is both projective and injective in $\Fr(\mcF)$. \end{prop}

\begin{proof} Let $u:X\to Y$ be an object in $\Fr(\mcF)$ and $\pi:u\to\id_T$ be an epimorphism. Then $\pi=[(f,g)]$ for some $f:X\to T$ and $g:Y\to T$ in $\mcF$ with $gu=f$. \begin{center}
    \begin{tikzcd}[sep=2cm]
     X \arrow[r,"u"] \arrow[d,"f",swap] & Y \arrow[d,"g"] \\
     T \arrow[r,"\id_T"] & T 
    \end{tikzcd}
 \end{center} Let $\psi:T\to C$ be a cofibre of $f$, so that $\psi f=0$. Then $[(\id_T,\psi)]:\id_T\to\psi$ is a morphism in $\Fr(\mcF)$ with $[(\id_T,\psi)]\pi=0$. \begin{center}
    \begin{tikzcd}[sep=2cm]
     X \arrow[r,"u"] \arrow[d,"f",swap] & Y \arrow[d,"g"] \\
     T \arrow[r,"\id_T"] \arrow[d,equal,swap] & T \arrow[d,"\psi"] \\
     T \arrow[r,"\psi"] & C
    \end{tikzcd}
 \end{center} Since $\pi$ is an epimorphism, this shows that $[(\id_T,\psi)]=0$ and therefore $\psi=0$. Since $f$ is a weak kernel of $\psi=0$, there exists a morphism $\ta:T\to X$ with $f\ta=\id_T$. \begin{center}
    \begin{tikzcd} [sep=2cm]
     T \arrow[d,dashed,"\ta",swap] \arrow[rd,"\id_T"] & {} & {} \\
     X \arrow[r,"f"] & T \arrow[r,"0"] & C \\ 
     \end{tikzcd}\end{center} We then obtain a right inverse $[(\ta,u\ta)]:\id_T\to u$ of $\pi$, so that $\pi$ is a split epimorphism. \begin{center}
    \begin{tikzcd}[sep=2cm]
     T \arrow[r,"\id_T"] \arrow[d,"\ta",swap] & T \arrow[d,"u\ta"] \\
     X \arrow[r,"u"] \arrow[d,"f",swap] & Y \arrow[d,"g"] \\
     T \arrow[r,"\id_T"] & T
    \end{tikzcd}
 \end{center} Since every epimorphism to $\id_T$ splits, and $\Fr(\mcF)$ is abelian, we conclude that $\id_T$ is projective in $\Fr(\mcF)$. A similar proof shows that $\id_T$ is injective in $\Fr(\mcF)$. \end{proof}
 
\begin{lem} \cite[Theorem 3.1]{Fre65} \label{freyd comp enough proj inj} $\Fr(\mcF)$ has enough projectives and enough injectives. \end{lem}

\begin{proof} Let $u:X\to Y$ be an object in $\Fr(\mcF)$. The following diagrams define morphisms $\id_X\to u$ and $u\to\id_Y$. 
\begin{center}
    \begin{tikzcd}[sep=2cm]
     X \arrow[r,"\id_X"] \arrow[d,equal,swap] & X \arrow[d,"u"] & {} & 
     X \arrow[r,"u"] \arrow[d,"u",swap] & Y \arrow[d,equal] \\
     X \arrow[r,"u"] & Y & {} & Y \arrow[r,"\id_Y"] & Y
    \end{tikzcd}
 \end{center} Since $X,Y\in\mcF$, we know from Lemma \ref{inj proj obj in freyd comp} that $\id_X$ is projective and $\id_Y$ is injective. Also since $\id_X$ is an epimorphism and $\id_Y$ is a monomorphism, we see from Lemma \ref{monos and epis preservation in freyd completion} that the above map $\id_X\to u$ is an epimorphism from a projective object, and the map $u\to\id_Y$ is a monomorphism to an injective object. \end{proof}
 
\begin{lem} \cite[Theorem 3.1]{Fre65} If \begin{center}
    \begin{tikzcd} [sep=2cm]
     X \arrow[r,"u"] & 
     Y \arrow[r,"v"] & 
     Z \arrow[r,"w"] & 
     X[1]\;\;\;.
    \end{tikzcd}
 \end{center} is exact in $\mcF$, then \begin{center}
    \begin{tikzcd} [sep=2cm]
     \cdots \arrow[r] & \id_X \arrow[r,"{[(u,u)]}"] & 
     \id_Y \arrow[r,"{[(v,v)]}"] & 
     \id_Z \arrow[r] & \cdots \;\;\;
    \end{tikzcd}
 \end{center} is exact in $\Fr(\mcF)$. \end{lem}
 
\begin{proof} We know from Lemma \ref{weak kernels implies kernels in freyd comp} that the kernel of $[(v,v)]$, and the cokernel of $[(u,u)]$, are \begin{center}
    \begin{tikzcd}[sep=2cm]
     X \arrow[r,"u"] \arrow[d,"u",swap] & Y \arrow[d,equal] & {} & 
     Y \arrow[r,"\id_Y"] \arrow[d,equal,swap] & Y \arrow[d,"v"] \\
     Y \arrow[r,"\id_Y"] & Y & {} & Y \arrow[r,"v"] & Z
    \end{tikzcd}
 \end{center} respectively. Then the image of $[(u,u)]$ is the kernel of the cokernel of $[(u,u)]$, which is $[(u,\id_Y)]$. Since $\im([(u,u)])=\ker([(v,v)])$, we see that the above sequence is exact in $\Fr(\mcF)$. \end{proof}
 
\begin{lem} \cite[Theorem 3.1]{Fre65} Let $u:X\to Y$ be an object in $\Fr(\mcF)$. Complete $u$ to an exact triangle \begin{center}
    \begin{tikzcd} [sep=2cm]
     X \arrow[r,"u"] & 
     Y \arrow[r,"v"] & 
     Z \arrow[r,"w"] & 
     X[1]\;\;\;.
    \end{tikzcd}
 \end{center} Then \begin{center}
    \begin{tikzcd} [sep=1.5cm]
    	 0 \arrow[r] &
    	 u \arrow[r] &
     \id_Y \arrow[r] & 
     \id_Z \arrow[r] & 
     \id_{X[1]} \arrow[r] &
     \cdots
    \end{tikzcd}
 \end{center} is an injective resolution of $u$ in $\Fr(\mcF)$ and \begin{center}
    \begin{tikzcd} [sep=1.5cm]
    	 \cdots \arrow[r] &
    	 \id_{Y[-1]} \arrow[r] &
     \id_{Z[-1]} \arrow[r] & 
     \id_X \arrow[r] & 
     u \arrow[r] &
     0
    \end{tikzcd}
 \end{center} is a projective resolution. \end{lem}

\begin{proof} We know that $u\to\id_Y$ is a monomorphism, so must have image $u$. One can check that the kernel of $\id_Y\to\id_Z$ is also $u$. This confirms exactness at $u$ and $\id_Y$. Then the whole sequence must be exact since the rest is part of the long exact sequence induced by \begin{center}
    \begin{tikzcd} [sep=2cm]
     X \arrow[r,"u"] & 
     Y \arrow[r,"v"] & 
     Z \arrow[r,"w"] & 
     X[1]\;\;\;.
    \end{tikzcd}
 \end{center} Since identity morphisms are injective, this is an injective resolution for $u$. We can similarly confirm that the claimed projective resolution is valid. \end{proof}
 
\begin{lem} Every monomorphism and every epimorphism in $\mcF$ is isomorphic (in $\Fr(\mcF)$) to an identity morphism in $\mcF$, and is therefore both projective and injective in $\Fr(\mcF)$. \end{lem}

\begin{proof} Let $u:X\to Y$ be a monomorphism in $\mcF$. By Proposition \ref{zeros in exact tr split monos epis}, we can complete $u$ to an exact triangle \begin{center}
    \begin{tikzcd} [sep=2cm]
     X \arrow[r,"u"] & 
     Y \arrow[r,"v"] & 
     Z \arrow[r,"0"] & 
     X[1] \;\;\;,
    \end{tikzcd}
 \end{center} with $u$ a split monomorphism and $v$ a split epimorphism, and there exists an isomorphism $\ta:Y\to X\vee Z$ for which the following diagram commutes. \begin{center}
    \begin{tikzcd} [sep=2cm] \label{split exact triangle diagram}
     X \arrow[r,"u"] \arrow[d,equal] & 
     Y \arrow[r,"v"] \arrow[d,"\ta",dashed] & 
     Z \arrow[r,"0"] \arrow[d,equal] & 
     X[1] \arrow[d,equal] \\
     X \arrow[r] &
     X\vee Z \arrow[r] &
     Z \arrow[r,"0"] &
     X[1] 
    \end{tikzcd}
 \end{center} Then $[(\id_X,\ta)]:u\to\io$ is an isomorphism, where $\io$ is the inclusion $X\incl X\vee Z$. One can check that \begin{center}
    \begin{tikzcd}[sep=2cm]
     X \arrow[r,"\io"] \arrow[d,equal,swap] & X\vee Z \arrow[d,"\pi"] \\
     X \arrow[r,"\id_X"] & X 
    \end{tikzcd}
 \end{center} represents an isomorphism $\al:\io_X\to\id_X$ in $\Fr(\mcF)$ whose inverse is represented by \begin{center}
    \begin{tikzcd}[sep=2cm]
     X \arrow[r,"\id_X"] \arrow[d,equal,swap] & X \arrow[d,"\io"] \\
     X \arrow[r,"\io"] & X\vee Z 
    \end{tikzcd}
 \end{center} So $u$ is isomorphic to $\id_X$. A similar argument shows that every epimorphism is isomorphic to an identity morphism. \end{proof}
 
\textbf{From now on, we shall assume that all coproducts in the ambient category $\mcC$ exist and that every cohomology functor on $\mcC$ is representable.}
 
\begin{lem} \cite[Proposition 1.6.8]{Nee01} \label{splitting of idempotents} If $X\in\mcF$ and $e:X\to X$ satisfies $e^2=e$, then there exists $Y\in\mcF$ and morphisms $f:X\to Y$ and $g:Y\to X$ with $gf=e$ and $fg=\id_Y$. \begin{center}
    \begin{tikzcd} [sep=2cm]
     {} & Y \arrow[rr,"\id_Y"] \arrow[rd,"g"] & {} & Y \arrow[rd,"g"] & {} \\
     X \arrow[rr,"e"] \arrow[ru,"f"] & {} & X \arrow[rr,"e"] \arrow[ru,"f"] & {} & X \\ 
     \end{tikzcd}\end{center}\end{lem}

\begin{proof} See Proposition 1.6.8 of \cite{Nee01}. \end{proof}

\begin{lem} Idempotents, and retracts of identity morphisms, are isomorphic to identity morphisms, and are therefore both projective and injective in $\Fr(\mcF)$. \end{lem}

\begin{proof} Let $e:X\to X$ be an idempotent. Then Lemma \ref{splitting of idempotents} shows that there exists $R\in\mcF$ and morphisms $\vphi:X\to R$, $\psi:R\to X$ with $e=\psi\vphi$ and $\vphi\psi=\id_R$. Then we can form a pair of inverse morphisms \begin{center}
    \begin{tikzcd}[sep=2cm]
     X \arrow[r,"\psi\vphi"] \arrow[d,"\vphi",swap] & X \arrow[d,"\vphi"] & {} & 
     R \arrow[r,"\id_R"] \arrow[d,"\psi",swap] & R \arrow[d,"\psi"] \\
     R \arrow[r,"\id_R"] & R & {} & X \arrow[r,"\psi\vphi"] & X
    \end{tikzcd}
 \end{center} to show that $e$ is isomorphic to $\id_R$. 
 
Now let $T\in\mcF$ and $u:X\to Y$ be a retract of $\id_T$. Then there exist morphisms $\al:u\to\id_T$ and $\bt:\id_T\to u$ in $\Fr(\mcF)$ with $\bt\al=\id_u$. Then there exist morphisms $f:Y\to T$ and $g:T\to X$ for which $\al=[(fu,f)]$ and $\bt=[(g,ug)]$. We can thus represent the composite $\bt\al=\id_u$ as follows. \begin{center}
    \begin{tikzcd}[sep=2cm]
     X \arrow[r,"u"] \arrow[d,"fu",swap] & Y \arrow[d,"f"] \\
     T \arrow[r,"\id_T"] \arrow[d,"g",swap] & T \arrow[d,"ug"] \\
     X \arrow[r,"u"] & Y
    \end{tikzcd}
 \end{center} Let $e=gfu:X\to X$. Then we see that $ue=u$, and we can thus form the following mutually inverse morphisms. \begin{center}
    \begin{tikzcd}[sep=2cm]
     X \arrow[r,"u"] \arrow[d,equal,swap] & Y \arrow[d,"gf"] & {} & 
     X \arrow[r,"e"] \arrow[d,equal,swap] & X \arrow[d,"u"] \\
     X \arrow[r,"e"] & X & {} & X \arrow[r,"u"] & Y
    \end{tikzcd}
 \end{center} So $u$ is isomorphic to $e$. Then $e^2=gfugfu=gfue=gfu=e$, so $e$ is an idempotent, and we have shown above that idempotents are isomorphic to identity morphisms, so $u$ must be isomorphic to an identity morphism. \end{proof}
 
\begin{lem} \cite[Theorem 3.1]{Fre65} \label{proj inj objects in freyd comp} Let $u:X\to Y$ be a morphism in $\mcF$. Then the following are equivalent. 

1. $u$ is an injective object in $\Fr(\mcF)$.

2. $u$ is a projective object in $\Fr(\mcF)$. 

3. $u$ is isomorphic to $\id_T$ for some $T\in\mcF$. \end{lem}

\begin{proof} If $u$ is a projective object in $\Fr(\mcF)$, then we know from Lemma \ref{freyd comp enough proj inj} that \begin{center}
    \begin{tikzcd}[sep=2cm]
     X \arrow[r,"\id_X"] \arrow[d,equal,swap] & X \arrow[d,"u"] \\
     X \arrow[r,"u"] & Y
    \end{tikzcd}
 \end{center} defines an epimorphism $\pi:\id_X\to u$. Since $u$ is projective, $\pi$ splits, so there exists a morphism $\al:u\to\id_X$ with $\pi\al=id_u$, so that $u$ is a retract of $\id_X$ and therefore isomorphic to $\id_T$ for some $T\in\mcF$. We know from Lemma \ref{inj proj obj in freyd comp} that every morphism isomorphic to an identity morphism is injective. If $u:X\to Y$ is injective, then there exists a monomorphism $\io:u\to\id_Y$, and this splits since $u$ is injective. So $u$ is a retract of $\id_Y$, which is projective, and therefore $u$ is projective. \end{proof}

\begin{definition} A \textit{Frobenius category} is an abelian category with enough projectives and enough injectives, where the injectives and projectives coincide.  \end{definition}

\begin{rmk} Frobenius categories are often defined more broadly than in the context of abelian categories (such as in \cite{HoZ12}), but we shall restrict our attention to the abelian case since $\Fr(\mcF)$ is always abelian. \end{rmk}

\begin{prop} \cite[Theorem 3.1]{Fre65} $\Fr(\mcF)$ is a Frobenius category whose projective-injective objects are those isomorphic to identity morphisms. \end{prop}

\begin{proof} We have already seen that $\Fr(\mcF)$ is an abelian category with enough projectives and enough injectives. Then the result follows from Lemma \ref{proj inj objects in freyd comp}. \end{proof}

\begin{prop} \cite[\S 8]{Fre65} If $u:X\to Y$ and $v:Z\to W$ are objects in $\Fr(\mcF)$, then $$\Ext_{\Fr(\mcF)}^{i+3}(u,v)=\Ext_{\Fr(\mcF)}^i(u,v[1])$$ for $i\geq 1$. \end{prop}

\begin{proof} Complete $v$ to an exact triangle \begin{center}
    \begin{tikzcd} [sep=2cm]
     Z \arrow[r,"v"] & 
     W \arrow[r,"\pi"] & 
     P \arrow[r,"\io"] & 
     Z[1] \;\;\;.
    \end{tikzcd}
 \end{center} Then \begin{center}
    \begin{tikzcd} [sep=1.5cm]
    	 0 \arrow[r] &
    	 v \arrow[r] &
     \id_W \arrow[r] & 
     \id_P \arrow[r] & 
     \id_{Z[1]} \arrow[r] &
     \cdots
    \end{tikzcd}
 \end{center} is an injective resolution for $v$ in $\Fr(\mcF)$. Since the suspension functor on $\Fr(\mcF)$ is exact, we see that \begin{center}
    \begin{tikzcd} [sep=1.3cm]
    	 0 \arrow[r] &
    	 v[1] \arrow[r] &
     \id_{W[1]} \arrow[r] & 
     \id_{P[1]} \arrow[r] & 
     \id_{Z[2]} \arrow[r] &
     \cdots
    \end{tikzcd}
 \end{center} is exact, and therefore an injective resolution for $v[1]$. To compute $\Ext_{\Fr(\mcF)}^{i+3}(u,v)$ and $\Ext_{\Fr(\mcF)}^i(u,v[1])$, we form the cochain complexes \begin{center}
    \begin{tikzcd} [sep=1cm]
    	 0 \arrow[r] &
     {[u,\id_W]} \arrow[r] & 
     {[u,\id_P]} \arrow[r] & 
     {[u,\id_{Z[1]}]} \arrow[r] &
     {[u,\id_{W[1]}]} \arrow[r] &
     \cdots
    \end{tikzcd}
 \end{center} and \begin{center}
    \begin{tikzcd} [sep=0.75cm]
    	 0 \arrow[r] &
     {[u,\id_{W[1]}]} \arrow[r] & 
     {[u,\id_{P[1]}]} \arrow[r] & 
     {[u,\id_{Z[2]}]} \arrow[r] &
     {[u,\id_{W[2]}]} \arrow[r] &
     \cdots \;\;\;,
    \end{tikzcd}
 \end{center} where here [-,-] denote hom-sets in $\Fr(\mcF)$. Then, for $i\geq 1$, we observe that the cohomology of the first cochain complex at $i+3$ is the same as the cohomology of the second cochain complex at $i$, and we are done. \end{proof}

\section{Stable Categories} \label{stable cat section}

In this section, we describe the stable category of a Frobenius category and its canonical triangulated structure, following \cite{HoZ12}. We fix a Frobenius category $\mcA$. 

\begin{definition} \cite{HoZ12} We say that a morphism $f:A\to B$ in $\mcA$ is \textit{projective-injective} if it factors through a projective-injective object. Given objects $A$ and $B$, it is routine to check that the projective-injective morphisms $A\to B$ form a subgroup of $\Hom(A,B)$. If $f:A\to B$ and $g:B\to C$ with either $f$ or $g$ projective-injective, then $gf$ is projective-injective. If $f,g:A\to B$, we write $f\sim g$ to indicate that $f-g$ is projective-injective. Then $\sim$ defines an equivalence relation on $\Hom(A,B)$. Since $\sim$ respects composition and additivity, we can define the \textit{stable category} of $\mcA$ to be $\St(\mcA)=\mcA/\sim$. There is an obvious quotient functor $\mcA\to\St(\mcA)$ sending each morphism to its equivalence class. 

For each object $A$ in $\mcA$, we can form a monomorphism $A\to I$ into a projective-injective object $I$, and then take its cokernel, which gives an object in $\St(\mcA)$. It follows from Schanuel's Lemma (Lemma 6 of \cite{HoZ12}) that this object is unique up to isomorphism in $\St(\mcA)$. Label it $T(A)$.

Let $f:A\to B$ be a morphism in $\mcA$. Let $A\to I$ and $B\to J$ be monomorphisms to projective-injective objects. Since $J$ is injective, we can find a morphism $g:I\to J$ making the left hand side of the following diagram commute.\begin{center}
    \begin{tikzcd} [sep=2cm]
     0 \arrow[r] \arrow[d] & 
     A \arrow[r] \arrow[d,"f"] & 
     I \arrow[r] \arrow[d,"g",dashed] & 
     T(A) \arrow[r] \arrow[d,"h",dashed] & 
     0 \arrow[d] \\
     0 \arrow[r] &
     B \arrow[r] &
     J \arrow[r] &
     T(B) \arrow[r] &
     0
    \end{tikzcd}
 \end{center}  

Then, given $g$, there exists a unique morphism $h:T(A)\to T(B)$ making the right hand side commute. One can check that this process gives a unique morphism $T(A)\to T(B)$ in $\St(\mcA)$ and a resulting additive functor $\mcA\to\St(\mcA)$. Since this sends projective-injective morphisms to zero morphisms, we obtain an additive functor $T:\St(\mcA)\to\St(\mcA)$.

We can similarly use epimorphisms from projective-injective objects to form another functor $S:\St(\mcA)\to\St(\mcA)$, and it turns out (Theorem 2 of \cite{HoZ12}) that these functors are inverses of each other, so that $T$ is an additive automorphism on $\St(\mcA)$. 

Given a morphism $f:A\to B$, and a monomorphism $A\to I$, where $I$ is projective-injective, we can construct a pushout $C$ of $f$ and $A\to I$. \begin{center}
    \begin{tikzcd} [sep=2cm]
     A \arrow[r] \arrow[d,"f"] & 
     I \arrow[r] \arrow[d] & 
     \cok(A\to I) \arrow[d,equal] \\
     B \arrow[r,"g"] &
     C \arrow[r,"h"] &
     \cok(A\to I)
    \end{tikzcd}
 \end{center} Then \begin{center}
    \begin{tikzcd} [sep=2cm]
     A \arrow[r,"{[f]}"] & 
     B \arrow[r,"{[g]}"] & 
     C \arrow[r,"{[h]}"] & 
     T(A)\;\;\;.
    \end{tikzcd}
 \end{center} is a triangle in $\St(\mcA)$. The triangles that are isomorphic to a triangle of the above form satisfy the axioms for a triangulation, so that $\St(\mcA)$ is a triangulated category (Theorem 4 of \cite{HoZ12}). See \cite{HoZ12} for the full details of these constructions. \end{definition}

\begin{rmk} \label{proj-inj become zero in stable cat} Let $A$ be a projective-injective object in a Frobenius category $\mcA$. Then $\id_A$ is a projective-injective morphism in $\mcA$, so that its equivalence class $[\id_A]$ is the zero morphism $A\to A$ in $\St(\mcA)$. So in $\St(\mcA)$, the identity on $A$ is the zero morphism. So $\Hom_{\St(\mcA)}(A,A)=0$ and we see that $A$ is the zero object in $\St(\mcA)$. So all projective-injective objects in $\mcA$ become zero objects in $\St(\mcA)$. \end{rmk}

\begin{rmk} We now return to the triangulated category $\mcF$ from \S \ref{freyd comp}. We know that $\Fr(\mcF)$ is a Frobenius category, so that $\St(\Fr(\mcF))$ is again a triangulated category. But for each $X\in\mcF$, we know that $\id_X$ is projective-injective in $\Fr(\mcF)$, so must be a zero object in $\St(\Fr(\mcF))$ by Remark \ref{proj-inj become zero in stable cat}. We conclude that $\mcF\to\St(\Fr(\mcF))$ is the zero functor. \end{rmk}

\section{Tensor Triangulated Categories} \label{tensor triang}

We now return to examining our fixed triangulated category $\mcC$ from \S\ref{freyd comp}, and introduce a compatible closed symmetric monoidal structure to $\mcC$, following Appendix A.2 of \cite{HPS97}. This provides a notion of tensor product with an appropriate 'tensor-hom adjunction'. We write the tensor product as $\w$ and call it the smash product to be consistent with the Spanier-Whitehead category. We shall introduce $i$-spheres for each $i\in\Z$, construct a graded-commutative ring $R$ and a functor $\pi_*:\mcC\to\Mod_R$. The action of this functor on the full subcategory $\mcF$ will be of particular importance in \S\ref{freyd semi freyd}.

\begin{definition} A \textit{closed symmetric monoidal structure} on $\mcC$ consists of a functor $\w:\mcC\x\mcC\to\mcC$ which we call the \textit{smash product}, an object $S$ which we call the \textit{$0$-sphere}, natural isomorphisms $$(X\w Y)\w Z\simeq X\w(Y\w Z),$$ $$S\w X\simeq X,$$ and one sending $$X\w Y\simeq Y\w X$$ for objects $X,Y,Z\in\mcC$ subject to various coherence conditions which we shall not list here. We also insist that for each object $Y\in\mcC$ that the functor $(-)\w Y:\mcC\to\mcC$ has a right adjoint $F(Y,-):\mcC\to\mcC$, so that we may form a functor $F:\mcC^{\text{op}}\to\mcC\to\mcC$ and a natural isomorphism $$[X\w Y,Z]\simeq[X,F(Y,Z)]$$ for objects $X,Y,Z\in\mcC$. There are also various requirements guaranteeing compatibility with the triangulated structure listed in Definition A.2.1 of \cite{HPS97}. From now on, we shall assume that $\mcC$ comes equipped with a closed symmetric monoidal structure compatible with its triangulation. \end{definition}

\begin{definition} \cite[Definition A.2.4]{HPS97} We write $D$ for the contravariant functor $F(-,S)$, and call this the \textit{Spanier-Whitehead dual} for $\mcC$. Note that $D$ is exact and takes coproducts to products. For each $X\in\mcC$, we have that $$[DX,DX]\simeq[X\w DX,S]\simeq[X,D^2X].$$ Denote the morphism $X\to D^2X$ corresponding to $\id_{DX}$ by $\chi_X$. Then $\chi:\id_{\mcC}\to D^2$ is a natural transformation. \end{definition}

\begin{theorem} \cite[Theorem A.2.5]{HPS97} We have the following.

1. If $X,Y\in\mcF$, then so are $X\w Y$ and $F(X,Y)$. In particular, $DX\in\mcF$.

2. $\chi_X:X\to D^2X$ is an isomorphism for all $X\in\mcF$. So $D$ restricts to a duality functor $\mcF\to\mcF$, and $DX=0$ if and only if $X=0$ for $X\in\mcF$. \end{theorem}

\begin{definition} \label{construction of stable homotopy ring} For each $i\in\Z$, define $S^i=S[i]$, which we call the \textit{$i$-sphere}. Then $S^i\w S^j\simeq S^{i+j}$ for all $i,j\in\Z$. Now define $\pi_i(X)=[S^i,X]$ for all $X\in\mcC$ and $i\in\Z$, the \textit{$i^{\text{th}}$ homotopy group} of $X$. This defines a graded abelian group $\pi_*(X)$ for all $X\in\mcC$. In particular, set $R=\pi_*(S)$.

For an object $X\in\mcC$, $i,j\in\Z$ and morphisms $S^i\to S$, $S^j\to X$ we can take the smash product to form a morphism $S^{i+j}\to X$. This gives a map $R_i\x\pi_j(X)\to \pi_{i+j}(X)$, and when $X=S$ we obtain a map $R_i\x R_j\to R_{i+j}$. We see that $R$ is a graded-commutative ring which we call the \textit{stable homotopy ring}, and that $\pi_*(X)$ is a graded $R$-module. Denote the category of graded $R$-modules by $\Mod_R$.

For each morphism $f:X\to Y$ in $\mcC$, $i\in\Z$, we set $\pi_i(f)=[S^i,f]$, producing a graded $R$-linear map $\pi_*(f):\pi_*(X)\to\pi_*(Y)$. We have now defined an additive functor $\pi_*:\mcC\to\Mod_R$. \end{definition}

\begin{definition} \cite[Definition 1.1.2]{HPS97} Let $X\in\mcC$. Using the adjunction of $\w$ and $F$, we obtain an isomorphism $$[DX,DX]\simeq[X\w DX,S].$$ Then $\id_{DX}$ corresponds to some morphism $$X\w DX\to S.$$ For each object $Y$ in $\mcC$, we take the smash product with $\id_Y$ to obtain a morphism $$X\w DX\w Y\to Y.$$ Then using the adjunction of $\w$ and $F$ again we obtain a morphism $$DX\w Y\to F(X,Y).$$ We say that $X$ is \textit{strongly dualisable} if the above morphism $$DX\w Y\to F(X,Y)$$ is an isomorphism for each $Y$. \end{definition}

\begin{exmp} $S$ is strongly dualisable. \end{exmp}

\begin{proof} For each object $A$, we have that $$[A,F(S,X)]\simeq[A\w S,X]\simeq[A,X]$$ by the adjunction of $\w$ and $F$. Then $F(S,X)=X$ by the Yoneda Lemma. Then $$F(S,S)\w X=S\w X=X=F(S,X),$$ so that $S$ is strongly dualisable. \end{proof}

\begin{definition} \cite[\S 4]{Fre65} We can extend the duality functor $D:\mcF\to\mcF$ to a duality functor $\Fr(\mcF)\to\Fr(\mcF)$ as follows. Let $u:X\to Y$ be an object in $\Fr(\mcF)$. Then $Du:DY\to DX$ is also an object in $\Fr(\mcF)$. Now let \begin{center}
    \begin{tikzcd}[sep=2cm]
     X \arrow[r,"u"] \arrow[d,"f",swap] & Y \arrow[d,"g"] \\
     Z \arrow[r,"v"] & W
    \end{tikzcd}
 \end{center} represent a morphism $\chi:u\to v$ in $\Fr(\mcF)$. Then we obtain a morphism \begin{center}
    \begin{tikzcd}[sep=2cm]
     DW \arrow[r,"Dv"] \arrow[d,"Dg",swap] & DZ \arrow[d,"Df"] \\
     DY \arrow[r,"Du"] & DX
    \end{tikzcd}
 \end{center} One can check that this construction respects equivalence, so that we obtain a map $$[(f,g)]\mapsto[(Dg,Df)]:\Hom_{\Fr(\mcF)}(u,v)\to\Hom_{\Fr(\mcF)}(Dv,Du).$$ This defines a functor $\Fr(\mcF)\to\Fr(\mcF)$ sending $\id_X$ to $\id_{DX}$ for all $X\in\mcF$.

Since $\chi:\id_{\mcF}\to D^2$ is a natural isomorphism, we have that for each object $u:X\to Y$ in $\Fr(\mcF)$, that \begin{center}
    \begin{tikzcd}[sep=2cm]
     X \arrow[r,"u"] \arrow[d,"\chi_X",swap] & Y \arrow[d,"\chi_Y"] \\
     D^2X \arrow[r,"D^2u"] & D^2Y
    \end{tikzcd}
 \end{center} represents an isomorphism $u\to D^2u$, which we label $\chi_u$. We thus obtain a natural isomorphism $\chi:\id_{\Fr(\mcF)}\to D^2$ with $\chi_{\id_X}=[(\chi_X,\chi_X)]$ for all $X\in\mcF$. \end{definition}

\begin{definition} We define the functor $\Fr(\mcF)\to\Mod_R$ as follows. We send each object $u:X\to Y$ in $\Fr(\mcF)$ to $\im(\pi_*(u))$ in $\Mod_R$. Represent a morphism $\chi$ in $\Fr(\mcF)$ by the commuting square \begin{center}
    \begin{tikzcd}[sep=2cm]
     X \arrow[r,"u"] \arrow[d,"f",swap] &
     Y \arrow[d,"g"] \\
     Z \arrow[r,"v"] & 
     W
    \end{tikzcd}
 \end{center} Then $\pi_*(g)(\im(\pi_*(u))\leq\im(\pi_*(v))$, so that we can form a morphism $\im(\pi_*(u))\to\im(\pi_*(v))$. We see that the composite $\mcF\to\Fr(\mcF)\to\Mod_R$ is the same as the inclusion $\mcF\to\Mod_R$, so that we may consider this functor to be an extension of $\pi_*$ to $\Fr(\mcF)$ and also label it $\pi_*$. \begin{center}
    \begin{tikzcd} [sep=1.5cm]
     {} & \Mod_R \\
     \mcF \arrow[r,hookrightarrow] \arrow[ru,"\pi_{*}"] & \Fr(\mcF) \arrow[u,dashed] \\ 
     \end{tikzcd}\end{center} \end{definition}
     
\section{Reflexive Modules} \label{reflexive graded modules}

In this section, we fix a graded-commutative ring $R$ and an injective graded $R$-module $\Ta$. Denote the category of graded $R$-modules by $\Mod_R$. We construct a full subcategory $\mcM$ of $\Mod_R$ and a duality functor $D:\mcM\to\mcM$. Once we have introduced some further assumptions on $\mcC$, specifically that it is a Freyd category, we will see in \S\ref{freyd semi freyd} that $R$ is self-injective, so that the results that follow will apply to the stable homotopy ring $R$ associated with $\mcC$.

However, we phrase this section in slightly broader generality to emphasise that the core arguments apply in other contexts that will be useful to us: they apply to $(A,\Ta)$, where $A$ is the Hahn ring and $\Ta$ is an injective $A$-module, in Chapter \ref{hahn ring chapter}, with further results in this special case presented in Chapter \ref{reflexivity}. They also apply to the infinite root algebra $P$ in Chapter \ref{infinite root algebra chapter}, which, as we will see, is self-injective. Finally, they apply to a complete local Noetherian ring together with the injective hull of its residue field (see, for example, Section 47.7 of \cite{StP}). Note that these final three cases are in the ungraded context.
     
\begin{definition} \cite[\S 2.2]{NaO04} Let $d\in\Z$. For each graded $R$-module $M$, we define the graded $R$-module $M(d)$ to be identical to $M$ as an $R$-module, but with a grading given by $M(d)_i=M_{i-d}$. For each morphism $\al:M\to N$, we define $\al(d):M(d)\to N(d)$ by setting $\al(d)_i=\al_{i-d}$. This defines an exact automorphism on $\Mod_R$ with inverse $M\mapsto M(-d)$. \end{definition}

\begin{definition} \cite[\S 2.4]{NaO04} Let $M$ and $N$ be graded $R$-modules, and set $$\Hom^*(M,N)=\bigoplus_i\Hom(M(i),N).$$ Since $$A_i\Hom(M(j),N)\subs\Hom(M(i+j),N),$$ $\Hom^*(M,N)$ defines a graded $A$-module. As with the ungraded Hom functor, we have $\Hom^*(R,M)\simeq M$ for all $M$. We have constructed a functor $$\Hom^*:\Mod_R^{\text{op}}\x\Mod_R\to\Mod_R.$$ \end{definition}

\begin{definition} \cite[\S 2]{ShS14} Set $D=\Hom^*(-,\Ta)$. For each graded $R$-module $M$, define the evaluation map $\chi_M:M\to D^2M$ by setting $$(\chi_M)(m)(\al)=\al(m)$$ for all $m\in M$ and $\al:M\to\Ta$. If $\al:M\to N$ is any morphism, then the following diagram commutes. \begin{center}
    \begin{tikzcd}[sep=2cm]
     M \arrow[r,"\al"] \arrow[d,"\chi_M",swap] &
     N \arrow[d,"\chi_N"] \\
     D^2M \arrow[r,"D^2\al"] & 
     D^2N \;\;\;.
    \end{tikzcd}
 \end{center} So $M\mapsto\chi_M$ defines a natural transformation $\chi:\id_{\Mod_R}\to D^2$. \end{definition}
 
\begin{prop} \cite[\S 2]{ShS14} \label{contrav internal hom exact} $D:\Mod_R\to\Mod_R$ is exact, and therefore switches kernels and cokernels, monomorphisms and epimorphisms. \end{prop}

\begin{proof} Since $\Mod_R$ is abelian and $\Ta$ is injective, $\Hom(-,\Ta)$ must be exact. One can check that the following diagram commutes for each $i$, where $\pi_i$ is the projection to the $i^{\text{th}}$ summand. \begin{center}
    \begin{tikzcd}[sep=2cm]
     \Mod_R \arrow[r,"D"] \arrow[d,"M\mapsto M(i)",swap] &
     \Mod_R \arrow[d,"\pi_i"] \\
     \Mod_R \arrow[r,"{\Hom(-,\Ta)}"] & 
     \text{Ab}
    \end{tikzcd}
 \end{center} So $\pi_iD$ is exact for all $i$ and then $D$ must be exact. The remaining claims follow easily from the fact that $D$ is exact. \end{proof}
 
\begin{prop} \label{identity for D and chi} If $M$ is a graded $R$-module, then $D(\chi_M)\chi_{DM}=\id_{DM}$. \end{prop}

\begin{definition} \cite[Definition 2.1]{ShS14} We say that a graded $R$-module $M$ is \textit{reflexive} if $\chi_M$ is an isomorphism. Denote by $\mcM$ the full subcategory of $\Mod_R$ consisting of the reflexive graded $R$-modules.\end{definition}

\begin{lem} \cite[Proposition 2.2]{ShS14} \label{if left and right are ref, then so is middle} Suppose that \begin{center}
    \begin{tikzcd} [sep=2cm]
     0 \arrow[r] & 
     L \arrow[r,"\al"] & 
     M \arrow[r,"\bt"] & 
     N \arrow[r] &
     0,
    \end{tikzcd}
   \end{center} is exact. If $L,N\in\mcM$, then $M\in\mcM$. \end{lem}
   
\begin{proof} This follows easily from the Snake Lemma. \end{proof}

\begin{prop} \cite[Proposition 2.2]{ShS14} \label{properties of reflexive modules general case} The following properties of reflexive graded $R$-modules hold. 

1. If $M,N\in\mcM$, then $M\oplus N\in\mcM$. 

2. If $N\leq M$, and both $N,M/N\in\mcM$, then $M\in\mcM$. 

3. The category $\mcM$ is closed under isomorphisms.

4. The category $\mcM$ is closed under kernels, cokernels and images.

5. If $M\in\mcM$, then $DM\in\mcM$. 

6. If $M\in\mcM$ and $i\in\Z$, then $M(i)\in\mcM$. \end{prop}

\begin{proof} 1. Apply Lemma \ref{if left and right are ref, then so is middle} to the split exact sequence \begin{center}
    \begin{tikzcd} [sep=1.5cm]
     0 \arrow[r] & 
     M \arrow[r] & 
     M\oplus N \arrow[r] & 
     N \arrow[r] &
     0.
    \end{tikzcd}
   \end{center} 
   
2. Apply Lemma \ref{if left and right are ref, then so is middle} to the short exact sequence \begin{center}
    \begin{tikzcd} [sep=1.5cm]
     0 \arrow[r] & 
     N \arrow[r] & 
     M \arrow[r] & 
     M/N \arrow[r] &
     0.
    \end{tikzcd}
   \end{center} 
   
3. This follows easily from the fact that $\chi$ is a natural transformation.

4. Let $\al:M\to N$ be a morphism in $\Mod_R$. Then since $\chi$ is a natural transformation, the following diagram commutes, and the rows are exact. \begin{center}
    \begin{tikzcd} [sep=2cm]
     0 \arrow[r] \arrow[d] & 
     \ker(\al) \arrow[r] \arrow[d,"\chi_{\ker(\al)}",swap] & 
     M \arrow[r,"\al"] \arrow[d,"\chi_M",swap] &
     N \arrow[d,"\chi_N",swap] \\
     0 \arrow[r] &
     D^2\ker(\al) \arrow[r,swap] &
     D^2M \arrow[r,"D^2\al",swap] &
     D^2N.
    \end{tikzcd}
   \end{center} Since $M,N\in\mcM$ we know that $\chi_M$ is an isomorphism and $\chi_N$ is a monomorphism, so that $\chi_{\ker(\al)}$ is an isomorphism by the Five Lemma and therefore $\ker(\al)\in\mcM$. A similar argument shows that $\cok(\al)\in\mcM$. Then $\im(\al)=\ker(N\to\cok(\al))\in\mcM$. 
   
5. If $M\in\mcM$, then we see from Proposition \ref{identity for D and chi} that $\chi_{DM}=D(\chi_M)^{-1}$ is an isomorphism, so that $DM\in\mcM$.

6. This follows from the fact that $\chi_{M(i)}=\chi_M(i)$. \end{proof}

\begin{theorem} \cite[Proposition 2.2]{ShS14} The category $\mcM$ is abelian, closed under isomorphisms, under $D$ and under $M\mapsto M(i)$ for each $i$. \end{theorem}

\begin{proof} This follows immediately from Proposition \ref{properties of reflexive modules general case}. \end{proof}

\section{Monogenic Stable Homotopy Categories} \label{monogenic}

We now introduce two further assumptions about $\mcC$ to arrive at the definition of a monogenic stable homotopy category, as first defined in \cite{HPS97}. 

\begin{definition} \cite[Definition 1.1.4]{HPS97} A \textit{monogenic stable homotopy category} is a tensor triangulated category $\mcC$ for which $S$ is small, $\loc(S)=\mcC$ (recall Definition \ref{thick localising subcat def}), arbitrary coproducts exists, and every cohomology functor is representable. 

Since we are already assuming the latter two properties, we are now simply introducing the requirements that $S$ is small and $\loc(S)=\mcC$. \end{definition}

\begin{lem} \label{monogenic properties} \cite[Lemma 1.4.5]{HPS97} We have the following. 

1. If $\eta$ is a natural transformation of (co)homology functors on $\mcC$ and $\eta_{S^i}$ is an isomorphism for all $i\in\Z$, then $\eta$ is a natural isomorphism. 

2. If $X\in\mcC$ and $\pi_*(X)=0$, then $X=0$. 

3. If $f:X\to Y$ is a morphism in $\mcC$ and $\pi_*(f)$ is an isomorphism, then $f$ is an isomorphism. \end{lem}

\begin{lem} Let $X\in\mcF$. Then $X[i]\in\mcF$ for $i\in\Z$. In particular, $S^i\in\mcF$ for all $i$, so that each $\pi_i$ is homological. \end{lem}

\begin{proof} Let $X\in\mcF$. Then \begin{center}
    \begin{tikzcd} [sep=2cm]
     X \arrow[r] & 
     0 \arrow[r] & 
     X \arrow[r,"\id_{X{[1]}}"] & 
     X[1]
    \end{tikzcd}
 \end{center} is exact, with $X,0\in\mcF$ and $\mcF$ thick, so that $X[1]\in\mcF$. Similarly, $X[-1]\in\mcF$ and we see by induction that $X[i]\in\mcF$ for all $i$. \end{proof}

\begin{lem} \label{pi star sends sphere to shifted ring} $\pi_*(S^i)\simeq R(i)$ for all $i$. \end{lem}

\begin{proof} For each $j$, we have that $$\pi_*(S^i)_j=[S^j,S^i]\simeq[S^{j-i},S]=R_{j-i}=R(i)_j.$$ \end{proof}

\begin{lem} \label{finite wedge of spheres gives iso} If $X$ is a finite wedge sum of spheres, then the induced map $$[X,Y]\to\Hom(\pi_*(X),\pi_*(Y))$$ is an isomorphism for all $Y\in\mcC$. \end{lem}

\begin{proof} Let $X=\bigvee_i S^{d_i}$. Since $[-,Y]$ is cohomological, we see that $$[X,Y]=\left[\bigvee_i S^{d_i},Y\right]\simeq\bigoplus_i[S^{d_i},Y]=\bigoplus_i\pi_{d_i}(Y)=\bigoplus_i\pi_*(Y)_{d_i}.$$ Then since $\Hom^*(R,M)\simeq M$ for each graded $R$-module $M$, we see that $$\bigoplus_i\pi_*(Y)_{d_i}\simeq\bigoplus_i\Hom^*(R,\pi_*(Y))_{d_i}=\bigoplus_i\Hom(R(d_i),\pi_*(Y)).$$ Then we use Lemma \ref{pi star sends sphere to shifted ring} to see that $$\bigoplus_i\Hom(R(d_i),\pi_*(Y))\simeq\Hom\left(\bigoplus_i R(d_i),\pi_*(Y)\right)\simeq\Hom\left(\bigoplus_i\pi_*(S^{d_i}),\pi_*(Y)\right).$$ Finally, the additivity of $\pi_*$ shows that $$\Hom\left(\bigoplus_i\pi_*(S^{d_i}),\pi_*(Y)\right)\simeq\Hom\left(\pi_*\left(\bigvee_i S^{d_i}\right),\pi_*(Y)\right)\simeq\Hom(\pi_*(X),\pi_*(Y)).$$ \end{proof}

\begin{lem} \cite[Lemma 2.1]{Hov07} \label{pi star gives isos} Let $Y\in\mcC$ with $\pi_*(Y)$ injective as a graded $R$-module. Then the induced map $$[X,Y]\to\Hom_R(\pi_*(X),\pi_*(Y))$$ is an isomorphism for all $X\in\mcC$. \end{lem}

\begin{proof} Since $\pi_*(Y)$ is injective, we know that both $[-,Y]$ and $\Hom_R(\pi_*(-),\pi_*(Y))$ are cohomology functors on $\mcC$. The functor $\pi_*$ induces a natural transformation $$\al:[-,Y]\to\Hom_R(\pi_*(-),\pi_*(Y)),$$ and we see from Lemma \ref{finite wedge of spheres gives iso} that $\al_{S^n}$ is an isomorphism for all $n\in\Z$. Then it follows from part 1 of Lemma \ref{monogenic properties} that $\al_X$ is an isomorphism for all $X\in\mcC$. \end{proof}

\begin{theorem} \cite[Theorem 2.1.3(d)]{HPS97} Let $\mcC$ be a monogenic stable homotopy category and $X\in\mcC$. Then the following are equivalent. 

1. $X\in\mcF$.

2. $X$ is strongly dualisable.

3. $X\in\text{thick}(S)$.

4. For any family of objects $\{Y_i\}_{i\in I}$ in $\mcC$, the natural map $$\bigvee_{i\in I}F(X,Y_i)\to F(X,\bigvee_{i\in I}Y_i)$$ is an isomorphism. \end{theorem}

\section{Freyd and semi-Freyd Categories} \label{freyd semi freyd}

We now introduce the notion of a Freyd category and a semi-Freyd category, following \cite{HPS97} and providing relatively straightforward axiomatic generalisations of results from \cite{Hov07}. Using the terminology of \cite{HPS97}, Freyd and semi-Freyd categories are special types of monogenic stable homotopy category. Semi-Freyd categories impose additional algebraic requirements on $R$, whereas a Freyd category is a semi-Freyd category with a much stronger assumption: that the restricted functor $\pi_*:\mcF\to\Mod_R$ is faithful.

\begin{definition} A \textit{semi-Freyd} category is a monogenic stable homotopy category $\mcC$ where $R_0$ is a complete local Noetherian ring, $R_i$ is finitely generated over $R_0$ for all $i\in\Z$, and $R_i=0$ for $i<0$. A \textit{Freyd category} is a semi-Freyd category for which the restricted functor $\pi_*:\mcF\to\Mod_R$ is faithful. 

We shall assume for the remainder of this section that $\mcC$ is a semi-Freyd category, and emphasise when results require $\mcC$ to be a Freyd category. \end{definition}

\begin{prop} \cite[Proposition 6.0.6]{HPS97} $[X,Y]$ is a finitely generated $R_0$-module for all $X,Y\in\mcF$. \end{prop}

\begin{proof} Let $\mcD$ be the full subcategory of $\mcC$ consisting of those $X\in\mcC$ for which $\pi_i(X)$ is finitely generated over $R_0$ for all $i\in\Z$. Since $\mcC$ is semi-Freyd, we know that $\pi_i(S)=R_i$ is finitely generated over $R_0$ for all $i$, so that $S\in\mcD$. 

Let \begin{center}
    \begin{tikzcd} [sep=2cm]
     X \arrow[r,"u"] & 
     Y \arrow[r,"v"] & 
     Z \arrow[r,"w"] & 
     X[1]
    \end{tikzcd}
 \end{center} be exact. Note that if $X\in\mcD$ and $n\in\Z$, then since $$\pi_i(X[n])=[S^i,X[n]]\simeq[S^{i-n},X]=\pi_{i-n}(X)$$ for all $i$, we see that $X[n]\in\mcD$. Let $X,Z\in\mcD$. Since each $\pi_i$ is homological, we know that \begin{center}
    \begin{tikzcd} [sep=1.3cm]
     \pi_i(Z[-1]) \arrow[r,"\al"] &
     \pi_i(X) \arrow[r,"\bt"] &
     \pi_i(Y) \arrow[r,"\gm"] & 
     \pi_i(Z) \arrow[r,"\dl"] & 
     \pi_i(X[1])
    \end{tikzcd}
 \end{center} is exact in $\Mod_{R_0}$. Then we can form a short exact sequence \begin{center}
    \begin{tikzcd} [sep=1.5cm]
     0 \arrow[r] &
     \pi_i(X)/\ker(\al) \arrow[r] &
     \pi_i(Y) \arrow[r] & 
     \im(\gm) \arrow[r] & 
     0
    \end{tikzcd}
 \end{center} Since $\pi_i(X)$ is finitely generated, we know that $\pi_i(X)/\ker(\al)$ is finitely generated. Since $R_0$ is Noetherian, $\im(\gm)\leq \pi_i(Z)$ and $\pi_i(Z)$ is finitely generated, we know that $\im(\gm)$ is finitely generated. So $\pi_i(Y)$ is an extension of finitely generated $R_0$-modules and therefore must be finitely generated. Since this holds for all $i$, we conclude that $Y\in\mcD$. Similarly, if $X,Y\in\mcD$, then $Z\in\mcD$, and if $Y,Z\in\mcD$, then $X\in\mcD$. 

If $Y\in\mcD$ and $X\in\mcC$ is a retract of $Y$, then $\pi_i(X)$ is a retract of $\pi_i(Y)$ for all $i$. Since $Y\in\mcD$, we know that $\pi_i(Y)$ is finitely generated and therefore $\pi_i(X)$ is finitely generated, so that $X\in\mcD$. We conclude that $\mcD$ is a thick subcategory of $\mcC$ containing $S$. So $D$ contains $\text{thick}(S)=\mcF$. So if $X\in\mcF$, then $\pi_i(X)$ is finitely generated over $R_0$ for all $i$. 

Now let $X,Y\in\mcF$. Then $$[X,Y]=[X[-1][1],Y]=[S^1\w X[-1],Y]=[S^1,F(X[-1],Y)]=\pi_1(F(X[-1],Y)).$$ Since $\mcF$ is closed under desuspensions and under $F$, we know that $F(X[-1],Y)\in\mcF$. Then $[X,Y]=\pi_1(F(X[-1],Y))$ is finitely generated over $R_0$ by the above argument, and we are done. \end{proof}

\begin{definition} Fix an injective hull $E_0$ of the $R_0$-module $R_0/\mfm$. Then the arguments of \S\ref{reflexive graded modules} apply to $(R_0,E_0)$ (see, for example, Section 47.7 of \cite{StP}). Since $E_0$ is injective and $\pi_0$ is homological on $\mcC$, we know that $\Hom_{R_0}(\pi_0(-),E_0)$ is cohomological on $\mcC$. Then since cohomological functors on $\mcC$ are representable, we can choose (and fix) an object $E\in\mcC$ and a natural isomorphism $$\eta:[-,E]\to\Hom_{R_0}(\pi_0(-),E_0).$$ For each $A\in\mcC$, we see that $$[A,F(F(A,E),E)]=[A\w F(A,E),E]=[F(A,E),F(A,E)]$$ and let $\sg_A:A\to F(F(A,E),E)$ be the morphism corresponding to $\id_{F(A,E)}$. \end{definition}

\begin{lem} \label{duality lemma for pistar inj} $\sg_A$ is an isomorphism for all $A\in\mcF$. \end{lem}

\begin{proof} Since $R_0$ is a complete local Noetherian ring, we know from Matlis duality (see, for example, Proposition 47.7.8 of the Stacks Project \cite{StP}) that the functor $\Hom_{R_0}(-,E_0)$ is an exact anti-equivalence between Noetherian and Artinian $R_0$-modules. Since $R_0$ is Noetherian, the Noetherian $R_0$-modules are simply the finitely generated ones. So the canonical map $$M\to\Hom_{R_0}(\Hom_{R_0}(M,E_0),E_0)$$ is an isomorphism for all finitely generated $R_0$-modules $M$ by Matlis duality.

Let $A\in\mcF$. Then $\pi_0(A)$ is a finitely generated $R_0$-module, so $$\pi_0(A)=\Hom_{R_0}(\Hom_{R_0}(\pi_0(A),E_0),E_0).$$ Then using the above natural transformation and the adjunction of $\w$ and $F$, we see that $$\Hom_{R_0}(\pi_0(A),E_0)=[A,E]=[S\w A,E]=[S,F(A,E)]=\pi_0(F(A,E)).$$ So $$\pi_0(A)=\Hom_{R_0}(\pi_0(F(A,E)),E_0)=[F(A,E),E]=[S\w F(A,E),E]=[S,F(F(A,E),E)],$$ which is equal to $$\pi_0(F(F(A,E)))$$ So $\pi_0(A)=\pi_0(F(F(A,E),E))$ and then $\pi_0(\sg_A)$ is an isomorphism. 

For each $i\in\Z$, since $A[-i]\in\mcF$ we know that $\pi_0(\sg_{A[-i]})$ is an isomorphism. Then from the axioms for a tensor triangulated category (see Definition A.2.1 of \cite{HPS97} for a list of the compatibility axioms), we see that $$F(F(A[-i],E),E)=F((F(A,E),E)[i])=F(F(A,E),E)[-i].$$ So $\pi_0(\sg_{A[-i]})=\pi_0(\sg_A[-i])=\pi_i(\sg_A)$, so that $\pi_i(\sg_A)$ is an isomorphism and finally $\pi_*(\sg_A)$ is an isomorphism. Then since $\mcC$ is monogenic, we know from part 3 of Lemma \ref{monogenic properties} that $\sg_A$ is an isomorphism. \end{proof}

\begin{definition} \cite[Definition 4.1.6]{HPS97} We say that a morphism $f:X\to Y$ in $\mcC$ is \textit{phantom} if for each $A\in\mcF$ and morphism $\vphi:A\to X$, we have that $f\vphi=0$. \end{definition}

\begin{definition} \cite[Definition 2.3.7]{HPS97} Let $H:\mcF\to\text{Ab}$ be an exact functor (i.e. it sends exact triangles in $\mcF$ to exact sequences of abelian groups). For each $X\in\mcC$, we construct a category $\Lmb(X)$ as follows. An object in $\Lmb(X)$ consists of an object $A\in\mcF$ and a morphism $A\to X$. A morphism in $\Lmb(X)$ from $A\to X$ to $B\to X$ is a morphism $A\to B$ in $\mcF$ making the following diagram commute. \begin{center}
    \begin{tikzcd} [sep=2cm]
     {} & X & {} \\
     A \arrow[rr] \arrow[ru] & {} & B \arrow[lu] \\ 
     \end{tikzcd}\end{center} We now define a functor $\Lmb(X)\to\text{Ab}$. It sends the object $A\to X$ in $\Lmb(X)$ to $H(A)$, and if we have a morphism $A\to B$ (going from $A\to X$ to $B\to X$ in $\Lmb(X)$), this is sent to $H(A\to B)$. This functor defines a diagram of abelian groups, and we write its colimit as $\widehat{H}(X)$. This leads to a functor $\widehat{H}:\mcC\to\text{Ab}$  extending $H$. \end{definition}
     
\begin{lem} \cite[Proposition 2.3.17]{HPS97} \label{homology iso lemma} Let $H$ be a homology functor on $\mcC$. Then the canonical map $\widehat{H}(X)\to H(X)$ is an isomorphism for all $X\in\mcC$. \end{lem}

\begin{lem} \label{no non-zero phantom maps} If $X\in\mcF$ and $\al:Y\to X$ is a phantom map, then $\al=0$. \end{lem}

\begin{proof} Let $T=F(X,E)$. Since $X\in\mcF$, we know from Lemma \ref{duality lemma for pistar inj} that $X=F(T,E)$, so that we have a phantom map $Y\to F(T,E)$. By the adjunction of $\w$ and $F$, our morphism $Y\to F(T,E)$ corresponds to a morphism $T\w Y\to E$. Then using the isomorphism $\sg_{T\w Y}$, this corresponds to an $R_0$-module homomorphism $\pi_0(T\w Y)\to E_0$. Since $\pi_0(T\w(-))$ is a homology functor on $\mcC$, we see from Lemma \ref{homology iso lemma} that $\pi_0(T\w Y)=\lim_{\Lambda(Y)}\pi_0(T\w Y_{\al})$. Since $Y_{\al}\in\mcF$ and $Y\to F(T,E)$ is a phantom map, we know that the composite $$Y_{\al}\to Y\to F(T,E)$$ is zero. Then $$\pi_0(T\w Y_{\al})\to\pi_0(T\w Y)\to\pi_0(E)=E_0$$ is zero for all $\al$, so that $\pi_0(T\w Y)\to E_0$ is zero. Then the map $T\w Y\to E$ is zero, so that the map $Y\to F(T,E)$ is zero, and finally the map $Y\to X$ is zero. \end{proof}

The following proposition is an axiomatic generalisation of \cite[Proposition 2.2]{Hov07}. 

\begin{prop} \cite[Proposition 2.2]{Hov07} \label{pi star injective} Let $\mcC$ be a Freyd category. Then $\pi_*(X)$ is injective as a graded $R$-module for all $X\in\mcF$. In particular, $R=\pi_*(S)$ is a self-injective graded ring. \end{prop}

\begin{proof} Let $X\in\mcF$ and $J$ be an injective hull of $\pi_*(X)$. Then the functor $\Hom(\pi_*(-),J)$ is cohomological, and therefore representable. So there exists $I\in\mcC$ and a natural isomorphism $$\eta:[-,I]\to\Hom_R(\pi_*(-),J).$$ Set $\psi=\eta_X^{-1}(\pi_*(X)\incl J):X\to I$. We can then choose $Y\in\mcC$ and form an exact triangle \begin{center}
    \begin{tikzcd} [sep=2cm]
     X \arrow[r,"\psi"] & 
     I \arrow[r] & 
     Y[1] \arrow[r] & 
     X[1] \;\;\;.
    \end{tikzcd}
 \end{center} Then \begin{center}
    \begin{tikzcd} [sep=2cm]
     Y \arrow[r] & 
     X \arrow[r,"\psi"] & 
     I \arrow[r] & 
     Y[1]
    \end{tikzcd}
 \end{center} is exact. Let $A\in\mcF$ and choose a morphism $A\to Y$. Since $Y\to X\to I$ is zero, the map $$A\to Y\to X\to I$$ is zero, and this is $[A,\psi](A\to Y\to X)$. So $A\to Y\to X$ is in the kernel of $[A,\psi]$. Let $f:A\to X$. Then since $\eta$ is a natural transformation, we see that \begin{center}
    \begin{tikzcd}[sep=3cm]
     {[X,I]} \arrow[r,"{[f,I]}"] \arrow[d,"\eta_X"] &
     {[A,I]} \arrow[d,"\eta_A"] \\
     {\Hom_{R}(\pi_{*}(X),J)} \arrow[r,"{\Hom_{R}(\pi_{*}(f),J)}"] & 
     {\Hom_{R}(\pi_{*}(A),J)}
    \end{tikzcd}
 \end{center} commutes. Then plugging in $\psi:X\to I$ into this diagram shows that the following diagram commutes. \begin{center}
    \begin{tikzcd}[sep=3cm]
     {[A,X]} \arrow[r,"{[A,\psi]}"] \arrow[d,"\pi_*"] &
     {[A,I]} \arrow[d,"\eta_A"] \\
     {\Hom_R(\pi_*(A),\pi_*(X))} \arrow[r,"{\Hom_R(\pi_*(A),\pi_*(X)\incl J)}"] & 
     {\Hom_R(\pi_*(A),J)}
    \end{tikzcd}
 \end{center} Since $A\in\mcF$ and $\mcC$ is a Freyd category, the vertical map $\pi_*$ is injective. Since $\io$ is a monomorphism, we know that $\al\mapsto\io\al$ is injective. So $[A,\psi]$ must be injective, giving that our map $A\to Y\to X$ is zero. We conclude that $Y\to X$ is a phantom map and therefore is zero by Lemma \ref{no non-zero phantom maps}. 

Then since \begin{center}
    \begin{tikzcd} [sep=2cm]
     Y \arrow[r] & 
     X \arrow[r,"\psi"] & 
     I \arrow[r] & 
     Y[1]
    \end{tikzcd}
 \end{center} is exact and $Y\to X$ is zero, we see from Lemma \ref{zeros in exact tr split monos epis} that $\psi$ is a split monomorphism, and then $\pi_*(\psi)$ is a split monomorphism. Since $$\pi_i(I)=[S^i,I]=\Hom_R(\pi_*(S^i),J)=\Hom_R(R(i),J)=\Hom^*(R,J)_i=J_i$$ for all $i\in\Z$, we know that $\pi_*(I)=J$ so that we have a split monomorphism $\pi_*(X)\to J$, and then $\pi_*(X)$ is a retract of $J$, which is injective. We conclude that $\pi_*(X)$ is injective. \end{proof}
 
This corollary is also an axiomatic generalisation of a result in \cite{Hov07}.
 
\begin{cor} \cite[Corollary 4.1]{Hov07} \label{pi star fully faithful} Let $\mcC$ be a Freyd category. Then $\pi_*:\mcF\to\Mod_R$ is fully faithful, and therefore an embedding of categories. \end{cor}

\begin{proof} Let $X,Y\in\mcF$. Then $\pi_*(Y)$ is injective by Proposition \ref{pi star injective}, so $$[X,Y]\to\Hom(\pi_*(X),\pi_*(Y))$$ is an isomorphism by Lemma \ref{pi star gives isos}. \end{proof}
 
\begin{lem} \cite[Lemma 12.3]{ShS14} Let \begin{center}
    \begin{tikzcd} [sep=2cm]
     X \arrow[r,"u"] & 
     Y \arrow[r,"v"] & 
     Z \arrow[r,"w"] & 
     X[1]
    \end{tikzcd}
 \end{center} be exact in $\mcF$. Then \begin{center}
    \begin{tikzcd} [sep=1.2cm]
    	 \cdots \arrow[r] &
     \pi_*(X) \arrow[r,"\pi_*(u)"] & 
     \pi_*(Y) \arrow[r,"\pi_*(v)"] & 
     \pi_*(Z) \arrow[r,"\pi_*(w)"] & 
     \pi_*(X[1]) \arrow[r] &
     \cdots
    \end{tikzcd}
 \end{center} is exact in $\Mod_R$. \end{lem}

\begin{proof} Let $i\in\Z$. Since \begin{center}
    \begin{tikzcd} [sep=2cm]
     X \arrow[r,"u"] & 
     Y \arrow[r,"v"] & 
     Z \arrow[r,"w"] & 
     X[1]
    \end{tikzcd}
 \end{center} is exact in $\mcF$ and $\pi_i$ is homological, we obtain the long exact sequence \begin{center}
    \begin{tikzcd} [sep=1.2cm]
    	 \cdots \arrow[r] &
     \pi_i(X) \arrow[r,"\pi_i(u)"] & 
     \pi_i(Y) \arrow[r,"\pi_i(v)"] & 
     \pi_i(Z) \arrow[r,"\pi_i(w)"] & 
     \pi_i(X[1]) \arrow[r] &
     \cdots \;\;\;.
    \end{tikzcd}
 \end{center} Then \begin{center}
    \begin{tikzcd} [sep=1.2cm]
    	 \cdots \arrow[r] &
     \pi_*(X) \arrow[r,"\pi_*(u)"] & 
     \pi_*(Y) \arrow[r,"\pi_*(v)"] & 
     \pi_*(Z) \arrow[r,"\pi_*(w)"] & 
     \pi_*(X[1]) \arrow[r] &
     \cdots \;\;\;.
    \end{tikzcd}
 \end{center} is exact in $\Mod_R$. \end{proof}

\begin{cor} \cite[Corollary 12.5]{ShS14} \label{fin gen implies free} Suppose that $\mcC$ is a Freyd category. If $X\in\mcF$ and $\pi_*(X)$ is finitely generated over $R$, then it is free. \end{cor}

\begin{proof} There exists a free and finitely generated graded $R$-module $F$ and an epimorphism $\al:F\to\pi_*(X)$. Then $F$ is isomorphic to $\pi_*(Y)$, where $Y$ is a finite wedge of spheres. So $Y\in\mcF$, and we have a morphism $\pi_*(Y)\to\pi_*(X)$. Then since $\mcC$ is a Freyd category, this is $\pi_*(f)$ for some $f:Y\to X$ in $\mcF$. Then $f$ is a split epimorphism, so $\al$ must be a split epimorphism. So $\pi_*(X)$ is a retract of $F$, which is free and therefore projective. So $\pi_*(X)$ must be projective. Since $R$ is local, $\pi_*(X)$ must therefore be free. \end{proof}

\begin{rmk} If $\mcC$ is a Freyd category, then since $R$ is self-injective, the arguments of \S\ref{reflexive graded modules} apply to $R$ in this case. We thus adopt the notation $D=\Hom^*(-,R)$ from that section. \end{rmk}

\begin{prop} \label{comparing d functors in f} Suppose that $\mcC$ is a Freyd category. Then $$\pi_*(DX)\simeq D(\pi_*(X))$$ for each $X\in\mcF$. \end{prop}

\begin{proof} Let $i\in\Z$. Then $$\pi_i(DX)=[S^i,DX]\simeq[D^2X,D(S^i)]\simeq[X,S^{-i}]\simeq\Hom(\pi_*(X),\pi_*(S^{-i})).$$ Then since $\pi_*(S^{-i})\simeq R(-i)$ (Lemma \ref{pi star sends sphere to shifted ring}), we have that $$\pi_i(DX)\simeq\Hom(\pi_*(X),R(-i))=\Hom(\pi_*(X)(i),R)=D(\pi_*(X))_i.$$ We conclude that $$\pi_*(DX)\simeq D(\pi_*(X)).$$ \end{proof}

\begin{prop} Let $M$ be a graded $R$-module. If $M$ is isomorphic to $\pi_*(X)$ for some $X\in\mcF$ and $DM=0$, then $M=0$. \end{prop}

\begin{proof} Suppose that $M$ is isomorphic to $\pi_*(X)$ for some $X\in\mcF$, and that $DM=0$. Then $\pi_*(DX)\simeq D(\pi_*(X))$ is isomorphic to $DM=0$ (by Proposition \ref{comparing d functors in f}), so that $DX=0$ (by part 2 of Lemma \ref{monogenic properties}) and then $X=0$. So $M$ is isomorphic to $\pi_*(X)=0$. \end{proof}

\begin{prop} $\im(\pi_*(u))\in\mcM$ for all $u\in\Fr(\mcF)$. \end{prop}

\begin{proof} Let $A\in\mcF$. We use Proposition \ref{comparing d functors in f} to see that $$\pi_*(A)\simeq\pi_*(D^2A)\simeq D(\pi_*(DA))\simeq D^2(\pi_*(A)),$$ so that $\pi_*(A)\in\mcM$. Now let $u:X\to Y$ be an object in $\Fr(\mcF)$. Then $\pi_*(X),\pi_*(Y)\in\mcM$, so that $\im(\pi_*(u))\in\mcM$. \end{proof}

\begin{conj} $\Fr(\mcF)\to\mcM$ is an equivalence of categories. \end{conj}

\begin{lem} (Graded Nakayama Lemma) \cite[Lemma 10.20.1]{StP} Let $M$ be a finitely generated graded $R$-module and $\{x_i\}$ be a finite set of homogeneous elements in $M$. Then $\{x_i\}$ is a minimal generating set for $M$ over $R$ if and only if $\{x_i+\mfm M\}$ is a basis for $M/\mfm M$ over $R/\mfm$. \end{lem}

\begin{proof} See Lemma 10.20.1 of the Stacks Project for a proof of the ungraded case. \end{proof}

\begin{lem} \cite[Section 4D]{Lam99} \label{all kernels fg for fin pres} Let $M$ be a graded $R$-module, and suppose that $F$ and $F'$ are free and finitely generated. Suppose that $F\to M$ and $F'\to M$ are epimorphisms, with $\ker(F\to M)$ finitely generated. Then $\ker(F'\to M)$ is finitely generated. \end{lem}
 
\begin{proof} Since $F$ is free, it is projective. Then since $F'\to M$ is an epimorphism, there exists a morphism $f:F\to F'$ fitting into the following diagram. \begin{center}
    \begin{tikzcd} [sep=2cm]
     0 \arrow[r] \arrow[d] & 
     \ker(F\to M) \arrow[r] \arrow[d,"f",dashed] & 
     F \arrow[r,"w"] \arrow[d,"g",dashed] & 
     M \arrow[r] \arrow[d,equal] &
     0 \arrow[d] \\
     0 \arrow[r] &
     \ker(F'\to M) \arrow[r] &
     F' \arrow[r,"w'"] &
     M \arrow[r] &
     0 
    \end{tikzcd}
 \end{center} Since $g(\ker(F\to M))\leq\ker(F'\to M)$, we can form a morphism $$f:\ker(F\to M)\to\ker(F'\to M)$$ fitting into the above diagram. Then \begin{center}
    \begin{tikzcd} [sep=1.5cm]
     0 \arrow[r] & 
     \im(f) \arrow[r] & 
     \ker(F'\to M) \arrow[r] & 
     \cok(f) \arrow[r]  &
     0
    \end{tikzcd}
 \end{center} is exact. Since $\ker(F\to M)$ is finitely generated, we know that $\im(f)$ is. Since $\id_M$ is an isomorphism, it follows from the Snake Lemma that $\cok(f)=\cok(g)$ is finitely generated, since $F'$ is. So $\ker(F'\to M)$ is an extension of finitely generated modules and therefore finitely generated. \end{proof}
 
\begin{definition} \cite[Definition 1.2]{ShS14} Let $R$ be a graded-commutative ring. We say that a graded $R$-module $M$ is \textit{finitely presented} if there exists a free and finitely generated graded $R$-module $F$ and an epimorphism $F\to M$ for which $\ker(F\to M)$ is finitely generated. We say that a graded-commutative ring $R$ is \textit{coherent} if all of its finitely generated graded ideals are finitely presented. In contrast, it is \textit{totally incoherent} if its only finitely presented graded ideals are $0$ and $R$. \end{definition}

\begin{lem} \cite[Proposition 12.6]{ShS14} \label{fin pres condition} Let $M$ be a graded $R$-module with $M_i=0$ for $i<0$. Then $M$ is finitely presented if and only if there exist finite wedges $X$ and $Y$ of non-negative spheres, and a morphism $f:X\to Y$, where $\cok(\pi_*(f))$ is isomorphic to $M$. \end{lem}
 
\begin{proof} First suppose that $M$ is finitely presented. Then there exists a free and finitely generated graded $R$-module $F$ and an epimorphism $\pi:F\to M$ for which $\ker(\pi)$ is finitely generated. Then Lemma \ref{all kernels fg for fin pres} shows that we may assume without loss of generality that $F$ maps a basis $\{e_i\}$ to a minimal generating set for $M$ over $R$. Since $M_i=0$ for $i<0$, the elements of the minimal generating set are homogeneous of non-negative degree, so that the basis elements in $F$ are homogeneous of non-negative degree. Then since $R_i=0$ for $i<0$, we see that $F_i=0$ for $i<0$. Then $\ker(\pi)_i=0$ for $i<0$ with $\ker(\pi)$ finitely generated, so we can similarly find a free and finitely generated graded $R$-module $F'$ and an epimorphism $\pi':F'\to\ker(\pi)$ with $(F')_i=0$ for $i<0$. 

We then form the composite $F'\to\ker(\pi)\to F$. Since every free and finitely generated graded $R$ module is isomorphic to a finite sum of shifted copies of $R$, $R(n)=\pi_*(S^n)$ for all $n$, and $\pi_*$ is additive, we see that $F$ and $F'$ must be isomorphic to $\pi_*(X)$ and $\pi_*(Y)$, where $X$ and $Y$ are wedges of spheres. Then since $F_i=(F')_i=0$ for $i<0$, these must all be non-negative spheres. Since $S\in\mcF$, and $\mcF$ is closed under finite smash products and coproducts, we see that $X,Y\in\mcF$. It then follows from Lemma \ref{finite wedge of spheres gives iso} that our morphism $F'\to F$ is $\pi_*(f)$ for some $f:X\to Y$ with $\im(\pi_*(f))=\ker(\pi)$. Then $M$ is isomorphic to $\cok(\pi_*(f))$. 

Conversely, if $M$ is isomorphic to $\cok(\pi_*(f))$, where $f:X\to Y$ is a morphism of finite wedges of non-negative spheres, then $\pi_*(Y)$ is free and finitely generated, $\pi_*(Y)\to M$ is an epimorphism and $\ker(\pi_*(Y)\to M)=\im(\pi_*(f))$, which is finitely generated since $\pi_*(X)$ is. So $M$ is finitely presented. \end{proof}

\begin{prop} \cite[Proposition 12.6]{ShS14} If $\mcC$ is a Freyd category, then $R$ is totally incoherent. \end{prop}

\begin{proof} Let $J$ be a finitely presented ideal in $R$. Since $R_i=0$ for $i<0$, we know that $J_i=0$ for $i<0$. Then we see from Lemma \ref{fin pres condition} that there exist finite wedges $X$ and $Y$ of non-negative spheres, and a morphism $f:X\to Y$, for which \begin{center}
    \begin{tikzcd} [sep=2cm]
     \pi_*(X) \arrow[r,"\pi_*(f)"] & 
     \pi_*(Y) \arrow[r] & 
     J \arrow[r] & 
     0
    \end{tikzcd}
 \end{center} is exact. We form the composite $\pi_*(Y)\to J\incl R=\pi_*(S)$ and then, using the fact that $\mcC$ is a Freyd category with $Y,S\in\mcF$, we can find $g:Y\to S$ for which $\pi_*(g):\pi_*(Y)\to R$ is the above composite, with \begin{center}
    \begin{tikzcd} [sep=2cm]
     \pi_*(X) \arrow[r,"\pi_*(f)"] & 
     \pi_*(Y) \arrow[r,"\pi_*(g)"] & 
     \pi_*(S) 
    \end{tikzcd}
 \end{center} exact and $\im(\pi_*(g))=J$. Since $X,Y\in\mcF$ and $\mcC$ is a Freyd category, we see that $gf=0$.

There exists an object $K\in\mcC$ for which \begin{center}
    \begin{tikzcd} [sep=2cm]
     Y \arrow[r] & 
     S \arrow[r] & 
     K[1] \arrow[r] & 
     Y[1]
    \end{tikzcd}
 \end{center} is exact. Then \begin{center}
    \begin{tikzcd} [sep=2cm]
     K \arrow[r] & 
     Y \arrow[r] & 
     S \arrow[r] & 
     K[1]
    \end{tikzcd}
 \end{center} is exact, so that \begin{center}
    \begin{tikzcd} [sep=2cm]
     S^{-1} \arrow[r,"d"] & 
     K \arrow[r,"i"] & 
     Y \arrow[r,"g"] & 
     S
    \end{tikzcd}
 \end{center} is exact. Since $gf=0$ and $i$ is a weak kernel of $g$, there exists a morphism $\td{f}:X\to K$ with $i\td{f}=f$. \begin{center}
    \begin{tikzcd} [sep=2cm]
    	 {} & X \arrow[rd,"f"] \arrow[d,"\td{f}",dashed,swap] & {} & {} \\
     S^{-1} \arrow[r,"d"] & 
     K \arrow[r,"i"] & 
     Y \arrow[r,"g"] & 
     S
    \end{tikzcd}
 \end{center} From $\td{f}$ and $d$, we form a morphism $\ta:X\vee S^{-1}\to K$. Apply $\pi_j$ to \begin{center}
    \begin{tikzcd} [sep=2cm]
     S^{-1} \arrow[r,"d"] & 
     K \arrow[r,"i"] & 
     Y \arrow[r,"g"] & 
     S
    \end{tikzcd}
 \end{center} to see that \begin{center}
    \begin{tikzcd} [sep=1.2cm]
    	 \cdots \arrow[r] &
     R_{j+1} \arrow[r,"\pi_j(d)"] & 
     \pi_j(K) \arrow[r,"\pi_j(i)"] & 
     \pi_j(Y) \arrow[r,"\pi_j(g)"] & 
     R_j \arrow[r] &
     \cdots 
    \end{tikzcd}
 \end{center} is exact. Then one can check that $\pi_j(\theta)$ is surjective, so that $\pi_*(K)=\im(\pi_*(\theta))$. Then since $\pi_*(X\vee S^{-1})$ is free and finitely generated, we know that $\pi_*(K)$ is finitely generated. Since $\mcF$ is thick and \begin{center}
    \begin{tikzcd} [sep=2cm]
     K \arrow[r,"i"] & 
     Y \arrow[r,"g"] & 
     S \arrow[r,"-d{[1]}"] & 
     K[1]
    \end{tikzcd}
 \end{center} is exact, we know that $K\in\mcF$. So $K\in\mcF$ with $\pi_*(K)$ finitely generated over $R$, therefore it must be free by Lemma \ref{fin gen implies free}. So $\pi_*(K)$ must be isomorphic to $\pi_*(Z)$ where $Z$ is a finite wedge of spheres. Since $Z,K\in\mcF$ and $\mcC$ is a Freyd category, there must be a unique morphism $Z\to K$ corresponding to this isomorphism. Then since $\pi_*(Z\to K)$ is an isomorphism, it follows from part 3 of Lemma \ref{monogenic properties} that this morphism $Z\to K$ must itself be an isomorphism. By appropriately adjusting the morphisms, we can thus replace $K$ with $Z$. \begin{center}
    \begin{tikzcd} [sep=2cm]
     S^{-1} \arrow[r,"d"] & 
     Z \arrow[r,"i"] & 
     Y \arrow[r,"g"] & 
     S
    \end{tikzcd}
 \end{center} Since $$\ker(\pi_j(d))=\im(\pi_{j+1}(g))=J_{j+1},$$ we obtain a monomorphism $R_{j+1}/J_{j+1}\to\pi_j(Z)$ and an exact sequence \begin{center}
    \begin{tikzcd} [sep=1.2cm]
    	 0 \arrow[r] &
     R_{j+1}/J_{j+1} \arrow[r] & 
     \pi_j(Z) \arrow[r,"\pi_j(i)"] & 
     \pi_j(Y) \arrow[r] & 
     J_j \arrow[r] &
     0 \;\;\;.
    \end{tikzcd}
 \end{center} So if $j<-1$, then $R_{j+1}=\pi_j(Y)=0$ and then $\pi_j(Z)=0$. So each sphere in the wedge sum $Z$ is a $k$-sphere for some $k\geq -1$. Let $e_i:S^{|e_i|}\to Z$ ($1\leq i\leq n$) be a homogeneous basis for $\pi_*(Z)$ over $R$. Since $R_i=0$ for $i<0$, if $|e_j|\geq 0$, then $e_j$ does not contribute to $\pi_{-1}(Z)$. If $|e_j|=-1$, then only $R_0$-multiples count. So $\pi_{-1}(Z)$ is a free $R_0$-module with those $e_i$ with $|e_i|=-1$ as a basis. The above exact sequence gives an isomorphism $R_0/J_0\to\pi_{-1}(Z)$. 

If $|e_i|\geq 0$ for all $i$, then $\pi_{-1}(Z)=0$ so that $J_0=R_0$. Then $1\in J_0$, so that $J=R$. Otherwise, $\pi_{-1}(Z)$ is a non-zero free $R_0$-module and then $\ann(\pi_{-1}(Z))=0$. Yet $J_0(R_0/J_0)=0$ so that $J_0\pi_{-1}(Z)=0$ and then $J_0\leq\ann(\pi_{-1}(Z))=0$. So $R_0=R_0/J_0$ is isomorphic to $\pi_{-1}(Z)$, giving that $\pi_{-1}(Z)$ has rank 1 over $R_0$. So there is a unique $i$ with $|e_i|=-1$ forming a basis for $\pi_{-1}(Z)$ over $R_0$, which we label $e_0$. Then $\pi_{-1}(d):R_0\to\pi_{-1}(Z)$ is an $R_0$-module isomorphism, so $\pi_{-1}(d)(1)=ae_0$ for some unique $a\in R_0$. 

If $a\in\mfm$, then for each $b\in R_0$ we have that $\pi_{-1}(d)(b)=b\pi_{-1}(d)(1)=baw_0\in\mfm w_0$. So $\im(\pi_{-1}(d))\leq\mfm w_0<\pi_{-1}(Z)$. Yet $\pi_{-1}(d)$ is an isomorphism, so $\im(\pi_{-1}(d))=\pi_{-1}(Z)$, a contradiction. So $a\notin\mfm$ and since $R_0$ is local, we must have $a\in R_0^\x$. To define a map $\pi_*(Z)\to\pi_*(S^{-1})$, we must choose where to send each $e_i$ in $\pi_*(Z)$ in a way that respects the gradings. Send $e_0$ to $a^{-1}\id_{S^{-1}}$ and $e_i$ to $0$ for all other $i$. Since $\mcC$ is a Freyd category and $Z,S^{-1}\in\mcF$, this gives a morphism $r:Z\to S^{-1}$ with $re_0=a^{-1}\id_{S^{-1}}$ and $re_i=0$ for all other $i$. Then $$\pi_*(rd)(\id_{S^{-1}})=\pi_*(r)(\pi_{-1}(d)(1))=\pi_{-1}(r)(ae_0)=are_0=aa^{-1}\id_{S^{-1}}=\id_{S^{-1}}.$$ So $\pi_*(rd)=\pi_*(\id_{S^{-1}})$ and then $rd=\id_{S^{-1}}$ since $S^{-1}\in\mcF$ and $\mcC$ is a Freyd category. We conclude that $d$ is a split monomorphism. Since $dg[-1]=0$ with $d$ a split monomorphism, we see that $g=g[-1][1]=0$ so that $J=\im(\pi_*(g))=0$. So the only finitely presented graded ideals in $R$ are $0$ and $R$, and we are done. \end{proof}

\section{The Spanier-Whitehead Category} \label{sw cat section}

We now construct the Spanier-Whitehead category $\mcF$. This section is largely based on \cite{Str20}. We write $I=[0,1]$. Given finite pointed CW complexes $X$ and $Y$, we write the reduced suspension of $X$ as $\Sg X$ and the pointed homotopy class of $X$ and $Y$ as $[X,Y]$.

Let $X$ and $Y$ be finite pointed CW complexes. We start by describing the group structure on $[\Sg X,Y]$. Identify $0$ and $1/2$ in $S^1$ to obtain the wedge sum $S^1\vee S^1$, and thus a quotient map $S^1\to S^1\vee S^1$. Taking the smash product of this map with $\id_X$, and using the distributivity of the smash product over the wedge sum, gives a map $\Sg X\to\Sg X\vee\Sg X$. There is an obvious map $Y\vee Y\to Y$. Given two maps $f,g:\Sg X\to Y$, we then form a map $f\vee g:\Sg X\vee\Sg X\to Y\vee Y$. Combining these together gives a map \begin{center}
    \begin{tikzcd} [sep=2cm]
     \Sg X \arrow[r] & 
     \Sg X\vee\Sg X \arrow[r,"f\vee g"] & 
     Y\vee Y \arrow[r] & 
     Y
    \end{tikzcd}
 \end{center} One can check that this construction respects homotopy, and that the resulting binary operation on $[\Sg X,Y]$ gives it the structure of a group. 

Since the reduced suspension of a finite pointed CW complex is again a finite pointed CW complex, for $n\geq 1$ we may replace $X$ with $\Sg^{n-1}X$ and $Y$ with $\Sg^n Y$ above; the resulting group structure on $[\Sg^n X,\Sg^n Y]$ is abelian for $n\geq 2$.

Taking the smash product of a map $f:X\to Y$ with $\id_{S^1}$ gives a map $\Sg f:\Sg X\to\Sg Y$. This construction also respects homotopy, so we obtain a map $[X,Y]\to[\Sg X,\Sg Y]$, and therefore maps $[\Sg^n X,\Sg^n Y]\to[\Sg^{n+1} X,\Sg^{n+1}Y]$ for all $n\geq 2$, all of which are group homomorphisms. We have thus constructed a diagram in the category of abelian groups; we denote its colimit by $[\Sg^\infty X,\Sg^\infty Y]$. We now use Freudenthal's Suspension Theorem to show that this colimit is attained at a finite stage.

\begin{theorem} (Freudenthal Suspension Theorem) \cite[Corollary 3.2.3]{Koc96} Let $X$ and $Y$ be finite pointed CW complexes and $n\geq 0$, where $X$ is $n$-connected and $\dim(Y)\leq 2n$. Then the map $$[X,Y]\to[\Sg X,\Sg Y]$$ induced by the suspension functor is an isomorphism. \end{theorem}

\begin{proof} See Corollary 3.2.3 of \cite{Koc96}. \end{proof}
 
\begin{cor} \label{freudenthal corollary} Let $X$ and $Y$ be finite pointed CW complexes. Then the group homomorphisms $$[\Sg^n X,\Sg^n Y]\to[\Sg^{n+1}X,\Sg^{n+1}Y]$$ eventually become isomorphisms. Thus we may identify $[\Sg^\infty X,\Sg^\infty Y]$ with $[\Sg^n X,\Sg^n Y]$ when $n$ is sufficiently large. \end{cor}

\begin{proof} Let $k\geq\dim(Y)+2$. Then $\Sg^k X$ and $\Sg^k Y$ are finite pointed CW complexes, where $\Sg^k X$ is $(k-1)$-connected and $\dim(\Sg^k Y)=k+\dim(Y)\leq 2(k-1)$. So $$[\Sg^k X,\Sg^k Y]\to[\Sg^{k+1} X,\Sg^{k+1} Y]$$ is an isomorphism.\end{proof}

Now choose $m,n\in\Z$. For sufficiently large $k$, we have that $\Sg^{k+m}X$ and $\Sg^{k+n} Y$ are finite pointed CW complexes, so that $[\Sg^\infty\Sg^{k+m}X,\Sg^\infty\Sg^{k+n} Y]$ is an abelian group, and is independent of $k$ by Corollary \ref{freudenthal corollary}; we label this $[\Sg^{\infty+m}X,\Sg^{\infty+n}Y]$.

The objects in the Spanier-Whitehead category $\mcF$ are pairs $(X,n)$, where $X$ is a finite pointed CW complex and $n\in\Z$. The morphisms from $(X,m)$ to $(Y,n)$ are given by $$\Hom_{\mcF}((X,m),(Y,n))=[\Sg^{\infty+m}X,\Sg^{\infty+n}Y].$$ To define composition in $\mcF$, we find sufficiently large $k$ and use the composition in the homotopy category of pointed spaces. The identity on $(X,n)$ is obtained from the identity on $\Sg^{k+n}X$ for sufficiently large $k$. The zero object is the zero space, together with $0\in\Z$. The wedge sum generalises to $\mcF$ by setting $$(X,m)\vee(Y,n)=(X\vee Y,m+n).$$ 

Given an object $(X,n)$, there is an obvious identification of $(\Sg X,n)$ with $(X,n+1)$. We define this object to be the reduced suspension $\Sg((X,n))$ of $(X,n)$ in $\mcF$, and we obtain an additive automorphism $\Sg$ on $\mcF$.

We now describe the canonical triangulated structure on $\mcF$. Let $(B,n)$ be an object in $\mcF$, and choose a subcomplex $A$ of $B$. Consider $A$ as a subspace of $I\w A$ via the map $a\mapsto[(1,a)]$. We then form an adjunction space $C$ from the inclusion $A\incl B$ and a canonical map $B\to C$. Quotienting out $B$ from $C$ identifies $(1,a)$ with $(0,a)$ for each $a\in A$, so we may identify $C/B$ with $S^1\w A=\Sg A$ and obtain a map $C\to\Sg A$. Applying $\Sg^{k+n}$ to these maps for sufficiently large $k$, we obtain a triangle \begin{center}
    \begin{tikzcd} [sep=2cm]
     (A,n) \arrow[r] & 
     (B,n) \arrow[r] & 
     (C,n) \arrow[r] & 
     \Sg((A,n))
    \end{tikzcd}
 \end{center} in $\mcF$.

\begin{prop} The class of all triangles in $\mcF$ that are isomorphic to one of the above form satisfies the axioms of a triangulated structure compatible with $\Sg$. \end{prop}

\begin{proof} See Theorem 7 of \cite{Mar83}. \end{proof}

Let $S=(S^0,0)$. In Definition \ref{construction of stable homotopy ring}, we saw how to construct the associated stable homotopy ring $R$ and functor $\pi_*:\mcF\to\Mod_R$ for a tensor triangulated category, so shall not repeat it here. We can now state Freyd's Generating Hypothesis (for $\mcF$).

\begin{conj} \label{original freyd gen hypo} \cite[\S 9]{Fre65} (Freyd's Generating Hypothesis) $\pi_*:\mcF\to\Mod_R$ is faithful. \end{conj}

\begin{rmk} We can adjust the definition of $\mcF$ slightly, tensoring the hom-sets with $\Z_p$, to obtain the \textit{$p$-completed Spanier-Whitehead category}. However, it is not obvious that this modified category has a canonical triangulation. \end{rmk}

\begin{theorem} There exists a monogenic stable homotopy category $\mcC$ whose full subcategory of strongly dualisable objects is $\mcF$. \end{theorem}

\begin{proof} See \cite{Mar83}. \end{proof}

\begin{rmk} \label{boardman cat remark} See \cite{Boa65} for details on a specific example of such a monogenic stable homotopy category, known as Boardman's stable homotopy category. \end{rmk}

\section{$p$-completion} \label{p completion}

In this section, we follow \cite{Str24} to establish that the $p$-completed Spanier-Whitehead category $\mcF_p$ is triangulated, and can be identified with the full subcategory of strongly dualisable objects in a semi-Freyd category. Let $\mcB$ be Boardman's stable homotopy category (see Remark \ref{boardman cat remark}), which contains the Spanier-Whitehead category $\mcF$ as its full subcategory of strongly dualisable objects. Also fix a prime $p$. Write $$\Z[1/p]=\left\{\frac{a}{p^k}:a\in\Z,k\geq 0\right\}$$ and $\Z/p^\infty$ for the quotient $\Z[1/p]/\Z$, which is the direct limit of \begin{center}
    \begin{tikzcd} [sep=1.5cm]
     0 \arrow[r] & 
     \Z/p \arrow[r] & 
     \Z/p^2 \arrow[r] & 
     \cdots \;\;\;,
    \end{tikzcd}
 \end{center} so that we have a short exact sequence \begin{center}
    \begin{tikzcd}
     0 \arrow[r] & 
     \Z \arrow[r] & 
     \Z{[1/p]} \arrow[r] & 
     \Z/p^{\infty} \arrow[r] &
     0 \;\;\;.
    \end{tikzcd}
 \end{center}

\begin{definition} \cite[Definitions 5.1(a) and 5.9(a)]{Str24} We say that an abelian group $A$ is \textit{Ext-$p$-complete} if $$\Hom(\Z[1/p],A)=\Ext(\Z[1/p],A)=0.$$ We say that an object $X\in\mcB$ is \textit{$p$-complete} if $\pi_n(X)=\Hom(S^n,X)$ is Ext-$p$-complete for all $n\in\Z$. Denote the full subcategory of $\mcB$ consisting of the $p$-complete objects by $\mcB_p$. This is a stable homotopy category but not monogenic (Remark 5.15 and Remark 6.2 of \cite{Str24}). Write $\mcB_p^d$ for the full subcategory of strongly dualisable objects in $\mcB_p$. \end{definition}

\begin{definition} \cite[Remark 3.5]{Str24} Let $A$ be an abelian group. Then there exist free abelian groups $P_0$ and $Q_0$ and maps $f_0:P_0\to Q_0$, $g_0:Q_0\to A$ with \begin{center}
    \begin{tikzcd} [sep=1.5cm]
     0 \arrow[r] & 
     P_0 \arrow[r,"f_0"] & 
     Q_0 \arrow[r,"g_0"] & 
     A \arrow[r] &
     0
    \end{tikzcd}
 \end{center} exact. Then $P_0=\bigoplus_{i\in I}\Z$ and $Q_0=\bigoplus_{j\in J}\Z$ for some $I,J$. Set $P=\bigvee_{i\in I}S$ and $Q=\bigvee_{j\in J}S$, so that $\pi_0(P)=P_0$ and $\pi_0(Q)=Q_0$. Then there exists a unique morphism $f:P\to Q$ with $\pi_0(f)=f_0$. Then there exists an object $M$, unique up to non-unique isomorphism, for which \begin{center}
    \begin{tikzcd} [sep=2cm]
     P \arrow[r,"f"] & 
     Q \arrow[r] & 
     M \arrow[r] & 
     \Sg P
    \end{tikzcd}
 \end{center} is exact. We write $SA=M$ and call $SA$ a \textit{Moore spectrum} of $A$. Write $S/p^\infty$ for a Moore spectrum of $\Z/p^\infty$ and $S^{-1}/p^\infty$ for $\Sg^{-1}(S/p^\infty)$. (Construction 3.13 of \cite{Str24}). We define the \textit{$p$-completion} $X_p$ of $X\in\mcB$ to be $F(S^{-1}/p^\infty,X)\in\mcB$, where $F$ is the internal Hom in $\mcB$. (Definition 5.11 of \cite{Str24}). \end{definition}

\begin{definition} \cite[Definition 6.1]{Str24} A \textit{strict module category} for $S_p$ consists of a monogenic stable homotopy category $\mcM_p$, an exact, strongly symmetric monoidal functor $P:\mcB\to\mcM_p$ and a right adjoint $U:\mcM_p\to\mcB$ such that $UP$ is naturally isomorphic to $S_p\w(-)$.  \end{definition}

\begin{theorem} \cite[Corollary 5.33 and Proposition 6.17]{Str24} Let $\mcM_p$ be a strict module category for $S_p$. The $p$-completed Spanier-Whitehead category $\mcF_p$, the strongly dualisable objects in $\mcB_p$, and the strongly dualisable objects in $\mcM_p$ are all equivalent. \end{theorem}

\begin{cor} \cite[Proposition 6.17]{Str24} $\mcF_p$ is triangulated, and can be identified with the full subcategory of strongly dualisable objects in a semi-Freyd category. \end{cor}

\begin{proof} $\mcB_p$ is a tensor triangulated category, so the full subcategory of its strongly dualisable objects is also a triangulated category. Since $\mcF_p$ is equivalent to this, it must be a triangulated category. Since $\mcM_p$ is a semi-Freyd category, $\mcF_p$ is equivalent to the full subcategory of strongly dualisable objects in a semi-Freyd category. \end{proof}

As demonstrated in Proposition 2.4 of \cite{Hov07}, the following conjecture is equivalent to the original statement of Freyd's Generating Hypothesis (Conjecture \ref{original freyd gen hypo}). 

\begin{conj} (Freyd's Generating Hypothesis for $\mcF_p$) The functor $\pi_*:\mcF_p\to\Mod_R$ is faithful (where $R$ is the stable homotopy ring associated with the semi-Freyd category $\mcM_p$ containing $\mcF_p$). Therefore, $\mcM_p$ is a Freyd category, with $\mcF_p$ as its full subcategory of strongly dualisable objects. \end{conj}

\newpage

\chapter{Multibasic Modules Over The Hahn Ring} \label{hahn ring chapter}

In Chapter \ref{freyd cats chapter}, we explored various algebraic properties of Freyd and semi-Freyd categories, particularly relating to the associated stable homotopy rings. In Chapter \ref{infinite root algebra chapter}, we will investigate the infinite root algebra $P$, a self-injective ring with no non-trivial finitely presented ideals, and how closely it resembles a stable homotopy ring for a Freyd category $\mcC$. We will see that its category of reflexive modules is abelian and has enough projectives and enough injectives, and obtain a complete description of the objects in this category up to isomorphism. We will also see that despite this, there can be no suitable embedding of a triangulated category into the injective reflexive $P$-modules. 

Since $P$ is not an integral domain, and it is a quotient of a ring $A$, we instead use this chapter and the next to develop the theory of reflexive modules over the Hahn ring $A$, with the vast majority of results being trivial to translate into corresponding results about $P$-modules. We focus on the ungraded context throughout the next three chapters, however one would expect that the results would generalise easily to the graded context when regarding $A$ and $P$ as trivially graded graded-commutative rings. 

The Hahn ring $A$ consists of the Hahn series with value group $\R$ and residue field being the field of order 2, that contain no terms with negative exponents. The definition of a Hahn series is due to Hans Hahn \cite{Hah95}. After giving some basic definitions, we shall construct the field $K$ of Hahn series and the Hahn ring $A$. We give criteria for flatness and injectivity of $A$-modules. Then we explore basic modules, and eventually multibasic modules, giving a classification of the $A$-submodules of a $K$-vector space. We finish the chapter by establishing that the decomposition of a multibasic module into a finite list of basic summands is unique up to the order of the summands. 

\section{Introduction}

We start by outlining some fundamental properties of well-ordered sets of real numbers, and use these to construct the Hahn field $K$ and the associated Hahn ring $A$. The below construction of the field of Hahn series $K$ mainly follows \cite{Pas77}.

\begin{definition} Denote the field of order $2$ by $\F$. We define the \textit{support} of a map $a:\R\to\F$ to be $$\supp(a)=\{q\in\R:a(q)=1\}.$$ Recall that a subset $A\subs\R$ is said to be \textit{well-ordered} if every non-empty subset $B$ of $A$ has a minimum $\min(B)$, or equivalently (see Lemma \ref{equiv condition for well-ordered} below), if every weakly decreasing sequence ($m\leq n$ implies $a_m\geq a_n$) in $A$ is eventually constant. Denote the set of all well-ordered subsets of $\R$ by $\mcW$. One can check that $\mcW$ is closed under subsets and finite unions.

A \textit{Hahn series} is a map $a:\R\to\F$ for which $\supp(a)\in\mcW$. We may think of $a$ as a formal power series $\sum_{q\in\R}a_qt^q$ in an indeterminate $t$, with exponents $q\in\R$ and coefficients $a_q\in\F$; then $t^q$ is realised as the map $\R\to\F$ given by $(t^q)_p=\dl_{pq}$. Denote the set of all Hahn series by $K$. It is routine to check that $K$ becomes an abelian group with respect to the addition $$(a+b)_q=a_q+b_q.$$ As we will show in Lemmas \ref{finitely many pairs in product of hahn series} and \ref{addition of well-ordered subsets forms commutative monoid}, the formula $$(ab)_q=\sum_{u+v=q}a_ub_v$$ gives a well-defined multiplication on $K$ giving it the structure of a commutative unital ring, and satisfying $$\supp(a+b)\subs\supp(a)\cup\supp(b)$$ and $$\supp(ab)\subs\supp(a)+\supp(b)$$ for all $a,b\in K$, where $$A+B=\{a+b:a\in A,\;b\in B\}$$ for $A,B\subs\R$. We will also see that $K$ is a field (Corollary \ref{hahn series form a field}). We set $\nu(0)=\infty$ and $$\nu(a)=\min(\supp(a))$$ for $a\in K^\x$; this defines a map $\nu:K\to\R\cup\{\infty\}$ which we call the \textit{valuation} on $K$. We have that $\nu(a)=\infty$ if and only if $a=0$, $\nu(ab)=\nu(a)+\nu(b)$, $$\nu(a+b)\geq\min\{\nu(a),\nu(b)\}$$ and $\nu(t^q)=q$ for all $q\in\R$. \end{definition}

\begin{definition} We set $$A=\{a\in K:\nu(a)\geq0\}=\{a\in K:\supp(a)\subs[0,\infty)\}$$ and $$I_{>0}=\{a\in K:\nu(a)>0\}=\{a\in K:\supp(a)\subs(0,\infty)\}.$$ Then $A$ is a unital subring of $K$ and therefore an integral domain, which we call the \textit{Hahn ring}. \end{definition}

\begin{lem} \label{lemma for invertible elements and expression of elements} We have that $$A^\x=A\tk I_{>0}=\{a\in K:\nu(a)=0\}.$$ This allows us to express each non-zero $a\in K$ in the form $a=t^{\nu(a)}u$ for unique $u\in A^\x$. \end{lem}

\begin{proof} It is routine to check the equivalent expressions of $A^\x$. If $a\in K$ is non-zero, then set $u=t^{-\nu(a)}a\in K$. Then since $\nu(u)=0$, we see that $u\in A^\x$ with $a=t^{\nu(a)}u$. Uniqueness of $u$ is clear. \end{proof}

\begin{definition} We write $N\leq M$ to indicate that $N$ is an $A$-submodule of $M$. We set $$I_q=t^qA=\{a\in K:\nu(a)\geq q\}=\{a\in K:\supp(a)\subs[q,\infty)\}$$ and $$I_{>q}=\bigcup_{r>q}I_r=\{a\in K:\nu(a)>q\}=\{a\in K:\supp(a)\subs(q,\infty)\}$$ for each $q\in\R$, and see that $I_q\leq K$ and $I_{>q}\leq K$. Of particular interest will be the quotients $K/I_{>0}$, $K/A$, and $A/I_{>0}$, which we label $\Ta$, $\Phi$, and $\F$, respectively. \end{definition}

\begin{prop} \label{submodules of hahn field} The $A$-submodules of $K$ are precisely $0$, $K$, $I_q$ and $I_{>q}$ for $q\in\R$. They are totally ordered by inclusion; if $p<q$, then $$0<I_{>q}<I_q<I_{>p}<I_p<K.$$ \end{prop}

\begin{proof} The only non-trivial part is to confirm that the above $A$-submodules are the only ones. To that end, let $V\leq K$ and set $$S=\{q\in\R:t^q\in V\}.$$ If $q\in S$ and $q'\geq q$, then $t^{q'}=t^{q'-q}t^q\in V$, since $t^{q'-q}\in A$ and $t^q\in V$. Then $q'\in S$, showing that $q\in S$ implies $[q,\infty)\subs S$. We see from this that either $S$ is empty, $S=\R$, $S=[q,\infty)$ or $S=(q,\infty)$ for some $q\in\R$. 

If $S$ is empty, let $a\in V$. If $a\neq 0$, then (as we will see from Lemma \ref{lemma for invertible elements and expression of elements}) there exists unique $u\in A^x$ for which $a=t^{\nu(a)}u$. Then $t^{\nu(a)}=au^{-1}\in V$, so that $\nu(a)\in S$. But $S$ is empty, a contradiction. So $a=0$, showing that $V=0$. 

If $S=\R$, then let $a\in K$. If $a=0$, then clearly $a\in V$. Otherwise, we can again use Lemma \ref{lemma for invertible elements and expression of elements} to see that $a=t^{\nu(a)}u$ for unique $u\in A^\x$. Then since $\nu(a)\in\R=S$, it must be the case that $a=t^{\nu(a)}u\in V$, and then $V=K$. Similarly, for $q\in\R$, $S=[q,\infty)$ implies $V=I_q$ and $S=(q,\infty)$ implies $V=I_{>q}$. \end{proof}

\begin{cor} \label{ideals in hahn ring} The ideals in $A$ are precisely $$0<I_{>q}<I_q<I_{>p}<I_p<I_{>0}<A,$$ for $0<p<q$. We see that $A$ has unique maximal ideal $I_{>0}$ and that all finitely generated ideals are principal. \end{cor}

\begin{rmk} Since every finitely generated ideal in $A$ is principal, $A$ is an example of a \textit{B\'ezout domain}. \end{rmk}

\begin{definition} \label{basic mb definition} We say that an $A$-module $M$ is \textit{basic} if it is non-zero, and whenever $x,y\in M$ we have that $x\in Ay$ or $y\in Ax$. We say that $M$ is \textit{multibasic} if it can be decomposed into a finite (possibly empty) direct sum of basic submodules. \end{definition}

\begin{exmp} $K$ is basic. \end{exmp}

\begin{proof} Let $a,b\in K$. If $a=0$ or $b=0$ then clearly $a\in Ab$ or $b\in Aa$, so we assume that $a,b\neq 0$. Then we will see from Lemma \ref{lemma for invertible elements and expression of elements} that $a=t^{\nu(a)}u$, $b=t^{\nu(b)}v$ for unique $u,v\in A^\x$. Assume without loss of generality that $\nu(a)\leq\nu(b)$, so that $t^{\nu(b)-\nu(a)}vu^{-1}\in A$. Then $$b=t^{\nu(b)}v=t^{\nu(b)-\nu(a)}vu^{-1}t^{\nu(a)}u=t^{\nu(b)-\nu(a)}vu^{-1}a\in Aa,$$ and we are done. \end{proof}

\begin{exmp} Let $M$ be an $A$-module that is either basic or zero, and $N\leq M$. Then $N$ and $M/N$ are either basic or zero. \end{exmp}

\begin{exmp} \label{quotient of submodules of hahn field is basic} If $V<U\leq K$, then $U/V$ is basic. For each $q\in\R$, multiplication by $t^q$ gives an isomorphism $U\to I_qU$ that induces an isomorphism $V\to I_qV$. So there is a unique isomorphism $U/V\to I_qU/I_qV$ for which the following diagram commutes. \begin{center}
    \begin{tikzcd}[sep=2cm]
     U \arrow[r,"\sim"] \arrow[d,swap] &
     I_qU \arrow[d] \\
     U/V \arrow[r,dashed] & 
     I_qU/I_qV
    \end{tikzcd}
 \end{center} \end{exmp}

\begin{exmp} An \textit{interval} is a non-empty subset $X$ of $\R$ where for each $x,y\in X$ with $x\leq y$, we have that $[x,y]\subs X$. The \textit{standard intervals} are as follows: $\R$, $[0,\infty)$, $(0,\infty)$, $(-\infty,0)$, $(-\infty,0]$, $\{0\}$, and $[0,q]$, $[0,q)$, $(0,q]$, $(0,q)$ for $q>0$. 

For each interval $X$, we make the following definitions. \begin{align} &\begin{aligned} L(X)=\{y\in\R:y\geq x\text{ for some }x\in X\} \end{aligned}\\ &\begin{aligned} R(X)=\{y\in\R:y>x\text{ for all }x\in X\} \end{aligned}\\ &\begin{aligned} Q(X)=\{a\in K:\supp(a)\subs L(X)\}/\{a\in K:\supp(a)\subs R(X)\} \end{aligned} \end{align} We can now use Example \ref{quotient of submodules of hahn field is basic} to see that $Q(X)$ is a basic $A$-module; we call the images of the standard intervals under $Q$ the \textit{standard basic modules}, which are listed below. \[ \begin{array}{|c|c|c|c|c|c|c|c|c|c|}\hline K      & A & I_{>0} & \Ta & \Phi & \F & A/I_q & A/I_{>q} & I_{>0}/I_q & I_{>0}/I_{>q} \\ \hline
    \end{array}
 \] The isomorphisms described in Example \ref{quotient of submodules of hahn field is basic} show that if $V<U\leq K$, then $U/V$ is canonically isomorphic to a standard basic module. For $q>0$, we say that the standard basic modules $A/I_q$, $A/I_{>q}$, $I_{>0}/I_q$ and $I_{>0}/I_{>q}$ have \textit{length} $q$. We say that $\F$ has length $0$. 

\end{exmp}

\begin{definition} \label{reflexivity intro definitions} Recall the constructions of \S\ref{reflexive graded modules}. For each $A$-module $M$, define the evaluation map $\chi_M:M\to D^2M$ by $$\chi_M(m)(\al)=\al(m).$$ We say that $M$ is \textit{$\Ta$-reflexive} if $\chi_M$ is an isomorphism, and write $\mcM$ for the full subcategory of $\Mod_A$ consisting of the $\Ta$-reflexive $A$-modules. \end{definition}

Our main theorem regarding $A$-modules is as follows.

\begin{theorem} An $A$-module $M$ is $\Ta$-reflexive if and only if it is multibasic. Every multibasic module has a unique decomposition into basic submodules, up to the ordering of the summands. \end{theorem}

\begin{proof} That multibasic modules are $\Ta$-reflexive is proved as Theorem \ref{multibasics are theta ref}; the converse is proved as Theorem \ref{theta-ref implies multibasic}. Uniqueness of decomposition is proved as Theorem \ref{uniqueness theorem for multibasics}. \end{proof}

Finally, we note the following corollary, which combines what we know about $\Ta$-reflexive modules and multibasic modules to give a more detailed picture of the category $\mcM$. 

\begin{cor} The category $\mcM$ of multibasic modules is abelian, has enough projectives and enough injectives. The projective objects are the flat multibasics and the injective objects are the injective multibasics. $\mcM$ is not a Frobenius category. Submodules, quotients and extensions of multibasic modules are multibasic. \end{cor}

\section{Constructions}

Recall that $\mcW$ is the set of all subsets $A\subs\R$ for which every non-empty subset $B\subs A$ has a minimum $\min(B)$.

\begin{lem} \label{well-ordered leq lemma} Let $A\subs\R$ and suppose that $\{x\in A:x\leq a\}\in\mcW$ for all $a\in A$. Then $A\in\mcW$. \end{lem}

\begin{proof} Let $B$ be a non-empty subset of $A$ and choose $b\in B$. Then $\{x\in B:x\leq b\}$ is a non-empty subset of $\{x\in A:x\leq b\}\in\mcW$. So $\{x\in B:x\leq b\}$ has a minimum $b_0$. If $b'\in B$, then either $b'\leq b$ so that $b'\geq b_0$, or $b'>b\geq b_0$. So $b_0$ is the minimum of $B$, proving that $A\in\mcW$. \end{proof}

\begin{prop} \cite[Lemma 9.3]{ShS14} \label{well ordered implies countable} Each $A\in\mcW$ is countable. \end{prop}

\begin{proof} Define $\sg:A\to A\cup\{\infty\}$ as follows. If $A$ has a maximum, then we set $\sg(\max(A))=\infty$. Whether or not a maximum exists, if $a\in A$ is not the maximum, then $\{x\in A:x>a\}$ is a non-empty subset of $A\in\mcW$, and we set $\sg(a)$ to be its minimum. Then for each $a\in A$, we have that $\sg(a)>a$, so that $\Q\cap(a,\sg(a))$ is a non-empty subset of $\Q$. We then choose $f(a)\in\Q\cap(a,\sg(a))$ for each $a\in A$, defining a map $f:A\to\Q$.

Then if $a,b\in A$ with $a<b$, we have that $f(a)<\sg(a)\leq b<f(b)$, showing that $f$ is strictly increasing, and therefore injective. We conclude that $A$ is countable. \end{proof}

\begin{lem} \label{equiv condition for well-ordered} Let $A\subs\R$. Then the following are equivalent. 

1. $A\in\mcW$ 

2. Every weakly decreasing sequence in $A$ is eventually constant. 

3. There are no strictly decreasing sequences in $A$. \end{lem}

\begin{proof} We first show that 1. implies 2. Suppose that $A\in\mcW$, and let $(a_n)$ be a weakly decreasing sequence in $A$. Suppose further that $a_n$ is not eventually constant. Then $(a_n)$ has a strictly decreasing subsequence $(a_{n_k})$. Since $\{a_{n_k}:k\geq 1\}$ is a non-empty subset of $A\in\mcW$, it must have a minimum $a_{n_K}$ for some $K$. Then $a_{n_{K+1}}<a_{n_K}$ since $(a_{n_k})$ is strictly decreasing, but $a_{n_{K+1}}\geq a_{n_K}$ since $a_{n_K}$ is the smallest term in the sequence, a contradiction. So $(a_n)$ must eventually be constant. So 1. implies 2.

If 2. holds, then any strictly decreasing sequence in $A$ is a weakly decreasing sequence and therefore must be constant, a contradiction. So 2. implies 3. 

Finally, we show that 3. implies 1. Suppose that there are no strictly decreasing sequences in $A$, and let $B$ be a non-empty subset of $A$. If $B$ has no minimum, then we can construct a strictly decreasing sequence in $B$, which of course is a strictly decreasing sequence in $A$, a contradiction. So $B$ has a minimum, showing that $A\in\mcW$. So 3. implies 1., and we are done. \end{proof}

\begin{lem} \label{sequence in well-ordered set has weakly increasing subsequence} \cite[Lemma 9.5]{ShS14} If $A\in\mcW$, then each sequence in $A$ has a weakly increasing subsequence. \end{lem}

\begin{proof} Let $(a_n)$ be a sequence in $A$. Since $A\in\mcW$, we know that $\{a_i:i\geq 1\}$ has a minimum. Since $\Z_{>0}\in\mcW$, we know that $$\{n\geq 1:a_n=\min\{a_i:i\geq 1\}\}$$ has a minimum, which we label $n_1$. Suppose that we have chosen $n_k$ for some $k\geq 1$. Then, as above, $\{a_i:i>n_k\}$ has a minimum and therefore $$\{n>n_k:a_n=\min\{a_i:i>n_k\}\}$$ has a minimum, which we label $n_{k+1}$. By induction, we have constructed a strictly increasing sequence $n_1<n_2<n_3<\cdots$ of positive integers and thus a subsequence $(a_{n_k})$ of $(a_n)$. 

Finally, we check that this subsequence is weakly increasing. Since $n_1\geq 1$, we have that $$a_{n_2}=\min\{a_i:i>n_1\}\geq\min\{a_i:i\geq 1\}=a_{n_1}.$$ For $k\geq 2$, we have that $n_k>n_{k-1}$, so that $$a_{n_{k+1}}=\min\{a_i:i>n_k\}\geq\min\{a_i:i>n_{k-1}\}=a_{n_k}.$$ \end{proof}

\begin{lem} \cite[Lemma 13.2.9]{Pas77} \label{finitely many pairs in product of hahn series} If $A,B\in\mcW$ and $q\in\R$, then there are only finitely many pairs $(a,b)\in A\x B$ for which $a+b=q$. \end{lem}

\begin{proof} It is clear that the set of all such pairs is in bijection with $\{a\in A:q-a\in B\}$, which is in $\mcW$ since $A\in\mcW$. If it is empty, then we are done. Otherwise, it is the image of a weakly increasing sequence $(a_n)$ by Proposition \ref{well ordered implies countable}. Then $(q-a_n)$ is a weakly decreasing sequence in $B$, and $B\in\mcW$, so $(q-a_n)$ (and therefore $(a_n)$) must be constant by Lemma \ref{equiv condition for well-ordered}. So $\{a\in A:q-a\in B\}$ is finite, and we are done. \end{proof}

\begin{lem} \cite{Pas77}*{Lemma 13.2.9} \label{addition of well-ordered subsets forms commutative monoid} The set $\mcW$ becomes a commutative monoid when equipped with the addition $(A,B)\mapsto A+B$ given by $$A+B=\{a+b:a\in A,\,b\in B\}.$$ We also write $nA=A+A+\cdots+A$ for $n\geq 1$ (not to be confused with $\{na:a\in A\}$, which is often denoted $nA$). \end{lem}

\begin{proof} Let $A,B\in\mcW$ and $(c_n)$ be a weakly decreasing sequence in $A+B$. Then there exist sequences $(a_n)$ in $A$ and $(b_n)$ in $B$ with $c_n=a_n+b_n$ for all $n$. Since $A,B\in\mcW$, we see from Lemma \ref{sequence in well-ordered set has weakly increasing subsequence} that there exists a strictly increasing sequence $(n_k)$ of positive integers for which $(a_{n_k})$ and $(b_{n_k})$ are both weakly increasing, so that $(c_{n_k})$ is both weakly increasing and weakly decreasing, so must be constant. 

Since $(c_n)$ is weakly decreasing and has a constant subsequence, it must eventually be constant. Then since every weakly decreasing sequence in $A+B$ is eventually constant, we see from Lemma \ref{equiv condition for well-ordered} that $A+B\in\mcW$. We see that $\mcW$ is closed under this addition of subsets, and it is now routine to check that $\mcW$ is a commutative monoid. \end{proof}

The following two lemmas are fundamental to establishing that $K$ is a field in Corollary \ref{hahn series form a field}.

\begin{lem} \cite[Lemma 13.2.10]{Pas77} \label{a tilde is well ordered} If $A\in\mcW$ and $A\subs(0,\infty)$, then $$\td{A}=\bigcup_{n=1}^\infty nA\in\mcW.$$ \end{lem}

\begin{proof} If $A$ is empty then this is trivial, so we assume that $A$ is not empty and therefore $\td{A}\subs(0,\infty)$ is also not empty. Let $w=\min(A)$, so that $nA\subs[nw,\infty)$ for all $n\geq 1$. Let $a\in\td{A}$, $B=\{x\in\td{A}:x\leq a\}$ and $S$ be a non-empty subset of $B$. If we can show that $S$ has a minimum, then $B\in\mcW$ and then it will follow from Lemma \ref{well-ordered leq lemma} that $\td{A}\in\mcW$. 

To that end, choose $n\geq 1$ with $nw>a$, and let $s\in S$. Then if $k\geq n$, we have that $$s\leq a<nw\leq kw,$$ so that $s\notin kA$. Since $s\in\td{A}$, this shows that $s\in A\cup\ldots\cup(n-1)A$, and we conclude that $$S\subs A\cup\cdots\cup(n-1)A.$$ Since $A\in\mcW$ and $\mcW$ is closed under finite unions, it follows from Lemma \ref{addition of well-ordered subsets forms commutative monoid} that $A\cup\cdots\cup(n-1)A\in\mcW$ so that $S$ must have a minimum, and we are done. \end{proof}

\begin{lem} \cite[Lemma 13.2.10]{Pas77} \label{tilde gives empty intersection} If $A\in\mcW$ and $A\subs(0,\infty)$, then $(n\td{A})_{n\geq 1}$ is a decreasing sequence of subsets of $\td{A}$ with empty intersection. \end{lem}

\begin{proof} If $A$ is empty, then this is trivial so we assume that it is not empty and then $\td{A}\subs(0,\infty)$ is not empty. Since $\td{A}+\td{A}\subs\td{A}$, we can see by induction that $(n\td{A})_{n\geq 1}$ is a decreasing sequence of subsets of $\td{A}$. Suppose that their intersection is non-empty and denote its minimum by $m$. For each $n\geq 2$, we choose $y_n\in(n-1)\td{A}$ and $z_n\in\td{A}$ with $m=y_n+z_n$. We see from Lemma \ref{finitely many pairs in product of hahn series} that $Y=\{y_n:n\geq 2\}$ is finite, and we choose $y\in Y$ and a subsequence $(y_{n_k})_{k\geq 1}$ of $(y_n)_{n\geq2}$ with $y_{n_k}=y$ for all $k$. For each $n\geq 2$, choose $k\geq 1$ with $n_k>n$. Then $$y=y_{n_k}\in(n_k-1)\td{A}\subs n\td{A},$$ so that $y\geq m$. But for any $k\geq 1$, we have that $y=y_{n_k}=m-z_{n_k}<\al$, a contradiction. We conclude that $\bigcap_{n=1}^\infty n\td{A}$ must be empty. \end{proof}

Having established various properties of well-ordered subsets of the real numbers, we can now proceed to construct the field of Hahn series, following \S 13 of \cite{Pas77}. 

\begin{lem} \cite[Theorem 13.2.11]{Pas77} \label{finite sum lemma} Let $a\in K$ be non-zero with $A=\supp(a)\subs(0,\infty)$. Then $\supp(a^n)\subs n\td{A}$ for all $n\geq 1$, and $\sum_{n=0}^\infty (a^n)_q\in\F$ for all $q\in\R$. \end{lem}

\begin{proof} If we set $A=\supp(a)$, then we see from Lemma \ref{a tilde is well ordered} and Lemma \ref{tilde gives empty intersection} that $\td{A}\in\mcW$ and $(n\td{A})_{n\geq 1}$ is a decreasing sequence of subsets of $\td{A}$ with empty intersection. Let $n\geq 1$ and $q\in\supp(a^n)$. Then $$(a^n)_q=\sum_{\sum_{i=1}^nu_i=q}a_{u_1}\cdots a_{u_n}=1,$$ so there exist $u_1,\ldots,u_n\in\supp(a)=A\subs\td{A}$ with $u_1+\cdots+u_n=q$, so that $q\in n\td{A}$. So $\supp(a^n)\subs n\td{A}$ for all $n\geq 1$. 

Now let $q\in\R$ and suppose that $q\in\supp(a^n)$ for infinitely many $n\geq 0$. Then $q\in n\td{A}$ for infinitely many $n\geq 1$, so that $q\in\bigcap_{n=1}^\infty n\td{A}=\emptyset$, a contradiction. So for each $q\in\R$, there are only finitely many $n\geq 0$ for which $q\in\supp(a^n)$, showing that $\sum_{n=0}^\infty (a^n)_q$ is in fact a finite sum of elements in $\F$. \end{proof}

\begin{lem} \cite[Theorem 13.2.11]{Pas77} \label{lemma for proof that hahn series form a field} If $a\in K$ with $A=\supp(a)\subs(0,\infty)$, then $1+a\in K^\x$. \end{lem}

\begin{proof} This is clearly true if $a=0$, so we assume that $a$ is non-zero. Then Lemma \ref{finite sum lemma} allows us to define $b:\R\to\F$ by $b_q=\sum_{n=0}^\infty(a^n)_q$. Suppose that $q\leq 0$. Then if $n\geq 1$, we know from Lemma \ref{finite sum lemma} that $\supp(a^n)\subs n\td{A}\subs(0,\infty)$, so that $(a_n)_q=0$. So $b_q=0$ for $q<0$ and $b_0=1$. Now suppose that $q>0$ with $q\notin\td{A}$. Then $(a^n)_q=0$ for $n\geq 1$, so that $b_q=1_q=0$. We see that $\supp(b)\subs\td{A}\cup\{0\}$ and since $\td{A}\in\mcW$, this shows that $b\in K$. We also have that $$\supp((1+a)b)\subs\td{A}\cup\{0\}$$ and $((1+a)b)_0=(1+a)_0b_0=1$; if we can show that $(1+a)b$ is zero on $\td{A}$, then $(1+a)b=1$, and we are done.

To that end, let $q\in\td{A}$. Then since $\bigcap_{n=1}^\infty n\td{A}$ is empty (Lemma \ref{tilde gives empty intersection}), there exists $n\geq 1$ for which $q\in n\td{A}\tk(n+1)\td{A}$. Set $c=b+\sum_{k=0}^n a^k\in K$. Then if $p\in\supp(c)$, we have that $c_p=\sum_{k=n+1}^\infty(a^k)_p=1$, giving that $p\in k\td{A}$ for some $k\geq n+1$ so that $p\in(n+1)\td{A}$. So $\supp(c)\subs(n+1)\td{A}$ and then $$\supp((1+a)c)\subs(n+1)\td{A}.$$ Then since $$(1+a)b=1+a^{n+1}+(1+a)c,$$ we see that $$\supp((1+a)b)\subs(n+1)\td{A}\cup\{0\}.$$ Since $q\notin (n+1)\td{A}$ and $q\neq 0$, we see that $((1+a)b)_q=0$, as required. \end{proof}

\begin{cor} \cite[Theorem 13.2.11]{Pas77} \label{hahn series form a field} $K$ is a field. \end{cor}

\begin{proof} Let $a\in K$ be non-zero and $b=1+t^{-\nu(a)}a\in K$. Then $\supp(b)\subs(0,\infty)$, so Lemma \ref{lemma for proof that hahn series form a field} gives that $1+b\in K^\x$. Then since $t^{\nu(a)}\in K^\x$, we have that $a=t^{\nu(a)}(1+b)\in K^\x$. \end{proof}

\begin{definition} For $U,V\leq K$, we define $$UV=\{ab:a\in U\text{, }b\in V\},$$ $$(U:V)=\{a\in K:aV\leq U\},$$ and $$V^\circ=(I_{>0}:V).$$ \end{definition}

These are tabulated as follows, with $U$ spanning the rows and $V$ spanning the columns.  \[ \begin{array}{|c|c|c|c|c|c|}\hline (U:V)      & 0 & K & I_q & I_{>q} \\ \hline 
     0         & K & 0 & 0      & 0  \\ \hline 
     K         & K & K & K      &  K  \\ \hline 
     I_p    & K & 0 & I_{p-q}      & I_{p-q}  \\ \hline 
     I_{>p}       & K & 0 & I_{>p-q}      & I_{p-q} \\ \hline
    \end{array}
 \]
 
\[ \begin{array}{|c|c|c|c|c|c|}\hline UV      & 0 & K & I_q & I_{>q} \\ \hline 
     0         & 0 & 0 & 0      & 0  \\ \hline 
     K         & 0 & K & K      &  K  \\ \hline 
     I_p    & 0 & K & I_{p+q}      & I_{>p+q}  \\ \hline 
     I_{>p}       & 0 & K & I_{>p+q}      & I_{>p+q} \\ \hline
    \end{array}
 \]
 
\[ \begin{array}{|c|c|c|c|c|c|}\hline V      & 0 & K & I_q & I_{>q} \\ \hline 
     V^\circ         & K & 0 & I_{>-q}      & I_{-q}  \\ \hline 
    \end{array}
 \]

Using these tables, we establish the following.

\begin{prop} \label{combining submods of k} If $U,V\leq K$, then $UV$, $(U:V)$ and $V^\circ$ are all $A$-submodules of $K$. The $A$-submodules of $K$ form a commutative monoid with identity $I_0=A$ with respect to the operation $(U,V)\mapsto UV$. Also $(V^\circ)^\circ=V$, $V\leq U$ implies $U^\circ\leq V^\circ$, and $$(U:V)=(VU^\circ)^\circ.$$ \end{prop}

\section{Flat and Injective Modules}

In this section, we establish criteria for flatness and injectivity of $A$-modules. We first recall Baer's Criterion for injectivity of modules.

\begin{theorem} \cite[Theorem 3.7]{Lam99} (Baer's Criterion) Let $M$ be a module over a commutative unital ring $R$. Then $M$ is injective if and only if for each ideal $L$ in $R$, we can extend every $R$-linear map $\al:L\to M$ to the whole of $R$. \end{theorem}

The following equivalent condition for flatness follows from Baer's Criterion and is adapted from the 'Modified Flatness Test' in \cite{Lam99} to the case of $A$-modules.

\begin{lem} \cite[4.12]{Lam99} \label{modified flatness test} $M$ is flat if and only if $(I_q\incl A)\tns\id_M$ is injective for $q>0$. \end{lem}

\begin{theorem} \label{flatness criterion} An $A$-module $M$ is flat if and only if $\ann(t,M)=0$. \end{theorem}

\begin{proof} For $q>0$, label the map $m\mapsto t^qm:M\to M$ by $\mu_q$. First we show that $M$ is flat if and only if $\mu_q$ is injective for all $q>0$. We already know (Lemma \ref{modified flatness test}) that $M$ is flat if and only if $$(I_q\incl A)\tns\id_M:I_q\tns M\to A\tns M$$ is injective for all $q>0$. We have an isomorphism $$m\mapsto 1\tns m:M\to A\tns M,$$ which we label $\ta$, and the isomorphism $a\mapsto t^qa:A\to I_q$ induces an isomorphism $$a\tns m\mapsto t^qa\tns m:A\tns M\to I_q\tns M.$$ Composing these two isomorphisms together, followed by $(I_q\incl A)\tns\id_M$, we obtain the map $$m\mapsto t^q\tns m:M\to A\tns M,$$ which is equal to $\ta\mu_q$. Since the first two maps are isomorphisms, we see that $M$ is flat if and only if $(I_q\incl A)\tns\id_M$ is injective if and only if $\ta\mu_q$ is injective if and only if $\mu_q$ is injective. 

Since $\ker(\mu_q)=\ann(t^q,M)$, we see that $M$ is flat if and only if $\ann(t^q,M)=0$ for all $q>0$. All that remains is to check that $\ann(t,M)=0$ implies $\ann(t^q,M)=0$ for all $q>0$. To that end, if $\ann(t,M)=0$, then we can use induction to show that $\ann(t^n,M)=0$ for all positive integers $n$. Then for any $q>0$, there exists a positive integer $n\geq q$ so that for each $m\in\ann(t^q,M)$ we have that $t^nm=t^{n-q}(t^qm)=0$, and then $m\in\ann(t^n,M)=0$, so that $\ann(t^q,M)=0$, as required. \end{proof}

Recall that $\Ta=K/I_{>0}$ and $\Phi=K/A$.

\begin{cor} \label{flatness corollary} Let $M$ be an $A$-module. Then the following are equivalent. 

1. $M$ is flat. 

2. The map $$m\mapsto 1\tns m:M\to K\tns M$$ is injective. 

3. $M$ is isomorphic to an $A$-submodule of a $K$-vector space. \end{cor}

\begin{proof} 1. implies 2. Since \begin{center}
    \begin{tikzcd}
     0 \arrow[r] & 
     A \arrow[r] & 
     K \arrow[r] & 
     \Phi \arrow[r] &
     0
    \end{tikzcd}
 \end{center} is exact and $M$ is flat, we know that \begin{center}
    \begin{tikzcd}
     0 \arrow[r] & 
     M \arrow[r] & 
     K\tns M \arrow[r] & 
     \Phi\tns M \arrow[r] &
     0
    \end{tikzcd}
 \end{center} is exact, and therefore $M\to K\tns M$ is injective. 

2. implies 3. If the map $M\to K\tns M$ is injective, then $M$ is isomorphic to a submodule of $K\tns M$, which is a $K$-vector space.

3. implies 1. By assumption, there exists a $K$-vector space $V$ and an injective map $\al:M\to V$. Let $m\in\ann(t,M)$. Then $t\al(m)=0$ with $\al(m)\in V$. Since $V$ is a $K$-vector space, we see that $\al(m)=t^{-1}(t\al(m))=0$, and then $m=0$ since $\al$ is injective. So $\ann(t,M)=0$ and we see from Theorem \ref{flatness criterion} that $M$ is flat. \end{proof}

\begin{rmk} We have seen in Corollary \ref{flatness corollary} that if $M$ is flat, then there exists a $K$-vector space $V$ and an injective map $\al:M\to V$. Suppose that $M$ is also finitely generated. Then $\im(\al)$ is finitely generated so that the $K$-linear subspace $V'$ of $V$ generated by $\im(\al)$ is finite-dimensional. Then $M$ must be isomorphic to an $A$-submodule of a \textbf{finite-dimensional} $K$-vector space. We will see later (Corollary \ref{flat fin gen implies free}) that if $M$ is flat and finitely generated, then it is in fact free. \end{rmk}

\begin{exmp} Let $V<U\leq K$. Then $U/V$ is flat if and only if $V=0$. \end{exmp}

\begin{proof} If $V=0$, then $U/V=U\leq K$ is flat by Corollary \ref{flatness corollary}. Conversely, if $U/V$ is flat, then $\ann(t,U/V)=0$ by Theorem \ref{flatness criterion}. Let $a\in U$. Then $ta\in V$ implies $a\in V$. If $V$ is non-zero, then $t^na\in V$ for some $n\geq 1$, so that $a\in V$ by induction, giving that $U=V$, a contradiction. So $V=0$. \end{proof}

We now use Baer's Criterion to derive an injectivity criterion specific to $A$-modules.

\begin{theorem} \label{injectivity criterion} An $A$-module $M$ is injective if and only if both of the following hold.

(I1) $tM=M$. 

(I2) For each morphism $\al:I_{>0}\to M$, there exists $m\in M$ for which $\al(a)=am$ for all $a\in I_{>0}$. \end{theorem}

\begin{proof} Suppose that $M$ is injective. For each $m\in M$, we use Baer's Criterion to extend the morphism $a\mapsto t^{-1}am:I_1\to M$ to a morphism $\bt:A\to M$ and we see that $m=t\bt(1)\in tM$, proving I1. We use Baer's Criterion again to extend each morphism $\al:I_{>0}\to M$ to some $\td{\al}:A\to M$ and see that $\al(a)=a\td{\al}(1)$ for all $a\in I_{>0}$, proving I2. 

Conversely, using I1, an induction argument shows that $t^nM=M$ for all integers $n\geq0$, and then it follows that $t^qM=M$ for all $q\geq 0$. Let $J$ be an ideal in $A$ and $\al:J\to M$ be any morphism. It is clear that $\al$ extends to the whole of $A$ if $J=0$ or $J=A$. If $J=I_q$ for some $q>0$, then $\al(t^q)\in M=t^qM$ so that $\al(t^q)=t^qm$ for some $m\in M$. Then $\al(a)=am$ for all $a\in J=I_q$, allowing us to extend $\al$ to the whole of $A$. 

Finally, suppose that $J=I_{>q}$ for some $q\geq 0$, and define the morphism $\bt:I_{>0}\to M$ by $\bt(a)=\al(t^qa)$. Then we see from I2 that there exists $m'\in M$ with $\bt(a)=am'$ for all $a\in I_{>0}$. Then $m'\in M=t^qM$ gives that $m'=t^qm$ for some $m\in M$, and we see that $\al(a)=am$ for all $a\in J$, allowing us again to extend $\al$ to the whole of $A$. Injectivity of $M$ then follows from Baer's Criterion. \end{proof}

\begin{exmp} \label{hahn field is injective module} $K$ is injective. \end{exmp}

\begin{proof} Clearly $tK=K$, so I1 holds. For each morphism $\al:I_{>0}\to K$, we set $m=t^{-1}\al(t)\in K$. If $q\geq 1$, then $$\al(t^q)=\al(t^{q-1}t)=t^{q-1}\al(t)=t^{q-1}tm=t^qm.$$ If $0<q<1$, then $$tm=\al(t)=\al(t^{1-q}t^q)=t^{1-q}\al(t^q),$$ so that $\al(t^q)=t^qm$ for all $q>0$. It is then clear that $\al(a)=am$ for all $a\in I_{>0}$, so that I2 holds and then $K$ is injective by Theorem \ref{injectivity criterion}. \end{proof}

We now consider two examples to show that neither of the conditions I1 and I2 in Theorem \ref{injectivity criterion} imply the other.

\begin{exmp} \label{infinite sum theta not injective} $\bigoplus_{i=0}^\infty\Ta$ satisfies I1 but not I2 of Theorem \ref{injectivity criterion}, and therefore is not injective. \end{exmp} 

\begin{proof} Set $M=\bigoplus_{i=0}^\infty\Ta$. If $m\in M$, then there exists a sequence $(a_i)$ in $K$ with $a_i=0$ for all but finitely many $i$ and with $$m=(a_i+I_{>0})=(t(t^{-1}a_i+I_{>0})),$$ so that $tM=M$; $I1$ holds. 

For each $i\geq 0$, define $\al_i:I_{>0}\to\Ta$ by $\al_i(a)=t^{-2^{-i}}a+I_{>0}$. Then $\ker(\al_i)=I_{2^{-i}}$ for all $i$, so that if $a\in I_{>0}$, then $\al_i(a)=0$ for sufficiently large $i$. This allows us to define $\al:I_{>0}\to M$ by $\al(a)=(\al_i(a))$. Suppose that there exists $m\in M$ with $\al(a)=am$ for all $a\in I_{>0}$. Then there exists a sequence $(a_i)$ in $K$ with $a_i=0$ for all but finitely many $i$ and for which $m=(a_i+I_{>0})$. So for each $q>0$, we have that $$(t^{q-2^{-i}}+I_{>0})=\al(t^q)=t^q(a_i+I_{>0})=(t^qa_i+I_{>0}).$$ So for all $q>0$ and integers $i\geq 0$, we have that $t^q(t^{-2^{-i}}+a_i)\in I_{>0}$ and then $\nu(t^{-2^{-i}}+a_i)>-q$. So for all $i$, we have that $t^{-2^{-i}}+a_i\in A$, and then $\nu(a_i)=-2^{-i}\leq 0$ gives that $a_i\notin I_{>0}$, a contradiction. So I2 fails. \end{proof}

\begin{exmp} $A$ satisfies I2 but not I1 of Theorem \ref{injectivity criterion}, and therefore is not injective. \end{exmp} 

\begin{proof} It is clear that I1 fails. Given a morphism $\al:I_{>0}\to A$, we can compose with the inclusion $A\incl K$ and use the injectivity of $K$ (Example \ref{hahn field is injective module}) to see that there exists $b\in K$ with $\al(a)=ab$ for all $a\in I_{>0}$. Then for all $q>0$, we have that $\al(t^q)=t^qa\in A$, so that $\nu(a)\geq -q$ for all $q>0$. Then $a\in A$, so that the second condition for injectivity in Theorem \ref{injectivity criterion} is satisfied. \end{proof}

After establishing an important connection between $\Ta$ and $\Phi$ in Proposition \ref{isomorphism lemma for phi injective}, we can prove that these are also injective. 

\begin{definition} Define $\zeta:K\to DK$ by $$\zeta(a)(b)=ab+I_{>0}.$$ \end{definition}

\begin{prop} \label{surjective lemma for injectivity of theta} The unique map $\Ta\to D(I_{>0})$ for which \begin{center}
    \begin{tikzcd}[sep=2cm]
     K \arrow[r,twoheadrightarrow] \arrow[d,"\zeta"swap] &
     \Ta \arrow[d,dashed] \\
     DK \arrow[r] & 
     D(I_{>0})
    \end{tikzcd}
 \end{center} commutes is surjective. \end{prop}

\begin{proof} Let $\al\in D(I_{>0})$. We want to show that there exists $a\in K$ with $\al(b)=ab+I_{>0}$ for all $b\in I_{>0}$.

To that end, we set $$X=\{a\in K:\supp(a)\subs(-\infty,0]\}.$$ For each $q>0$, we define $\rho(q)$ to be the unique $a\in X$ for which $\al(t^q)=a+I_{>0}$. If $r>s>0$, then $\supp(\rho(r)-t^{r-s}\rho(s))$ is a subset of $$\supp(\rho(r))\cup\supp(t^{r-s}\rho(s))\subs(-\infty,0]\cup(-\infty,r-s]\subs(-\infty,r-s].$$ Then $$\rho(r)+I_{>0}=\al(t^r)=t^{r-s}\al(t^s)=t^{r-s}\rho(s)+I_{>0},$$ so that $\rho(r)-t^{r-s}\rho(s)\in I_{>0}$ and then $$\supp(\rho(r)-t^{r-s}\rho(s))\subs(0,\infty).$$ We conclude that $$\supp(\rho(r)-t^{r-s}\rho(s))\subs(0,r-s].$$ 

For $q>0$, set $\xi(q)=t^{-q}\rho(q)\in K$. Then (still with $r>s>0$) we see that $$\xi(r)-\xi(s)=t^{-r}\rho(r)-t^{-s}\rho(s)=t^{-r}(\rho(r)-t^{r-s}\rho(s)),$$ so that $$\supp(\xi(r)-\xi(s))\subs\{-r\}+(0,r-s]=(-r,-s].$$ Then checking the cases $r/2>s$, $r/2<s$ and $r/2=s$ shows that $$\supp(\xi(r/2)-\xi(s))\subs(-\max\{r/2,s\},-\min\{r/2,s\}]$$ for all $r,s>0$. Finally, if $r\geq s>0$, then since $-r\leq-\max\{r/2,s\}$ we see that $(\xi(r/2)-\xi(s))_{-r}=0$ and then \begin{equation} \label{xi property} \xi(r/2)_{-r}=\xi(s)_{-r}.
\end{equation} 

Define $a:\R\to\F$ by $a_r=\xi(-r/2)_r$ for $r<0$ and $a_r=0$ for $r\geq0$. Let $r>0$ and $q\leq -r$. Since $q<0$, we see that $a_q=\xi(-q/2)_q$. Since $-q\geq r>0$, we see from \eqref{xi property} that $\xi(-q/2)_q=\xi(r)_q$. So $a_q=\xi(r)_q$ and then $q\notin\supp(a-\xi(r))$. Finally, we see that \begin{equation} \label{supp xi subset} \supp(a-\xi(r))\subs(-r,\infty) \end{equation} for all $r>0$. 

Let $r\in\supp(a)$ and $Q=\{q\in\supp(a):q\leq r\}$. Since $r\in\supp(a)$, we know that $r<0$, and then we see from \eqref{supp xi subset} that $$\supp(a-\xi(-r))\subs(r,\infty).$$ So if $q\in Q$, then $\xi(-r)_q=a_q=1$ and then $q\in\supp(\xi(-r))$. We conclude that $Q\subs\supp(\xi(-r))$, so that $Q\in\mcW$. It then follows from Lemma \ref{well-ordered leq lemma} that $\supp(a)\in\mcW$, so that $a\in K$. 

For each $q>0$, we have that $$\al(t^q)=t^{q/2}\al(t^{q/2})=t^{q/2}\rho(q/2)+I_{>0}=t^q\xi(q/2)+I_{>0},$$ and then we use \eqref{supp xi subset} to see that $$\supp(t^qa-t^q\xi(q/2))\subs\{q\}+(-q/2,\infty)=(q/2,\infty)\subs(0,\infty).$$ So $$\al(t^q)=t^q\xi(q/2)+I_{>0}=t^qa+I_{>0},$$ and we conclude that $\al(b)=ab+I_{>0}$ for all $b\in I_{>0}$, as required. \end{proof}

\begin{prop} \label{theta is injective} $\Ta$ is injective, and therefore $D$ is exact. \end{prop}

\begin{proof} We use Theorem \ref{injectivity criterion}. Clearly I1 holds. Let $\al:I_{>0}\to\Ta$. Then we see from Proposition \ref{surjective lemma for injectivity of theta} that there exists $a\in K$ for which $\al(b)=ab+I_{>0}$ for all $b\in I_{>0}$, so that I2 holds, and we are done. \end{proof}

\begin{prop} \label{isomorphism lemma for phi injective} The map $\Phi\to D(I_{>0})$ for which \begin{center}
    \begin{tikzcd}[sep=2cm]
     K \arrow[r,twoheadrightarrow] \arrow[d,"\zeta"swap] &
     \Phi \arrow[d,dashed] \\
     DK \arrow[r] & 
     D(I_{>0})
    \end{tikzcd}
 \end{center} commutes is an isomorphism. \end{prop} 

\begin{proof} It is clear that the surjective map defined in Proposition \ref{surjective lemma for injectivity of theta} has kernel $A/I_{>0}\simeq\F$. The result then follows from the first and third isomorphism theorems. \end{proof}

\begin{prop} $\Phi$ is injective. \end{prop}

\begin{proof} We use Proposition \ref{isomorphism lemma for phi injective} and the tensor-hom adjunction to see that $$\Hom(M,\Phi)=\Hom(M,D(I_{>0}))=D(I_{>0}\tns M).$$ Since $I_{>0}$ is a submodule of $K$, it is flat (Corollary \ref{flatness corollary}) so that $I_{>0}\tns(-)$ is exact, and we have just seen that $D$ is exact (Example \ref{theta is injective}). So $\Hom(-,\Phi)=D(I_{>0}\tns(-))$ is exact and therefore $\Phi$ is injective. \end{proof}

\begin{exmp} Let $V<U\leq K$. Then $U/V$ is injective if and only if $U=K$. \end{exmp}

\begin{proof} If $U=K$, then $U/V=K/V$ is isomorphic to $K$, $\Ta$ or $\Phi$ and therefore must be injective. Conversely, if $U/V$ is injective, then $t(U/V)=U/V$ by Theorem \ref{injectivity criterion}. If $U\neq K$, then $tU<U$ so that $t(U/V)<U/V$, a contradiction. We conclude that $U/V$ is injective. \end{proof}

We now introduce injective hulls and use these to classify the injective $A$-modules.

\begin{definition} A \textit{uniserial ring} is a commutative unital ring $R$ (not necessarily an integral domain) whose ideals are totally ordered by inclusion. Let $R$ be a uniserial ring. An $R$-module $I$ is said to be \textit{indecomposable} if it has no non-trivial splitting into submodules. We say that an extension $I$ of $M$ is \textit{essential} if whenever $E\leq I$ and $E\cap M=0$, we have that $E=0$. An \textit{injective hull} of an $R$-module $M$ is an injective essential extension of $M$. Injective hulls are unique up to isomorphism fixing $M$. \end{definition}

\begin{exmp} \label{injective hulls of sb mods} Let $V<U\leq K$. Then $K/V$ is an injective hull of $U/V$. So every basic module has injective hull $K$, $\Ta$ or $\Phi$. \end{exmp}

\begin{proof} $K/V$ is clearly an injective extension of $U/V$. Let $W/V$ be a submodule of $K/V$, where $V\leq W\leq K$ and $W/V\cap U/V=0$. If $U\leq W$, then $U/V\leq W/V$, so that $U/V=0$, a contradiction. So $W\leq U$ gives that $W/V=W/V\cap U/V=0$, proving that $K/V$ is an essential extension of $U/V$. \end{proof}

The following example concerns finding the injective hull of $\bigoplus_{i=0}^\infty\Ta$. It is inconclusive, but rules out two obvious cases to try.

\begin{exmp} We have seen in Example \ref{infinite sum theta not injective} that $L=\bigoplus_{i=0}^\infty\Ta$ is not injective. Can we find an explicit description for the injective hull of $L$? The larger module $N=\prod_{i=0}^\infty\Ta$ is a product of injective modules and therefore injective. Let $$\Dl=\{(x,x,x,\ldots):x\in\Ta\}.$$ Then $\Dl\leq N$ and $\Dl\cap L=0$, so that $N$ is not an essential extension of $L$ and therefore not an injective hull of $L$. 

Now define $$M=\{x\in N:t^qx=0\text{ for some }q\geq 0\}.$$ Then $M\leq N$ and $t^q(t^{-i}+I_{>0})_{i\geq 0}$ is non-zero for all $q\geq 0$, so $M<N$. Also, $t^2(t^{-1/i}+I_{>0})_{i\geq 0}=0$, so $(t^{-1/i}+I_{>0})\in M$ but clearly not in $L$, giving that $L<M<N$. If $x\in M$, then since $N$ is injective we have that $x\in N=tN$ so that $x=tx'$ for some $x'\in N$. But $x\in M$, so there exists $q\geq 0$ for which $t^qx=t^{q+1}x'=0$, giving that $x'\in M$ and then $tM=M$. 

Let $\al:I_{>0}\to M$. Then since $N$ is injective, we compose with the inclusion $M\incl N$ to see that there exists $x\in N$ with $\al(a)=ax$ for all $a\in I_{>0}$. Since $tx=\al(t)\in M$, there exists $q\geq 0$ with $t^q(tx)=t^{q+1}x=0$ and therefore $x\in M$, proving injectivity of $M$. 

Is $M$ the injective hull of $L$? For any $x\in\Dl$, there exists $q\geq 0$ with $t^qx=0$, so that $\Dl\leq M$. But $\Dl\cap L=0$, so $M$ is not an essential extension of $L$ and therefore not injective hull of $L$. \end{exmp}

\begin{theorem} \cite[Theorem 4.4]{FuS85} \label{classification of injectives thm} Let $R$ be a uniserial ring. The indecomposable injective $R$-modules are precisely the injective hulls of the cyclic $R$-modules. The injective $R$-modules are precisely the injective hulls of (possibly infinite) direct sums of the indecomposable injectives. \end{theorem}

\begin{prop} \label{classification of injective a mods} The indecomposable injective $A$-modules are $K$, $\Ta$ and $\Phi$. The injective $A$-modules are the injective hulls of direct sums of copies of $K$, $\Ta$ and $\Phi$. \end{prop}

\begin{proof} This follows easily from Theorem \ref{classification of injectives thm} and Example \ref{injective hulls of sb mods}. \end{proof}

Finally, we use the fact that $A$ is local to note the very simple classification of projective $A$-modules.

\begin{theorem} \cite[Theorem 2]{Kap58} \label{projective iff free} An $A$-module is projective if and only if it is free. \end{theorem}

In Chapter \ref{reflexivity}, we will obtain simple classifications of the injective objects and the projective objects in the full subcategory $\mcM$. We will use the fact that $\Ta$ and $\Phi$ are injective $A$-modules to show in Chapter \ref{infinite root algebra chapter} that $A/I_{>1}$ and $A/I_1$ are injective modules over the infinite root algebra $P=A/I_{>1}$, so that, in particular, $P$ is a self-injective ring.

\section{Basic Modules} 

In this section, we discuss basic modules and see that each of these is isomorphic to a subquotient of $K$ (Theorem \ref{basic iso to standard basic}), in particular to a standard basic module. We shall see in Theorem \ref{uniqueness theorem for multibasics} that this standard basic module is unique, providing a classification of basic modules up to isomorphism. 

We derive formulae, and provide tables, for the action of the functor $D$, and all Ext and Tor modules, on the standard basic modules, demonstrating these all of these functors send basic modules either to other basic modules or to zero. Finally, we find projective and injective resolutions for the standard basic modules.

\begin{definition} (Reminder of Definition \ref{basic mb definition}) We say that an $A$-module $M$ is \textit{basic} if it is non-zero, and for all $x,y\in M$ we have that $x\in Ay$ or $y\in Ax$. \end{definition}

The following lemma will be used to show that every basic module is isomorphic to a standard basic module in Theorem \ref{basic iso to standard basic}.
 
\begin{lem} \label{basic module induces map lemma} For each basic module $M$, there exists a unique map $$\psi:(M\tk\{0\})\x(M\tk\{0\})\to\R$$ satisfying $\psi(ax,x)=\nu(a)$ and $\psi(x,ax)=-\nu(a)$ for all $a\in A$ and $x\in M$ with $ax\neq 0$. Furthermore, we have that $$\psi(x,y)+\psi(y,z)=\psi(x,z)$$ and $\psi(y,x)=-\psi(x,y)$ for all non-zero $x,y,z\in M$. Finally, $\psi(y,x)\geq 0$ if and only if $y\in Ax$. \end{lem}

\begin{proof} Let $z,w\in M$ be non-zero. Since $M$ is basic, we know that $z=aw$ or $w=az$ for some non-zero $a\in A$. This shows that each element of the domain of $\psi$ has the form $(ax,x)$ or $(x,ax)$ for some $a\in A$ and $x\in M$ with $ax\neq 0$. We may assume without loss of generality that any element of the domain takes the first of these forms. 

First of all, we show that $\nu(a)$ does not depend on the choice of $a$ and $x$. To that end, let it be expressed as $(ax,x)=(by,y)$ for some $a,b\in A$, $x,y\in M$ with $ax,by\neq 0$. Since $A$ is basic, we may assume without loss of generality that $b\in Aa$, so that $b=ac$ for some non-zero $c\in A$. Since $x=y$ and $ax=by$, we see that $a(1+c)x=(a+b)x=0$. If $c\in I_{>0}$, then $1+c\in A^\x$, so that $ax=0$, a contradiction. So $c\in A^\x$ and then $\nu(a)=\nu(a)+\nu(c)=\nu(b)$. 

All that remains to construct $\psi$ is to check that it gives consistent values if an element of the domain can be expressed in both forms. To that end, suppose that our element is $(ax,x)=(y,by)$ for some $a,b\in A$ and $x,y\in M$ with $ax=by\neq 0$. Then $(1+ab)x=0$, so if $ab\in I_{>0}$, we have that $1+ab\in A^\x$ and then $x=0$, a contradiction. So $ab\in A^\x$, giving that $a,b\in A^\x$. Then $\nu(a)=-\nu(b)=0$, as required. Uniqueness is clear. 

If $x,y,z\in M$ are non-zero, then without loss of generality we may assume that $x=ay$ and $y=bz$ for some non-zero $a,b\in A$. We can then easily show that $$\psi(x,y)+\psi(y,z)=\psi(x,z)$$ and that $\psi(y,x)=-\psi(x,y)$. The final part is routine. \end{proof}

\begin{definition} Let $M$ be a basic $A$-module and $\psi$ be the induced map for $M$ described in Lemma \ref{basic module induces map lemma}. For each non-zero $m_0\in M$, we set $\xi:M\tk\{0\}\to\R$ by $\xi(m)=\psi(m,m_0)$, and $X=\im(\xi)$. \end{definition}

\begin{lem} \label{smallest el case basic is st basic} Let $M$ be a basic module. Then for each non-zero $m_0\in M$, the associated subset $X\subs\R$ is an interval. \end{lem}

\begin{proof} Let $x\leq x'\leq x''$ with $x,x''\in X$, and let $m'=t^{x'-x}m$. We want to show that $m'$ is non-zero and $\xi(m')=x'$, so that $x'\in X$ and then $X$ must be an interval. 

Let $x=\xi(m)$ and $x''=\xi(m'')$ for some non-zero $m,m''\in M$. We see that $$\psi(m'',m)=\psi(m'',m_0)-\psi(m,m_0)=\xi(m'')-\xi(m)=x''-x\geq 0,$$ so that $m''=am$ for some non-zero $a\in A$. Then $$x''-x=\psi(m'',m)=\psi(am,m)=\nu(a).$$ Since $a=t^{\nu(a)}u$ for unique $u\in A^\x$, we see that $$m''=am=t^{\nu(a)}um=t^{x''-x}um=t^{x''-x'}um'.$$ Then since $m''\neq 0$, we see that $m'\neq 0$. Finally, $$x'=x'-x+\xi(m)=\psi(t^{x'-x}m,m)+\psi(m,m_0)=\psi(m',m_0)=\xi(m'),$$ as required. \end{proof}

\begin{rmk} If $m_0'\in M$ is non-zero, then its associated interval $X'$ is simply $X+\{\psi(m_0,m_0')\}$. So $X'$ has a smallest element if and only if $X$ does. \end{rmk} 

We now consider whether or not the associated intervals for a basic module have a smallest element.

\begin{lem} \label{no smallest el case basic is st basic} Let $M$ be a basic module, and $m_0\in M$ be non-zero. Then the associated interval $X$ has a smallest element if and only if $M$ is cyclic. \end{lem}

\begin{proof} Let $m_0\in M$ be non-zero with its associated interval $X$ having a smallest element $\min(X)$. Then $\min(X)=\xi(m_1)$ for some non-zero $m_1\in M$. Let $m\in M$ be non-zero, and then $$\psi(m,m_1)=\psi(m,m_0)+\psi(m_0,m_1)=\psi(m,m_0)-\psi(m_1,m_0)=\xi(m)-\xi(m_1)\geq 0,$$ so that $m\in Am_1$. So $M=Am_1$ is cyclic. 

Conversely, suppose that $M$ is cyclic. Then $M=Am$ for some non-zero $m\in M$, and $\xi(m)\in X$. If $x\in X$, then $x=\xi(m')$ for some non-zero $m'\in M$. Then $$x-\xi(m)=\psi(m',m_0)-\psi(m_0,m)=\psi(m',m)\geq 0$$ since $m'\in M=Am$. So $\xi(m)$ must be the smallest element of $X$. \end{proof}

\begin{lem} Let $M$ be a basic module, and $m_0\in M$ be non-zero. If the associated interval $X$ has no smallest element, then $M$ is isomorphic to a subquotient of $K$ that is not finitely generated. \end{lem}

\begin{proof} Let $m_0\in M$ be non-zero with its associated interval $X$ having no smallest element. We see that $0\in X$ and $\inf(X)\in[-\infty,0)$. Since $X$ is an interval, this shows that $(\inf(X),0]\subs X$. If $\inf(X)=-\infty$ we set $x_n=-n$ for all $n\geq 0$, and if $\inf(X)\in(-\infty,0)$ we set $x_0=0$ and $x_n=\inf(X)(1-2^{-n})$ for $n\geq 1$. In both cases, $(x_n)$ is a strictly decreasing sequence in $X$ with $x_0=0$ and $x_n\to\inf(X)$. 

Note that $\xi(m_0)=x_0$ and suppose that for some $n\geq0$ we have chosen non-zero $m_n\in M$ with $\xi(m_n)=x_n$. Then $x_{n+1}=\xi(m')$ for some non-zero $m'\in M$, and we see that $$\psi(m',m_n)=\psi(m',m_0)-\psi(m_n,m_0)=\xi(m')-\xi(m_n)=x_{n+1}-x_n<0.$$ Since $M$ is basic, this shows that $m_n=am'$ for some non-zero $a\in A$. Then $$x_n-x_{n+1}=\psi(m_n,m')=\psi(am',m')=\nu(a)$$ and finally $$t^{x_n-x_{n+1}}m'=t^{\nu(a)}m'=uam'=um_n$$ for some $w\in A^\x$, and therefore must be non-zero. Then since $\psi(t^{x_n-x_{n+1}}m',m_n)$ is equal to $$\psi(t^{x_n-x_{n+1}}m',m')+\psi(m',m_n)=x_n-x_{n+1}+\psi(m',am')=\nu(a)-\nu(a)=0,$$ we see that $t^{x_n-x_{n+1}}um'=m_n$ for some $u\in A^\x$, and we set $m_{n+1}=um'$. Then $$\xi(m_{n+1})=\xi(um')=\psi(um',m')+\psi(m',m_0)=\nu(u)+\xi(m')=x_{n+1}$$ and $t^{x_n-x_{n+1}}m_{n+1}=m_n$. By induction, we have constructed a sequence $(m_n)_{n\geq0}$ of non-zero elements in $M$, starting with the $m_0$ fixed earlier, and satisfying $\xi(m_n)=x_n$ and \begin{equation} \label{induction step basic iso st basic proof} t^{x_n-x_{n'}}m_{n'}=m_n \end{equation} for all $0\leq n\leq n'$. 

Recall that $$L(X)=\{y\in\R:y\geq x\text{ for some }x\in X\},$$ and set $$V=\{a\in K:\supp(a)\subs L(X)\}.$$ It is routine to check that $A\leq V\leq K$. We now construct a map $\al:V\to M$. Let $a\in V$, and suppose that $\nu(a)\leq x_n$ for all $n\geq 0$. Taking limits, we see that $\nu(a)\leq\inf(X)<0$. So $a\neq 0$ and then $a=t^{\nu(a)}u$ for unique $u\in A^\x$. So $t^{\nu(a)}\in V$, which gives that $\nu(a)\in L(X)$. So $\nu(a)\geq x$ for some $x\in X$. Since $X$ has no minimum, $x>\inf(X)\geq\nu(a)\geq x$, a contradiction. So there exists $n\geq 0$ with $x_n<\nu(a)$, where $\nu(t^{-x_n}a)=\nu(a)-x_n>0$ for any such $n$. If $n,n'\geq 0$ with $x_n,x_{n'}<\nu(a)$, then without loss of generality we may assume that $n\leq n'$, so that $$t^{-x_n}am_n=(t^{-x_n}a)(t^{x_n-x_{n'}}m_{n'})=t^{-x_{n'}}am_{n'}$$ by \eqref{induction step basic iso st basic proof}. We have now defined a map $\al:V\to M$ given by $\al(a)=t^{-x_n}am_n$ for any $n\geq 0$ with $x_n<\nu(a)$. Let $a,b\in V$. It is clear that $\al(a+b)=\al(a)+\al(b)$ when $a=0$ or $b=0$. If $a,b\neq 0$, then choose $n,n'\geq 0$ with $x_n<\nu(a)$ and $x_{n'}<\nu(b)$. Without loss of generality, we may assume that $n\leq n'$. Then $$\nu(a+b)\geq\min\{\nu(a),\nu(b)\}>\min\{x_n,x_{n'}\}=x_{n'},$$ so that $$\al(a+b)=t^{-x_{n'}}(a+b)m_{n'}=t^{-x_{n'}}am_{n'}+t^{-x_{n'}}bm_{n'}=\al(a)+\al(b).$$ It is now easy to see that $\al$ is a morphism; we want to show that it is surjective so that $M\simeq V/\ker(\al)$, a subquotient of $K$. Since $X$ has no smallest element, $L(X)$ is either $\R$ or $(q,\infty)$ for some $q\in\R$, so that $V$ is either $K$ or $I_{>q}$. So $M$ will not be finitely generated. 

Let $m\in M$. We want to find $a\in V$ and $n\geq 0$ with $x_n<\nu(a)$ and $t^{-x_n}am_n=m$. The $m=0$ case is clear. If $m$ is non-zero, then $\xi(m)\in X$ and then $\xi(m)>\inf(X)$. Since $x_n\to\inf(X)$, there exists $n\geq 0$ with $x_n<\xi(m)$. Since $$\psi(m,m_n)=\psi(m,m_0)-\psi(m_n,m_0)=\xi(m)-\xi(m_n)=\xi(m)-x_n>0,$$ we know that $m=bm_n$ for some non-zero $b\in A$. Set $a=bt^{x_n}$. Since $$\supp(a)=\{x_n\}+\supp(b)\subs[x_n,\infty)\subs L(X),$$ we see that $a\in V$. Then since $$x_{n+1}<x_n\leq x_n+\nu(b)=\nu(a),$$ we see that $$m=bm_n=b(t^{x_n-x_{n+1}}m_{n+1})=t^{-x_{n+1}}am_{n+1},$$ as required. \end{proof}

\begin{theorem} \label{basic iso to standard basic} Every basic module is isomorphic to a standard basic module. \end{theorem}

\begin{proof} Let $M$ be a basic module. If the associated intervals for $M$ have smallest elements, then $M$ is cyclic by Lemma \ref{smallest el case basic is st basic}, and therefore isomorphic to a finitely generated standard basic module. Otherwise, it is isomorphic to a subquotient of $K$ that is not finitely generated by Lemma \ref{no smallest el case basic is st basic}. \end{proof}

\begin{lem} \label{hahn field is purely complete} For each sequence $(a_n)$ in $K$ with $a_{n+1}-a_n\in I_n$ for all $n$, there exists $a\in K$ with $a_n-a\in I_n$ for all $n$. \end{lem}

\begin{proof} Let $m,n\geq 1$. Then $$a_m-a_n=(a_{\max\{m,n\}}-a_{\max\{m,n\}-1})+\cdots+(a_{\min\{m,n\}+1}-a_{\min\{m,n\}})\in I_{\min\{m,n\}},$$ so that $\supp(a_m-a_n)\subs[\min\{m,n\},\infty)$. If $q<\min\{m,n\}$, then $(a_m)_q=(a_n)_q$. This allows us to define $a:\R\to\F$ by $a_q=(a_n)_q$ for any positive integer $n>q$. 

Let $r\in\supp(a)$ and $Q=\{q\in\supp(a):q\leq r\}$. Choose a positive integer $n>r$. If $q\in Q$, then $n>q$ so that $(a_n)_q=a_q=1$. This shows that $Q\subs\supp(a_n)$ and therefore that $Q\in\mcW$. Then it follows from Lemma \ref{well-ordered leq lemma} that $\supp(a)\in\mcW$ so that $a\in K$.

Let $n\geq 1$ and $q\in\R$ with $q<n$. Then $(a_n)_q=a_q$ gives that $q\notin\supp(a_n-a)$ so that $\supp(a_n-a)\subs[n,\infty)$, and finally $a_n-a\in I_n$. \end{proof}

\begin{lem} \label{hahn field is purely complete infinitely generated case} For each sequence $(a_n)$ in $K$ with $a_{n+1}-a_n\in I_{>n}$ for all $n$, there exists $a\in K$ with $a_n-a\in I_{>n}$ for all $n$. \end{lem}

\begin{proof} Using Lemma \ref{hahn field is purely complete}, we see that there exists $a\in K$ with $a_n-a\in I_n$ for all $n$. Then $$a_n-a=(a_n-a_{n+1})+(a_{n+1}-a)\in I_{>n}$$ for all $n$. \end{proof}

\begin{lem} \label{lemma for hahn field is reflexive} $\zeta$ is an isomorphism with $D(\zeta)\chi_K=\zeta$. \begin{center}
    \begin{tikzcd} [sep=2cm]
     {} & D^2K \arrow[rd,"D(\zeta)"] & {} \\
     K \arrow[rr,"\zeta"] \arrow[ru,"\chi_K"] & {} & DK \\ 
     \end{tikzcd}\end{center} \end{lem}

\begin{proof} If $a,b\in K$, then $$(D(\zeta)\chi_K)(a)(b)=(\chi_K(a)\zeta)(b)=\chi_K(a)(\zeta(b))=\zeta(b)(a)=\zeta(a)(b),$$ proving the identity. If $a\in K$ and $\zeta(a)=0$, then $\zeta(a)(b)=ab+I_{>0}=I_{>0}$ for all $b\in K$, so that $aK\leq I_{>0}$ and then $a\in K^\circ=0$, proving injectivity.

All that remains to prove that $\zeta$ is surjective. Let $\al\in DK$ and choose a sequence $(a_n)$ in $K$ with $\al(t^{-n})=a_n+I_{>0}$ for all $n$. Then $$ta_{n+1}+I_{>0}=t\al(t^{-n-1})=\al(t^{-n})=a_n+I_{>0},$$ so that $$t^{n+1}a_{n+1}-t^na_n=t^n(ta_{n+1}-a_n)\in t^nI_{>0}=I_{>n},$$ for all $n$. Then we see from Lemma \ref{hahn field is purely complete infinitely generated case} that there exists $a\in K$ with $t^na_n-a\in I_{>n}$ for all $n$,and then $a_n-t^{-n}a\in t^{-n}I_{>n}=I_{>0}$ so that $a_n+I_{>0}=t^{-n}a+I_{>0}$. If $b\in K$, then $b=t^{-n}c$ for some $n\geq 1$ and $c\in A$. Then $$\al(b)=\al(t^{-n}c)=c\al(t^{-n})=c(t^{-n}a+I_{>0})=ab+I_{>0}=\zeta(a)(b)$$ for all $b\in K$. Finally, we have that $\al=\zeta(a)$, proving surjectivity. \end{proof}

Now we derive formulae to describe the actions of various functors on the standard basic modules.

\begin{prop} Let $U\leq K$. Then the unique map $K/U^\circ\to DU$ for which \begin{center}
    \begin{tikzcd}[sep=2cm]
     K \arrow[r,"\zeta"] \arrow[d,twoheadrightarrow,swap] &
     DK \arrow[d,twoheadrightarrow] \\
     K/U^\circ \arrow[r,dashed] & 
     DU
    \end{tikzcd}
 \end{center} commutes is an isomorphism. \end{prop}
 
\begin{proof} We know from Lemma \ref{lemma for hahn field is reflexive} that $\zeta$ is an isomorphism, and since $D$ is exact, the map $DK\to DU$ is surjective. So the composite $K\to DK\to DU$ is surjective, and it is routine to check that its kernel is $U^\circ$. \end{proof}

\begin{prop} \label{lemma for iso for d of basic module} Let $V<U\leq K$. Then the unique map $V^\circ\to D(U/V)$ for which \begin{center}
    \begin{tikzcd}[sep=2cm]
     V^\circ \arrow[r,rightarrowtail] \arrow[d,dashed] &
     K \arrow[d,twoheadrightarrow] \\
     D(U/V) \arrow[r,rightarrowtail] & 
     DU
    \end{tikzcd}
 \end{center} commutes is surjective with kernel $U^\circ$, so that \begin{center}
    \begin{tikzcd}
     0 \arrow[r] & 
     U^\circ \arrow[r] & 
     V^\circ \arrow[r] & 
     D(U/V) \arrow[r] &
     0
    \end{tikzcd}
 \end{center} is exact and $D(U/V)\simeq V^\circ/U^\circ$. \end{prop}

\begin{proof} Since $D$ is exact, we use Proposition \ref{lemma for iso for d of basic module} to show that the following diagram commutes, with all three columns and the bottom two rows forming short exact sequences. \begin{center}
    \begin{tikzcd}[sep=2cm]
     U^\circ \arrow[r,rightarrowtail] \arrow[d,rightarrowtail] &
     V^\circ \arrow[r,twoheadrightarrow] \arrow[d,rightarrowtail] &
     D(U/V) \arrow[d,rightarrowtail] \\
     U^\circ \arrow[r,rightarrowtail] \arrow[d,twoheadrightarrow] & 
     K \arrow[r,twoheadrightarrow] \arrow[d,twoheadrightarrow] &
     DU \arrow[d,twoheadrightarrow] \\
     0 \arrow[r,rightarrowtail] &
     DV \arrow[r,twoheadrightarrow] &
     DV
    \end{tikzcd}
   \end{center} Then it follows from the Nine Lemma that the top row forms a short exact sequence. \end{proof}
   
\begin{exmp} \label{tensor of k and k} The map $$a\tns b\to ab:K\tns K\to K$$ is an isomorphism. \end{exmp}

\begin{proof} If $q>0$, then $$1\tns t^{-q}=t^qt^{-q}\tns t^{-q}=t^{-q}\tns 1=t^{-q}(1\tns 1).$$ This shows that $1\tns 1$ forms a basis for $K\tns K$ over $K$, and then it is routine to check that the specified map is an isomorphism. \end{proof}
   
\begin{lem} \label{tensor of submods of k} Let $L_0,L_1\leq K$. Then there exists a unique isomorphism $L_0\tns L_1\to L_0L_1$ for which \begin{center}
    \begin{tikzcd}[sep=2cm]
     L_0\tns L_1 \arrow[r,dashed] \arrow[d,rightarrowtail] &
     L_0L_1 \arrow[d,rightarrowtail] \\
     K\tns K \arrow[r,"\sim"] & 
     K
    \end{tikzcd}
 \end{center} commutes. \end{lem}
 
\begin{proof} The vertical map on left is the composite \begin{center}
    \begin{tikzcd}
     L_0\tns L_1 \arrow[r] & 
     K\tns L_1 \arrow[r] & 
     K\tns K \;\;\;,
    \end{tikzcd}
 \end{center} which is injective since $L_1$ and $K$ are flat. Since the map $K\tns K\to K$ from Example \ref{tensor of k and k} is an isomorphism mapping $L_0\tns L_1$ perfectly onto $L_0L_1$, we obtain the required isomorphism $L_0\tns L_1\to L_0L_1$. \end{proof}

\begin{prop} \label{formulas for tensor and tor} Let $M_i=L_i/N_i$ with $N_i<L_i\leq K$ for $i=0,1$. Then there exists a unique isomorphism making the following diagram commute. \begin{center}
    \begin{tikzcd}[sep=2cm]
     L_0\tns L_1 \arrow[r,"\sim"] \arrow[d,twoheadrightarrow] &
     L_0L_1 \arrow[d,twoheadrightarrow] \\
     M_0\tns M_1 \arrow[r,dashed] & 
     \frac{L_0L_1}{L_0N_1+N_0L_1}
    \end{tikzcd}
 \end{center} Also, $$\Tor(M_0,M_1)\simeq\frac{N_0L_1\cap L_0N_1}{N_0N_1}$$ and $\Tor_i(M_0,M_1)=0$ for $i\geq 2$. \end{prop}

\begin{proof} Since, for each $i$, \begin{center}
    \begin{tikzcd}
     0 \arrow[r] & 
     N_i \arrow[r] & 
     L_i \arrow[r] &
     M_i \arrow[r] &
     0 \;\;\;,
    \end{tikzcd}
 \end{center} is exact with $L_i$ and $N_i$ flat, we tensor these short exact sequences with $L_i$, $N_i$ and $M_i$ to obtain the following commutative diagram, with the first two rows and the first two columns forming short exact sequences, and the final row and final column forming right exact sequences. \begin{center}
    \begin{tikzcd}[sep=2cm]
     N_0\tns N_1 \arrow[r,rightarrowtail] \arrow[d,rightarrowtail] &
     N_0\tns L_1 \arrow[r,twoheadrightarrow] \arrow[d,rightarrowtail] &
     N_0\tns M_1 \arrow[d] \\
     L_0\tns N_1 \arrow[r,rightarrowtail] \arrow[d,twoheadrightarrow] & 
     L_0\tns L_1 \arrow[r,twoheadrightarrow] \arrow[d,twoheadrightarrow] &
     L_0\tns M_1 \arrow[d,twoheadrightarrow] \\
     M_0\tns N_1 \arrow[r] &
     M_0\tns L_1 \arrow[r,twoheadrightarrow] &
     M_0\tns M_1
    \end{tikzcd}
   \end{center} 
   
Let $$U=\im(L_0\tns N_1\to L_0\tns L_1)+\im(N_0\tns L_1\to L_0\tns L_1).$$ It is clear that $$U\leq\ker(L_0\tns L_1\to M_0\tns M_1),$$ giving a map \begin{equation} \label{map for tensor product formula} (x\tns y)+U\mapsto(x+N_0)\tns(y+N_1):\frac{L_0\tns L_1}{U}\to M_0\tns M_1. \end{equation} We now construct an inverse for this map. Let $x,x'\in L_0$, $y,y'\in L_1$ with $x-x'\in N_0$ and $y-y'\in N_1$. Then $$x\tns y-x'\tns y'=x\tns(y-y')+(x-x')\tns y\in U,$$ and one can check that the map $$(x+N_0)\tns(y+N_1)\mapsto(x\tns y)+U:M_0\tns M_1\to\frac{L_0\tns L_1}{U}$$ is an inverse for \eqref{map for tensor product formula}. So we see that $$M_0\tns M_1\simeq\frac{L_0\tns L_1}{U}.$$ Using Lemma \ref{tensor of submods of k}, we see that the isomorphism $L_0\tns L_1\to L_0L_1$ maps $\ker(L_0\tns L_1\to M_0\tns M_1)$ perfectly onto $L_0N_1+N_0L_1$, and therefore induces an isomorphism $$\frac{L_0\tns L_1}{U}\simeq\frac{L_0L_1}{L_0N_1+N_0L_1},$$ as required. Since \begin{center}
    \begin{tikzcd}
     0 \arrow[r] & 
     N_0 \arrow[r] & 
     L_0 \arrow[r] &
     M_0 \arrow[r] &
     0
    \end{tikzcd}
 \end{center} is exact and $\Tor(L_0,M_1)=0$ (because $L_0$ is flat), we obtain from the long exact sequence of Tor modules that \begin{center}
    \begin{tikzcd}
     0 \arrow[r] & 
     \Tor(M_0,M_1) \arrow[r] & 
     N_0\tns M_1 \arrow[r] &
     L_0\tns M_1 \arrow[r] &
     M_0\tns M_1 \arrow[r] &
     0
    \end{tikzcd}
 \end{center} is exact. So $$\Tor(M_0,M_1)\simeq\ker(N_0\tns M_1\to L_0\tns M_1).$$ Expressing $N_0\tns M_1$ and $L_0\tns M_1$ as standard basic modules, we see that $\Tor(M_0,M_1)$ is isomorphic to the kernel of the unique map for which the following diagram commutes. \begin{center}
    \begin{tikzcd}[sep=2cm]
     N_0\tns M_1 \arrow[r] \arrow[d] &
     L_0\tns M_1 \arrow[d] \\
     \frac{N_0L_1}{N_0N_1} \arrow[r,dashed] & 
     \frac{L_0L_1}{L_0N_1}
    \end{tikzcd}
 \end{center} We also obtain from the long exact sequences of Tor modules that \begin{center}
    \begin{tikzcd}
     \Tor_i(L_0,M_1) \arrow[r] & 
     \Tor_i(M_0,M_1) \arrow[r] &
     \Tor_{i-1}(N_0,M_1)
    \end{tikzcd}
 \end{center} is exact for $i\geq 2$. Then since $L_0$ and $N_0$ are flat and $i,i-1\geq 1$, we know that $$\Tor_i(L_0,M_1)=\Tor_{i-1}(N_0,M_1)=0$$ and then $\Tor_i(M_0,M_1)=0$, as required. \end{proof}
   
\begin{prop} \label{formulas for hom and ext} Let $M_i=L_i/N_i$ with $N_i<L_i\leq K$ for $i=0,1$. There is an isomorphism $$\frac{(L_1:L_0)\cap(N_1:N_0)}{(N_1:L_0)}\to\Hom(M_0,M_1)$$ given by $$a+(N_1:L_0)\mapsto(b+N_0\mapsto ab+N_1).$$ Also, $$\Ext(M_0,M_1)\simeq\frac{(L_1:N_0)}{(L_1:L_0)+(N_1:N_0)}$$ and $\Ext^i(M_0,M_1)=0$ for $i\geq 2$. \end{prop}

\begin{proof} Using the formulas for $D$ and $(-)\tns(-)$ acting on basic modules, together with the tensor-hom adjunction, we establish isomorphisms $$\Hom(M_0,M_1)\simeq\Hom\left(M_0,D\left(\frac{N_1^\circ}{L_1^\circ}\right)\right)\simeq D\left(M_0\tns\frac{N_1^\circ}{L_1^\circ}\right)\simeq D\left(\frac{L_0N_1^\circ}{L_0L_1^\circ+N_0N_1^\circ}\right).$$ Using the formula for $D$, we have another isomorphism from the last of these to \begin{equation} \label{expression for hom formula} \frac{(L_0L_1^\circ+N_0N_1^\circ)^\circ}{(L_0N_1^\circ)^\circ}. \end{equation} Using Proposition \ref{combining submods of k}, we see that if $L_0L_1^\circ\leq N_0N_1^\circ$, then $$(N_1:N_0)=(N_0N_1^\circ)^\circ\leq(L_0L_1^\circ)^\circ=(L_1:L_0).$$ so the expression \eqref{expression for hom formula} is equal to $$\frac{(N_0N_1^\circ)^\circ}{(L_0N_1^\circ)^\circ}=\frac{(N_1:N_0)}{(N_1:L_0)}=\frac{(L_1:L_0)\cap(N_1:N_0)}{(N_1:L_0)}.$$ Otherwise, $N_0N_1^\circ\leq L_0L_1^\circ$ and we use a similar argument to arrive at the same formula. One can check that the composite isomorphism is the required one.
     
Now for the Ext formula. Since \begin{center}
    \begin{tikzcd}
     0 \arrow[r] & 
     M_1 \arrow[r] & 
     K/N_1 \arrow[r] &
     K/L_1 \arrow[r] &
     0
    \end{tikzcd}
 \end{center} is an injective resolution for $M_1$, we see that $\Ext(M_0,M_1)$ is the cohomology of the cochain complex \begin{equation} \label{cochain complex ext formula} 
    \begin{tikzcd}
     0 \arrow[r] & 
     \Hom(M_0,K/N_1) \arrow[r] & 
     \Hom(M_0,K/L_1) \arrow[r] &
     0
    \end{tikzcd}
     \end{equation} at $\Hom(M_0,K/L_1)$ i.e. $$\Ext(M_0,M_1)\simeq\cok(\Hom(M_0,K/N_1\twoheadrightarrow K/L_1)).$$  Then expressing $\Hom(M_0,K/N_1)$ and $\Hom(M_0,K/L_1)$ as standard basic modules, we see that $\Ext(M_0,M_1)$ is isomorphic to the cokernel of the unique map making the following diagram commute. \begin{center}
    \begin{tikzcd}[sep=2cm]
     \Hom(M_0,K/N_1) \arrow[r] \arrow[d] &
     \Hom(M_0,K/L_1) \arrow[d] \\
     \frac{(N_1:N_0)}{(N_1:L_0)} \arrow[r,dashed] & 
     \frac{(L_1:N_0)}{(L_1:N_0)}
    \end{tikzcd}
 \end{center} The image of this map is $$\frac{(L_1:L_0)+(N_1:N_0)}{(L_1:L_0)},$$ so that $$\Ext(M_0,M_1)\simeq\frac{(L_1:N_0)}{(L_1:L_0)+(N_1:N_0)},$$ as required. For $i\geq 2$, $\Ext^i(M_0,M_1)$ is the cohomology of the cochain complex \eqref{cochain complex ext formula} at $0$, so must be $0$. \end{proof}

\begin{prop} \label{functor tables} By applying $D$, $(-)\tns(-)$, $\Hom(-,-)$, $\Tor(-,-)$ and $\Ext(-,-)$ to standard basic modules using the formulas in Propositions \ref{formulas for tensor and tor} and \ref{formulas for hom and ext}, we obtain the following tables. We set $r=\min\{p,q\}$. Also \[  X_{pq}= \left\{
\begin{array}{ll}
      I_p & p\leq q \\
      I_{>q} & p>q
\end{array} 
\right. \] and \[ Y_{pq}=  \left\{
\begin{array}{ll}
      I_{>0} & p\leq q \\
      A & p>q
\end{array} 
\right. \] Note that these are related by $I_rX_{pq}^\circ=Y_{pq}$ and $I_rY_{pq}^\circ=X_{pq}$.

\[ \begin{array}{|c|c|c|c|c|c|c|c|c|c|c|}\hline
     M       & K & A & I_{>0} & \Ta & \Phi & \F & A/I_q & A/I_{>q} & I_{>0}/I_q & I_{>0}/I_{>q} \\ \hline 
     DM    & K & \Ta & \Phi      & A   & I_{>0} & \F & I_{>0}/I_{>q}     & A/I_{>q}        & I_{>0}/I_q          & A/I_q             \\ \hline
    \end{array}
 \]

 \[ \begin{array}{|c|c|c|c|c|c|c|c|c|c|c|}\hline M\tns N      & K & A & I_{>0} & \Ta & \Phi & \F & A/I_q & A/I_{>q} & I_{>0}/I_q & I_{>0}/I_{>q} \\ \hline 
     K         & K & K & K      & 0   & 0  & 0 & 0     & 0        & 0          & 0             \\ \hline 
     A         & K & A & I_{>0}      & \Ta   & \Phi & \F  & A/I_q     & A/I_{>q}        & I_{>0}/I_q          & I_{>0}/I_{>q}             \\ \hline 
     I_{>0}    & K & I_{>0} & I_{>0}      & \Ta   & \Ta  & 0  & I_{>0}/I_{>q}     & I_{>0}/I_{>q}        & I_{>0}/I_{>q}          & I_{>0}/I_{>q}             \\ \hline 
     \Ta       & 0 & \Ta & \Ta      & 0   & 0  & 0 & 0     & 0        & 0          & 0 \\ \hline 
     \Phi       & 0 & \Phi & \Ta      & 0   & 0  &  0   & 0   & 0   & 0        & 0           \\ \hline
     \F      &  0  &  \F  &  0  &  0  &  0  &  \F  &  \F  &  \F  &  0  &  0  \\ \hline
     A/I_p 	  & 0 & A/I_	p & I_{>0}/I_{>p}      &	 0	&	0  & \F &  A/I_r    &    A/X_{pq}     &     I_{>0}/X_{qp}      &	 I_{>0}/I_{>r}									\\ \hline A/I_{>p}
     & 0 & A/I_{>p} & I_{>0}/I_{>p}      &	  0	&	0 & \F &  A/X_{qp}    &    A/I_{>r}     &     I_{>0}/X_{qp} &	  I_{>0}/I_{>r}												 \\ \hline I_{>0}/I_p 
    & 0 & I_{>0}/I_p & I_{>0}/I_{>p}      &	  0 	&	0  & 0 &  I_{>0}/X_{pq}    &    I_{>0}/X_{pq} &   I_{>0}/I_{>r} &	 I_{>0}/I_{>r}													\\ \hline   I_{>0}/I_{>p}
    & 0 & I_{>0}/I_{>p} & I_{>0}/I_{>p}      &	  0	&	0  & 0 &   I_{>0}/I_{>r}  &   I_{>0}/I_{>r}     &    I_{>0}/I_{>r}    &	  I_{>0}/I_{>r}													 \\ \hline      
    \end{array} \]
    
 \[ \begin{array}{|c|c|c|c|c|c|c|c|c|c|c|}\hline
\Hom(M,N)      & K & A & I_{>0} & \Ta & \Phi & \F & A/I_q & A/I_{>q} & I_{>0}/I_q & I_{>0}/I_{>q} \\ \hline 
     K         & K & 0 & 0      & K   & K  & 0  & 0     & 0        & 0          & 0             \\ \hline 
     A         & K & A & I_{>0}      & \Ta   & \Phi  & \F & A/I_q     & A/I_{>q}        & I_{>0}/I_q          & I_{>0}/I_{>q}             \\ \hline 
     I_{>0}    & K & A & A      & \Phi   & \Phi & 0  & A/I_q     & A/I_q        & A/I_q          & A/I_q             \\ \hline 
     \Ta       & 0 & 0 & 0      & A   & A & 0  & 0     & 0        & 0          & 0             \\ \hline 
     \Phi       & 0 & 0 & 0      & I_{>0}   & A & 0  & 0   & 0   & 0        & 0           \\ \hline
     \F     &  0  &  0  &  0  &  \F  &  0 & \F &  0  &  \F  &  0  &  \F  \\ \hline
     A/I_p 	  & 0 & 0 & 0      &	  I_{>0}/I_{>p}	&	A/I_p  & \F &  A/I_r    &   Y_{pq}/I_{>r}     &     Y_{qp}/I_r      &	  I_{>0}/I_{>r}										\\ \hline A/I_{>p}
     & 0 & 0 & 0      &	  A/I_{>p}	&	A/I_p & \F &  A/I_r    &    A/I_{>r}     &     Y_{qp}/I_r      &	  Y_{qp}/I_{>r}														 \\ \hline I_{>0}/I_p 
    & 0 & 0 & 0      &	  I_{>0}/I_p 	&	A/I_p & 0 &  A/I_r    &    Y_{pq}/I_r     &     A/I_r      &	  Y_{pq}/I_r 														\\ \hline   I_{>0}/I_{>p}
    & 0 & 0 & 0      &	  A/I_p	&	A/I_p & 0 &  A/I_r    &    A/I_r     &     A/I_r      &	  A/I_r 														 \\ \hline      
    \end{array}
 \]
 
 \[ \begin{array}{|c|c|c|c|c|c|c|c|c|c|c|}\hline \Tor(M,N) & K & A & I_{>0} & \Ta & \Phi & \F & A/I_q & A/I_{>q} & I_{>0}/I_q & I_{>0}/I_{>q} \\ \hline

K & 0 & 0 & 0 & 0 & 0 & 0 & 0 & 0 & 0 & 0 \\ \hline

A & 0 & 0 & 0 & 0 & 0 & 0 & 0 & 0 & 0 & 0 \\ \hline

I_{>0} & 0 & 0 & 0 & 0 & 0 & 0 & 0 & 0 & 0 & 0 \\ \hline 

     \Ta     & 0 & 0 & 0        & \Ta   & \Ta  &  0   & I_{>0}/I_{>q}     & I_{>0}/I_{>q}        & I_{>0}/I_{>q} & I_{>0}/I_{>q} \\ \hline
      
     \Phi      & 0 & 0 & 0       & \Ta   & \Phi  &  \F   & A/I_q     & A/I_{>q}        & I_{>0}/I_q          & I_{>0}/I_{>q}           \\ \hline
     
     \F      &  0  &  0  &  0  &  0  &  \F  &  0  &  \F  &  0  &  \F  &  0  \\ \hline
     A/I_p 	 & 0 & 0 & 0      &	 I_{>0}/I_{>p} 	&  A/I_p  &  \F  &  A/I_r    &    Y_{pq}/I_{>r}     &     Y_{qp}/I_r      &	 I_{>0}/I_{>r}									\\ \hline 
     A/I_{>p}
     & 0 & 0 & 0    &	  I_{>0}/I_{>p}	&	A/I_{>p} &  0  &  Y_{qp}/I_{>r}    &     I_{>0}/I_{>r}    &     Y_{qp}/I_{>r} &	 I_{>0}/I_{>r}												 \\ \hline 
         I_{>0}/I_p  & 0 & 0 & 0 
        &	  I_{>0}/I_{>p} 	&	I_{>0}/I_p  &  \F  &  Y_{pq}/I_r    &    Y_{pq}/I_{>r} &   I_{>0}/I_r &	 I_{>0}/I_{>r}													\\ \hline   
        I_{>0}/I_{>p} & 0 & 0 & 0 
        &	  I_{>0}/I_{>p}	&	I_{>0}/I_{>p}  &  0  &   I_{>0}/I_{>r}  &   I_{>0}/I_{>r}     &    I_{>0}/I_{>r}    &	  I_{>0}/I_{>r}													 \\ \hline
              
    \end{array} \]
    
\[ \begin{array}{|c|c|c|c|c|c|c|c|c|c|c|}\hline \Ext(M,N)      & K & A & I_{>0} & \Ta & \Phi & \F & A/I_q & A/I_{>q} & I_{>0}/I_q & I_{>0}/I_{>q} \\ \hline 
     K         & 0 & 0 & 0      & 0   & 0  & 0  & 0     & 0        & 0          & 0             \\ \hline 
     A         & 0 & 0 & 0      & 0   & 0  & 0  & 0     & 0        & 0          & 0             \\ \hline 
     I_{>0}    & 0 & 0 & 0      & 0   & 0  & 0  & 0     & 0        & 0          & 0             \\ \hline 
     \Ta       & 0 & A & A      & 0   & 0  & 0  & A/I_q     & A/I_q        & A/I_q          & A/I_q \\ \hline 
     \Phi       & 0 & A & I_{>0}      & 0   & 0  & \F  & A/I_q     & A/I_{>q}        & I_{>0}/I_q          & I_{>0}/I_{>q}           \\ \hline
     \F      &  0  &  0  &  \F  &  0  &  0  &  0  &  0  &  0  &  \F  &  \F  \\ \hline
     A/I_p 	  & 0 & A/I_p & I_{>0}/I_{>p}      &	 0	&  0  & \F &  A/I_r    &    A/X_{pq}     &     I_{>0}/X_{qp}      &	 I_{>0}/I_{>r}									\\ \hline A/I_{>p}
     & 0 & A/I_p & A/I_{>p}      &	  0	&	0 & 0 &  A/I_r    &    A/I_r     &     A/X_{qp} &	 A/X_{qp}												 \\ \hline I_{>0}/I_p 
    & 0 & A/I_p & I_{>0}/I_p      &	  0 	&	0 & \F  &  A/I_r    &    A/X_{pq} &   I_{>0}/I_r &	 I_{>0}/X_{pq}													\\ \hline   I_{>0}/I_{>p}
    & 0 & A/I_p & A/I_p      &	  0	&	0 & 0  &   A/I_r  &   A/I_r     &    A/I_r    &	  A/I_r													 \\ \hline      
    \end{array} \]
    
\end{prop}

\begin{rmk} \label{multibasics in category z rmk} We note the following from Proposition \ref{functor tables}. If $M$ is basic, then $K\tns M\cong K$ if and only if $M$ is flat. Also $K\tns DM\cong K$ if and only if $M$ is injective. So if $M$ is multibasic, then the dimensions of $K\tns M$ and $K\tns DM$ over $K$ are the number of flat summands and injective summands, respectively. So $K\tns M=K\tns DM=0$ if and only if $M$ has no flat or injective summands i.e. no copies of $K$, $A$, $I_{>0}$, $\Ta$ or $\Phi$. We conclude that the multibasic modules in the subcategory $\mcZ$ are precisely those whose summands come from $\F$, $A/I_q$, $A/I_{>q}$, $I_{>0}/I_q$, $I_{>0}/I_{>q}$ for $q>0$. \end{rmk}

\begin{exmp} \label{first proj res exmps} Set $P=\bigoplus_{i=0}^\infty A$. Define $\phi:P\to K$ by $\phi(e_i)=t^{-i}$, $\psi:P\to I_{>0}$ by $\psi(e_i)=t^{2^{-i}}$, $\al:P\to P$ by $\al(e_i)=e_i+te_{i+1}$ and $\bt:P\to P$ by $\bt(e_i)=e_i+t^{2^{-i-1}}e_{i+1}$ for all $i$. Then \begin{center}
    \begin{tikzcd}
     0 \arrow[r] & 
     P \arrow[r,"\al"] & 
     P \arrow[r,"\phi"] &
     K \arrow[r] &
     0
    \end{tikzcd}
 \end{center} and \begin{center}
    \begin{tikzcd}
     0 \arrow[r] & 
     P \arrow[r,"\bt"] & 
     P \arrow[r,"\psi"] &
     I_{>0} \arrow[r] &
     0
    \end{tikzcd}
 \end{center} are exact. \end{exmp}

\begin{proof} If $\sum_ia_ie_i\in\ker(\al)$, then $\sum_ia_i(e_i+te_{i+1})=0$, so that $a_0=0$ and $a_i=ta_{i-1}$ for $i\geq 1$. So $a_i=t^ia_0=0$ for all $i$ and then $\sum_ia_ie_i=0$, and we conclude that $\al$ is injective. 

If $\sum_ia_ie_i\in\ker(\bt)$, then $\sum_ia_i(e_i+t^{2^{-i-1}}e_{i+1})=0$, so that $a_0=0$ and $a_i=t^{2^{-i}}a_{i-1}$ for $i\geq 1$. So $a_i=t^{1-2^{-i}}a_0=0$ for all $i$ and then $\sum_ia_ie_i=0$, and we conclude that $\bt$ is injective.

If $a\in K$, then there exists $i\geq 0$ with $\nu(a)\geq -i$. So $a=t^{-i}b$ for some $b\in A$, and $a=b\phi(e_i)=\phi(be_i)\in\im(\phi)$, showing that $\phi$ is surjective. 

If $a\in I_{>0}$, then there exists $i\geq 0$ with $\nu(a)\geq 2^{-i}$. So $a=t^{2^{-i}}b$ for some $b\in A$. Then $a=b\psi(e_i)=\psi(be_i)\in\im(\psi)$, showing that $\psi$ is surjective. 

For each $i$, we have that $$\phi(\al(e_i))=\phi(e_i+te_{i+1})=t^{-i}+tt^{-i-1}=0.$$ So $\phi(\al(x))=0$ for all $x\in P$, giving that $\im(\al)\leq\ker(\phi)$. Similarly, $$\psi(\bt(e_i))=\psi(e_i+t^{2^{-i-1}}e_{i+1})=t^{2^{-i}}+t^{2^{-i-1}}t^{2^{-i-1}}=0$$ for all $i$, so that $\im(\bt)\leq\ker(\psi)$. 

Now let $\sum_ia_ie_i\in\ker(\phi)$. Then $\sum_ia_it^{-i}=0$. If $a_i=0$ for all $i$, then clearly $\sum_ia_ie_i\in\im(\al)$. Otherwise, let $n\geq 0$ be largest such that $a_n\neq 0$. Then $$a_0+a_1t^{-1}+\cdots+a_nt^{-n}=0$$ so that $$a_0t^n+a_1t^{n-1}+\cdots+a_n=0.$$ Let $$b_i=a_i+a_{i-1}t+\cdots+a_0t^i$$ for all $i$. Then $b_{i+1}=a_{i+1}+b_it$ for all $i\geq 0$. We see that $b_n=0$, and then since $a_i=0$ for $i>n$ we see that $b_i=a_i+b_{i-1}t=b_{i-1}t$ for $i>n$ so that $b_i=b_nt^{i-n}=0$ for $i\geq n$ and finally $\sum_ib_ie_i\in P$. Then applying $\al$ to this shows that $$\sum_ia_ie_i=\al(\sum_ib_ie_i)\in\im(\al),$$ and we conclude that $\im(\al)=\ker(\phi)$, so that \begin{center}
    \begin{tikzcd}
     0 \arrow[r] & 
     P \arrow[r,"\al"] & 
     P \arrow[r,"\phi"] &
     K \arrow[r] &
     0
    \end{tikzcd}
 \end{center} is exact.

Now let $\sum_ia_ie_i\in\ker(\psi)$. Then $\sum_ia_it^{2^{-i}}=0$. If $a_i=0$ for all $i$, then clearly $\sum_ia_ie_i\in\im(\bt)$. Otherwise, let $n\geq 0$ be largest such that $a_n\neq 0$. Then $$ta_0+t^{2^{-1}}a_1+\cdots+t^{2^{-n}}a_n=0.$$ Noting that $$2^{-i}-2^{-n}=2^{-i-1}+\cdots+2^{-n}$$ for $i<n$, we multiply through by $t^{-2^{-n}}$ to see that $$a_n+t^{2^{-n}}a_{n-1}+\cdots+t^{2^{-1}+\cdots+2^{-n}}a_0=0.$$ Set $$b_i=a_i+t^{2^{-i}}a_{i-1}+\cdots+t^{2^{-1}+\cdots+2^{-i}}a_0$$ for $i\geq 0$. Then $a_0=b_0$ and $a_i=b_i+t^{2^{-i}}b_{i-1}$ for $i\geq 1$. We see that $b_n=0$, and since $a_i=0$ for $i>n$, we also see that $$b_i=t^{2^{-i}}b_{i-1}=\cdots=t^{2^{-i}+\cdots+2^{-n-1}}b_n=0$$ for all $i>n$, giving that $\sum_ib_ie_i\in P$. Since $a_0=b_0$ and $a_i=b_i+t^{2^{-i}}b_{i-1}$ for $i\geq 1$, we see that $\sum_ia_ie_i=\bt(\sum_ib_ie_i)\in\im(\bt)$, so that $\ker(\psi)=\im(\bt)$, and we conclude that \begin{center}
    \begin{tikzcd}
     0 \arrow[r] & 
     P \arrow[r,"\bt"] & 
     P \arrow[r,"\psi"] &
     I_{>0} \arrow[r] &
     0
    \end{tikzcd}
 \end{center} is exact.
\end{proof}

\begin{exmp} Using the notation from Example \ref{first proj res exmps}, and defining $\gm:A\to P$ by $\gm(a)=ae_0$, we have that $$0\to P \xrightarrow{\renewcommand{\arraystretch}{1}\bbm\bt\\ 0\ebm} P\oplus P \xrightarrow{\renewcommand{\arraystretch}{1}\bbm \gm|_{I_{>0}}\psi & \al\ebm} P\to\Ta\to 0$$ is a projective resolution for $\Ta$. \end{exmp}

\begin{proof} Since $\bt$ is injective, it is clear that $\renewcommand{\arraystretch}{1}\bbm\bt\\ 0\ebm$ is injective, and it has image $\im(\bt)\oplus 0=\ker(\psi)\oplus 0$, which is clearly contained within the kernel of the next map. Let $(a,b)\in\ker\left(\renewcommand{\arraystretch}{1}\bbm \gm|_{I_{>0}}\psi & \al\ebm\right)$, with $a=\sum_i a_ie_i$ and $b=\sum_i b_ie_i$. Then $\gm(\psi(a))=\al(b)$, and we see that $b_i=t^i\psi(a)$ for all $i$. Since $b\in P$, when $i$ is sufficiently large we have $t^i\psi(a)=b_i=0$, so that $\psi(a)=0$ and then $b=0$ and $a\in\im(\bt)$. We conclude that the sequence is exact at $P\oplus P$. 

If $a,b\in P$, then $$\phi(\gm(\psi(a))+\phi(\al(b))=\phi(\gm(\psi(a))=\psi(a)\in I_{>0},$$ so that $$\im\left(\renewcommand{\arraystretch}{1}\bbm \gm|_{I_{>0}}\psi & \al\ebm\right)\leq\ker(P\to\Ta).$$ If $a\in\ker(P\to\Ta)$, then $\phi(a)\in I_{>0}=\im(\psi)$, so that $$\phi(a)=\psi(c)=\phi(\gm(\psi(c)))$$ for some $c\in P$. Then $a-\gm(\psi(c))\in\ker(\phi)=\im(\al)$, so that $a\in\im\left(\renewcommand{\arraystretch}{1}\bbm \gm|_{I_{>0}}\psi & \al\ebm\right)$. Finally, the map $P\to\Ta$ is surjective, and we are done. \end{proof}

\begin{exmp} Using the notation from Example \ref{first proj res exmps}, we have that $$0\to P \xrightarrow{\renewcommand{\arraystretch}{1}\bbm\bt\\ t^q\ebm} P\oplus P \xrightarrow{\renewcommand{\arraystretch}{1}\bbm t^q & \bt\ebm} P\to I_{>0}/I_{>q}\to 0$$ is a projective resolution for $I_{>0}/I_{>q}$. \end{exmp}

\begin{proof} It is routine to check that $\renewcommand{\arraystretch}{1}\bbm\bt\\ t^q\ebm$ is injective with its image contained in $\ker\left(\renewcommand{\arraystretch}{1}\bbm t^q & \bt\ebm\right)$. Now let $(a,b)\in\ker\left(\renewcommand{\arraystretch}{1}\bbm t^q & \bt\ebm\right)$ with $a=\sum_i a_ib_i$ and $b=\sum_ib_ie_i$. Then $\bt(b)=t^qa$, and we see that $t^qa_0=b_0$, and $t^qa_i=b_i+b_{i-1}t^{2^{-i}}$ for $i\geq 1$. From this, we see by induction that $t^{-q}b_i\in A$ for all $i$. Let $b'=\sum_it^{-q}b_ie_i\in P$, and then $b=t^qb'$ with $t^q\bt(b')=t^qa$. So $a=\bt(b')$ and then $(a,b)\in\im\left(\renewcommand{\arraystretch}{1}\bbm\bt\\ t^q\ebm\right)$.

It is routine to check that $$\im(\renewcommand{\arraystretch}{1}\bbm t^q & \bt\ebm)\leq\ker(P\to I_{>0}/I_{>q}).$$ Now let $a\in\ker(P\to I_{>0}/I_{>q})$. Then $$\psi(a)\in I_{>q}=t^qI_{>0}=t^q\im(\psi),$$ so that $\psi(a)=t^q\psi(c)=\psi(t^qc)$ for some $c\in P$. Then $a-t^qc\in\ker(\psi)=\im(\bt)$, so that $a=t^qc+\bt(b)\in\im\left(\renewcommand{\arraystretch}{1}\bbm t^q & \bt\ebm\right)$ for some $b\in P$. Then since the map $P\to I_{>0}/I_{>q}$ is surjective, we are done. \end{proof}

The projective resolutions for the remaining standard basic modules can be obtained similarly to the above examples, and are tabulated below. We use the notation from Example \ref{first proj res exmps}, and for $q>0$ define $j_q:A\to P$ by $j_q(a)=t^{q-2^{-k}}ae_k$, where $k$ is the smallest non-negative integer for which $2^{-k}\leq q$.

  \[ \setlength{\tabcolsep}{20pt}
     \renewcommand{\arraystretch}{2}
     \begin{array}{|l|l||l|l|}\hline 
      M & \text{projective resolution of } M & 
      M & \text{projective resolution of } M \\ \hline
      K             & P \xrightarrow{\al} P \xrightarrow{\phi} K & \F & P\xrightarrow{\bt} P\xrightarrow{\io\psi} A\to \F \\ \hline
      A        & A\xrightarrow{\id_A} A &
      A/I_q         & A\xrightarrow{t^q} A \to A/I_q \\ \hline
      I_{>0}             & P \xrightarrow{\bt} P \xrightarrow{\psi} I_{>0} &
      A/I_{>q}      & P \xrightarrow{\bt} P \xrightarrow{t^q\psi} A \\ \hline
      \Ta          & P \xrightarrow{\renewcommand{\arraystretch}{1}\bbm\bt\\ 0\ebm} P\oplus P \xrightarrow{\renewcommand{\arraystretch}{1}\bbm \gm|_{I_{>0}}\psi & \al\ebm} P &
      I_{>0}/I_q    & A\oplus P \xrightarrow{\renewcommand{\arraystretch}{1}\bbm j_q & \bt \ebm} P \\ \hline 
      \Phi          & A\oplus P \xrightarrow{\renewcommand{\arraystretch}{1}\bbm \gm & \al\ebm} P &
      I_{>0}/I_{>q} & P \xrightarrow{\renewcommand{\arraystretch}{1}\bbm\bt\\ t^q\ebm} P\oplus P \xrightarrow{\renewcommand{\arraystretch}{1}\bbm t^q & \bt\ebm} P \\ \hline
    \end{array}
  \]
  
\begin{rmk} Since $A$ is local and $D(\bigoplus_{i\in I}A)\simeq\prod_{i\in I}\Ta$ for any set $I$, we see that $D$ sends projective objects to injective objects. Then since $D$ is exact, we can apply it a projective resolution of $M$ to obtain an injective resolution of $DM$. This allows us to find injective resolutions for the standard basic modules using the above table of projective resolutions. Also note that \begin{center}
    \begin{tikzcd}
     0 \arrow[r] & 
     U/V \arrow[r] & 
     K/V \arrow[r] &
     K/U \arrow[r] &
     0
    \end{tikzcd}
 \end{center} is an injective resolution of $U/V$ for any $V<U\leq K$. \end{rmk}

\section{Multibasic Modules} 

We now discuss multibasic modules, and will show in Theorem \ref{theta-ref implies multibasic} that these are simply the objects of $\mcM$ i.e. the $\Ta$-reflexive modules. Using the properties of $\mcM$, we will then see in Corollary \ref{multibasics ab and closed under exts subs quos} that $\mcM$ is abelian, and is closed under submodules, quotients and extensions. 

In this section, we classify the $A$-submodules of finite-dimensional $K$-vector spaces, showing that they are all flat and multibasic. We also investigate $\Ext(M,N)$, where $M$ and $N$ are multibasic and at least one of them is finite. Firstly, we recall the definition of a multibasic module.

\begin{definition} We say that an $A$-module $M$ is \textit{multibasic} if it can be decomposed as the direct sum of a finite (possibly empty) list of basic submodules. It is clear that the category of multibasic modules is closed under finite direct sums and isomorphisms. 

We say that a multibasic module $M$ \textit{involves} a collection $\mcS$ of standard basic modules if it has a finite decomposition into basic submodules, each of which is isomorphic to a standard basic module in $\mcS$. \end{definition}

\begin{prop} If $M$ and $N$ are multibasic, then so are $M\tns N$, $\Hom(M,N)$, $\Tor(M,N)$ and $\Ext(M,N)$. In particular, $DM=\Hom(M,\Ta)$ is multibasic. \end{prop}

\begin{proof} If $M$ and $N$ are basic, then it is clear from Proposition \ref{functor tables} that all four of these modules are either basic or zero. If $M$ and $N$ are multibasic, then since all four of these functors are additive and multibasic modules are finite direct sums of basic modules, all four modules must be multibasic. \end{proof}

\begin{lem} \label{extensions flat multibasic split lemma} Let $N$ be multibasic. Then $N$ is flat if and only if for every multibasic $L$, we have that every short exact sequence \begin{center}
    \begin{tikzcd}
     0 \arrow[r] & 
     L \arrow[r] & 
     M \arrow[r] &
     N \arrow[r] &
     0
    \end{tikzcd}
 \end{center} splits. \end{lem}

\begin{proof} First suppose that $N$ is flat, and let \begin{equation} \label{short exact flat mb splitting seq}
    \begin{tikzcd}
     0 \arrow[r] & 
     L \arrow[r] & 
     M \arrow[r] &
     N \arrow[r] &
     0
    \end{tikzcd}
 \end{equation} be exact with $L$ multibasic. Since $N$ is flat, it involves only $K$, $A$ and $I_{>0}$, the flat standard basic modules. Then we see from Proposition \ref{functor tables} that $\Ext(V,X)=0$ when $V$ and $X$ are basic and $V$ is flat. Since $\Ext$ is additive, we see that $\Ext(N,L)=0$ so that \eqref{short exact flat mb splitting seq} splits. 
 
Conversely, suppose that for every multibasic $L$, that every short exact sequence \begin{center}
    \begin{tikzcd}
     0 \arrow[r] & 
     L \arrow[r] & 
     M \arrow[r] &
     N \arrow[r] &
     0
    \end{tikzcd}
 \end{center} splits. Then $\Ext(N,L)=0$ for every multibasic $L$, so that $\Ext(N,I_{>0})=0$. This shows that $\Ext(N_i,I_{>0})=0$ for every summand $N_i$ in $N$ and we see from Proposition \ref{functor tables} that $N_i\in\{K,A,I_{>0}\}$, so that $N_i$ is flat for all $i$ and then $N$ is flat. \end{proof}

\begin{lem} \label{extending morphisms on inf gen to ones on fin gen lemma} If $F$ is free and finitely generated and $q\in\R$, then each morphism $\al:I_{>q}\to F$ can be extended to $I_q$. \end{lem}

\begin{proof} Suppose that this holds in the case $F=A$, and let $\al:I_{>q}\to F$. Then for each $i$, we have maps $\pi_i\al:I_{>q}\to A$, which extend to maps $\bt_i:I_q\to A$. Then set $\bt=\sum_i\io_i\bt_i:I_q\to F$, so that $\bt$ extends $\al$. So we can assume without loss of generality that $F=A$. Using the obvious isomorphisms $I_{>0}\to I_{>q}$ and $I_q\to A$, we can also assume that $q=0$, so all that remains is to show that each map $I_{>0}\to A$ can be extended to $A$.

If $0<p<1$, then $$\nu(\al(t))=\nu(t^p\al(t^{1-p}))=p+\nu(\al(t^{1-p}))\geq p,$$ so that $\nu(\al(t))\geq 1$ and then $\al(t)\in I_1$. So $\al(t)=ta$ for some $a\in A$. So for $b\in I_{>0}$, we have that $t(\al(b)-ab)=0$ and then multiplying by $t^{-1}$, we see that $\al(b)=ab$. We conclude that $\al$ can be extended to $A$. \end{proof}

We separate out three claims within the following proof by induction to keep track of various assumptions being made at different stages.

\begin{prop} \label{classification of submodules of free fin gen} Let $F$ be a free $A$-module of rank $n$. Then for each $M\leq F$, there exists a basis $\{u_1,\ldots,u_n\}$ for $F$ over $A$ and ideals $L_1\leq\cdots\leq L_n$ in $A$ for which $M=\bigoplus_i L_iu_i$. In particular, we have that $M\simeq\bigoplus_iL_i$ and $F/M\simeq\bigoplus_iA/L_i$ are both multibasic, with the former flat and the latter finitely generated. Furthermore, every finitely generated module is multibasic. \end{prop}

\begin{proof} We proceed by induction on $n$. It is clear for $n=0$; we now let $n\geq 1$ and assume that it holds for $n-1$. Let $F$ be a free $A$-module of rank $n$, and let $M\leq F$. We must show that there exists a basis $\{u_1,\ldots,u_n\}$ for $F$ and ideals $L_1\leq\ldots\leq L_n$ in $A$ for which $M=\bigoplus_i L_iu_i$. Choose a basis $\{v_1,\ldots,v_n\}$ for $F$ with dual basis $\{\pi_1,\ldots,\pi_n\}$ for $\Hom(F,A)$. 

We can reduce to the case where $$\pi_1(M)\leq\cdots\leq\pi_n(M).$$

\noindent\rule{\textwidth}{0.5pt}

Suppose that for each $N\leq F$ with $$\pi_1(N)\leq\cdots\leq\pi_n(N)$$ that there exists a basis $\{w_1,\ldots,w_n\}$ for $F$ and ideals $L_1\leq\cdots\leq L_n$ in $A$ for which $N=\bigoplus_i L_iw_i$. Since the ideals in $A$ are ordered by inclusion, there exists $\sigma\in S_n$ for which $$\pi_{\sigma(1)}(M)\leq\cdots\leq\pi_{\sigma(n)}(M).$$ Let $\td{\sigma}$ be the automorphism of $F$ sending $v_i\mapsto v_{\sigma(i)}$. Let $m\in M$. Then $m=\sum_j\lmb_j v_j$ for unique $\lmb_j\in A$, and we see that $$\pi_{\sigma(i)}(m)=\pi_{\sigma(i)}\left(\sum_j\lmb_j v_j\right)=\lmb_{\sigma(i)}$$ and $$\pi_i(\td{\sigma}^{-1}(m))=\pi_i(\sum_j\lmb_j\td{\sigma}^{-1}(v_j))=\pi_i\left(\sum_j\lmb_j v_{\sigma^{-1}(j)}\right)=\lmb_{\sigma(i)}.$$ This shows that $\pi_{\sigma(i)}(M)=\pi_i(\td{\sigma}^{-1}(M))$ for all $i$, so that $$\pi_1(\td{\sigma}^{-1}(M))\leq\cdots\leq\pi_n(\td{\sigma}^{-1}(M)).$$ Then, by our assumption in this claim, there exists a basis $\{w_1,\ldots,w_n\}$ for $F$ and ideals $L_1\leq\cdots\leq L_n$ in $A$ for which $\td{\sigma}^{-1}(M)=\bigoplus_iL_iw_i$. Setting $u_i=\td{\sigma}(w_i)$, we see that $\{u_1,\ldots,u_n\}$ is a basis for $F$. If $m\in M$, then $\td{\sigma}^{-1}(m)\in\td{\sigma}^{-1}(M)=\bigoplus_i L_i w_i$. So $\td{\sigma}^{-1}(m)=\sum_j \lmb_j w_j$ for unique $\lmb_j\in L_j$. Then $$m=\sum_j\lmb_j\td{\sigma}(w_j)=\sum_j\lmb_j u_j\in\bigoplus_i L_iu_i.$$ If $m\in\bigoplus_i L_iu_i$, then $m=\sum_j\lmb_j u_j=\td{\sigma}(\sum_j\lmb_j w_j)$ for unique $\lmb_j\in L_j$. Then $\sum_j\lmb_j w_j\in\bigoplus_i L_iw_i=\td{\sigma}^{-1}(M)$, so that $m=\td{\sigma}(\sum_j\lmb_j w_j)\in M$. This establishes that $M=\bigoplus_i L_iu_i$, as required. 

\noindent\rule{\textwidth}{0.5pt}

So we shall assume without loss of generality that $$\pi_1(M)\leq\cdots\leq\pi_n(M).$$ Set $F'=A\{v_1,\ldots,v_{n-1}\}$ and $M'=M\cap F'$. Since $F'$ is a free $A$-module of rank $n-1$ and $M'\leq F'$, the induction hypothesis shows that there exists a basis $\{u_1,\ldots,u_{n-1}\}$ for $F'$ with dual basis $\{\phi_1,\ldots,\phi_{n-1}\}$ for $\Hom(F',A)$, and ideals $L_1\leq\cdots\leq L_{n-1}$ in $A$, for which $M'=\bigoplus_{i<n}L_iu_i$. Set $L_n=\pi_n(M)$. 

We now prove that $M\leq L_nF$ and $L_{n-1}\leq L_n$.

\noindent\rule{\textwidth}{0.5pt}

Let $m\in M$, so that $m=\sum_i\pi_i(m)v_i$. If $\pi_i(m)=0$ for all $i$, then $m=0\in L_nF$. Otherwise, we set $$q=\min\{\nu(\pi_i(m)):\pi_i(m)\neq 0\}.$$ Then for each $i$, we have that $\pi_i(m)=t^qa_i$ for some $a_i\in A$. There exists $i$ with $\pi_i(m)\neq 0$ and $\nu(\pi_i(m))=q$, so that $a_i\in A^\x$. Then $$t^qa_i=\pi_i(m)\in\pi_i(M)\leq\pi_n(M)=L_n,$$ so that $t^q\in L_n$. We see that $m=t^q(\sum_i a_iv_i)\in L_nF$, so that $M\leq L_nF$.

If $a\in L_{n-1}$, then $au_{n-1}\in\bigoplus_{j<n} L_ju_j=M'\leq M$, so that $\pi_i(au_{n-1})\in L_n$ for $i<n$, and then $$a=\phi_{n-1}(au_{n-1})=\phi_{n-1}\left(\sum_{i<n}\pi_i(au_{n-1})v_i\right)=\sum_{i<n}\pi_i(au_{n-1})\phi_{n-1}(v_i)\in L_n.$$ So $L_{n-1}\leq L_n$, as required.

\noindent\rule{\textwidth}{0.5pt}

All that now remains is to show that there exists $u_n\in F$ for which $\{u_1,\ldots,u_n\}$ is a basis for $F$ and $M=\bigoplus_i L_iu_i$. In fact, we actually only need to prove that there exists $u_n\in F$ with $\pi_n(u_n)=1$ and $L_nu_n\leq M$, as we now demonstrate.

\noindent\rule{\textwidth}{0.5pt}

Suppose that there exists $u_n\in F$ with $\pi_n(u_n)=1$ and $L_nu_n\leq M$. If $i<n$, then $u_i\in F'$, so must be a linear combination of $v_i$ for $i<n$. This shows that $\pi_n(u_i)=0$ for $i<n$. If $\sum_i\lmb_iu_i=0$ for some $\lmb_i\in A$, then applying $\pi_n$ shows that $\lmb_n=0$. Then since $\{u_1,\ldots,u_{n-1}\}$ is linearly independent, we see that $\lmb_i=0$ for all $i$ and then $\{u_1,\ldots,u_n\}$ is linearly independent. If $x\in F$, then $x-\pi_n(x)u_n\in F'$ and therefore is spanned by $u_i$ for $i<n$. This shows that $\{u_1,\ldots,u_n\}$ spans $F$ and is therefore a basis for $F$. 

If $m\in M$, then $m=\sum_i\lmb_iu_i$ for unique $\lmb_i\in A$. Applying $\pi_n$, we see that $\lmb_n=\pi_n(m)\in L_n$, so that $\lmb_nu_n\in L_nu_n\leq M$. We see that $$\sum_{i<n}\lmb_iu_i=m-\lmb_nu_n\in M'=\bigoplus_{i<n}L_iu_i,$$ and therefore $\lmb_i\in L_i$ for all $i$, giving that $M\leq\bigoplus_iL_iu_i$. If $m\in\bigoplus_i L_iu_i$, then $m=\sum_i\lmb_iu_i$ for unique $\lmb_i\in L_i$. Then $\lmb_n\in L_n$ gives that $\lmb_n u_n\in M$, and $$\sum_{i<n}\lmb_iu_i\in\bigoplus_{i<n}L_iu_i=M'\leq M.$$ We conclude that $M=\bigoplus_i L_iu_i$, as required.

\noindent\rule{\textwidth}{0.5pt}

We shall now find $u_n\in F$ for which $\pi_n(u_n)=1$ and $L_nu_n\leq M$. 

Restricting $\pi_n$ gives a surjective map $\pi:M\to L_n$. Then $\ker(\pi)\leq F'$ and $F'$ is free and finitely generated of rank $n-1$, so it follows from the induction hypothesis that $\ker(\pi)$ is multibasic. Since $\pi$ is surjective and $L_n$ is flat and multibasic, we know from Lemma \ref{extensions flat multibasic split lemma} that the short exact sequence \begin{center}
    \begin{tikzcd}
     0 \arrow[r] & 
     \ker(\pi) \arrow[r] & 
     M \arrow[r] &
     L_n \arrow[r] &
     0
    \end{tikzcd}
 \end{center} splits, and then there exists a morphism $\al:L_n\to M$ with $\pi\al=\id_{L_n}$. 

The $L_n=0$ case is clear. If $L_n$ is non-zero, then there exists $q\geq 0$ for which $L_n$ is either $I_q$ or $I_{>q}$. First of all, suppose that $L_n=I_q$. Then $t^q\in I_q=L_n$, so that $\al(t^q)\in M\leq L_nF=I_qF$. So $\al(t^q)=t^qu_n$ for some $u_n\in F$. Applying $\pi_n$ to both sides, we see that $t^q=t^q\pi_n(u_n)$ and therefore $\pi_n(u_n)=1$. Each element of $L_nu_n$ has the form $t^qau_n=a\al(t^q)=\al(t^qa)\in M$ for some $a\in A$, so $L_nu_n\leq M$.

Finally, we turn our attention to the case $L_n=I_{>q}$. Since $\{v_i\}$ is a basis for $F$, one can check that $\{t^qv_i\}$ is a basis for $I_qF$, so that $I_qF$ is also free and finitely generated. Then since $M\leq L_nF=I_{>q}F\leq I_qF$, Lemma \ref{extending morphisms on inf gen to ones on fin gen lemma} shows that there exists a morphism $\bt:I_q\to I_qF$ with $\bt(a)=\al(a)$ for all $a\in I_{>q}=L_n$. Denote the composite map \begin{center}
    \begin{tikzcd}
     I_q \arrow[r,"\bt"] & 
     I_q	F \arrow[r,"\pi_n"] & 
     A 
    \end{tikzcd}
 \end{center} by $\gm$. Then $$\gm(t^r)=\pi_n(\bt(t^r))=\pi_n(\al(t^r))=\pi(\al(t^r))=t^r$$ for all $r>q$. Then $t\gm(t^q)=\gm(t^{q+1})=t^{q+1}$ gives that $\gm(t^q)=t^q$ and therefore $\gm(a)=a$ for all $a\in I_q$. Since $\bt(t^q)\in I_qF$, we can choose $u_n\in F$ for which $\bt(t^q)=t^qu_n$. Applying $\pi_n$ to both sides, we see that $t^q=t^q\pi_n(u_n)$, so that $\pi_n(u_n)=1$. We now know that $\bt(a)=au_n$ for all $a\in I_q$, so each element of $L_nu_n$ has the form $$au_n=\bt(a)=\al(a)\in M$$ for some $a\in L_n=I_{>q}$. We conclude that $L_nu_n\leq M$. Then the result follows by induction. \end{proof}

\begin{cor} \label{flat fin gen implies free} Every flat and finitely generated module is free. \end{cor}

\begin{proof} Let $M$ be flat and finitely generated. Since $M$ is finitely generated, it must be multibasic with only finitely generated summands by Proposition \ref{classification of submodules of free fin gen}. Since $M$ is also flat, all of its summands must be both flat and finitely generated. The only flat and finitely generated standard basic module is $A$, so $M$ is a finite direct sum of copies of $A$ and must therefore be free. \end{proof}

\begin{rmk} An alternative proof of Corollary \ref{flat fin gen implies free} is as follows. Suppose that $M$ is flat and finitely generated. Since every non-zero finitely generated ideal in $A$ is of the form $I_q$ ($q\geq 0$), $I_{-q}$ is a fractional ideal and $I_qI_{-q}=A$, we see that every non-zero finitely generated ideal in $A$ is invertible. Since $M$ is finitely generated, it then follows from Theorem 1 of \cite{Kap52} that $M/T$ is free, where $T$ is the torsion submodule of $M$. Since $M$ is flat, it must be torsion-free, so that $T=0$ and then $M$ is free.

Also note that Corollary \ref{flat fin gen implies free} is equivalent to the statement that every flat and finitely generated module is finitely presented. If we assume the latter, and let $M$ be flat and finitely generated, then it is flat and finitely presented, and therefore projective by Theorem 3.56 of \cite{Rot09}. Then $M$ must be free by Theorem \ref{projective iff free}. The converse is clear. \end{rmk}

\begin{theorem} \label{classification of submodules of fin dim k-vector spaces} Let $V$ be an $n$-dimensional $K$-vector space, regarded as an $A$-module. Then for each $M\leq V$ there exists a basis $\{u_1,\ldots,u_n\}$ for $V$ over $K$, and $A$-modules $L_1\leq\ldots\leq L_n\leq K$ with $L_i\in\{K,A,I_{>0},0\}$, for which $M=\bigoplus_i L_iu_i$. In particular, we have that $M\simeq\bigoplus_i L_i$ and $V/M\simeq\bigoplus_i K/L_i$ are both multibasic, with the former flat and the latter injective. \end{theorem}

\begin{proof} We proceed by induction on $n$. It is clear for $n=0$; we now let $n\geq 1$ and assume that it holds for $n-1$. Let $V$ be an $n$-dimensional $K$-vector space and $M$ be an $A$-submodule of $V$. We must show that there exists a basis $\{u_1,\ldots,u_n\}$ for $V$ over $K$, and $A$-submodules $L_1\leq\cdots\leq L_n\leq K$ with $L_i\in\{K,A,I_{>0},0\}$, for which $M=\bigoplus_iL_iu_i$. Choose a basis $\{v_1,\ldots,v_n\}$ for $V$ over $K$ with dual basis $\{\pi_1,\ldots,\pi_n\}$ for $\Hom_K(V,K)$ over $K$. 

We first reduce to the case where $$\pi_1(M)\leq\cdots\leq\pi_n(M).$$ 

\noindent\rule{\textwidth}{0.5pt}

Suppose that for every $N\leq V$ with $$\pi_1(N)\leq\cdots\leq\pi_n(N)$$ there exists a basis $\{w_1,\ldots,w_n\}$ for $V$ over $K$ and $L_1\leq\cdots\leq L_n\leq K$ with $L_i\in\{K,A,I_{>0},0\}$ and $N=\bigoplus_i L_iw_i$. Then there exists $\sigma\in S_n$ for which $$\pi_{\sigma(1)}(M)\leq\cdots\leq\pi_{\sigma(n)}(M),$$ and we define $\td{\sigma}$ to be the $K$-linear automorphism of $V$ sending $v_i\mapsto v_{\sigma(i)}$. Then, similarly to the case in Proposition \ref{classification of submodules of free fin gen}, we see that $\pi_{\sigma(i)}(M)=\pi_i(\td{\sigma}^{-1}(M))$ for all $i$, so that $\td{\sigma}^{-1}(M)\leq V$ and $$\pi_1(\td{\sigma}^{-1}(M))\leq\cdots\leq\pi_n(\td{\sigma}^{-1}(M)),$$ and then it follows from our assumption above (not to be confused with the induction hypothesis) that there exists a basis $\{w_1,\ldots,w_n\}$ for $V$ over $K$, and $L_1\leq\cdots\leq L_n\leq K$ satisfying $L_i\in\{K,A,I_{>0},0\}$, and $\td{\sigma}^{-1}(M)=\bigoplus_i L_iw_i$. Set $u_i=\td{\sigma}(w_i)$ for all $i$. Then $\{u_1,\ldots,u_n\}$ is a basis for $V$ over $K$, and one can check that $M=\bigoplus_i L_iu_i$, as required.

\noindent\rule{\textwidth}{0.5pt}

So we shall assume without loss of generality that $$\pi_1(M)\leq\cdots\leq\pi_n(M).$$

First of all, we consider the case where $\pi_n(M)=K$. Set $V'=K\{v_1,\ldots,v_{n-1}\}$ and $M'=M\cap V'$, and write the restriction $\pi_n|_M:M\to K$ as $\pi$. Since $V'$ has dimension $n-1$ over $K$ and $\ker(\pi)\leq V'$, we see from the induction hypothesis that $\ker(\pi)$ is multibasic. Since $\pi$ is surjective, we have a short exact sequence \begin{center}
    \begin{tikzcd}
     0 \arrow[r] & 
     \ker(\pi) \arrow[r] & 
     M \arrow[r,"\pi"] &
     K \arrow[r] &
     0 \;\;\;,
    \end{tikzcd}
 \end{center} which, since $K$ is flat, splits by Lemma \ref{extensions flat multibasic split lemma}. So there exists $\al:K\to M$ with $\pi\al=\id_K$. Since $M\leq V$ and $V$ is a $K$-vector space, we can multiply elements in $M$ by scalars in $K$. If $q>0$, then $$\al(1)=\al(t^qt^{-q})=t^q\al(t^{-q}).$$ We can then multiply through by $t^{-q}\in K$ to see that $\al(t^{-q})=t^{-q}\al(1)$. From this, it is easy to see that $\al(a)=a\al(1)$ for all $a\in K$, and therefore that $K\al(1)\leq M$. Since $V'$ is a $K$-vector space of dimension $n-1$ and $M'\leq V'$, we see from the induction hypothesis that there exists a basis $\{u_1,\ldots,u_{n-1}\}$ for $V'$ over $K$, and $A$-submodules $L_1\leq\cdots\leq L_{n-1}\leq K$ with $L_i\in\{K,A,I_{>0},0\}$, for which $M'=\bigoplus_{i<n}L_iu_i$. 

If we set $u_n=\al(1)\in V$, then (similarly to the proof of Proposition \ref{classification of submodules of free fin gen}) we see that $\pi_n(u_n)=1$ so that $\{u_1,\ldots,u_n\}$ is a basis for $V$ over $K$. If we also set $L_n=K$, then $L_1\leq\cdots\leq L_n\leq K$ with each $L_i\in\{K,A,I_{>0},0\}$. If $m\in M$, then $m=\sum_i\lmb_iu_i$ for unique $\lmb_i\in A$, and we see that $\lmb_nu_n\in K\al(1)\leq M$, so that $$\sum_{i<n}\lmb_iu_i=m-\lmb_nu_n\in M'=\bigoplus_{i<n}L_iu_i$$ and therefore $\lmb_i\in L_i$ for all $i$, so $M\leq\bigoplus_i L_iu_i$. The converse is clear, so $M=\bigoplus_iL_iu_i$, as required. 

Now we turn to the case where $\pi_n(M)<K$, so that there exists $q\geq 0$ for which $\pi_n(M)\leq t^{-q}A$. Since $\{v_1,\ldots,v_n\}$ is a basis for $V$ over $K$, we see that $\{t^{-q}v_1,\ldots,t^{-q}v_n\}$ is linearly independent over $A$ and therefore is a basis for $F=A\{t^{-q}v_1,\ldots,t^{-q}v_n\}$. If $m\in M$, then $$\pi_i(m)\in\pi_i(M)\leq\pi_n(M)\leq t^{-q}A,$$ so that $\pi_i(m)=t^{-q}a_i$ for some $a_i\in A$. Then $$m=\sum_i\pi_i(m)v_i=\sum_i a_i(t^{-q}v_i)\in F.$$ So $F$ is a free $A$-module of rank $n$ containing $M$, and therefore there must exist a basis $\{w_1,\ldots,w_n\}$ for $F$ over $A$ and ideals $L_1\leq\cdots\leq L_n$ in $A$ for which $M=\bigoplus_iL_iw_i$ by Proposition \ref{classification of submodules of free fin gen}. 

If $L_i=0$, we set $u_i=w_i\in V$ and $N_i=0$; if $L_i=I_r$ for some $r\geq 0$, we set $u_i=t^rw_i\in V$ and $N_i=A$; if $L_i=I_{>r}$ for some $r\geq 0$, we set $u_i=t^rw_i\in V$ and $N_i=I_{>0}$. Let $\lmb_1,\ldots,\lmb_n\in K$ with $\sum_i\lmb_iu_i=0$. Set $\mu_i=\lmb_i$ if $L_i=0$ and $\mu_i=t^r\lmb_i$ if $L_i$ is non-zero. Then $\lmb_iu_i=\mu_iw_i$ for all $i$ and one can check that there exists $s\geq 0$ for which $t^s\mu_i\in A$ for all $i$. So $\sum_i\mu_iw_i=0$ and then $\sum_i t^s\mu_iw_i=0$ with $t^s\mu_i\in A$ and $\{w_1,...,w_n\}$ linearly independent over $A$. So $t^s\mu_i=0$ for all $i$ and therefore $\lmb_i=0$ for all $i$, so $\{u_1,\ldots,u_n\}$ is linearly independent over $K$. 

Each $x\in V$ takes the form $x=\sum_it^q\pi_i(x)t^{-q}v_i$. Then each $t^{-q}v_i\in F$, and $\{w_1,\ldots,w_n\}$ is a basis for $F$ over $A$, so $t^{-q}v_i$ can be expressed as an $A$-linear combination of the $w_i$, each of which is a $K$-multiple of $u_i$. So $\{u_1,\ldots,u_n\}$ $K$-spans $V$ and therefore is a basis for $V$ over $K$. It is clear that $N_iu_i=L_iw_i$ for each $i$, so that $M=\bigoplus_iN_iu_i$ with $N_i\in\{A,I_{>0},0\}$ for each $i$.

Since $N_1,\ldots,N_n\leq K$, there exists $\sigma\in S_n$ for which $$N_{\sigma(1)}\leq\cdots\leq N_{\sigma(n)}.$$ Then $\{u_{\sigma(1)},\ldots,u_{\sigma(n)}\}$ is a basis for $V$ over $K$, $N_{\sigma(1)}\leq\cdots\leq N_{\sigma(n)}\leq K$ with $N_{\sigma(i)}\in\{K,A,I_{>0},0\}$ and $M=\bigoplus_iN_{\sigma(i)}u_{\sigma(i)}$. It is routine to check that the final maps are well-defined isomorphisms. \end{proof}

The following four propositions lay the groundwork for establishing that the category of multibasics is closed under extensions. We start by examining $\Ext(M,N)$ where either $M$ or $N$ is finite. Half of the work is covered in the four propositions below, but we defer the dual results to Corollary \ref{extensions of multibasics finite cases}, where we use the fact that $M$ is multibasic if and only if $DM$ is. These results regarding extensions will be useful in proving that every $\Ta$-reflexive module is multibasic, since for $M\in\mcM$, we will see that $\F\tns M$ and $\Tor(\F,M)$ are finite ((Proposition \ref{f tns m is finite} and Theorem \ref{tor(f2,m) is finite}). We will then use the fact that $\mcM$ is closed under extensions to see that the multibasics are closed under extensions.

\begin{prop} \label{extensions f over multibasic} Let \begin{center}
    \begin{tikzcd}
     0 \arrow[r] & 
     L \arrow[r,"\al"] & 
     M \arrow[r,"\bt"] &
     \F \arrow[r] &
     0
    \end{tikzcd}
 \end{center} be exact. If $L$ is multibasic, then so is $M$. \end{prop}

\begin{proof} Let $L=\bigoplus_{i=1}^m U_i/V_i$, where $V_i<U_i\leq K$, and let $r$ be the number of summands with $U_i=I_{>0}$. Rearrange so that $U_i=I_{>0}$, $U_j=A$ and $U_k=K$ implies $i<j<k$. Then rearrange so that $V_1\geq\cdots\geq V_r$. Now let $I=\bigoplus_i K/V_i$ and $J=\bigoplus_i K/U_i$, so that \begin{center}
    \begin{tikzcd}
     0 \arrow[r] & 
     L \arrow[r,"\io"] & 
     I \arrow[r,"\dl"] &
     J \arrow[r] &
     0
    \end{tikzcd}
 \end{center} is an injective resolution for $L$. Note that since $U_i\in\{K,A,I_{>0}\}$ we know that $K/U_i\in\{0,\Ta,\Phi\}$. So $J$ involves only copies of $\Ta$ and $\Phi$, with all copies of $\Ta$ to the left of all copies of $\Phi$. Since $I$ is injective, there exists a map $\rho:M\to I$ for which the following diagram commutes. \begin{center}
    \begin{tikzcd} [sep=2cm]
     L \arrow[r,"\al",rightarrowtail] \arrow[d,"\io",swap] & 
     M \arrow[ld,"\rho",dashed] \\
     I 
     \end{tikzcd}\end{center} If $\bt=0$, then $\al$ is an isomorphism so that $M\leq L$ is multibasic. Otherwise, we can choose $m\in M$ with $\bt(m)=1+I_{>0}$. Then $\dl(\rho(m))\in J$. If $m'\in M$ with $\bt(m')=1+I_{>0}$, then $m-m'\in\ker(\bt)=\im(\al)$ so that $m-m'=\al(x)$ for some $x\in L$. Then $$\dl(\rho(m))-\dl(\rho(m'))=\dl(\rho(\al(x)))=\dl(\io(x))=0,$$ so that this element in $J$ is independent of the choice of $m\in M$ with $\bt(m)=1+I_{>0}$. If $a\in I_{>0}$, then $\bt(am)=I_{>0}$, so $am\in\ker(\bt)=\im(\al)$. So $am=\al(x)$ for some $x\in L$ and then $$a\dl(\rho(m))=\dl(\rho(am))=\dl(\rho(\al(x)))=\dl(\io(x))=0.$$ So our element in $J$ is contained in $\ann(I_{>0},J)$. We define $\gm:\F\to J$ by sending $1+I_{>0}$ to this element. One can then check that $\gm$ is the unique morphism for which the right-hand square of the following diagram commutes, and therefore forms a pullback square. \begin{center}
    \begin{tikzcd} [sep=2cm]
     0 \arrow[r] \arrow[d] &
     L \arrow[r,"\al"] \arrow[d,equal] &
     M \arrow[r,"\bt"] \arrow[d,"\rho"] & 
     \F \arrow[r] \arrow[d,"\gm",dashed] &
     0 \arrow[d] \\
     0 \arrow[r] &
     L \arrow[r,"\io"] & 
     I \arrow[r,"\dl"] &
     J \arrow[r] &
     0
     \end{tikzcd}\end{center} 
     
If $L=0$, then $M=\F$ is multibasic. Also, if $\gm=0$, then $\dl\rho=0$ so that $\rho(m)\in L$ and then $\al$ is a split monomorphism, which gives that $M=\F\oplus L$ is multibasic. So we assume that $L$ and $\gm$ are both non-zero. 

Note that $\Hom(\F,\Ta)=\F$ and $\Hom(\F,\Phi)=0$, so $\gm$ is determined by which components mapping into $\Ta$ are zero, and which are inclusions. Set $$\gm(1+I_{>0})=(\gm_1+I_{>0},\ldots,\gm_r+I_{>0},0+A,\ldots,0+A,0,\ldots,0),$$ where each $\gm_i\in\{0,1\}$. Let the integer $s$ be maximal such that $\gm_s=1$. Set $$m=(1+I_{>0},\gm_1+V_1,\ldots,\gm_r+V_r,0+V_{r+1},\ldots,0+V_m)\in\F\oplus I$$ so that $\ann(m)=V_s$ and $I_{>0}m\leq L$. Now define $$X:=\{(a+I_{>0},x)\in\F\oplus I:\gm(a+I_{>0})=\dl(x)\},$$ so that $m\in X\tk L$. One can check that $$X=A\{L\cup\{m\}\}=\{x+\lmb m:x\in L,\lmb\in\{0,1\}\}.$$ The map $$a\mapsto am:A\to X$$ has kernel $\ann(m)$, so we get a monomorphism $A/V_s\to X$. We then get a surjective map $\phi:\frac{A}{\ann(m)}\oplus L\to X$ given by $$(a+\ann(m),x)\mapsto am+x$$ and thus an isomorphism $$(a+\ann(m),x)+\ker(\phi)\mapsto am+x:\frac{\frac{A}{\ann(m)}\oplus L}{\ker(\phi)}\to X.$$ One can check that $$\ker(\phi)=\{(a+\ann(m),am):a\in I_{>0}\}.$$ Define $\psi:\frac{A}{\ann(m)}\oplus L\to\frac{A}{\ann(m)}\oplus L$ as follows. It sends $$(a+\ann(m),b_1+V_1,\ldots,b_m+V_m)$$ to $$(a+b_s+\ann(m),b_1+\gm_1b_s+V_1,\ldots,b_{s-1}+\gm_{s-1}b_s+V_{s-1},b_s+V_s,b_m+V_m).$$ One can check that $\psi$ is an automorphism with $\psi(\ker(\phi))=I_{>0}/V_s$, a summand in $L$. Then $\frac{L}{I_{>0}/V_s}$ is multibasic. So $\psi$ induces an isomorphism $$\frac{\frac{A}{\ann(m)}\oplus L}{\ker(\phi)}\to\frac{A}{\ann(m)}\oplus\frac{L}{I_{>0}/V_s},$$ giving that $$X\simeq\frac{A}{V_s}\oplus\frac{L}{I_{>0}/V_s}$$ is multibasic. Since $$x\mapsto(0+I_{>0},x):L\to X$$ is injective and has cokernel $\{L,m+L\}\simeq\F$, we can form a short exact sequence $$0\to L\to X\to\F\to 0.$$ Let $\rho'$ to be the restriction of the projection $\F\oplus I\to I$ to $X$. Then $\rho'\al'=\io$ and $\gm\bt'=\dl\rho'$, so that $(X,\rho',\bt')$ also forms a pullback of $\gm$ and $\dl$. We conclude that $$M\simeq X\simeq\frac{A}{V_s}\oplus\frac{L}{I_{>0}/V_s}$$ is multibasic. \end{proof}

\begin{prop} \label{extension finite on right} Let \begin{center}
    \begin{tikzcd}
     0 \arrow[r] & 
     L \arrow[r,"\al"] & 
     M \arrow[r,"\bt"] &
     \F^n \arrow[r] &
     0
    \end{tikzcd}
 \end{center} be exact. If $L$ is multibasic, then so is $M$. \end{prop}

\begin{proof} We proceed by induction on $n$. If $n=0$, then clearly $M\simeq L$ is multibasic. The $n=1$ case was covered in Proposition \ref{extensions f over multibasic}. Let $n\geq 1$ and assume that this holds in the $n-1$ case. Let $\pi:\F^n\to\F$ be the projection to the $n^{\text{th}}$ coordinate, and set $M'=\ker(\pi\bt)$. Since $\pi\bt$ is surjective, we know that \begin{center}
    \begin{tikzcd}
     0 \arrow[r] & 
     M' \arrow[r] & 
     M \arrow[r,"\pi\bt"] &
     \F \arrow[r] &
     0
    \end{tikzcd}
 \end{center} is exact. If $x\in L$, then $$(\pi\bt)(\al(x))=\pi(\bt(\al(x)))=\pi(0)=0$$ so that $\al(x)\in\ker(\pi\bt)=M'$ and then $\im(\al)\leq M'$; denote the corestriction of $\al$ to $M'$ by $\al'$. Since $\al$ is injective, we know that $\al'$ is injective. Let $\io_\F:\F^{n-1}\to\F^n$ be given by $$\io_\F(x_1,\ldots,x_{n-1})=(x_1,\ldots,x_{n-1},0);$$ this is the categorical kernel of $\pi$. Then since $$(\pi(\bt\io))(m)=(\pi\bt)(m)=0$$ for all $m\in M'$, we know that $\pi(\bt\io)=0$, so by the universal property of the kernel there must be a unique map $\bt':M'\to\F^{n-1}$ for which the following diagram commutes. \begin{center}
    \begin{tikzcd}[sep=2cm]
     M' \arrow[d,"\bt'",dashed,swap] \arrow[rd,"\bt\io"] &
     {} & {} \\
     \F^{n-1} \arrow[r,"\io_\F"] & 
     \F^n \arrow[r,"\pi"] & 
     \F
    \end{tikzcd}
 \end{center}

Let $x=(a_1,\ldots,a_{n-1})\in\F^{n-1}$. Then $(x,0)\in\F^n=\im(\bt)$, so that $(x,0)=\bt(m)$ for some $m\in M$. Then $$(\pi\bt)(m)=\pi(x,0)=0,$$ so that $m\in M'$. Then $$(\bt'(m),0)=\io_\F(\bt'(m))=\bt(\io(m))=\bt(m)=(x,0),$$ so that $x=\bt'(m)\in\im(\bt')$, proving that $\bt'$ is surjective. 

Let $m\in\im(\al')$. Then $$m=\al'(x)=\al(x)\in\im(\al)=\ker(\bt)$$ for some $x\in L$. Then $(\pi\bt)(m)=0$ so that $m\in M'$. So $$(\io_{\F}\bt')(m)=(\bt\io)(m)=\bt(m)=0$$ gives that $\bt'(m)\in\ker(\io_{\F})=0$ and then $m\in\ker(\bt')$. So $\im(\al')\leq\ker(\bt')$. If $m\in\ker(\bt')$, then $$\bt(m)=(\bt\io)(m)=\io_{\F}(\bt'(m))=0,$$ so that $m\in\ker(\bt)=\im(\al)=\im(\al')$. Then $\im(\al')=\ker(\bt')$, and we have shown that \begin{center}
    \begin{tikzcd}
     0 \arrow[r] & 
     L \arrow[r,"\al'"] & 
     M' \arrow[r,"\bt'"] &
     \F^{n-1} \arrow[r] &
     0
    \end{tikzcd}
 \end{center} is exact. Since $L$ is multibasic, it follows from the induction hypothesis that $M'$ must be multibasic. Then since \begin{center}
    \begin{tikzcd}
     0 \arrow[r] & 
     M' \arrow[r] & 
     M \arrow[r,"\pi\bt"] &
     \F \arrow[r] &
     0
    \end{tikzcd}
 \end{center} is exact and $M'$ is multibasic, we see from Proposition \ref{extensions f over multibasic} that $M$ must be multibasic, as required. \end{proof}

\begin{prop} \label{extensions f on left} Let \begin{center}
    \begin{tikzcd}
     0 \arrow[r] & 
     \F \arrow[r,"\al"] & 
     M \arrow[r,"\bt"] &
     N \arrow[r] &
     0
    \end{tikzcd}
 \end{center} be exact. If $M$ multibasic, then so is $N$. \end{prop}

\begin{proof} We use automorphisms of $M$ to reduce $M$ to a basic module. Let $\psi:M\to M$ be an automorphism, and $M_1$ be the direct sum of a minimal collection of basic summands in $M$ for which $\psi(\im(\al))\leq M_1$. Then $M=M_0\oplus M_1$ for some multibasic $M_0$. Corestricting $\psi\al$, we form the short exact sequence \begin{center}
    \begin{tikzcd}
     0 \arrow[r] & 
     \F \arrow[r,"\psi\al"] & 
     M_1 \arrow[r] &
     \cok(\F\to M_1) \arrow[r] &
     0 \;\;\;.
    \end{tikzcd}
 \end{center} Then since \begin{center}
    \begin{tikzcd}
     0 \arrow[r] & 
     \F \arrow[r,"\psi\al"] & 
     M \arrow[r] &
     M_0\oplus\cok(\F\to M_1) \arrow[r] &
     0
    \end{tikzcd}
 \end{center} is exact, we see that $$N\simeq\cok(\al)\simeq\cok(\psi\al)\simeq M_0\oplus\cok(\F\to M_1),$$ so that we only need to show that $\cok(\F\to M_1)$ is multibasic (due to the minimality of $M_1$, the maps from $\F$ to each component in $M_1$ are non-zero, so each basic summand in $M_1$ must have denominator $I_{>q}$ for some $q>0$). One can check that $$(a_i+I_{>q_i})\mapsto(a_1+I_{>q_1},a_2+t^{q_2-q_1}a_1+I_{>q_2},\ldots,a_r+t^{q_r-q_1}a_1+I_{>q_r})$$ defines an automorphism of the summands in $M$ of the type $I_{>0}/I_{>q}$ that reduces these summands to the left-most summand; it can then be extended to an automorphism of $M$ that fixes the summands of other types. Using this approach, we can reduce to the case where $M$ is basic, in which case it is clear that $N$ is multibasic. \end{proof}

\begin{prop} \label{l finite m multibasic implies n multibasic} Let \begin{center}
    \begin{tikzcd}
     0 \arrow[r] & 
     \F^n \arrow[r,"\al"] & 
     M \arrow[r,"\bt"] &
     N \arrow[r] &
     0
    \end{tikzcd}
 \end{center} be exact. If $M$ is multibasic, then so is $N$. \end{prop}

\begin{proof} We proceed by induction on $n$. If $n=0$, then $N\simeq M$ is multibasic. If $n=1$, then $N$ is multibasic by Proposition \ref{extensions f on left}. Let $n\geq 2$ and assume that this proposition holds for $n-1$.

Identifying $\{(0,\ldots,0,x_i):x_i\in\F\}\leq\F^n$ with $\F$, we define $\al':\F^{n-1}\to M/\al(\F)$ to be the unique map for which the left-hand corner of the following diagram commutes. \begin{equation} \label{nine lemma diagram}
    \begin{tikzcd}[sep=2cm]
    	 \F \arrow[r,"\sim"] \arrow[d,rightarrowtail] &
     \al(\F) \arrow[r] \arrow[d,rightarrowtail] &
     0 \arrow[d] \\
     \F^n \arrow[r,"\al",rightarrowtail] \arrow[d,twoheadrightarrow] & 
     M \arrow[r,"\bt",twoheadrightarrow] \arrow[d,twoheadrightarrow] &
     N \arrow[d,equal] \\
     \F^{n-1} \arrow[r,"\al'",dashed,rightarrowtail] &
     M/\al(\F) \arrow[r,"\bt'",twoheadrightarrow] &
     N
    \end{tikzcd}
 \end{equation} In this diagram, $\bt'$ is the composite \begin{center}
    \begin{tikzcd}
     M/\al(\F) \arrow[r,twoheadrightarrow] & 
     \cok(\al) \arrow[r,"\sim"] & 
     N \;\;\;.
    \end{tikzcd}
 \end{center} One can check that the whole diagram commutes, with all columns and the top two rows forming short exact sequences. It then follows from the Nine Lemma that the bottom row also forms a short exact sequence. Since \begin{center}
    \begin{tikzcd}
     0 \arrow[r] & 
     \F \arrow[r] & 
     M \arrow[r] &
     M/\al(\F) \arrow[r] &
     0
    \end{tikzcd}
 \end{center} is exact and $M$ is multibasic, we also see from Proposition \ref{extensions f on left} that $M/\al(\F)$ is multibasic. Then we use the induction hypothesis on the short exact sequence formed from the bottom row of \eqref{nine lemma diagram} to see that $N$ is multibasic, as required. \end{proof}

\section{Uniqueness of Decomposition}

In this short section, we prove that every multibasic module has a unique finite decomposition into basic submodules, up to the ordering of the summands. For this purpose, we briefly develop the theory of multibasic invariants, and construct ones that determine the number of summands in the decomposition of a multibasic module that are isomorphic to a given standard basic module.

\begin{definition} A \textit{multibasic invariant} is a map $\xi$ that associates an integer $\xi_M\in\Z$ to each multibasic module $M$, with the following properties. 

1. $\xi_M=\xi_N$ when $M$ and $N$ are isomorphic. 

2. $\xi_{M\oplus N}=\xi_M+\xi_N$ for all $M$ and $N$. 

Note that since every multibasic module is isomorphic to a finite direct sum of standard basic modules, the isomorphism classes of multibasic modules must form a set. Since multibasic invariants, by definition, respect isomorphisms, we can think of a multibasic invariant as a function on the set of isomorphism classes of multibasic modules. The set of multibasic invariants has an abelian group structure given by $$(\xi+\tau)_M=\xi_M+\tau_M.$$ We say that a multibasic invariant $\xi$ \textit{counts} a collection $\mcS$ of standard basic modules if $\xi_V=1$ for all $V\in\mcS$ and $\xi_V=0$ for all standard basic $V$ not in $\mcS$. It is then immediate that $\xi_M$ is the number of basic summands in $M$ that are in the collection $\mcS$. Since a multibasic invariant is uniquely determined by its action on the standard basic modules, any multibasic invariant counting a collection $\mcS$ of standard basic modules must be unique. \end{definition}

\begin{exmp} For multibasic $M$, set $f_M=\dim_K(K\tns M)$ and $g_M=\dim_{\F}(\F\tns M)$. These define multibasic invariants $f$ and $g$. \end{exmp}

\begin{exmp} If $\xi$ is a multibasic invariant and $E:\Mod_A\to\Mod_A$ is an additive functor (covariant or contravariant) that sends multibasics to multibasics, then $M\mapsto\xi_{EM}$ is also a multibasic invariant. \end{exmp}

\begin{exmp} Since $D:\Mod_A\to\Mod_A$ is an additive contravariant functor that sends multibasics to multibasics, we know that $M\mapsto f_{DM}$ and $M\mapsto g_{DM}$ are multibasic invariants. We similarly see that, for each $q\geq 0$, that $M\mapsto g_{I_qM}$ defines a multibasic invariant. \end{exmp}

The proof of the following lemma is routine.

\begin{lem} \label{automatic duals for mb invariants} Let $\mcS$ be a collection of standard basic modules and $\xi$ be a multibasic invariant. Then $\xi$ counts $\mcS$ if and only if $M\mapsto\xi_{DM}$ counts $D\mcS$, the duals of those in $\mcS$. \end{lem}

\begin{exmp} We now tabulate the values of some of the above examples of multibasic invariants on the standard basic modules. Note that the results in the final four columns are independent of the choice of $q>0$. We will use these examples to construct multibasic invariants that count each standard basic module in Lemma \ref{mb mod counting each sb mod}.

\[ \begin{array}{|c|c|c|c|c|c|c|c|c|c|c|}\hline
     M       & K & A & I_{>0} & \Ta & \Phi & \F & A/I_q & A/I_{>q} & I_{>0}/I_q & I_{>0}/I_{>q} \\ \hline 
     f_M    & 1 & 1 & 1      & 0   & 0  & 0  & 0     & 0        & 0          & 0             \\ \hline 
     f_{DM}   & 1 & 0 & 0      & 1   & 1  & 0  & 0     & 0        & 0          & 0             \\ \hline 
     g_M    & 0 & 1 & 0      & 0   & 0 & 1   & 1     & 1        & 0          & 0             \\ \hline 
     g_{DM}   & 0 & 0 & 0      & 1   & 0  & 1  & 0     & 1        & 0          & 1             \\ \hline 
    \end{array}
 \] 
 
\end{exmp}

\begin{rmk} We see from the above table that if $M$ is multibasic, then $f_M$ counts the flat summands, $f_{DM}$ counts the injective summands, $g_M$ counts the finitely generated summands and $g_{DM}$ counts the duals of finitely generated summands. \end{rmk}

We now use $f$ and $g$ to construct some more sophisticated multibasic invariants.

\begin{definition} For piecewise-constant maps $f:[0,\infty)\to\Z$, define $\delta^+(f):(0,\infty)\to\Z$ and $\delta^-(f):[0,\infty)\to\Z$ by $$\delta^+(f)(q)=\lim_{p\to q^-}f(p)-f(q)$$ and $$\delta^-(f)(q)=f(q)-\lim_{p\to q^+}f(p).$$ \end{definition}

\begin{exmp} Let $\xi=\{\xi_q\}_{q\geq 0}$ be a family of multibasic invariants, where for each multibasic $M$, the map $\xi_M:[0,\infty)\to\Z$ given by $(\xi_M)(q)=(\xi_q)_M$ is piecewise constant. Then $$M\mapsto\dl^+(\xi_M)(q)$$ is a multibasic invariant for $q>0$, and $$M\mapsto\dl^-(\xi_M)(q)$$ is a multibasic invariant for $q\geq 0$. We label these $\dl_q^+(\xi)$ and $\dl_q^-(\xi)$, respectively. Also $$\left(\lim_{q\to\infty}\xi_q\right)_M=\lim_{q\to\infty}(\xi_q)_M$$ defines a multibasic invariant. \end{exmp}

\begin{exmp} For each $q\geq 0$ and multibasic $M$, define $(\eta_q)_M=g_{I_qM}$. Then $\eta=\{\eta_q\}_{q\geq 0}$ is a family of multibasic invariants, where the maps $\eta_M$ are piecewise constant for each $M$. So $\dl_q^+(\eta)$ for $q>0$, $\dl_q^-(\eta)$ for $q\geq 0$, and $\lim_{q\to\infty}\eta_q$ are all multibasic invariants. \end{exmp}

\begin{lem} \label{mb mod counting each sb mod} For each standard basic module $V$, there exists a unique multibasic invariant $\Psi(V)$ that counts $V$. \end{lem}

\begin{proof} It is routine to check that $\Psi(A)=\lim_{q\to\infty}\eta_q$, $\Psi(\F)=\dl_0^-(\eta)$, $\Psi(A/I_q)=\dl_q^+(\eta)$ and $\Psi(A/I_{>q})=\dl_q^-(\eta)$ for each $q>0$. We also see that $$\Psi(A)_M+\Psi(I_{>0})_M=\Psi(A)_{\Hom(I_{>0},M)},$$ giving $\Psi(I_{>0})$. Then $$f=\Psi(K)+\Psi(A)+\Psi(I_{>0})$$ gives $\Psi(K)$. Using Lemma \ref{automatic duals for mb invariants}, we immediately can identify $\Psi(\Ta)$, $\Psi(\Phi)$ and $\Psi(I_{>0}/I_{>q})$ for $q>0$. All that remains is $I_{>0}/I_q$ for $q>0$. But $\Psi(A/I_q)_{\Hom(I_{>0},M)}$ counts the summands of length $q$, and we already know $\Psi(A/I_q)$, $\Psi(A/I_{>q})$ and $\Psi(I_{>0}/I_{>q})$, giving $\Psi(I_{>0}/I_q)$. \end{proof}

We now list the multibasic invariants counting each of the standard basic modules.

\begin{align} 
&\begin{aligned} \Psi(A)=\lim_{q\to\infty}\eta_q, \end{aligned} \\
&\begin{aligned} \Psi(\F)=\dl_0^-(\eta), \end{aligned} \\
&\begin{aligned} \Psi(A/I_q)=\dl_q^+(\eta), \end{aligned} \\
&\begin{aligned} \Psi(A/I_{>q})=\dl_q^-(\eta), \end{aligned} \\ 
&\begin{aligned} \Psi(\Ta)_M=\Psi(A)_{DM}, \end{aligned} \\
&\begin{aligned} \Psi(I_{>0}/I_{>q})_M=\Psi(A/I_q)_{DM} \end{aligned} \\ 
&\begin{aligned} \Psi(I_{>0})_M=\Psi(A)_{\Hom(I_{>0},M)}-\Psi(A)_M, \end{aligned} \\ 
&\begin{aligned} \Psi(K)=f-\Psi(A)-\Psi(I_{>0}) \end{aligned} \\ 
&\begin{aligned} \Psi(\Phi)_M=\Psi(I_{>0})_{DM}, \end{aligned} \\ 
&\begin{aligned} \Psi(I_{>0}/I_q)_M=\Psi(A/I_q)_{\Hom(I_{>0},M)}-\Psi(A/I_q)_M-\Psi(A/I_{>q})_M-\Psi(I_{>0}/I_{>q})_M.  \end{aligned} \end{align} 

Finally, we prove uniqueness of decomposition. 

\begin{theorem} \label{uniqueness theorem for multibasics} Each multibasic module has a unique decomposition into standard basic modules, up to the ordering of the summands.\footnote{Alternatively, this can be approached by applying the Krull-Schmidt Theorem to finite direct sums of indecomposable pure injectives. See Section 3.3 of \cite{BeG99} for details.} \end{theorem}

\begin{proof} Let $M$ be multibasic, with decompositions $\bigoplus_iM_i$ and $\bigoplus_jN_j$, each $M_i$, $N_j$ basic. For each standard basic $V$, we have that $\Psi(V)_{\bigoplus_i M_i}=\Psi(V)_{\bigoplus_j N_j}$, so that both decompositions have the same number of copies of $V$. So the number of copies of $V$ appearing in $M$ is independent of the choice of decomposition. Since this applies to all standard basic $V$, it shows that all decompositions have exactly the same summands, but possibly in a different order. \end{proof}

\newpage

\chapter{Theta-Reflexivity} \label{reflexivity} In this chapter, we continue our investigation of $A$-modules by focusing on the category $\mcM$ of $\Ta$-reflexives. We show that $\mcM$ is abelian, closed under extensions, submodules and quotients. We also see that $DM$ is multibasic if and only if $M$ is. Our main result is that $\mcM$ is in fact the full subcategory of multibasic modules, providing a complete description of the objects in $\mcM$ up to isomorphism.

\section{The Category of Theta-Reflexives}

In this section, we establish properties of the category $\mcM$ of $\Ta$-reflexive modules. In particular, we will see that $\mcM$ is abelian and contains the multibasic modules.

\begin{lem} \label{chi defines natural transformation on a modules relation chi and functor} $M\mapsto\chi_M$ defines a natural transformation $\chi:\id_{\Mod_A}\to D^2$. For each $A$-module $M$, we have that $$D(\chi_M)\chi_{DM}=\id_{DM}.$$ \end{lem}

\begin{proof} We have a morphism $\chi_M:M\to D^2M$ for each $A$-module $M$. Let $\al:M\to N$. Then for all $m\in M$, we have that $$(D^2\al\chi_M)(m)=(D^2\al)(\chi_M(m))=D(D\al)(\chi_M(m)).$$ For $m\in M$ and $f\in DN$, we then have that $D(D\al)(\chi_M(m))(f)$ is equal to $$(\chi_M(m)(D\al))(f)=\chi_M(m)((D\al)(f))=\chi_M(m)(f\al)=(f\al)(m)=f(\al(m)).$$ We also have $$(\chi_N\al)(m)(f)=\chi_N(\al(m))(f)=f(\al(m)),$$ so that $D^2\al\chi_M=\chi_N\al$, confirming that $M\mapsto\chi_M$ defines a natural transformation. \begin{center}
    \begin{tikzcd}[sep=2cm]
     M \arrow[r,"\al"] \arrow[d,"\chi_M",swap] &
     N \arrow[d,"\chi_N"] \\
     D^2M \arrow[r,"D^2\al"] & 
     D^2N
    \end{tikzcd}
 \end{center}

For the final part, let $f\in DM$. Then $$(D(\chi_M)\chi_{DM})(f)=D(\chi_M)(\chi_{DM}(f))=\chi_{DM}(f)\chi_M.$$ Then for all $m\in M$, we have that $$(\chi_{DM}(f)\chi_M)(m)=\chi_{DM}(f)(\chi_M(m))=\chi_M(m)(f)=f(m)=\id_{DM}(f)(m),$$ and we conclude that $$D(\chi_M)\chi_{DM}=\id_{DM}.$$ \end{proof}

\begin{lem} \label{theta is cogenerator} Let $M$ be an $A$-module. Then $DM=0$ if and only if $M=0$. \end{lem}

\begin{proof} It is clear that $M=0$ implies $DM=0$. Conversely, suppose that $DM=0$. Let $m\in M$ and $\al:Am\to\Ta$. Since $\Ta$ is injective, we can extend $\al$ to $M$. But since $DM=0$, this extension is the zero map, forcing $\al=0$. So $D(A/\ann(m))\simeq D(Am)=0$. But $A/\ann(m)$ is a subquotient of $K$, so it follows from Proposition \ref{lemma for iso for d of basic module} that $A/\ann(m)=0$ and then $\ann(m)=A$ gives that $m=0$ and finally $M=0$. \end{proof}

\begin{lem} \label{hom functor d switches injective and surjective} $D$ switches kernels and cokernels, injective maps and surjective maps. \end{lem}

\begin{proof} That $D$ switches kernels and cokernels follows easily from the fact that $D$ is exact. That $D$ switches injective maps and surjective maps then follows from Lemma \ref{theta is cogenerator}. \end{proof}

\begin{lem} \label{chi is mono equivalent cok condition} For each $A$-module $M$, we have that $\chi_M$ is injective, so that $M\in\mcM$ if and only if $\cok(\chi_M)=0$. \end{lem}

\begin{proof} Let $m\in\ker(\chi_M)$ and $\al\in D(Am)$. Since $\Ta$ is injective, we can extend $\al$ to a morphism $\td{\al}\in DM$. Since $m\in\ker(\chi_M)$, we see that $$\al(m)=\td{\al}(m)=\chi_M(m)(\td{\al})=0$$ and therefore $\al=0$. This shows that $D(Am)=0$, so that $Am=0$ by Lemma \ref{theta is cogenerator} and then $m=0$. So $\ker(\chi_M)=0$ and we conclude that $\chi_M$ is injective.

Finally, since $\chi_M$ is automatically injective, we see that it is an isomorphism if and only if it is surjective, which is equivalent to $\cok(\chi_M)=0$. \end{proof}

\begin{prop} \label{d is faithful} $D$ is faithful. Furthermore, if $M,N\in\mcM$ then the induced map $$\Hom(M,N)\to\Hom(DN,DM)$$ is an isomorphism. \end{prop}

\begin{proof} Let $\al:M\to N$ with $D\al=0$. Let $m\in M$ and $f\in D(A\al(m))$. Since $\Ta$ is injective, we can extend $f$ to a morphism $\td{f}\in DN$. Then $$f(\al(m))=\td{f}(\al(m))=(\td{f}\al)(m)=(D\al(\td{f}))(m)=0$$ so that $f=0$. Then $D(A\al(m))=0$ so that $A\al(m)=0$ and finally $\al(m)=0$. Since $m\in M$ was chosen arbitrarily, we see that $\al=0$ and therefore that $\al\mapsto D\al$ is injective for all $M$ and $N$, so that $D$ is faithful. 

Now assume that $M,N\in\mcM$, let $\bt:DN\to DM$ and set $\gm=\chi_N^{-1}D\bt\chi_M:M\to N$. We now use Lemma \ref{chi defines natural transformation on a modules relation chi and functor} to show that $$D(\chi_N^{-1}D\bt\chi_M)=D(\chi_M)D^2\bt D(\chi_N)^{-1}=D(\chi_M)D^2\bt\chi_{DN}=D(\chi_M)\chi_{DM}\bt=\bt,$$ proving that $\al\mapsto D\al$ is an isomorphism. \end{proof}

\begin{rmk} Note that the above proof that the induced map $\Hom(M,N)\to\Hom(DN,DM)$ is an isomorphism also works without the assumption that $M\in\mcM$. \end{rmk}

\begin{lem} \label{induced short exact sequence of cokernels} Each short exact sequence \begin{center} \begin{tikzcd}
     0 \arrow[r] & 
     L \arrow[r,"\vphantom{\bt}{\al}"] & 
     M \arrow[r,"\bt"] & 
     N \arrow[r] &
     0
    \end{tikzcd}
   \end{center} of $A$-modules induces a short exact sequence \begin{center} \begin{tikzcd}
     0 \arrow[r] &
     \cok(\chi_L) \arrow[r] & 
     \cok(\chi_M) \arrow[r] & 
     \cok(\chi_N) \arrow[r] &
     0
    \end{tikzcd}
   .\end{center} \end{lem}
   
\begin{proof} Since $M\mapsto\chi_M$ is a natural transformation, the following diagram commutes. \begin{center}
    \begin{tikzcd}[sep=2cm]
     0 \arrow[r] \arrow[d] &
     L \arrow[r,"\vphantom{\bt}{\al}"] \arrow[d,rightarrowtail,"\chi_L"'] &
     M \arrow[r,"\bt"] \arrow[d,rightarrowtail,"\chi_M"'] &
     N \arrow[d,rightarrowtail,"\chi_N",swap] \arrow[r] &
     0 \arrow[d] \\
     0 \arrow[r] &
     D^2L \arrow[r,"\vphantom{D^2\bt}{D^2\al}"] & 
     D^2M \arrow[r,"D^2\bt"] &
     D^2N \arrow[r] &
     0
    \end{tikzcd}
   \end{center} Since $D^2$ is exact, both rows are exact. We know from Lemma \ref{chi is mono equivalent cok condition} that the vertical maps are injective, and the result then follows from the Snake Lemma. \end{proof}

We shall make ample use of the following properties of the category $\mcM$ of $\Ta$-reflexive $A$-modules.

\begin{theorem} \label{properties of reflexive modules} The category $\mcM$ is abelian and is closed under submodules, quotients, extensions and isomorphisms, where $M\in\mcM$ if and only if $DM\in\mcM$. The functor $D$ restricts to a duality functor on $\mcM$, and $\chi$ restricts to a natural isomorphism $\id_{\mcM}\to D^2$. \end{theorem}

\begin{proof} If \begin{center}
    \begin{tikzcd}
     0 \arrow[r] & 
     L \arrow[r] & 
     M \arrow[r] &
     N \arrow[r] &
     0
    \end{tikzcd}
 \end{center} is exact, then it follows from Lemma \ref{induced short exact sequence of cokernels} that \begin{center}
    \begin{tikzcd}
     0 \arrow[r] & 
     \cok(\chi_L) \arrow[r] & 
     \cok(\chi_M) \arrow[r] &
     \cok(\chi_N) \arrow[r] &
     0
    \end{tikzcd}
 \end{center} is exact. So $\cok(\chi_M)=0$ if and only if $\cok(\chi_L)=\cok(\chi_N)=0$. We know from Lemma \ref{chi is mono equivalent cok condition} that $\cok(\chi_Z)=0$ if and only if $Z\in\mcM$, so we conclude that $M\in\mcM$ if and only if $L,N\in\mcM$. This shows that $\mcM$ is closed under submodules, quotients and extensions. Since finite direct sums are split extensions, kernels are submodules and cokernels are quotients, we see that $\mcM$ is abelian. 

If $M\in\mcM$ and $N$ is isomorphic to $M$ then there exists an isomorphism $\al:M\to N$ with $D^2\al\chi_M=\chi_N\al$ since $\chi$ is a natural transformation. \begin{center}
    \begin{tikzcd}[sep=2cm]
     M \arrow[r,"\al"] \arrow[d,"\chi_M",swap] &
     N \arrow[d,"\chi_N"] \\
     D^2M \arrow[r,"D^2\al"] & 
     D^2N
    \end{tikzcd}
 \end{center} So $D^2\al$, $\chi_M$ and $\al$ are isomorphisms, giving that $\chi_N$ also is and therefore $N\in\mcM$, so that $\mcM$ is closed under isomorphisms.

If $M\in\mcM$, then $\chi_{DM}=D(\chi_M)^{-1}$ is an isomorphism so that $DM\in\mcM$. If $DM\in\mcM$, then we apply the previous sentence to $DM$ to see that $D^2M\in\mcM$. Then, since $\mcM$ is closed under submodules and $\im(\chi_M)\leq D^2M$, we see that $\im(\chi_M)\in\mcM$. Finally, since $\chi_M$ is injective and $\mcM$ is closed under isomorphisms, we see that $M\simeq\im(\chi_M)\in\mcM$.

The remaining claims follow from Proposition \ref{d is faithful}. \end{proof}

\begin{lem} \label{mono/epi in theta-refs iff inj/surj} A morphism in $\mcM$ is a monomorphism (epimorphism) in $\mcM$ if and only if it is injective (surjective). \end{lem}

\begin{proof} It is clear that injective (surjective) morphisms in $\mcM$ are monomorphisms (epimorphisms) in $\mcM$. Let $\io:M\to N$ be a monomorphism in $\mcM$. Since $M,N\in\mcM$, we know from Theorem \ref{properties of reflexive modules} that $\ker(\io)\in\mcM$. Since $\io$ is a monomorphism in $\mcM$, and composing $\io$ with the inclusion $\ker(\io)\incl M$ and the zero map $\ker(\io)\to M$ both give zero, we see that these maps must be equal, and therefore $\ker(\io)=0$. So $\io$ is injective (a monomorphism in $\Mod_A$). A similar argument shows that epimorphisms in $\mcM$ are surjective. \end{proof}

\begin{prop} \label{d functor switches inj and proj in theta ref} Let $M\in\mcM$. Then $M$ is injective (projective) in $\mcM$ if and only if $DM$ is projective (injective) in $\mcM$. \end{prop}

\begin{proof} Let $P$ be projective in $\mcM$, $\io:M\to N$ be a monomorphism in $\mcM$ and $\al:M\to DP$. Then $D\io:DN\to DM$ is an epimorphism in $\mcM$ by Lemma \ref{mono/epi in theta-refs iff inj/surj}, so there exists a morphism $\bt:P\to DN$ in $\mcM$ for which the following diagram commutes. \begin{center}
    \begin{tikzcd} [sep=2cm]
     {} & DN \arrow[d,"D\io",twoheadrightarrow] \\
     P \arrow[r,"(D\al)\chi_P",swap] \arrow[ru,"\bt",dashed] & DM
     \end{tikzcd}\end{center} Then we use Lemma \ref{chi defines natural transformation on a modules relation chi and functor} to see that $(D\bt)\chi_N\io=(D\bt)(D^2\io)\chi_M$ is equal to $$D((D\io)\bt)\chi_M=D((D\al)\chi_P)\chi_M=D(\chi_P)D^2\al\chi_M=D(\chi_P)\chi_{DP}\al=\al.$$ \begin{center}
    \begin{tikzcd} [sep=2cm]
     M \arrow[d,"\al",swap] \arrow[r,"\io",rightarrowtail] & N \arrow[ld,"(D\bt)\chi_N",dashed] \\
     DP 
     \end{tikzcd}\end{center} We conclude that $DP$ is injective in $\mcM$. A similar argument shows that if $I$ is injective in $\mcM$ then $DI$ is projective in $\mcM$. For the converse, if $DP$ is injective in $\mcM$, then, from the argument above, $P\simeq D^2P$ is projective in $\mcM$, with a similar argument showing that if $DI$ is projective in $\mcM$, then $I$ is injective in $\mcM$. \end{proof}

\begin{prop} \label{inj proj and theta ref implies inj proj in theta ref} If an object in $\mcM$ is injective (projective) in $\Mod_A$, then it is injective (projective) in $\mcM$. \end{prop}

\begin{proof} Let $I\in\mcM$ be injective as an $A$-module. If $\io:M\to N$ is a monomorphism in $\mcM$, and $\al:M\to I$ is a morphism in $\mcM$, then we know from Proposition \ref{mono/epi in theta-refs iff inj/surj} that $\io$ is injective, and therefore a monomorphism in $\Mod_A$. Since $I$ is injective in $\Mod_A$, there exists a morphism $\bt:N\to I$ with $\bt\io=\al$. \begin{center}
    \begin{tikzcd} [sep=2cm]
     M \arrow[d,"\al",swap] \arrow[r,"\io",rightarrowtail] & N \arrow[ld,"\bt",dashed] \\
     I 
     \end{tikzcd}\end{center} But since $\mcM$ is full, $\bt$ is also a morphism in $\mcM$, so $I$ must be injective in $\mcM$. A similar argument shows that if an object in $\mcM$ is projective as an $A$-module, then it must be a projective object in $\mcM$. \end{proof}
     
\begin{prop} \label{inj in t ref implies inj} If a multibasic module is injective in $\mcM$, then it is injective as an $A$-module. \end{prop}

\begin{proof} Let $\bigoplus_i U_i/V_i$ be multibasic and injective in $\mcM$. Then the short exact sequence \begin{center}
    \begin{tikzcd}
     0 \arrow[r] & 
     \bigoplus_i U_i/V_i \arrow[r] & 
     \bigoplus_i K/V_i \arrow[r] &
     \bigoplus_i K/U_i \arrow[r] &
     0
    \end{tikzcd}
 \end{center} in $\mcM$ must split, so that $$\bigoplus_i K/V_i=\left(\bigoplus_i U_i/V_i\right)\oplus\left(\bigoplus_i K/U_i\right).$$ We know that $\bigoplus_i K/V_i$ is injective, so $\bigoplus_i U_i/V_i$ must also be injective. \end{proof}

\begin{lem} \label{hahn field is reflexive} $K\in\mcM$. \end{lem}

\begin{proof} We know from Lemma \ref{lemma for hahn field is reflexive} that $\zeta$ is an isomorphism, and therefore that $D\zeta$ is an isomorphism. Then applying $(D\zeta)^{-1}$ on the left of $D(\zeta)\chi_K=\zeta$ shows that $\chi_K=(D\zeta)^{-1}\zeta$ is an isomorphism, and therefore that $K\in\mcM$. \end{proof}

\begin{theorem} \label{multibasics are theta ref} Every multibasic module is contained in $\mcM$. \end{theorem}

\begin{proof} We know that $K\in\mcM$ (Lemma \ref{hahn field is reflexive}). If $V<U\leq K$, then since $\mcM$ is closed under submodules and quotients (Theorem \ref{properties of reflexive modules}) we know that $U/V\in\mcM$. If $M$ is basic, then it is isomorphic to $U/V$ for some $V<U\leq K$ (Theorem \ref{basic iso to standard basic}), so that $M\in\mcM$ follows from the fact that $\mcM$ is closed under isomorphisms. 

Finally, if $M$ is multibasic, then $M=\bigoplus_i M_i$, where each $M_i$ is basic and therefore contained in $\mcM$. It then follows from the fact that $\mcM$ is abelian (and therefore additive) that $M=\bigoplus_i M_i\in\mcM$, and we are done. \end{proof}

\begin{exmp} Let $M$ be injective and multibasic, or equivalently, a multibasic module whose basic summands are all injective, namely $K$, $\Ta$ and $\Phi$. Then $M$ is an injective $A$-module contained in $\mcM$ by Theorem \ref{multibasics are theta ref}, so must be an injective object in $\mcM$ by Proposition \ref{inj proj and theta ref implies inj proj in theta ref}. 

Now let $M$ be flat and multibasic, or equivalently, a multibasic module whose basic summands are all flat, namely $K$, $A$ and $I_{>0}$. Then since $DK=K$, $DA=\Ta$ and $D(I_{>0})=\Phi$, we see that $DM$ is injective and multibasic, and therefore an injective object in $\mcM$ by the above paragraph. Then $M$ is a projective object in $\mcM$ by Proposition \ref{d functor switches inj and proj in theta ref}. We will extend this to cover all flat modules contained in $\mcM$ in Corollary \ref{flat and reflexive implies projective object}. \end{exmp}

\begin{cor} \label{dm mb iff m mb} $M$ is multibasic if and only if $DM$ is. \end{cor}

\begin{proof} We already know that if $M$ is multibasic, then so is $DM$. If $DM$ is multibasic, then $DM\in\mcM$ so that $M\in\mcM$ and then $M\simeq D^2M$. Since $DM$ is multibasic, we know that $D^2M$ is multibasic and therefore $M$ is. \end{proof}

\begin{cor} \label{extensions of multibasics finite cases} Let \begin{center}
    \begin{tikzcd}
     0 \arrow[r] & 
     L \arrow[r] & 
     M \arrow[r] &
     N \arrow[r] &
     0
    \end{tikzcd}
 \end{center} be exact. If $L$ is finite, then $M$ is multibasic if and only if $N$ is multibasic. If $N$ is finite, then $M$ is multibasic if and only if $L$ is multibasic. \end{cor}

\begin{proof} If $L$ is finite and $M$ is multibasic, then $N$ is multibasic by Proposition \ref{l finite m multibasic implies n multibasic}. If $N$ is finite and $L$ is multibasic, then $M$ is multibasic by Proposition \ref{extension finite on right}.

Since $D$ is exact, we know that \begin{center}
    \begin{tikzcd}
     0 \arrow[r] & 
     DN \arrow[r] & 
     DM\arrow[r] &
     DL \arrow[r] &
     0
    \end{tikzcd}
 \end{center} is exact. If $L$ is finite and $N$ is multibasic, then $DL$ is finite and $DN$ is multibasic, so that $DM$ is multibasic by Proposition \ref{extension finite on right} and then $M$ is multibasic by Corollary \ref{dm mb iff m mb}. If $N$ is finite and $M$ is multibasic, then $DN$ is finite and $DM$ is multibasic, so that $DL$ is multibasic by Proposition \ref{l finite m multibasic implies n multibasic} and then $L$ is multibasic by Corollary \ref{dm mb iff m mb}. \end{proof}
 
\begin{cor} \label{flat and reflexive implies projective object} If $P\in\mcM$ is flat, then it is a projective object in $\mcM$. \end{cor}

\begin{proof} Let $\pi:M\to N$ be an epimorphism in $\mcM$ and $\al:P\to N$ be any morphism in $\mcM$. The map $$(x,f)\mapsto f(\al(x)):P\x DN\to\Ta$$ is bilinear, so there exists a unique morphism $\al':P\tns DN\to\Ta$ satisfying $$\al'(x\tns f)=f(\al(x))$$ for all $x\in P$ and $f\in DN$. We know from Lemma \ref{mono/epi in theta-refs iff inj/surj} that $\pi$ is surjective, so that $D\pi:DN\to DM$ is injective. Then since $P$ is flat, we know that $$\id_P\tns D\pi:P\tns DN\to P\tns DM$$ is injective. Since $\Ta$ is injective, there exists a morphism $\bt'':P\tns DM\to\Ta$ making the following diagram commute. \begin{center}
    \begin{tikzcd} [sep=2cm]
     P\tns DN \arrow[d,"\al'",swap] \arrow[r,"\id_P\tns D\pi",rightarrowtail] & P\tns DM \arrow[ld,"\bt''",dashed] \\
     \Ta 
     \end{tikzcd}\end{center} Then $$\bt''(x\tns(f\pi))=(\bt''(\id_P\tns D\pi))(x\tns f)=\al'(x\tns f)=f(\al(x))$$ for all $x\in P$ and $f\in DN$. Define $\bt':P\to D^2M$ by $$\bt'(x)(f)=\bt''(x\tns f)$$ for all $x\in P$ and $f\in DM$. Then $$\bt'(x)(f\pi)=\bt''(x\tns f\pi)=f(\al(x))$$ for all $x\in P$ and $f\in DN$. Noting that $M\in\mcM$, we set $\bt=\chi_M^{-1}\bt':P\to M$. Let $x\in P$ and $f\in DN$. Then $\bt'(x)\in D^2M$ and $M\in\mcM$, so $\bt'(x)=\chi_M(m)$ for some $m\in M$. This gives that $$(\pi\bt)(x)=\pi(\chi_M^{-1}(\bt'(x)))=\pi(m).$$ Then $$f(\pi(m))=\chi_M(m)(f\pi)=\bt'(x)(f\pi)=f(\al(x)).$$ So $$\chi_N(\al(x))(f)=f(\al(x))=f(\pi(m))=\chi_N(\pi(m))(f)$$ for all $f\in DN$ gives that $\chi_N(\al(x))=\chi_N(\pi(m))$. Then since $\chi_N$ is injective, we have that $$(\pi\bt)(x)=\pi(m)=\al(x)$$ for all $x\in P$, and we conclude that $\pi\bt=\al$. \begin{center}
    \begin{tikzcd} [sep=2cm]
     {} & M \arrow[d,"\pi",twoheadrightarrow] \\
     P \arrow[r,"\al",swap] \arrow[ru,"\bt",dashed] & N
     \end{tikzcd}\end{center} 

Since the subcategory $\mcM$ is full, $\bt$ must be contained in $\mcM$ and we conclude that $P$ is projective in $\mcM$. \end{proof}

We now present a partial converse to Corollary \ref{flat and reflexive implies projective object}.

\begin{prop} \label{proj in t ref implies flat} If a multibasic module is projective in $\mcM$, then it is flat. \end{prop}

\begin{proof} Let $M=\bigoplus_i U_i/V_i$ be multibasic and projective in $\mcM$. Then since \begin{center}
    \begin{tikzcd}
     0 \arrow[r] & 
     \bigoplus_i V_i \arrow[r] & 
     \bigoplus_i U_i \arrow[r] &
     \bigoplus_i U_i/V_i \arrow[r] &
     0
    \end{tikzcd}
 \end{center} is exact in $\mcM$, it must split so that $$\bigoplus_i U_i=\left(\bigoplus_i V_i\right)\oplus\left(\bigoplus_i U_i/V_i\right).$$ Then since $\bigoplus_i U_i$ is flat, we know that $\bigoplus_i U_i/V_i$ is flat. \end{proof}
 
\begin{prop} \label{for mb inj iff flat} Let $M$ be multibasic. Then $M$ is injective (flat) if and only if $DM$ is flat (injective). \end{prop}

\begin{proof} Let $M$ be multibasic. If $M$ is injective, then since it is contained in $\mcM$, it must be injective in $\mcM$. Then $DM$ is multibasic and projective in $\mcM$, so must be flat. Conversely, if $M$ is flat, then it must be projective in $\mcM$, so that $DM$ is injective in $\mcM$ and therefore an injective $A$-module. 

If $DM$ is flat, then since $M\in\mcM$ we can apply the above paragraph to see that $M\simeq D^2M$ is injective. Similarly, if $DM$ is injective, then $M\simeq D^2M$ is flat by applying the above paragraph. \end{proof}

\begin{exmp} $I_{>0}$ is flat and contained in $\mcM$, so is a projective object in $\mcM$. But if $I_{>0}$ is a projective $A$-module, then it must be free, and therefore isomorphic to copies of $A$. Since $I_{>0}\in\mcM$, it must be isomorphic to finitely many copies of $A$ and therefore be finitely generated, a contradiction. So $I_{>0}$ is a projective object in $\mcM$, but not a projective $A$-module. \end{exmp}

\section{Finiteness} \label{finiteness}

Now that we have established various properties of the category $\mcM$, we work towards proving that various modules have finite dimension over the fields $K$ or $\F$. Lemma \ref{direct sum is reflexive implies finitely many summands} is of fundamental importance to this goal, demonstrating that any infinite direct sum of non-zero $A$-modules cannot be contained in $\mcM$.

\begin{lem} \label{direct sum is reflexive implies finitely many summands} If $(M_i)_{i\in I}$ is a family of $A$-modules for which $\bigoplus_{i\in I}M_i\in\mcM$, then $M_i=0$ for all but finitely many $i$. \end{lem}

\begin{proof} Let $M=\bigoplus_i M_i$. We construct a family of maps $\al_i:M_i\to\Ta$ as follows. For each $i\in I$, if $M_i\neq 0$, then $DM_i\neq 0$, so we can choose non-zero $\al_i\in DM_i$. If $M_i=0$, set $\al_i=0$. This induces a map $\al:DM$. If $M_i=0$, then $\al_i=0$, so $\al(M_i)=0$. If $M_i\neq 0$, then $\al_i\neq 0$, so we can choose $m\in M_i$ with $\al_i(m)\neq 0$. Then $\al(m)=\al_i(m)\neq 0$, so that $\al(M_i)\neq 0$. So $M_i=0$ if and only if $\al(M_i)=0$, and we need to show that $\al(M_i)=0$ for all but finitely many $i$.

To that end, let $L$ be the submodule of $DM$ consisting of those $\gm:M\to\Ta$ for which $\gm(M_i)=0$ for all but finitely many $i$, and let $\bt\in D(DM/L)$. Since $M\in\mcM$, the following composite map \begin{center}
    \begin{tikzcd}
     DM \arrow[r,twoheadrightarrow] & 
     DM/L \arrow[r,"\bt"] & 
     \Ta
    \end{tikzcd}
 \end{center} is $\chi_M(m)$ for some $m\in M$. Let $i\in I$ and $\vphi\in DM_i$. Then $\vphi\pi_i\in DM$ and $(\vphi\pi_i)(M_j)=0$ for $j\neq i$, so $\vphi\pi_i\in L$. We then see that $\chi_{M_i}(\pi_i(m))(\vphi)$ is equal to $$\vphi(\pi_i(m))=(\vphi\pi_i)(m)=\chi_M(m)(\vphi\pi_i)=\bt(\vphi\pi_i+L)=0.$$ Since $\vphi\in DM_i$ was chosen arbitrarily, we conclude that $\chi_{M_i}(\pi_i(m))=0$. Then since $\chi_{M_i}$ is injective, we see that $\pi_i(m)=0$. Since $i\in I$ was chosen arbitrarily, this gives that $m=0$, so that $\bt=0$. Since $\bt\in D(DM/L)$ was chosen arbitrarily, we see that $D(DM/L)=0$, so that $DM/L=0$ and then $L=DM$. We conclude that $\al(M_i)=0$ for all but finitely many $i$, and therefore that $M_i=0$ for all but finitely many $i$, as required. \end{proof}
 
\begin{cor} An injective $A$-module is $\Ta$-reflexive if and only if it is multibasic involving only $K$, $\Ta$ and $\Phi$. A projective $A$-module is $\Ta$-reflexive if and only if it is free and finitely generated. \end{cor}

\begin{proof} By Proposition \ref{classification of injective a mods}, an injective $A$-module is either a finite direct sum of indecomposable injectives, which are both $\Ta$-reflexive and multibasic, or the injective hull of an infinite direct sum, which are neither. 

Since $A$ is local, the projective $A$-modules are simply the free ones. Those of finite rank are both $\Ta$-reflexive, and free and finitely generated. Those of infinite rank are neither. \end{proof}

Note that for each $A$-module $M$, the $A$-multiplication on $K\tns M$ extends uniquely to a $K$-multiplication satisfying $a(b\tns m)=(ab)\tns m$ for all $a,b\in K$ and $m\in M$. 

\begin{prop} \label{k tns m has fin dem over k} If $M\in\mcM$, then $K\tns M$ has finite dimension over $K$. \end{prop} 

\begin{proof} Since $\{1\tns m:m\in M\}$ $K$-spans $K\tns M$, there must be a subset $\{m_i\}_{i\in I}$ of $M$ for which $\{1\tns m_i\}_{i\in I}$ is a basis for $K\tns M$ over $K$. The multiplication maps $a\mapsto am_i:A\to M$ induce a map $\al:\bigoplus_{i\in I}A\to M$. If $a\in\bigoplus_{i\in I}A$ with $\al(a)=0$, then $$1\tns\al(a)=\sum_{i\in I}a_i(1\tns m_i)=0.$$ Then since $\{1\tns m_i\}$ is linearly independent over $K$, we see that $a_i=0$ for all $i$ and then $a=0$, so that $\al$ is injective. Then since $M\in\mcM$ we see from Theorem \ref{properties of reflexive modules} that $\bigoplus_{i\in I}A\in\mcM$, so that $I$ is finite by Lemma \ref{direct sum is reflexive implies finitely many summands}, and therefore $K\tns M$ has finite dimension over $K$. \end{proof}

\begin{lem} \label{f tensor m as quotient} For each $A$-module $M$, the sequence \begin{center}
    \begin{tikzcd}
     0 \arrow[r] &
     I_{>0}M \arrow[r] & 
     M \arrow[r] & 
     \F\tns M \arrow[r] &
     0
    \end{tikzcd}
 \end{center} is exact. \end{lem}
 
\begin{proof} Since \begin{center}
    \begin{tikzcd}
     0 \arrow[r] &
     I_{>0} \arrow[r] & 
     A \arrow[r] & 
     \F \arrow[r] &
     0
    \end{tikzcd}
 \end{center} is exact, we use the long exact sequence for Tor modules to see that \begin{center}
    \begin{tikzcd}
     I_{>0}\tns M \arrow[r] & 
     M \arrow[r] & 
     \F\tns M \arrow[r] &
     0
    \end{tikzcd}
 \end{center} is exact, and therefore $$\F\tns M\simeq\cok(I_{>0}\tns M\to M).$$ Every element of $I_{>0}\tns M$ can be expressed in the form $t^q\tns m$ for some $q>0$ and $m\in M$, and we use this to see that $\im(I_{>0}\tns M\to M)=I_{>0}M$. \end{proof}

\begin{prop} \label{f tns m is finite} If $M\in\mcM$, then $\F\tns M$ has finite dimension over $\F$. \end{prop}

\begin{proof} We choose points $\{m_i\}_{i\in I}$ in $M$ for which $\{m_i+I_{>0}M\}_{i\in I}$ forms a basis for $M/I_{>0}M$ over $\F$. Then $M/I_{>0}M$ is isomorphic to $\bigoplus_{i\in I}\F$ both over $A$ and over $\F$. Since $M\in\mcM$, we know that $M/I_{>0}M\in\mcM$, so that $\bigoplus_{i\in I}\F\in\mcM$ and then $I$ is finite by Lemma \ref{direct sum is reflexive implies finitely many summands}. So $\F\tns M\simeq M/I_{>0}M$ has finite dimension over $\F$ by Lemma \ref{f tensor m as quotient}. \end{proof}

\begin{lem} \label{ext is zero flat tref} Let $M$ and $N$ be $A$-modules. Then $$\Ext^i(M,DN)\simeq D\Tor_i(M,N).$$ \end{lem}

\begin{proof} Let \begin{center}
    \begin{tikzcd}
     \cdots \arrow[r] &
     P_2 \arrow[r] & 
     P_1 \arrow[r] & 
     P_0 \arrow[r] &
     M \arrow[r] &
     0
    \end{tikzcd}
 \end{center} be a projective resolution for $M$. This gives the cochain complex \begin{center}
    \begin{tikzcd}
     0 \arrow[r] &
     \Hom(P_0,DN) \arrow[r] & 
     \Hom(P_1,DN) \arrow[r] & 
     \Hom(P_2,DN) \arrow[r] &
     \cdots
    \end{tikzcd}
 \end{center} and we see that $\Ext^i(M,DN)$ is the cohomology of this complex at $\Hom(P_i,DN)$. Using the tensor-hom adjunction, this is the same as the cohomology of the complex \begin{center}
    \begin{tikzcd}
     0 \arrow[r] &
     D(P_0\tns N) \arrow[r] & 
     D(P_1\tns N) \arrow[r] & 
     D(P_2\tns N) \arrow[r] &
     \cdots
    \end{tikzcd}
 \end{center} at $D(P_i\tns N)$. Since exact functors commute with cohomology, this is the result of applying $D$ to the homology of \begin{center}
    \begin{tikzcd}
     \cdots \arrow[r] &
     P_2\tns N \arrow[r] & 
     P_1\tns N \arrow[r] & 
     P_0\tns N \arrow[r] &
     0
    \end{tikzcd}
 \end{center} at $P_i\tns N$, which is just $\Tor_i(M,N)$. \end{proof}

\begin{cor} \label{hom/ext in terms of tns/tor} If $N\in\mcM$, then $$\Ext^i(M,N)\simeq D\Tor_i(M,DN).$$ \end{cor} 

\begin{rmk} \label{ext in terms of tor remark} We see from Corollary \ref{hom/ext in terms of tns/tor} that the action of the functors $\Ext^i$ on modules in $\mcM$ can be understood purely in terms of the $\Tor_i$ functors and $D$. In particular, we have that $\Hom(\F,M)\simeq D(\F\tns DM)$. But since $\F\tns DM$ is unnaturally isomorphic to $\F^n$ for some $n$, we see that $\Hom(\F,M)\simeq D(\F\tns DM)$ is unnaturally isomorphic to $\F\tns DM$. Similarly, $\Ext(\F,M)\simeq D(\Tor(\F,DM))$ is unnaturally isomorphic to $\Tor(\F,DM)$. \end{rmk}

\begin{cor} If $M\in\mcM$, then $K\tns M$, $\Hom(K,M)$, $\Hom(M,K)$ have finite dimension over $K$, and $\F\tns M$, $\Hom(\F,M)$ and $\Hom(M,\F)$ are finite.  \end{cor}

\begin{proof} That $K\tns M$ has finite dimension over $K$ is Proposition \ref{k tns m has fin dem over k}. Then since $DK\simeq K$, we use the tensor-hom adjunction to see that $$\Hom(M,K)\simeq\Hom(M,DK)\simeq D(K\tns M).$$ Then since $K\tns M$ is unnaturally isomorphic to $K^n$ for some $n$, we see that $D(K\tns M)$ is unnaturally isomorphic to $K\tns M$, and therefore $\Hom(M,K)$ has the same finite dimension over $K$ as $K\tns M$. Then since $M\in\mcM$, we see that $$\Hom(K,M)\simeq\Hom(DM,DK)\simeq\Hom(DM,K)$$ has finite dimension over $K$, since $DM\in\mcM$. 

That $\F\tns M$ is finite is Proposition \ref{f tns m is finite}. Similarly to the case for $K$, we use the fact that $D\F\simeq\F$ to see that $\Hom(M,\F)\simeq D(\F\tns M)$ has the same finite dimension over $\F$ as $\F\tns M$. Then $\Hom(\F,M)\simeq\Hom(DM,\F)$ also must be finite. \end{proof}

\begin{cor} If $M$ is flat and $i\geq 1$, then $\Ext^i(M,DN)=0$. If $N\in\mcM$, then $\Ext^i(M,N)=0$. \end{cor}

We now work towards proving that $\Tor(\F,M)$ is finite whenever $M\in\mcM$ (Theorem \ref{tor(f2,m) is finite}), which will be of fundamental importance when proving that every object in $\mcM$ is multibasic in Theorem \ref{theta-ref implies multibasic}. We start by demonstrating that $\Tor(\F,M)$ is a colimit of a diagram $T(M)$ of subquotients of $M$. Then when $M\in\mcM$, this diagram is attained at a finite stage, and all objects in the diagram are finite, so that $\Tor(\F,M)$ must itself be finite. 

\begin{definition} Let $M$ be an $A$-module. For each $k\geq 0$, set $u_k=t^{2^{-k}}\in I_{>0}$, also denote the map $$m\mapsto u_km:M\to M$$ by $u_k$, and note that $u_{k+1}^2=u_k$. Then set $T_k'(M)=\ann(u_k,M)$ and $$T_k''(M)=\{m\in M:u_k\tns m=0\}.$$ Since $$T_k''(M)\leq T_k'(M)\leq M,$$ we can set $$T_k(M)=T_k'(M)/T_k''(M).$$ Now assume $k\geq 1$. Since $u_k$ maps $T_{k-1}'(M)$ into $T_k'(M)$ and $T_{k-1}''(M)$ into $T_k''(M)$, we can restrict $u_k$ to obtain a map $T_{k-1}'(M)\to T_k'(M)$, and then a map $T_{k-1}(M)\to T_k(M)$, which we label $\al_k$. We have thus constructed a diagram \begin{center}
    \begin{tikzcd}
     T_0(M) \arrow[r,"\al_1"] &
     T_1(M) \arrow[r,"\al_2"] & 
     T_2(M) \arrow[r,"\al_3"] & 
     \cdots
    \end{tikzcd}
 \end{center} of $A$-modules, which we label $T(M)$. 
\end{definition}

\begin{lem} \label{tor is colimit} $\Tor(\F,M)$ is the colimit of the diagram $T(M)$. Also, if $t^qm=0$ for some $q\in[0,2^{-k})$, then $m\in T_k''(M)$. \end{lem}

\begin{proof} Since \begin{center}
    \begin{tikzcd}
     0 \arrow[r] &
     I_{>0} \arrow[r] & 
     A \arrow[r] & 
     \F \arrow[r] &
     0
    \end{tikzcd}
 \end{center} is exact and $A$ is flat, we see from the long exact sequence of Tor modules that \begin{center}
    \begin{tikzcd}
     0 \arrow[r] &
     \Tor(\F,M) \arrow[r] & 
     I_{>0}\tns M \arrow[r] & 
     M \arrow[r] &
     \F\tns M \arrow[r] &
     0
    \end{tikzcd}
 \end{center} is exact. Noting that $\im(I_{>0}\tns M\to M)=I_{>0}M$, we label the corestriction of the map $I_{>0}\tns M\to M$ to its image by $\mu$. Then \begin{center}
    \begin{tikzcd}
     0 \arrow[r] &
     \Tor(\F,M) \arrow[r] & 
     I_{>0}\tns M \arrow[r,"\mu"] & 
     I_{>0}M \arrow[r] &
     0
    \end{tikzcd}
 \end{center} is exact, and we see that $\Tor(\F,M)\simeq\ker(\mu)$.

For each $k\geq 0$, the map $m\mapsto u_k\tns m:M\to I_{>0}\tns M$ sends $T_k'(M)$ into $\ker(\mu)$, and restricts to a map $T_k'(M)\to\ker(\mu)$ which has kernel $T_k''(M)$. We thus obtain an injective map $T_k(M)\to\ker(\mu)$ sending $m+T_k''(M)$ to $u_k\tns m$, which we label $i_k$.

Note that for $j\geq k$, we have that $$u_{k+1}\cdots u_j=t^{2^{-k-1}}\cdots t^{2^{-j}}=t^{2^{-k-1}+\cdots+2^{-j}}=t^{2^{-k}-2^{-j}}=u_ku_j^{-1}.$$ Then $$u_k\tns m=u_{k+1}\tns u_{k+1}m=\cdots=u_j\tns(u_{k+1}\cdots u_jm)=u_j\tns(u_ku_j^{-1}m),$$ so that the following diagram commutes. \begin{center}
    \begin{tikzcd} [sep=2cm]
     {} & {} & {} & \ker(\mu) \\
     T_k(M) \arrow[r,"\al_{k+1}"] \arrow[rrru,"i_k"] & T_{k+1}(M) \arrow[r,"\al_{k+2}"] & \cdots \arrow[r,"\al_j"] & T_j(M) \arrow[u,"i_j"]\\ 
     \end{tikzcd}\end{center} We conclude that $$(i_k:T_k(M)\to\ker(\mu))_{k\geq 0}$$ is a co-cone for the diagram $T(M)$. 

Let $$(\td{i_k}:T_k(M)\to\td{M})_{k\geq 0}$$ be another co-cone for this diagram. Each element of $\ker(\mu)$ has the form $$u_k\tns m=i_k(m+T_k''(M))$$ for some $k\geq 0$ and $m\in T_k'(M)$. We want to construct a well-defined map $\xi:\ker(\mu)\to\td{M}$ with $\xi i_k=\td{i_k}$ for all $k$; it must not depend on the choice of $k$ or $m$. To demonstrate that this is possible, suppose that $k_1,k_2\geq 0$, $m_1\in T_{k_1}'(M)$, $m_2\in T_{k_2}'(M)$ with $$i_{k_1}(m_1+T_{k_1}''(M))=i_{k_2}(m_2+T_{k_2}''(M)).$$ Assume without loss of generality that $k_1\leq k_2$. Then since the $\td{i_k}$ form a co-cone, we see that $$\td{i_{k_1}}(m_1+T_{k_1}''(M))=\td{i_{k_2}}(u_{k_1+1}\cdots u_{k_2}m_1+T_{k_2}''(M))=\td{i_{k_2}}(u_{k_1}u_{k_2}^{-1}m_1+T_{k_2}''(M)).$$ Since $u_{k_1}\tns m_1=u_{k_2}\tns m_2$, we have that $$u_{k_2}\tns(u_{k_1}u_{k_2}^{-1}m_1-m_2)=0$$ and then $$u_{k_1}u_{k_2}^{-1}m_1-m_2\in T_{k_2}''(M).$$ Then $$\td{i_{k_1}}(m_1+T_{k_1}''(M))=\td{i_{k_2}}(m_2+T_{k_2}''(M)),$$ as required. So we can define $\xi:\ker(\mu)\to\td{M}$ by $$\xi(i_k(m+T_k''(M)))=\td{i_k}(m+T_k''(M))$$ for all $k\geq 0$ and $m\in T_k'(M)$. It is clear that this is the unique morphism $\ker(\mu)\to\td{M}$ for which the following diagram commutes. \begin{center}
    \begin{tikzcd} [sep=2cm]
     {} & \td{M} \\
     T_k(M) \arrow[r,"i_k"] \arrow[ru,"\td{i_k}"] & \ker(\mu) \arrow[u,"\xi",dashed,swap] 
     \end{tikzcd}\end{center} We conclude that $$(i_k:T_k(M)\to\ker(\mu))_{k\geq 0}$$ is a colimit for the above diagram. Since $\Tor(\F,M)\simeq\ker(\mu)$, $\Tor(\F,M)$ is the colimit of $T(M)$. 
     
Finally, suppose that $m\in M$ and $t^qm=0$ for some $q\in[0,2^{-k})$. Then we can choose $j\geq k$ with $2^{-k}-2^{-j}\geq q$, so that $$u_k\tns m=u_j\tns(u_ku_j^{-1}m)=u_j\tns(t^{2^{-k}-2^{-j}}m)=0$$ and then $m\in T_k''(M)$. \end{proof}

\begin{lem} \label{diagram has injective maps} Each $\al_k$ is injective, and $T_k(M)\simeq\F\tns T_k(M)$ for each $k$. \end{lem}

\begin{proof} Each element of $\ker(\al_k)$ has the form $m+T_{k-1}''(M)$ for some $m\in T_{k-1}'(M)$ with $u_km\in T_k''(M)$. Then $$u_{k-1}\tns m=u_k^2\tns m=u_k\tns(u_km)=0,$$ so that $m\in T_{k-1}''(M)$. We conclude that $\al_k$ is injective.

If $a\in I_{>0}$ and $m\in T_k'(M)$, then $a=t^qb$ for some $q>0$ and $b\in A$, and $u_km=0$. If $q\geq 2^{-k}$, then $$am=bt^{q-2^{-k}}(u_km)=0\in T_k''(M).$$ Otherwise, $q<2^{-k}$, so that $2^{-k}-q\in[0,2^{-k})$ and $$t^{2^{-k}-q}(am)=t^{2^{-k}-q}(t^qbm)=b(u_km)=0,$$ so that $am\in T_k''(M)$. We conclude that $I_{>0}T_k'(M)\leq T_k''(M)$. Then each element of $I_{>0}T_k(M)$ has the form $$a(m+T_k''(M))=am+T_k''(M)=T_k''(M)$$ for some $a\in I_{>0}$ and $m\in T_k'(M)$, so that $I_{>0}T_k(M)=0$. Then $$T_k(M)=T_k(M)/I_{>0}T_k(M)\simeq\F\tns T_k(M).$$ \end{proof}

\begin{lem} Set $Q_k(M)=\cok(\al_k)$ and $$Z_k(M)=T_k''(M)+u_k(T_{k-1}'(M)).$$ If $M\in\mcM$, then $T_k(M)$ and $Q_k(M)\simeq T_k'(M)/Z_k(M)$ are both finite. \end{lem}

\begin{proof} If $M\in\mcM$, then $T_k'(M)\leq M$ shows that $T_k'(M)\in\mcM$. Then $$T_k(M)=T_k'(M)/T_k''(M)\in\mcM,$$ so that $T_k(M)\simeq\F\tns T_k(M)$ is finite by Lemma \ref{diagram has injective maps}. Then $$Q_k(M)=\cok(\al_k)=T_k(M)/\im(\al_k)$$ is a quotient of a finite module, and therefore must be finite. Finally, we have that $\im(\al_k)=Z_k(M)/T_k''(M)$ and therefore $Q_k(M)\simeq T_k'(M)/Z_k(M)$. \end{proof}

We now use bases for $T_k'(M)/Z_k(M)$ for each $k$ to construct an $A$-module $P$ and a morphism $p:P\to M$. We will prove by induction that $p$ is injective, so that if $M\in\mcM$, then $P\in\mcM$. This will eventually lead to the conclusion that our diagram $T(M)$ is attained at a finite stage. 

\begin{definition} Let $M\in\mcM$, and for each $k\geq 1$, let $\{m_{jk}:j\in J_k\}$ be a finite subset of $T_k'(M)$ for which $\{m_{jk}+Z_k(M):j\in J_k\}$ forms a basis for $T_k'(M)/Z_k(M)$ over $\F$. We construct an $A$-module $P$ and a morphism $p:P\to M$ as follows.

Let $k\geq 1$. Then set $B_k=A/I_{2^{-k}}$ and $P_k=\bigoplus_{j\in J_k}B_k$. For each $j\in J_k$, multiplying by $m_{jk}$ sends $I_{2^{-k}}$ to $0$, giving a map $B_k\to T_k'(M)$. Combining these maps for each $j$ gives a map $p_k:P_k\to T_k'(M)$ with $$p_k\left(\sum_{j\in J_k}(a_j+I_{2^{-k}})\right)=\sum_{j\in J_k}a_jm_{jk}.$$

Now set $P_{\leq k}=\bigoplus_{i=1}^k P_i$. For each $i$, we have the composite \begin{center}
    \begin{tikzcd}
     P_i \arrow[r,"p_i"] & 
     T_i'(M) \arrow[r,hookrightarrow] & 
     M \;\;\;.
    \end{tikzcd}
 \end{center} Combining these maps together gives a map $p_{\leq k}:P_{\leq k}\to M$ with $$p_{\leq k}\left(\sum_{i\leq k}x_i\right)=\sum_{i\leq k}p_i(x_i).$$ 

Now set $P=\bigoplus_{k=1}^\infty P_k$. If $x\in P$, then $x\in P_{\leq k}$ for some $k\geq 1$ and then $p_{\leq k}(x)\in M$. If $x\in P_k$ and $x\in P_{k'}$, then without loss of generality we may assume that $k\leq k'$, and it is clear that $p_{\leq k}(x)=p_{\leq {k'}}(x)$. This allows us to define $p:P\to M$ by $p(x)=p_{\leq k}(x)$ for any $k\geq 1$ and $x\in P_{\leq k}$. \end{definition}

\begin{lem} \label{technical lemma for tor f m finite} Let $k\geq 1$ and $x\in P_k$ be non-zero. Then $x=t^qy$ for some $q\in[0,2^{-k})$ and $y\in P_k$ with $p_k(y)\notin Z_k(M)$. \end{lem}

\begin{proof} We know that $x=\sum_{j\in J_k}(a_j+I_{2^{-k}})$, where $a_j\in A$ are non-zero. Since $J_k$ is finite and non-empty, we can set $q=\min\{\nu(a_j):j\in J_k\}\in[0,2^{-k})$. Then $a_j+I_{2^{-k}}=t^q(b_j+I_{2^{-k}})$ for some non-zero $b_j\in A$, for each $j$. Set $y=\sum_{j\in J_k}(b_j+I_{2^{-k}})\in P_k$. Then $x=t^qy$ with $y\neq 0$.

Suppose that $y\in I_{>0}P_k$. Then $y=t^py'$ for some $p>0$ and $y'\in P_k$. Then $y'=\sum_{j\in J_k}(c_j+I_{2^{-k}})$ for some non-zero $c_j\in A$. If $p\geq 2^{-k}$, then $y=t^{p-2^{-k}}(u_ky')$. Since $u_kB_k=0$, we see that $u_kP_k=0$, so that $u_ky'=0$ and then $y=0$, a contradiction. So $p\in(0,2^{-k})$. Since there exists $j\in J_k$ with $\nu(a_j)=q$, we know that $a_j=t^qe_j$ for some $e_j\in A^\x$. Then $t^q(b_j-e_j)\in I_{2^{-k}}$, so that $b_j-e_j\in I_{2^{-k}}-q\in I_{>0}$. If $b_j\in I_{>0}$, then $e_j\in I_{>0}$, a contradiction, so $b_j\in A^\x$. Then $b_j+I_{2^{-k}}=t^p(c_j+I_{2^{-k}})$, so that $b_j-t^pc_j\in I_{2^{-k}}\leq I_{>0}$. But $t^pc_j\in I_p\leq I_{>0}$, giving that $b_j\in I_{>0}$ and then $b_j\notin A^\x$, a contradiction. So $y\notin I_{>0}P_k$.

Let $j\in J_k$. If $s\in(0,2^{-k})$, then $t^{2^{-k}-s}(t^sm_{jk})=u_km_{jk}=0$ since $m_{jk}\in T_k'(M)$. Then since $2^{-k}-s\in[0,2^{-k})$, we have that $t^s m_{jk}\in T_k''(M)\leq Z_k(M)$. If $s\geq 2^{-k}$, then $t^sm_{jk}=t^{s-2^{-k}}(u_km_{jk})=0\in Z_k(M)$. So $t^sm_{jk}\in Z_k(M)$ for all $s>0$, and then $I_{>0}m_{jk}\leq Z_k(M)$ for all $j\in J_k$. Then $$p_k(y)+Z_k(M)=\sum_{j\in J_k}b_jm_{jk}+Z_k(M)=\sum_{j\in J_k}\lmb_j m_{jk}+Z_k(M),$$ where $\lmb_j\in\{0,1\}$ for all $j$ and $\lmb_j=1$ for some $j$ ($\lmb_j=1$ if and only if $b_j\in A^\x$). Since $\{m_{jk}+Z_k(M):j\in J_k\}$ forms a basis for $T_k'(M)/Z_k(M)$ over $\F$, we know that $\sum_{j\in J_k}\lmb_j m_{jk}+Z_k(M)\neq Z_k(M)$, so that $p_k(y)\notin Z_k(M)$. 
\end{proof}

\begin{lem} $p_1$ is injective. \end{lem}

\begin{proof} Suppose not. Then there exists non-zero $x\in P_1$ with $p_1(x)=0$. So $x=t^qy$ for some $q\in[0,\frac{1}{2})$ and $y\in P_1$ with $p_1(y)\notin Z_1(M)$ by Lemma \ref{technical lemma for tor f m finite}. Then $p_1(y)\notin T_1''(M)$, but $$t^qp_1(y)=p_1(x)=0$$ so that $p_1(y)\in T_1''(M)$, a contradiction. \end{proof}

\begin{lem} $p_{\leq k}$ is injective for all $k\geq 1$. \end{lem}

\begin{proof} We proceed by induction on $k$. Since $p_{\leq 1}=p_1$, the $k=1$ case clearly holds. Now let $k\geq 2$ and assume that $p_{\leq k-1}$ is injective. Let $x\in\ker(p_{\leq k})$, so that $$p_{\leq k}(x)=\sum_{i\leq k}p_i(x_i)=0.$$

First of all, suppose that $x_k$ is non-zero. Then $x_k=t^qy$ for some $q\in[0,2^{-k})$ and $y\in P_k$ with $p_k(y)\notin Z_k(M)$ by Lemma \ref{technical lemma for tor f m finite}. Set $r=2^{-k}-q>0$ and $$x'=\sum_{i<k}t^rx_i\in P_{\leq k-1}.$$ Then $$t^rx_k=t^{r+q}y=u_ky\in u_kP_k=0.$$ Then since $x'\in P_{\leq k-1}$, we see that $$p_{\leq k-1}(x')=\sum_{i<k}p_i((x')_i)=\sum_{i\leq k}p_i((x')_i)=t^r\sum_{i\leq k}p_i(x_i)=t^rp_{\leq k}(x)=0.$$ So $x'\in\ker(p_{\leq k-1})$, and $p_{\leq k-1}$ is injective by the induction hypothesis, giving that $x'=0$. So $t^rx_i=0$ for all $i\leq k$. 

For $i<k$, note that $$2^{-i}-r\geq 2^{1-k}-r=2^{-k}+(2^{-k}-r)=2^{-k}+q,$$ so that $t^{2^{-i}-r}=t^qu_kb_i$ for some $b_i\in A$. Let $x_i=\sum_{j\in J_i}(a_j+I_{2^{-i}})$ for some $a_j\in A$. Then $$\sum_{j\in J_i}(t^ra_j+I_{2^{-i}})=t^rx_i=0.$$ So for all $j\in J_i$, we have that $t^ra_j\in I_{2^{-i}}$, so that $\nu(a_j)\geq 2^{-i}-r$, and then $a_j\in I_{2^{-i}-r}$, so that $a_j=t^{2^{-i}-r}c_j$ for some $c_j\in A$. 

Let $z_i=b_i\sum_{j\in J_i}(c_j+I_{2^{-i}})\in P_i$. Then $$x_i=\sum_{j\in J_i}(t^{2^{-i}-r}c_j+I_{2^{-i}})=t^{2^{-i}-r}\sum_{j\in J_i}(c_j+I_{2^{-i}})=t^qu_kz_i.$$ Set $w=\sum_{i<k}p_i(z_i)\in M$ and $b=p_k(y)+u_kw\in M$. If $i<k$, then $$u_{k-1}z_i=u_k^2z_i=t^{r+q}u_kz_i=t^rt^qu_kz_i=t^rx_i=0,$$ so that $$u_{k-1}w=u_{k-1}\sum_{i<k}p_i(z_i)=\sum_{i<k}p_i(u_{k-1}z_i)=0$$ and then $w\in T_{k-1}'(M)$. Then $$t^qu_kw=\sum_{i<k}p_i(t^qu_kz_i)=\sum_{i<k}p_i(x_i)=p_{\leq k}(x)-p_k(x_k)=p_k(x_k)=t^qp_k(y).$$ Then $$t^qb=t^q(p_k(y)+u_kw)=0,$$ with $q\in[0,2^{-k})$, so that $b\in T_k''(M)$. So $$p_k(y)=b+u_kw\in T_k''(M)+u_k(T_{k-1}'(M))=Z_k(M),$$ a contradiction. So $x_k=0$. 

Then $x\in P_{\leq k-1}$, so that $p_{\leq k-1}(x)=p_{\leq k}(x)=0$, and since $p_{\leq k-1}$ is injective by the induction hypothesis, we see that $x=0$. So $\ker(p_{\leq k})=0$ and we conclude that $p_{\leq k}$ is injective. The result then follows by induction. \end{proof}

\begin{lem} \label{colimit attained at finite stage for tor} If $M\in\mcM$, then the colimit of the diagram $T(M)$ is attained at a finite stage. \end{lem}

\begin{proof} Since $M\in\mcM$, for each $k$ we have that $T_k'(M)/Z_k(M)$ is finite, so there exists a finite subset $\{m_{jk}:j\in J_k\}$ for which $\{m_{jk}+Z_k(M):j\in J_k\}$ is a basis for $T_k'(M)/Z_k(M)$ over $\F$. This induces a map $p:P\to M$, where $p_{\leq k}$ is injective for all $k\geq 1$. Let $x\in\ker(p)$. Then $p(x)=p_{\leq k}(x)=0$ for some $k\geq 1$ with $x\in P_{\leq k}$. So $x\in\ker(p_{\leq k})=0$ and therefore $p$ is injective. Since $M\in\mcM$, this shows that $P\in\mcM$ and therefore $P_k=\bigoplus_{j\in J_k}B_k=0$ for sufficiently large $k$. Since $B_k$ is non-zero for all $k$, it must be the case that $J_k$ is empty for sufficiently large $k$. But since $J_k$ indexes a basis for $T_k'(M)/Z_k(M)$ over $\F$, we see that $J_k$ is empty if and only if $$\cok(\al_k)=Q_k(M)\simeq T_k'(M)/Z_k(M)=0$$ if and only if $\al_k$ is surjective. But we already know that $\al_k$ is injective for all $k$. So $\al_k$ is an isomorphism for sufficiently large $k$, so that $T_k(M)$ is a colimit of $T(M)$ for sufficiently large $k$. \end{proof}

\begin{theorem} \label{tor(f2,m) is finite} If $M\in\mcM$, then $\Tor(\F,M)$ is finite. \end{theorem}

\begin{proof} We know that $\Tor(\F,M)$ is the colimit of the diagram $T(M)$, which is attained at a finite stage. So $\Tor(\F,M)\simeq T_k(M)$ for sufficiently large $k$. Since $T_k(M)$ is finite, we are done. \end{proof}

\begin{cor} \label{ext functors are finite} If $M\in\mcM$, then $\Ext(\F,M)$ and $\Ext(M,\F)$ are finite. \end{cor}

\begin{proof} We use Corollary \ref{hom/ext in terms of tns/tor} to see that $\Ext(\F,M)\simeq D(\Tor(\F,DM))$, which is unnaturally isomorphic to $\Tor(\F,DM)$, which is finite since $DM\in\mcM$. Also $\Ext(M,\F)\simeq D(\Tor(\F,M))$ is unnaturally isomorphic to $\Tor(\F,M)$, which is finite.  \end{proof}

\begin{cor} \label{m tref implies i gr 0 tns m tref} If $M\in\mcM$, then so is $I_{>0}\tns M$. \end{cor}

\begin{proof} As demonstrated in the proof of Lemma \ref{tor is colimit}, we can form a short exact sequence \begin{center}
    \begin{tikzcd}
     0 \arrow[r] &
     \Tor(\F,M) \arrow[r] & 
     I_{>0}\tns M \arrow[r] & 
     I_{>0}M \arrow[r] &
     0 \;\;\;.
    \end{tikzcd}
 \end{center} Since $M\in\mcM$, we know from Lemma \ref{tor(f2,m) is finite} that $\Tor(\F,M)$ is finite and therefore contained in $\mcM$. Also $I_{>0}M\leq M$, so $I_{>0}M\in\mcM$. Then since $\mcM$ is closed under extensions, we conclude that $I_{>0}\tns M\in\mcM$. \end{proof}

\begin{lem} If $M\in\mcM$, then $$D(I_{>0}\tns M)\simeq\Hom(I_{>0},DM),$$ $$I_{>0}\tns M\simeq D\Hom(I_{>0},DM),$$ and $$I_{>0}\tns DM\simeq D\Hom(I_{>0},M).$$ \end{lem}

\begin{proof} It follows from the tensor-hom adjunction that $D(I_{>0}\tns M)\simeq\Hom(I_{>0},DM)$. Applying $D$ to this we see that $D^2(I_{>0}\tns M)\simeq D\Hom(I_{>0},DM)$. Then since $M\in\mcM$, we know from Corollary \ref{m tref implies i gr 0 tns m tref} that $I_{>0}\tns M\in\mcM$. So $$I_{>0}\tns M\simeq D^2(I_{>0}\tns M)\simeq D\Hom(I_{>0},DM).$$ We then apply this to $DM\in\mcM$ to see that $$I_{>0}\tns DM\simeq D\Hom(I_{>0},D^2M)\simeq D\Hom(I_{>0},M).$$ \end{proof}

\begin{prop} \label{isomorphism lemma for D of quotient f2 thing} If $M\in\mcM$, then $$D\left(\frac{M}{\ann(I_{>0},M)}\right)\simeq I_{>0}DM.$$ \end{prop}

\begin{proof} Since $\ker(M\to\Hom(I_{>0},M))=\ann(I_{>0},M)$ we see that $$D\left(\frac{M}{\ann(I_{>0},M)}\right)\simeq D(\im(M\to\Hom(I_{>0},M)))\simeq\im(D(M\to\Hom(I_{>0},M))).$$ Then since \begin{center}
    \begin{tikzcd} [sep=2cm]
     I_{>0}\tns DM \arrow[r] \arrow[d] & DM \arrow[d,equal] \\
     D\Hom(I_{>0},M) \arrow[r] & DM
     \end{tikzcd}\end{center} commutes, with both vertical maps being isomorphisms, we see that the horizontal maps have the same image. So $$D\left(\frac{M}{\ann(I_{>0},M)}\right)\simeq\im(I_{>0}\tns DM\to DM)=I_{>0}DM.$$ \end{proof}

\begin{cor} Let $M,N\in\mcM$ with at least one of these modules being multibasic. Then $M\tns N$, $\Hom(M,N)$, $\Tor(M,N)$ and $\Ext(M,N)$ are all contained in $\mcM$, and $$\Tor_i(M,N)=\Ext^i(M,N)=0$$ for $i\geq 2$. \end{cor}

\begin{proof} Since the functors $\Tor_i$ and $\Ext^i$ are additive, and every basic module is isomorphic to a subquotient of $K$, we only need to show that if $M\in\mcM$ and $V<U\leq K$, then $M\tns(U/V)$, $\Hom(M,U/V)$, $\Hom(U/V,M)$, $\Tor(M,U/V)$, $\Ext(M,U/V)$, $\Ext(U/V,M)$ are contained in $\mcM$, and that $$\Tor_i(M,U/V)=\Ext^i(M,U/V)=\Ext^i(U/V,M)=0$$ for $i\geq 2$. By expressing $\Ext^i$ in terms of $\Tor_i$ and $D$, we can reduce this further, so that we only need to show that $M\tns(U/V),\Tor(M,U/V)\in\mcM$ and that $\Tor^i(M,U/V)=0$ for $i\geq 2$. Then since \begin{center}
    \begin{tikzcd}
     0 \arrow[r] &
     V \arrow[r] & 
     U \arrow[r] & 
     U/V \arrow[r] &
     0
    \end{tikzcd}
 \end{center} is exact, we use the long exact sequence of Tor modules, and the fact that $U$ and $V$ are flat, to see that \begin{center}
    \begin{tikzcd}
     0 \arrow[r] &
     \Tor(U/V,M) \arrow[r] & 
     V\tns M \arrow[r] & 
     U\tns M \arrow[r] &
     (U/V)\tns M \arrow[r] &
     0
    \end{tikzcd}
 \end{center} and \begin{center}
    \begin{tikzcd}
     0 \arrow[r] &
     \Tor_i(U/V,M) \arrow[r] & 
     0
    \end{tikzcd}
 \end{center} are exact for $i\geq 2$. So $\Tor_i(U/V,M)=0$ for $i\geq 2$, $\Tor(U/V,M)\simeq\ker(V\tns M\to U\tns M)$ and $(U/V)\tns M\simeq\cok(V\tns M\to U\tns M)$, reducing to simply showing that $U\tns M\in\mcM$ when $M\in\mcM$ and $U$ is flat and basic. So $U$ is isomorphic to one of $K$, $A$, $I_{>0}$. The $U=A$ case is clear. $K\tns M\in\mcM$ since it has finite dimension over $K$, covering the $U=K$ case. All that remains is to show that $I_{>0}\tns M\in\mcM$, which is covered by Lemma \ref{m tref implies i gr 0 tns m tref}. \end{proof}

\section{Theta-Reflexive Modules are Multibasic} \label{decomp of theta ref section}

In this section, we prove that every $\Ta$-reflexive module can be decomposed into a finite list of basic submodules i.e. that it is multibasic, and therefore that $\mcM$ is simply the category of multibasic modules. We begin by decomposing a $\Ta$-reflexive module into a multibasic summand and a summand in the subcategory $\mcZ$, whose definition we now recall.

\begin{definition} (Reminder of Definition \ref{reflexivity intro definitions}) $\mcZ$ is the full subcategory of $\mcM$ consisting of those $M$ for which $K\tns M=K\tns DM=0$. \end{definition}

\begin{lem} \label{related modules that are zero} $\mcZ$ is closed under submodules, quotients, extensions and the functor $D$. \end{lem}

\begin{proof} Let \begin{center}
    \begin{tikzcd}
     0 \arrow[r] &
     L \arrow[r] & 
     M \arrow[r] & 
     N \arrow[r] &
     0
    \end{tikzcd}
 \end{center} be exact. Since $K$ is flat, $K\tns(-)$ is exact. Also $D$ is exact, so \begin{center}
    \begin{tikzcd}
     0 \arrow[r] &
     K\tns L \arrow[r] & 
     K\tns M \arrow[r] & 
     K\tns N \arrow[r] &
     0
    \end{tikzcd}
 \end{center} and \begin{center}
    \begin{tikzcd}
     0 \arrow[r] &
     K\tns DN \arrow[r] & 
     K\tns DM \arrow[r] & 
     K\tns DL \arrow[r] &
     0
    \end{tikzcd}
 \end{center} are both exact. So $K\tns M=0$ if and only if $K\tns L=K\tns N=0$, and $K\tns DM=0$ if and only if $K\tns DL=K\tns DN=0$. So $M\in\mcZ$ if and only if $L,N\in\mcZ$. Proving the remaining claim is routine. \end{proof}

\begin{lem} \label{splitting reflexive mod k tns m is zero} Each $M\in\mcM$ splits as $M=M_0\oplus M_1$, with $M_0$ multibasic, $M_1\in\mcM$ and $K\tns M_1=0$. \end{lem}

\begin{proof} Let $\al:M\to K\tns M$ be given by $\al(m)=1\tns m$. Since $M\in\mcM$, we know that $K\tns M$ has finite dimension over $K$, so that $\im(\al)\leq K\tns M$ must be flat and multibasic, and therefore a projective object in $\mcM$. Since $\ker(\al)\in\mcM$ also, we see that \begin{center}
    \begin{tikzcd}
     0 \arrow[r] &
     \ker(\al) \arrow[r] & 
     M \arrow[r] & 
     \im(\al) \arrow[r] &
     0
    \end{tikzcd}
 \end{center} is a split exact sequence in $\mcM$. Finally, since $a\tns m=a\al(m)=0$ for all $a\in K$ and $m\in\ker(\al)$, we see that $K\tns\ker(\al)=0$. \end{proof}

\begin{lem} \label{splits with K tns M and K tns DM zero} Each $M\in\mcM$ splits as $M=M_0\oplus M_1$, where $M_0$ is multibasic and $M_1\in\mcZ$. \end{lem}

\begin{proof} We see from Lemma \ref{splitting reflexive mod k tns m is zero} that $M$ splits as $M=M_0\oplus M_1$, with $M_0$ multibasic, $M_1\in\mcM$ and $K\tns M_1=0$; applying it to $DM_1$ shows that $DM_1=N_0\oplus N_1$, with $N_0$ multibasic, $N_1\in\mcM$ and $K\tns N_1=0$. Then applying $D$ to this splitting shows that $M_1=M_2\oplus M_3$, where $M_2=DN_0$ and $M_3=DN_1$. Since $K$ is flat, the inclusion $M_3\incl M_1$ induces an injective map $K\tns M_3\to K\tns M_1=0$, confirming that $K\tns M_3=0$. Then since $N_1\in\mcM$, we see that $K\tns DM_3=K\tns N_1=0$, so that $M_3\in\mcZ$. Finally, $M=(M_0\oplus M_2)\oplus M_3$, where $M_0\oplus M_2$ is multibasic and $M_3\in\mcZ$. \end{proof}

\begin{lem} \label{existence of max basic submods} Every non-zero $A$-module has a maximal basic submodule. \end{lem}

\begin{proof} Let $M$ be a non-zero $A$-module and $\mcP$ be the set of its basic submodules. Since $M$ is non-zero, we know that $Am\in\mcP$ for some non-zero $m\in M$, so $\mcP$ cannot be empty. It is clear that $\mcP$ is a partially ordered set with respect to inclusion. If $\{M_i\}_{i\in I}$ is a chain in $\mcP$, then one can directly use the definition of a basic module to check that $\bigcup_{i\in I}M_i\in\mcP$ and is an upper bound for the $M_i$. The result then follows from Zorn's Lemma. \end{proof}

\begin{lem} \label{maximal basic submod is injective} If $M\in\mcM$ and $\F\tns M=\Ext(\F,M)=0$, then every maximal basic submodule of $M$ is injective. \end{lem}

\begin{proof} Let $N$ be a maximal basic submodule of $M$. Since $N$ is basic, there exist $V<U\leq K$ and injective map $\al:U/V\to M$ with $\im(\al)=N$ and where $U$ is $K$, $A$ or $I_{>0}$. We want to exclude the cases $U=A$ and $U=I_{>0}$, so that $U=K$ and then $U/V=K/V$ is injective.

Suppose that $U=A$. Since $M/I_{>0}M\simeq\F\tns M=0$, we see that $M=I_{>0}M$. Then $\al(1+V)=t^qm$ for some $q>0$ and $m\in M$. If $a,b\in A$ and $t^{-q}a+V=t^{-q}b+V$, then $t^{-q}(a-b)\in V$, so that $\al(t^{-q}(a-b)+V)=0$. But since $V<U=A$, we also have that $$\al(t^{-q}(a-b)+V)=t^{-q}(a-b)\al(1+V)=am-bm,$$ so that $am=bm$. This allows us to define $\al':I_{-q}/V\to M$ by $\al'(t^{-q}a+V)=am$ for all $a\in A$. $I_{-q}/V$ properly extends $U/V$, and if $a\in A=U$, we have that $$\al'(a+V)=\al'(t^{-q}t^qa+V)=t^qam=a\al(1+V)=\al(a+V).$$ So $\al'$ extends $\al$, giving that $N\leq\im(\al')$. Since $I_{-q}/V$ is basic, so is $\im(\al')\leq M$. Since $N$ is maximal, this gives that $N=\im(\al')$. 

Then $m=\al'(t^{-q}+V)\in\im(\al')=\im(\al)$, so that $$m=\al(a+V)=a\al(1+V)=t^qam$$ for some $a\in A$. If $V$ is non-zero, then $I_{nq}\leq V$ for some $n\geq 1$. Then since $(t^qa)^n\in I_{nq}\leq V$, we have that $m=(t^qa)^{n+1}m=t^q((t^qa)^na)m=\al((t^qa)^na+V)=0$, so that $\al=0$ and then $N=0$. But $N$ is basic, a contradiction, so $V=0$. Then $\al(1+V)=t^qm=t^{2q}am=\al(t^qa+V)$ with $\al$ injective, so that $1+t^qa\in V=0$. So $t^{-q}=a\in A$, giving that $q\leq 0$, a contradiction. So $U\neq A$.

Now suppose that $U=I_{>0}$. Since $\Ext(\F,M)=0$ and \begin{center}
    \begin{tikzcd}
     0 \arrow[r] &
     I_{>0} \arrow[r] & 
     A \arrow[r] & 
     \F \arrow[r] &
     0 \;\;\;,
    \end{tikzcd}
 \end{center} is exact, we can use one of the long exact sequences for Ext modules to see that \begin{center}
    \begin{tikzcd}
     M \arrow[r] & 
     \Hom(I_{>0},M) \arrow[r] &
     0
    \end{tikzcd}
 \end{center} is exact i.e. that the map $M\to\Hom(I_{>0},M)$ is surjective. So there exists $m\in M$ for which $\al(a+V)=am$ for all $a\in I_{>0}$. Then $\ann(m)=V$. Then we define $\al':A/V\to M$ by $\al'(a+V)=am$ for all $a\in A$, which clearly extends $\al$, so that $N\leq\im(\al')$. Since $\im(\al')\leq M$ is basic and $N$ is maximal, we see that $N=\im(\al')$. Then $m=\al'(1+V)\in\im(\al')=\im(\al)$, so that $m=\al(a+V)=am$ for some non-zero $a\in I_{>0}$. If $V$ is non-zero, then there exists $n\geq 1$ with $I_{n\nu(a)}\leq V$. So $a^n\in I_{n\nu(a)}$ and then $m=a^nm=0$, so that $\al=0$ implies $N=0$, a contradiction. So $V=0$, and then $\al(a+V)=am=a^2m=\al(a^2+V)$ gives that $a+V=a^2+V$ since $\al$ is injective. So $a+a^2\in V=0$ gives that $a=a^2$ and then $\nu(a)=0$. But $a\in I_{>0}$, a contradiction. We conclude that $U\neq I_{>0}$. 
 
The only remaining possibility is that $U=K$, as required. \end{proof}

\begin{cor} \label{conditions for multibasic} If $M\in\mcM$, then we have the following. 

1. If $$\F\tns M=\Ext(\F,M)=0,$$ then $M$ is injective and multibasic. 

2. If $$\Tor(\F,M)=\Hom(\F,M)=0,$$ then $M$ is flat and multibasic. 

3. If $$\F\tns M=\Ext(\F,M)=\Tor(\F,M)=\Hom(\F,M)=0,$$ then $M$ is flat, injective and multibasic, and therefore has finite dimension over $K$. \end{cor}

\begin{proof} 1. If $M=0$, then clearly it is injective and multibasic. Otherwise, it contains a maximal basic submodule $M_0$ (Lemma \ref{existence of max basic submods}), which must be injective by Lemma \ref{maximal basic submod is injective}. So the short exact sequence \begin{center}
    \begin{tikzcd}
     0 \arrow[r] &
     M_0 \arrow[r] & 
     M \arrow[r] & 
     M/M_0 \arrow[r] &
     0
    \end{tikzcd}
 \end{center} splits, giving that $M=M_0\oplus M_1$ for some $M_1$. Since the tensor functor and Ext functor are additive, we see that $M_1\in\mcM$ with $$\F\tns M_1=\Ext(\F,M_1)=0.$$ So we can repeat this process with $M_1$ to see that either $M_1=0$ or $M_1=M_2\oplus M_3$, where $M_2$ is injective and basic. This process continues splitting off injective basic submodules until the remaining summand is zero. If this process never terminates, then $M$ would contain an infinite direct sum of basic modules that is contained in $\mcM$, contradicting Lemma \ref{direct sum is reflexive implies finitely many summands}. So this process must eventually terminate, giving that $M$ is a finite direct sum of injective basic modules, and therefore must be injective and multibasic.

2. We see from Remark \ref{ext in terms of tor remark} that $DM\in\mcM$ with $$\F\tns DM=\Ext(\F,DM)=0,$$ so we can apply the above argument to see that $DM$ is injective and multibasic, so that $M$ is flat and multibasic by Proposition \ref{for mb inj iff flat}.

3. Using 1. and 2., we see that $M$ must be flat, injective and multibasic. So all of its summands must be flat and injective basic modules, all of which are isomorphic to $K$. So $M$ must have finite dimension over $K$. \end{proof}

\begin{lem} \label{reflexive in Z is extension of multibasics} If $M\in\mcZ$, then there exists finitely generated $L$ and multibasic $N$ with $\F\tns N=0$ forming a short exact sequence \begin{center}
    \begin{tikzcd}
     0 \arrow[r] &
     L \arrow[r] & 
     M \arrow[r] & 
     N \arrow[r] &
     0 \;\;\;.
    \end{tikzcd}
 \end{center} \end{lem}

\begin{proof} We know that $\F\tns M$ has finite dimension over $\F$ (Proposition \ref{f tns m is finite}). Since all elements of $\F\tns M$ take the form $1\tns m$ for some $m\in M$, we can choose a finite subset $\{m_i\}_{i\in I}$ of $M$ for which $\{1\tns m_i\}_{i\in I}$ forms a basis for $\F\tns M$ over $\F$. Let $L$ be the submodule of $M$ generated by $\{m_i\}_{i\in I}$, and let $N=M/L$. Then \begin{center}
    \begin{tikzcd}
     0 \arrow[r] &
     L \arrow[r] & 
     M \arrow[r] & 
     N \arrow[r] &
     0
    \end{tikzcd}
 \end{center} is exact, with $L$ finitely generated. Each element of $\F\tns N$ takes the form $1\tns(m+L)$ for some $m\in M$. Then $1\tns m\in\F\tns M$, so that $$1\tns m=\sum_j \left(1\tns m_j\right)=1\tns\left(\sum_j m_j\right),$$ where $j$ runs through a subset of $I$. Let $x=\sum_jm_j\in L$, so that $1\tns(m-x)=0$ and then $$m-x\in\ker(M\to\F\tns M)=\im(I_{>0}\tns M\to M)=I_{>0}M.$$ So $m=x+t^qm'$ for some $q>0$ and $m'\in M$. Then $$1\tns(m+L)=1\tns(t^qm'+L)=t^q\tns(m'+L)=0,$$ and we conclude that $\F\tns N=0$. Since $N\in\mcM$, we know from Corollary \ref{ext functors are finite} that $\Ext(\F,N)$ is finite. Since \begin{center}
    \begin{tikzcd}
     N \arrow[r] & 
     \Hom(I_{>0},N) \arrow[r] & 
     \Ext(\F,N) \arrow[r] &
     0
    \end{tikzcd}
 \end{center} is exact, we see that $\Ext(\F,N)\simeq\cok(N\to\Hom(I_{>0},N))$. So the latter of these is finite, and then we can choose a finite set of maps $\al_i:I_{>0}\to N$ for which $\{\al_i+\im(N\to\Hom(I_{>0},N))\}_{i\in I}$ forms a basis for $\cok(N\to\Hom(I_{>0},N))$ over $\F$. These maps induce a map $\al:\bigoplus_{i\in I}I_{>0}\to N$, and the quotient map $N\to\cok(\al)$ induces a surjective map $0=\F\tns N\to\F\tns\cok(\al)$, so that $\F\tns\cok(\al)=0$. Let $\bt:I_{>0}\to\cok(\al)$. Since $I_{>0}$ is projective, there exists $\gm:I_{>0}\to N$ for which the following diagram commutes. \begin{center}
    \begin{tikzcd} [sep=2cm]
     {} & N \arrow[d,twoheadrightarrow] \\
     I_{>0} \arrow[r,"\bt"] \arrow[ru,"\gm",dashed] & \cok(\al)
     \end{tikzcd}\end{center} Then we use our basis for $\cok(N\to\Hom(I_{>0},N))$ to see that $$\gm+\im(N\to\Hom(I_{>0},N))=\sum_i\al_i+\im(N\to\Hom(I_{>0},N)),$$ where $i$ runs through a subset of the $I$. Setting $$\phi=\gm-\sum_i\al_i\in\im(N\to\Hom(I_{>0},N)),$$ we see that there exists $y\in N$ for which $\phi(a)=ay$ for some $y\in N$. Since the composite \begin{center}
    \begin{tikzcd}
     I_{>0} \arrow[r,"\al_i"] & 
     N \arrow[r,twoheadrightarrow] & 
     \cok(\al)
    \end{tikzcd}
 \end{center} is zero for all $i$, we see that $\bt$ is the composite \begin{center}
    \begin{tikzcd}
     I_{>0} \arrow[r,"\phi"] & 
     N \arrow[r,twoheadrightarrow] & 
     \cok(\al) \;\;\;.
    \end{tikzcd}
 \end{center} So $\bt(a)=a(y+\im(\al))$ for all $a\in I_{>0}$, giving that $\bt\in\im(\cok(\al)\to\Hom(I_{>0},\cok(\al)))$. We conclude that $\cok(\cok(\al)\to\Hom(I_{>0},\cok(\al)))=0$. Since \begin{center}
    \begin{tikzcd}
     \cok(\al) \arrow[r] & 
     \Hom(I_{>0},\cok(\al)) \arrow[r] & 
     \Ext(\F,\cok(\al)) \arrow[r] &
     0
    \end{tikzcd}
 \end{center} is exact, we see that $\Ext(\F,\cok(\al))=0$. Since $\cok(\al)\in\mcM$ and $$\F\tns\cok(\al)=\Ext(\F,\cok(\al))=0,$$ we see from Corollary \ref{conditions for multibasic} that $\cok(\al)$ is injective and multibasic. But $M\in\mcZ$, so $\cok(\al)\in\mcZ$ and then $K\tns D(\cok(\al))=0$. Since $\cok(\al)$ is injective and multibasic, we know that $D(\cok(\al))$ is flat, so that the canonical map $D(\cok(\al))\to K\tns D(\cok(\al))=0$ is injective, giving that $\cok(\al)=0$. So $\al$ is surjective, giving an isomorphism $\frac{\bigoplus_{i\in I}I_{>0}}{\ker(\al)}\to N$ and then a monomorphism $N\to\frac{\bigoplus_{i\in I}A}{\ker(\al)}$, which has cokernel $\bigoplus_{i\in I}\F$. So we have a short exact sequence  \begin{center}
    \begin{tikzcd}
    	 0 \arrow[r] &
     N \arrow[r] & 
     \frac{\bigoplus_{i\in I}A}{\ker(\al)} \arrow[r] & 
     \bigoplus_{i\in I}\F \arrow[r] &
     0
    \end{tikzcd}
 \end{center} Then since $\bigoplus_{i\in I}\F$ is finite and $\frac{\bigoplus_{i\in I}A}{\ker(\al)}$ is finitely generated and therefore multibasic, we know from Corollary \ref{extensions of multibasics finite cases} that $N$ is multibasic, as required. \end{proof}

\begin{lem} \label{hom i grt zero fin gen} If $M\in\mcZ$, then $\Hom(I_{>0},M)$ is finitely generated. \end{lem}

\begin{proof} We know from Lemma \ref{reflexive in Z is extension of multibasics} that there exist multibasic $L$ and $N$ with \begin{center}
    \begin{tikzcd}
     0 \arrow[r] &
     L \arrow[r] & 
     M \arrow[r] & 
     N \arrow[r] &
     0
    \end{tikzcd}
 \end{center} exact. Since $I_{>0}$ is flat and contained in $\mcM$, it is a projective object in $\mcM$, and therefore the functor $\Hom_{\mcM}(I_{>0},-):\mcM\to\mcM$ is exact. Since $\Hom_{\mcM}(I_{>0},X)=\Hom(I_{>0},X)$ for all $X\in\mcM$, we see that \begin{equation} \label{hom exact sequence}
    \begin{tikzcd}
     0 \arrow[r] &
     \Hom(I_{>0},L) \arrow[r] & 
     \Hom(I_{>0},M) \arrow[r] & 
     \Hom(I_{>0},N) \arrow[r] &
     0
    \end{tikzcd}
 \end{equation} is exact in $\Mod_A$. We know from Lemma \ref{related modules that are zero} that $L,N\in\mcZ$. Then from Proposition \ref{functor tables} and Remark \ref{multibasics in category z rmk}, we see that $\Hom(I_{>0},L)$ and $\Hom(I_{>0},N)$ must be finitely generated. Since an extension of finitely generated modules is finitely generated, we see from \eqref{hom exact sequence} that $\Hom(I_{>0},M)$ is finitely generated. \end{proof}

\begin{prop} If $M\in\mcM$, then $I_{>0}\tns M$ is multibasic. \end{prop}

\begin{proof} First of all, suppose that this holds for all $M\in\mcZ$. We know from Lemma \ref{splits with K tns M and K tns DM zero} that $M=M_0\oplus M_1$, where $M_0$ is multibasic and $M_1\in\mcZ$. Then $I_{>0}\tns M_1$ is multibasic by our assumption, and since $M_0$ is multibasic we know that $I_{>0}\tns M_0$ is multibasic. So $I_{>0}\tns M=I_{>0}\tns M_0\oplus I_{>0}\tns M_1$ is multibasic. So we may assume without loss of generality that $M\in\mcZ$. 

Since $M\in\mcZ$, we know that $DM\in\mcZ$ so that $\Hom(I_{>0},DM)$ is finitely generated and therefore multibasic by Lemma \ref{hom i grt zero fin gen}. Then since $$D(I_{>0}\tns M)\simeq\Hom(I_{>0},DM),$$ we see that $D(I_{>0}\tns M)$ is multibasic and therefore $I_{>0}\tns M$ is multibasic. \end{proof}

\begin{theorem} \label{theta-ref implies multibasic} If $M\in\mcM$, then $M$ is multibasic. \end{theorem}

\begin{proof} Since \begin{center}
    \begin{tikzcd}
     0 \arrow[r] &
     I_{>0} \arrow[r] & 
     A \arrow[r] & 
     \F \arrow[r] &
     0
    \end{tikzcd}
 \end{center} is exact, we know that \begin{center}
    \begin{tikzcd}
     0 \arrow[r] &
     \Tor(\F,M) \arrow[r] & 
     I_{>0}\tns M \arrow[r] & 
     M \arrow[r] &
     \F\tns M \arrow[r] &
     0
    \end{tikzcd}
 \end{center} is exact. Since $\im(I_{>0}\tns M\to M)=I_{>0}M$, we know that \begin{center}
    \begin{tikzcd}
     0 \arrow[r] &
     \Tor(\F,M) \arrow[r] & 
     I_{>0}\tns M \arrow[r] & 
     I_{>0}M \arrow[r] &
     0
    \end{tikzcd}
 \end{center} and \begin{center}
    \begin{tikzcd}
     0 \arrow[r] &
     I_{>0}M \arrow[r] & 
     M \arrow[r] & 
     \F\tns M \arrow[r] &
     0
    \end{tikzcd}
 \end{center} are exact. Since $\Tor(\F,M)$ is finite and $I_{>0}\tns M$ multibasic, we see from Corollary \ref{extensions of multibasics finite cases} that $I_{>0}M$ is multibasic. Then since $I_{>0}M$ is multibasic and $\F\tns M$ is finite, we see again from Corollary \ref{extensions of multibasics finite cases} that $M$ is multibasic. \end{proof}

\begin{cor} \label{multibasics ab and closed under exts subs quos} The category of multibasic modules is abelian and is closed under submodules, quotients and extensions. \end{cor}

\begin{cor} $\mcM$ has enough projectives and enough injectives. \end{cor}

\begin{proof} Let $M\in\mcM$. Then $M$ is multibasic, so $M=\bigoplus_i U_i/V_i$ with $V_i<U_i\leq K$. Then the inclusion $\bigoplus_i U_i/V_i\to\bigoplus_i K/V_i$ is a monomorphism to an injective object, and the quotient map $\bigoplus_i U_i\to\bigoplus_i U_i/V_i$ is an epimorphism from a projective object. \end{proof}

\begin{exmp} \label{proj inj res within t ref} Let $V_i<U_i\leq K$. Then \begin{center}
    \begin{tikzcd}
     0 \arrow[r] &
     \bigoplus_i U_i/V_i \arrow[r] & 
     \bigoplus_i K/V_i \arrow[r] & 
     \bigoplus_i K/U_i \arrow[r] &
     0
    \end{tikzcd}
 \end{center} and \begin{center}
    \begin{tikzcd}
     0 \arrow[r] &
     \bigoplus_i V_i \arrow[r] & 
     \bigoplus_i U_i \arrow[r] & 
     \bigoplus_i U_i/V_i \arrow[r] &
     0
    \end{tikzcd}
 \end{center} are injective and projective resolutions for the multibasic module $\bigoplus_i U_i/V_i$ in $\mcM$. Note that the injective resolution in $\mcM$ is also an injective resolution of $A$-modules. \end{exmp}
 
\begin{cor} $\mcM$ has global dimension 1. \end{cor}

\begin{cor} The injective objects in $\mcM$ are the injective and multibasic modules. The projective objects in $\mcM$ are the flat and multibasic modules. For $M\in\mcM$, we have that $M$ is injective if and only if $DM$ is flat. $\mcM$ is not a Frobenius category. \end{cor}

\begin{proof} We know that injective modules contained in $\mcM$ are injective objects in $\mcM$. Conversely, suppose that $M$ is an injective object in $\mcM$. Then it is multibasic, and must also be injective by Proposition \ref{inj in t ref implies inj}. The classification of projective objects in $\mcM$ is proved similarly using Proposition \ref{proj in t ref implies flat}. That $M$ is injective if and only if $DM$ is flat then follows from Proposition \ref{for mb inj iff flat}. Since the flat multibasics and the injective modules do not coincide ($A$ is projective but not injective, $\Ta$ is injective but not projective), the projective and injective objects in $\mcM$ do not coincide, so $\mcM$ is not a Frobenius category. \end{proof}

\begin{conj} Let $n(X)$ be the number of basic summands in a given module $X\in\mcM$, and \begin{center}
    \begin{tikzcd}
     0 \arrow[r] &
     L \arrow[r,"\al"] & 
     M\arrow[r,"\bt"] & 
     N \arrow[r] &
     0
    \end{tikzcd}
 \end{center} be exact, with $L,N\in\mcZ$. Then $n(M)\leq n(L)+n(N)$. \end{conj}

\chapter{The Infinite Root Algebra} \label{infinite root algebra chapter}

In Chapters \ref{hahn ring chapter} and \ref{reflexivity}, we introduced the Hahn ring $A$ and developed the theory of multibasic $A$-modules. We then developed the theory of $\Ta$-reflexive $A$-modules, demonstrating that the full subcategory $\mcM$ of $\Mod_A$ consisting of the $\Ta$-reflexive modules is simply the category of multibasic modules, and that each $M\in\mcM$ has a unique decomposition into basic summands up to the ordering of the summands. 

In this chapter, we turn our attention to the infinite root algebra $P=A/I_{>1}$, where we use our knowledge of $A$-modules to derive corresponding results for $P$-modules. Our interest in the infinite root algebra stems from the fact that it is self-injective (Proposition \ref{p and q injective}) and has no non-trivial finitely presented ideals (Proposition \ref{infinite root algebra totally incoherent}), similar to the properties for the associated stable homotopy ring of a Freyd category. We consider the ungraded context, however given that $P$ is a trivially graded graded-commutative ring one would expect that these results generalise easily to the graded context. Since $P$ is self-injective and has no non-trivial finitely presented ideals, and the category $\mcM'$ of $\Ta$-reflexive $P$-modules is abelian with enough projectives and enough injectives, it is plausible that there could be an ungraded triangulation structure for $P$ (see Definition \ref{def ungraded triang structure}). However, as we will see in Theorem \ref{no ungraded triang structure for inf root alg}, no such ungraded triangulation structure for $P$ exists that is contained within $\mcM'$. 

The infinite root algebra is constructed from the field of Hahn series with value group $\R$ and residue field $\F$, and was first described in \cite[\S 9]{ShS14}. Hahn series were first described by Hans Hahn \cite{Hah95}. The infinite root algebra construction can be generalised by replacing $\F$ with an arbitrary field, as in \cite{ShS14}, but we restrict to the case of $\F$ to keep the coefficients simple. The exposition of the basic properties of the infinite root algebra follows \cite{ShS14}.

\section{The Infinite Root Algebra}

\begin{definition} We define the \textit{infinite root algebra} to be the quotient ring $P=A/I_{>1}$. Note that $(t+I_{>1})^2=I_{>1}$, so that $P$ is not an integral domain. \end{definition}

\begin{definition} We say that an $A$-module $M$ is \textit{truncated} if $I_{>1}M=0$. It is clear that a $P$-module is the same as a truncated $A$-module, with the $P$-multiplication given by $(a+I_{>1})m=am$ for all $a\in A$ and $m\in M$. A morphism of $P$-modules is the same as a morphism of truncated $A$-modules. We thus identify $\Mod_P$ with the full subcategory of truncated $A$-modules in $\Mod_A$. We set $Q=A/I_1$, which is truncated and therefore a $P$-module. \end{definition}

\begin{prop} \cite[Proposition 9.21]{ShS14} \label{infinite root algebra totally incoherent} $P$ has no non-trivial finitely presented ideals. \end{prop}

\begin{proof} Let $J$ be a non-zero finitely presented ideal in $P$. Then $J$ is finitely generated, so $J=I_q/I_{>1}$ for some $0\leq q\leq 1$. Let $\pi:P\to I_q/I_{>1}$ be given by $\pi(a+I_{>1})=t^qa+I_{>1}$ for all $a\in A$. This is clearly surjective, and since $J$ is finitely presented, $\ker(\pi)$ must be finitely generated. But $a+I_{>1}\in\ker(\pi)$ if and only if $t^qa\in I_{>1}$ if and only if $a\in I_{>1-q}$ if and only if $a+I_{>1}\in I_{>1-q}/I_{>1}$. So $\ker(\pi)=I_{>1-q}/I_{>1}\simeq I_{>0}/I_{>q}$ is finitely generated and therefore $q=0$, giving that $J=P$, as required. \end{proof}

\begin{prop} An $A$-module $M$ is truncated if and only if $DM$ is. \end{prop}

\begin{proof} Suppose that $M$ is truncated. Then if $q>1$ and $\al\in DM$, then $$(t^q\al)(m)=\al(t^qm)\in\al(I_{>1}M)=0$$ for all $m\in M$. So $t^q\al=0$ and we see that $DM$ is truncated. 

Conversely, if $DM$ is truncated we apply the above paragraph to see that $D^2M$ is truncated. Then the injective map $\chi_M:M\to D^2M$ induces an injective map $I_{>1}M\to I_{>1}D^2M$, so that $M$ is truncated. \end{proof}

We can thus regard $D$ as a functor $\Mod_P\to\Mod_P$. 

\begin{prop} Let $\gm:P\to\Ta$ be given by $\gm(a+I_{>1})=t^{-1}a+I_{>0}$ for all $a\in A$, and for each truncated $A$-module $M$, set $\eta_M=\Hom(M,\gm)$. Then $\eta:D\to\Hom_P(-,P)$ is a natural isomorphism. So $\Hom_P(M,P)=0$ implies $M=0$. \end{prop}

\begin{proof} It is clear that $\eta$ defines a natural transformation; we just need to show that $\eta_M$ is an isomorphism. Since $\gm$ is injective, it is clear that $\eta_M$ is injective. Now let $\al\in DM$. Since $I_{>1}M=0$, we see that $I_{>1}\im(\al)=0$. Since $\im(\al)\leq\Ta$, we know that $\im(\al)=J/I_{>0}$ for some $I_{>0}\leq J\leq K$. So $I_{>1}(J/I_{>0})=0$, giving that $I_{>1}J\leq I_{>0}$. Then $J\leq I_{>1}^\circ=I_{-1}$, so that $\im(\al)\leq I_{-1}/I_{>0}$. Then we can form the composite \begin{center}
    \begin{tikzcd}
     M \arrow[r,"\al"] & 
     I_{-1}/I_{>0} \arrow[r,"\sim"] & 
     P
    \end{tikzcd}
 \end{center} and see that $\Hom(M,\gm)$ sends this composite to $\al$, so that $\Hom(M,\gm)$ is an isomorphism, as required. \end{proof}
 
We shall write $\Hom_P(-,P)$ as $D$ for the remainder of this chapter.

\begin{definition} For each $A$-module $M$, we define $$\ann(I_{>1},M)=\{m\in M:I_{>1}m=0\},$$ which is a truncated $A$-submodule of $M$ and therefore a $P$-module. \end{definition} 

The following proposition allows us to convert injective $A$-modules into injective $P$-modules, which we will use to prove the self-injectivity of $P$.

\begin{prop} \label{injective p mods from injective a mods} If $J$ is an injective $A$-module, then $\ann(I_{>1},J)$ is an injective $P$-module. \end{prop}

\begin{proof} Let $\io:M\to N$ be a monomorphism of truncated $A$-modules and $\al:M\to\ann(I_{>1},J)$. One can check that $\ker(\io)$ is truncated, and then since $\io$ is a monomorphism of truncated $A$-modules, $\ker(\io)=0$ so that $\io$ is injective, and therefore a monomorphism in $\Mod_A$. Then since $J$ is an injective $A$-module, there must exist a map $\td{\al}:N\to J$ for which the following diagram commutes.   \begin{center}
    \begin{tikzcd}[sep=2cm]
     M \arrow[r,"\io",rightarrowtail] \arrow[d,"\al"] & 
     N \arrow[ldd,"\td{\al}",dashed] \arrow[ld,"\bt",dashed] \\ 
     \ann(I_{>1},J) \arrow[d,rightarrowtail] \\ 
     J 
     \end{tikzcd}
     \end{center} If $q>1$ and  $y\in N$, then $$t^q\td{\al}(y)=\td{\al}(t^qy)\in\td{\al}(I_{>1}N)=0$$ since $N$ is truncated. So $\im(\td{\al})\leq\ann(I_{>1},J)$ and we can corestrict $\td{\al}$ to form a map $\bt:N\to\ann(I_{>1},J)$ with $\bt\io=\al$, as shown in the above diagram. We conclude that $\ann(I_{>1},J)$ is an injective object in the category of truncated $A$-modules, and therefore must be an injective $P$-module. \end{proof}
     
\begin{prop} \label{p and q injective} $P$ and $Q$ are injective $P$-modules. In particular, $P$ is a self-injective ring. \end{prop}
     
\begin{proof} Since $\Ta$ and $\Phi$ are injective $A$-modules, we use Proposition \ref{injective p mods from injective a mods} to see that $$\ann(I_{>1},\Ta)=(I_{>0}:I_{>1})/I_{>0}=I_{-1}/I_0\simeq A/I_{>1}=P$$ and $$\ann(I_{>1},\Phi)=(A:I_{>1})/A=I_{-1}/A\simeq A/I_1=Q$$ are injective $P$-modules. \end{proof}

For a more direct proof that $P$ is self-injective that does not involve $A$-modules, see Proposition 9.20 of \cite{ShS14} (the original source of this result). Note that since $P$ is self-injective, the arguments of \S\ref{reflexive graded modules} apply here. We now use Theorem \ref{classification of injectives thm} to classify the injective $P$-modules.

\begin{theorem} \label{classification of injective p mods} The indecomposable injective $P$-modules are $P$ and $Q$. The injective $P$-modules are the finite direct sums of copies of $P$ and $Q$, and the injective hulls of infinite direct sums of copies of $P$ and $Q$. \end{theorem}

\begin{proof} Since $A$ is a uniserial ring, it is clear that $P$ also is. Its cyclic modules are $A/I_{>q}$ for $0\leq q\leq 1$ and $A/I_q$ for $0<q\leq 1$, which have injective hulls $P$ and $Q$, respectively. So the indecomposable injective $P$-modules are precisely $P$ and $Q$ by Theorem \ref{classification of injectives thm}. Then the injective $P$-modules are precisely the injective hulls of the direct sums of copies of $P$ and $Q$, again by Theorem \ref{classification of injectives thm}. Of course, the injective hull of such a \textbf{finite} direct sum is just the direct sum itself. Then the only remaining case is that of the injective hulls of \textbf{infinite} direct sums of copies of $P$ and $Q$. \end{proof}

\begin{definition} We say that a $P$-module $M$ is \textit{basic} if it is non-zero and $x,y\in M$ implies that $x\in Py$ or $y\in Px$. It is \textit{multibasic} if it has a finite decomposition into basic submodules. It is \textit{$\Ta$-reflexive} if the canonical map $M\to D^2M$ is an isomorphism. We denote the category of $\Ta$-reflexive $P$-modules by $\mcM'$. It is routine to check that these definitions are equivalent to $M$ having the corresponding properties as an $A$-module. \end{definition}

The following theorem follows immediately from the equivalence of the above definitions for $A$-modules and $P$-modules.

\begin{theorem} \label{t ref iff mb for p mods} A $P$-module is $\Ta$-reflexive if and only if it is multibasic. Every $M\in\mcM'$ has a unique finite decomposition into basic submodules, up to the ordering of the summands. \end{theorem}

\begin{cor} $\mcM'$ is an abelian category with enough projectives and enough injectives. It contains $P$ and is closed under submodules, quotients and isomorphisms. The functor $D$ restricts to a duality functor on $\mcM'$. \end{cor}

\begin{proof} This follows easily from Theorem \ref{t ref iff mb for p mods} and the corresponding properties for $\mcM$. \end{proof}

\begin{exmp} One noteworthy difference between $\mcM$ and $\mcM'$ is that $\mcM'$ is not closed under extensions, as can be seen from the fact that \begin{center}
    \begin{tikzcd} [sep=1.5cm]
     0 \arrow[r] & 
     Q \arrow[r] & 
     A/I_2 \arrow[r] & 
     Q \arrow[r] &
     0 
     \end{tikzcd}
     \end{center} is exact, with $Q\in\mcM'$, yet $A/I_2\in\mcM$ but not in $\mcM'$ since $I_{>1}(A/I_2)=I_{>1}/I_2\neq 0$. \end{exmp}

\begin{exmp} \label{simple inj res for p mods} Note the following short exact sequences. \begin{center}
    \begin{tikzcd} [sep=1.5cm]
     0 \arrow[r] & 
     \F \arrow[r] & 
     P \arrow[r] & 
     Q \arrow[r] &
     0 
     \end{tikzcd}
     \end{center}
     \begin{center}
    \begin{tikzcd} [sep=1.5cm]
     0 \arrow[r] & 
     I_{>0}/I_1 \arrow[r] & 
     Q \arrow[r] & 
     \F \arrow[r] &
     0 
     \end{tikzcd}
     \end{center}
     \begin{center}
    \begin{tikzcd} [sep=1.5cm]
     0 \arrow[r] & 
     I_{>0}/I_{>1} \arrow[r] & 
     P \arrow[r] & 
     \F \arrow[r] &
     0 
     \end{tikzcd}
     \end{center} We then factor through $\F$ to obtain maps $\al:Q\to P$ and $\bt:P\to P$. We have already found an injective resolution for $\F$; we use these maps between $P$ and $Q$ to form injective resolutions \begin{center}
    \begin{tikzcd} [sep=1.5cm]
     0 \arrow[r] & 
     I_{>0}/I_1 \arrow[r] & 
     Q \arrow[r,"\al"] & 
     P \arrow[r] &
     Q \arrow[r] &
     0 
     \end{tikzcd}
     \end{center} and 
     \begin{center}
    \begin{tikzcd} [sep=1.5cm]
     0 \arrow[r] & 
     I_{>0}/I_{>1} \arrow[r] & 
     P \arrow[r,"\bt"] & 
     P \arrow[r] &
     Q \arrow[r] &
     0  \;\;\;.
     \end{tikzcd}
     \end{center}\end{exmp}

\begin{exmp} \label{more inj res for p mods} For $0<q<1$, we note the following short exact sequences. \begin{center}
    \begin{tikzcd} [sep=1.5cm]
     0 \arrow[r] & 
     A/I_q \arrow[r] & 
     Q \arrow[r] & 
     A/I_{1-q} \arrow[r] &
     0 
     \end{tikzcd}
     \end{center}
     \begin{center}
     \begin{tikzcd} [sep=1.5cm]
     0 \arrow[r] & 
     A/I_{>q} \arrow[r] & 
     P \arrow[r] & 
     A/I_{1-q} \arrow[r] &
     0 
     \end{tikzcd}
     \end{center}
     \begin{center}
     \begin{tikzcd} [sep=1.5cm]
     0 \arrow[r] & 
     I_{>0}/I_q \arrow[r] & 
     Q \arrow[r] & 
     A/I_{>1-q} \arrow[r] &
     0 
     \end{tikzcd}
     \end{center}
     \begin{center}
     \begin{tikzcd} [sep=1.5cm]
     0 \arrow[r] & 
     I_{>0}/I_{>q} \arrow[r] & 
     P \arrow[r] & 
     A/I_{>1-q} \arrow[r] &
     0 
     \end{tikzcd}
     \end{center} We then form composites $\al_q$, $\bt_q$, $\gm_q$ and $\dl_q$ given by \begin{center}
    \begin{tikzcd} [sep=1.5cm]
     Q \arrow[r,twoheadrightarrow] & 
     A/I_q \arrow[r,rightarrowtail] & 
     Q
     \end{tikzcd}
     \end{center} \begin{center}
    \begin{tikzcd} [sep=1.5cm]
     P \arrow[r,twoheadrightarrow] & 
     A/I_q \arrow[r,rightarrowtail] & 
     Q
     \end{tikzcd}
     \end{center}
     \begin{center}
    \begin{tikzcd} [sep=1.5cm]
     Q \arrow[r,twoheadrightarrow] & 
     A/I_{>q} \arrow[r,rightarrowtail] & 
     P
     \end{tikzcd}
     \end{center}
     \begin{center}
    \begin{tikzcd} [sep=1.5cm]
     P \arrow[r,twoheadrightarrow] & 
     A/I_{>q} \arrow[r,rightarrowtail] & 
     P
     \end{tikzcd}
     \end{center} respectively, and see that $\im(\al_q)=\ker(\al_{1-q})$ and $\im(\al_{1-q})=\ker(\al_q)$. So \begin{center}
    \begin{tikzcd} [sep=1.5cm]
     Q \arrow[r,"\al_{1-q}"] & 
     Q \arrow[r,"\al_q"] & 
     Q \arrow[r,"\al_{1-q}"] & 
     Q \arrow[r,"\al_q"] &
     Q \arrow[r,"\al_{1-q}"] & 
     \cdots
     \end{tikzcd}
     \end{center} is exact. Using the above short exact sequences as starting points, together with the above maps between $P$ and $Q$, we can form injective resolutions for $A/I_q$, $A/I_{>q}$, $I_{>0}/I_q$ and $I_{>0}/I_{>q}$ of the form \begin{center}
    \begin{tikzcd} [sep=1.5cm]
     0 \arrow[r] & 
     M \arrow[r] & 
     I_0 \arrow[r] & 
     I_1 \arrow[r] &
     Q \arrow[r,"\al_{1-q}"] & 
     Q \arrow[r,"\al_q"] &
     \cdots \;\;\;,
     \end{tikzcd}
     \end{center} where $I_0$ and $I_1$ are either $P$ or $Q$, depending on the type of standard basic module. \end{exmp}
     
\begin{exmp} \label{inj res summary for p mods} Let $M$ be a basic $P$-module. If it has length $0$ or $1$, then its injective dimension is no greater than $2$. Otherwise, it has an injective resolution of the form \begin{center}
    \begin{tikzcd} [sep=1.3cm]
     0 \arrow[r] & 
     M \arrow[r] & 
     I_0 \arrow[r] & 
     I_1 \arrow[r] &
     Q \arrow[r,"\al"] &
     Q \arrow[r,"\bt"] &
     Q \arrow[r,"\al"] &
     \cdots
     \end{tikzcd}
     \end{center} for some maps $\al,\bt:Q\to Q$. \end{exmp}
     
\begin{proof} This follows easily from Examples \ref{simple inj res for p mods} and \ref{more inj res for p mods}. \end{proof}

\begin{prop} $\Ext_{\mcM'}^{i+2}(M,N)=\Ext_{\mcM'}^i(M,N)$ for $i\geq 3$. \end{prop}

\begin{proof} First suppose that $N$ is basic. Using Example \ref{inj res summary for p mods}, we can form an injective resolution for $N$, and thus an eventually 2-periodic cochain complex \begin{center}
    \begin{tikzcd} [sep=1.3cm]
     0 \arrow[r] & 
     I_0 \arrow[r,] & 
     I_1 \arrow[r] &
     C \arrow[r,"\al"] &
     C \arrow[r,"\bt"] &
     C \arrow[r,"\al"] &
     \cdots
     \end{tikzcd}
     \end{center} giving $\Ext_{\mcM'}^j(M,N)$ for $j\geq 0$. The result then follows from the 2-periodicity of this cochain complex and the additivity of the $\Ext_{\mcM'}^i$ functors. \end{proof}

\begin{theorem} There does not exist a triangulated category $\mcC$ with $\Fr(\mcC)\simeq\mcM'$. \end{theorem}

\begin{proof} Suppose that such a $\mcC$ does exist. Then $\mcM'\simeq\Fr(\mcC)$ is a Frobenius category. Since $Q$ is injective in $\mcM'$, it must also be projective in $\mcM'$. Then since \begin{center}
    \begin{tikzcd}[sep=1.5cm]
     0 \arrow[r] & 
     \F \arrow[r] & 
     P \arrow[r] & 
     Q \arrow[r] &
     0 
     \end{tikzcd}
     \end{center} is exact in $\mcM'$, this sequence splits, so that $P$ is isomorphic to $\F\oplus Q$, both as a $P$-module and as an $A$-module, contradicting the uniqueness of decomposition of multibasic $A$-modules. \end{proof}

\begin{prop} \label{inj t ref p mod implies fin gen} Every injective and $\Ta$-reflexive $P$-module is isomorphic to one of the form $P^m\oplus Q^n$ for some $m,n\geq 0$, and therefore finitely generated. \end{prop}

\begin{proof} Let $M$ be injective and $\Ta$-reflexive. Then by Theorem \ref{classification of injective p mods}, it must be an injective hull of $\left(\bigoplus_{i\in I}P\right)\oplus\left(\bigoplus_{j\in J}Q\right)$. Then, as a submodule, $\left(\bigoplus_{i\in I}P\right)\oplus\left(\bigoplus_{j\in J}Q\right)$ must itself be $\Ta$-reflexive, so that both $I$ and $J$ are finite by Lemma \ref{direct sum is reflexive implies finitely many summands}. So $M$ is isomorphic to $P^m\oplus Q^n$ for some $m,n$, and therefore finitely generated. \end{proof}

\begin{definition} \cite[Definition 12.2]{ShS14} \label{def ungraded triang structure} An \textit{ungraded triangulation structure} for an (ungraded) self-injective ring $R$ consists of 

1. A full additive subcategory $\mcC$ of $\Mod_R$ consisting solely of injective modules and containing $R$. 

2. A triangulated structure (in the sense of Definition \ref{definition of triangulated structure}) on $\mcC$ compatible with the identity functor. \end{definition}

\begin{lem} \cite[Lemma 12.3]{ShS14} \label{exact triangles are exact sequences} Let $(\mcC,\Dl)$ be an ungraded triangulation structure for an ungraded self-injective ring $R$. Then exact triangles in $\mcC$ are exact sequences in $\Mod_R$. \end{lem}

\begin{proof} Let \begin{center}
    \begin{tikzcd}[sep=1.5cm]
     X \arrow[r] & 
     Y \arrow[r] & 
     Z \arrow[r] & 
     X[1] 
     \end{tikzcd}
     \end{center} be exact in $\mcC$. Since $R\in\mcC$, $\mcC$ is full, and $\Hom_{\mcC}(R,-)$ is exact, we see that \begin{center}
    \begin{tikzcd}[sep=1.5cm]
     \Hom(R,X) \arrow[r,"{\Hom(R,u)}"] & 
     \Hom(R,Y) \arrow[r,"{\Hom(R,v)}"] & 
     \Hom(R,Z) \arrow[r,"{\Hom(R,w)}"] & 
     \Hom(R,X[1]) 
     \end{tikzcd}
     \end{center} is exact in $\Mod_R$. Then since $\Hom(R,M)\simeq M$ for all $M$, we see that \begin{center}
    \begin{tikzcd}[sep=1.5cm]
     X \arrow[r,"u"] & 
     Y \arrow[r,"v"] & 
     Z \arrow[r,"w"] & 
     X[1] 
     \end{tikzcd}
     \end{center} is exact in $\Mod_R$, as required. \end{proof}

\begin{theorem} \label{no ungraded triang structure for inf root alg} Let $(\mcC,\Dl)$ be an ungraded triangulation structure for $P$. Then $\mcC$ is not contained in $\mcM'$. \end{theorem}

\begin{proof} Suppose that $\mcC\leq\mcM'$. Let $v:P\to P$ be given by $v(a+I_{>1})=ta+I_{>1}$. Since $P\in\mcC$ and $\mcC$ is full, this is a morphism in $\mcC$. We can use it to form an exact triangle \begin{center}
    \begin{tikzcd}[sep=2cm]
     X \arrow[r,"u"] & 
     P \arrow[r,"v"] & 
     P \arrow[r,"w"] & 
     X 
     \end{tikzcd}
     \end{center} in $\mcC$, which must be exact in $\Mod_P$ by Lemma \ref{exact triangles are exact sequences}. So $\im(u)=\ker(v)=I_{>0}/I_{>1}$. But since $X\in\mcC$, we know that $X$ is injective, and by our assumption, it must also be in $\mcM'$. So $X$ must be finitely generated by Proposition \ref{inj t ref p mod implies fin gen}. But then $I_{>0}/I_{>1}=\im(u)$ is finitely generated, a contradiction. \end{proof}

\newpage

\chapter{Appendix: Reflexive Modules Over Discrete Valuation Rings}

In this appendix, we very briefly summarise some key results in the theory of $\Ta$-reflexive and multibasic modules over a discrete valuation ring with a sequential completeness condition. The author investigated this case as a warm-up exercise for the more involved case in Chapters \ref{hahn ring chapter} and \ref{reflexivity}. As the proofs are all essentially just simplified versions of the corresponding results in Chapters \ref{hahn ring chapter} and \ref{reflexivity}, we shall simply document the results here. 

As with the Hahn ring and the infinite root algebra, the $\Ta$-reflexives and the multibasics are the same. However, over a discrete valuation ring the theory is considerably simplified by the fact that all of the ideals are principal. Many of the arguments of Section \ref{reflexive graded modules} apply here, so we shall focus on summarising the differences between this case and that of the Hahn ring in Chapters \ref{hahn ring chapter} and \ref{reflexivity}.

\section{Reflexive Modules Over Discrete Valuation Rings}

Fix a discrete valuation ring $A$ and a generator $p$ of its unique non-zero maximal ideal. Denote the field of fractions of $A$ by $K$, which we think of as an $A$-module. Each $a\in K^\x$ can be expressed uniquely in the form $p^ku$ for $k\in\Z$ and $u\in A^\x$. Setting $\nu(p^ku)=k$ and $\nu(0)=\infty$ defines a valuation $\nu$ on $K$. Set $$I_k=p^kA=\{a\in K:\nu(a)\geq k\}$$ for each $k\in\Z$. 

\textbf{We shall assume that for every sequence $(a_n)$ in $K$ with $a_{n+1}-a_n\in I_n$ for all $n$, that there exists $a\in K$ with $a_n-a\in I_n$ for all $n$.} 

Set $\Ta=K/I_1$ and $D=\Hom(-,\Ta)$. As with the Hahn ring, $\Ta$ is injective, so that $D$ is exact. Also $DM=0$ if and only if $M=0$. For each $A$-module $M$, we denote the evaluation map $$m\mapsto(\al\mapsto\al(m)):M\to D^2M$$ by $\chi_M$ and say that $M$ is \textit{$\Ta$-reflexive} if $\chi_M$ is an isomorphism. We write $\mcM$ to denote the full subcategory of $\Mod_A$ consisting of the $\Ta$-reflexive modules.

\begin{prop} \label{properties of theta-reflexive modules} $\mcM$ is abelian and is closed under submodules, quotients, extensions and isomorphisms. $D$ restricts to a duality functor $\mcM\to\mcM$, where $DM\in\mcM$ implies $M\in\mcM$. \end{prop}

\begin{exmp} \label{k tns m has finite dimension over k} If $M\in\mcM$, then $K\tns M$ and $\Hom(M,K)$ have finite dimension over $K$. Also $A/I_1\tns M=M/I_1M$ has finite dimension over the field $A/I_1$. \end{exmp}

\begin{definition} We call a non-zero $A$-module $M$ \textit{basic} if for all $x,y\in M$, we have that $x\in Ay$ or $y\in Ax$, and \textit{multibasic} if it is a direct sum of a (possibly empty) finite list of basic modules. Every basic module is isomorphic to $K$, $A$, $\Ta$ or $A/I_n$ for some $n\geq 1$. We call this smaller collection of basic modules the \textit{standard basic modules}. Denote the field $A/I_1$ by $\F$. \end{definition}

\begin{exmp} The following two tables demonstrate how $\Hom(M,N)$ and $M\tns N$ act on the standard basic modules, where $M$ varies along the rows and $N$ varies along the columns. Below these are the values for several invariants of the standard basic modules.

\[ \begin{array}{|c|c|c|c|c|}                           \hline
     \Hom(M,N) & K & A & \Ta     & A/I_n           \\  \hline
     K           & K & 0 & K       & 0               \\  \hline
     A           & K & A & \Ta     & A/I_n           \\  \hline
     \Ta         & 0 & 0 & A       & 0               \\  \hline 
     A/I_m       & 0 & 0 & A/I_m   & A/I_{\min\{m,n\}} \\  \hline 
    \end{array} \]  
    
    \[ \begin{array}{|c|c|c|c|c|}                             \hline
     M\tns N & K & A     & \Ta & A/I_n             \\  \hline
     K       & K & K     & 0   & 0                 \\  \hline
     A       & K & A     & \Ta & A/I_n             \\  \hline
     \Ta     & 0 & \Ta   & 0   & 0                 \\  \hline 
     A/I_m   & 0 & A/I_m & 0   & A/I_{\min\{m,n\}} \\  \hline 
    \end{array} \]

\[ \begin{array}{|c|c|c|c|c|}\hline
     M                            & K & A & \Ta & A/I_n \\ \hline
     \dim_K(K\tns M)              & 1 & 1 & 0 & 0       \\ \hline
     \dim_K(K\tns DM)             & 1 & 0 & 1 & 0       \\ \hline
     \dim_{\F}(\F\tns M)    & 0 & 1 & 0 & 1       \\ \hline 
     \dim_{\F}(\F\tns DM)   & 0 & 0 & 1 & 1       \\ \hline
     \ann(M)                      & 0 & 0 & 0 & I_n     \\ \hline 
    \end{array} \] \end{exmp}

\begin{prop} \label{basic module isomorphic to unique standard basic module} Every basic $A$-module is isomorphic to a unique standard basic $A$-module. \end{prop}

\begin{theorem} The decomposition of a multibasic module into standard basic modules is unique up to the order of the summands. \end{theorem}

\begin{prop} If $M$ and $N$ are multibasic, then so are $\Hom(M,N)$ and $M\tns N$. \end{prop}

\begin{lem} \label{theta-reflexive and zero tensor with dm implies fin gen} If $M\in\mcM$ and $K\tns DM=0$, then $M$ is finitely generated. \end{lem}

\begin{lem} \label{if k tns m is zero then m is sum of torsion basics} If $M\in\mcM$ and $K\tns M=K\tns DM=0$, then $M$ is a finite direct sum of modules of the form $A/I_n$. \end{lem}

\begin{theorem} $M\in\mcM$ if and only if it is multibasic. \end{theorem}

\end{document}